\magnification=\magstephalf

\input amstex
\documentstyle{amsppt}
\NoBlackBoxes
\hoffset 1.25truecm
\hsize=12.4 cm
\vsize=19.7 cm
\TagsOnRight

\font\bbb=bbold10 scaled\magstephalf
\font\bbt=bbold8

\font\tenscr=rsfs10 
\font\sevenscr=rsfs7 
\font\fivescr=rsfs5 
\skewchar\tenscr='177 \skewchar\sevenscr='177 \skewchar\fivescr='177
\newfam\scrfam \textfont\scrfam=\tenscr \scriptfont\scrfam=\sevenscr
\scriptscriptfont\scrfam=\fivescr
\def\scr{\fam\scrfam}

\def\az{\alpha}  \def\bz{\beta}
    \def\dz{\delta}
\def\ez{\eta}    \def\fz{\varphi}
\def\gz{\gamma}  \def\kz{\kappa}
\def\lz{\lambda} \def\mz{\mu}
\def\nz{\nu}     \def\oz{\omega}
\def\pz{\pi}

\def\tz{\tau}        \def\uz{\theta}
\def\vz{\varepsilon} \def\xz{\xi}
\def\zz{\zeta}

\def\ggz{\Gamma}  \def\ooz{\Omega}

\def\kq{{\scr K}}

\def\flsh{\flushpar}

\def\qd{\quad}
\def\qqd{\qquad}

\def\defn{Definition}
\def\lmm{Lemma}
\def\prp{Proposition}
\def\thm{Theorem}
\def\crl{Corollary}
\def\rmk{Remark}
\def\prf{Proof}
\def\xmp{Example}

\def\d{\text{\rm d}}

\def\supp{\text{\rm supp}\,}
\def\le{\leqslant}
\def\ge{\geqslant}

\def\dbl{\text{\bbb{1}}}

 \topmatter
\topinsert \vskip-0.65truecm \captionwidth{10 truecm}Front. Math. China 5:3 (2010), 379--515
\endinsert
\title{Speed of stability for birth--death processes}\endtitle
\author{Mu-Fa Chen}\endauthor
\affil{(Beijing Normal University)\\
        May 25, 2007---November 10, 2009}\endaffil
\thanks {Research supported in part by the Creative Research Group Fund of the
         National Natural Science Foundation of China (No. 10721091),
         by the ``985'' project from the Ministry of Education in China.}\endthanks
\address {\rm \flsh School of Mathematical Sciences, Beijing Normal University,
Laboratory of Mathematics and Complex Systems (Beijing Normal University),
Ministry of Education, Beijing 100875,
    The People's Republic of China.\newline E-mail: mfchen\@bnu.edu.cn\newline Home page:
    http://www.bnu.edu.cn/\~{}chenmf/main$\_$eng.htm} \endaddress
\abstract{This paper is a continuation of the
study on the stability speed for Markov processes. It
extends the previous study of the ergodic convergence speed to
the non-ergodic one, in which the processes are even allowed to be explosive
or to have general killings. At the beginning stage, this paper is
concentrated on the birth-death processes. According to the
classification of the boundaries, there are four cases plus one
having general killings. In each case, some dual variational
formulas for the convergence rate are presented, from which the
criterion for the positivity of the rate and an approximating
procedure of estimating the rate are deduced. As the first step of
the approximation, the ratio of the resulting bounds is usually no
more than 2. The criteria as well as basic estimates for more
general types of stability are also presented. Even though
the paper contributes mainly to the non-ergodic case, there are some
improvements in the ergodic one. To illustrate the power of the
results, a large number of examples are included. }\endabstract
\endtopmatter

\document

\head{Contents}\endhead
{\roster
\item "\S\, 1." Introduction \dotfill 2
\item "\S\, 2." Absorbing (Dirichlet) boundary at infinity: dual variational formulas \dotfill 7
\item "\S\, 3." Absorbing (Dirichlet) boundary at infinity: criterion, approximating procedure and examples\dotfill 27
\item "\S\, 4." Absorbing (Dirichlet) boundary at origin and reflecting (Neumann) boun\-dary at infinity\dotfill 39
\item "\S\, 5." Dual approach\dotfill  50
\item "\S\, 6." Reflecting (Neumann) boundaries at origin and infinity\,(ergodic case) \dotfill 59
\item "\S\, 7." Bilateral absorbing (Dirichlet) boundaries \dotfill 71
\item "\S\, 8." Criteria for Poincar\'e-type inequalities\dotfill 92
\item "\S\, 9." General killing\dotfill 99
\item "\S\, 10." Notes\dotfill 121
\item "{}" References \dotfill 130
\endroster}

\head{1. Introduction}
\endhead

Consider a birth--death process on the nonnegative integers ${\Bbb
Z}_+$ with birth rates $b_n>0\,(n\ge 0)$ and death rates
$a_n>0\,(n\ge 1)$. Define
$$\mu_0^{}=1,\quad
\mu_n= \frac{b_0\cdots b_{n-1}}{a_1\cdots a_n},\qquad n\ge  1. \tag 1.1$$
We say that the birth--death process is nonexplosive if the following
Dobrushin's uniqueness criterion holds:
$$\sum_{k=0}^ \infty \frac{1}{b_k \mu_k} \sum_{i=0}^k \mu_i
\bigg[=\sum_{i=0}^\infty  \mu_i\sum_{k=i}^\infty \frac{1}{b_k \mu_k}\bigg]
 = \infty \tag 1.2$$
(cf. Dobrushin (1952), or Wang and Yang (1992, \crl\;5.2.1), or [10; \crl\;3.18]).
This implies a useful condition that
$$\sum_{k=0}^ \infty \bigg(\frac{1}{b_k \mu_k}+ \mu_k\bigg)
= \infty.\tag 1.3
$$
When $\sum_0^\infty \mu_k<\infty$, each of (1.2) and (1.3) is equivalent to
the recurrent condition: $\sum_0^\infty (b_n\mu_n)^{-1}=\infty$.
Otherwise, (1.3) cannot imply (1.2)
since one can easily construct a counterexample so that
$\sum_0^\infty \mu_k=\infty$ but
$$\sum_{i=0}^\infty \mu_i \sum_{k=i}^\infty \frac{1}{b_k \mu_k}<\infty.$$
Thus, under (1.3), the process may not be unique.

It is well known that for a birth--death process, the transition probabilities
$(p_{ij}(t))$ satisfy
$$\lim_{t\to\infty} p_{ij}(t)=:\pi_j\ge 0 \tag 1.4$$
for all $i, j\in {\Bbb Z}_+$. We are now interested in the exponential convergence rate
$$\az^*=\sup\Big\{\az:\; |p_{ij}(t)-\pz_j|=O\big(\exp[-\az\, t]\big)\; \text{as} \;
t\to\infty \;\text{for all} \; i,j\in E \Big\}.\tag 1.5$$
In the ergodic case \big(i.e., $\lim_{t\to\infty} p_{ij}(t)>0$
for all $i,\, j$\big), we have  $Z:=\sum_{j=0}^\infty\mu_j<\infty$ and then $\pi_j:=\mu_j/Z>0$
for all $j\ge 0$. In this case, the problem has been well studied, see, for instance, van Doorn (1981; 2002),
Zeifman (1991), Kijima (1997), [2, 12], and the references therein. The problem becomes
trivial in the zero-recurrent case for general irreducible Markov chains, since we have on the one hand
$\pi_j=0$ for all $j$, and on the other hand,
$\int_0^\infty p_{ii}(t)\d t=\infty$ for all $i.$
Hence, the exponential decay can only happen in the transient case:
$$\sum_{n=0}^\infty\frac{1}{b_n\mu_n}<\infty. \tag 1.6$$
Since the process is $\mu$-symmetric: $\mu_i p_{ij}(t)=\mu_j p_{ji}(t)$ for all $i$, $j$ and $t$, it is natural,
as we did in the ergodic case, to use the $L^2$-theory. As usual, denote by $\|\cdot\|$ and $(\cdot, \cdot)$
the norm and the inner product on the real Hilbert space $L^2(\mu)$, respectively. Let
$$\kq=\{f: f \text{ has finite support}\}. \tag 1.7$$
Define $$
D(f)=\sum_{i\ge 0} \mu_i b_i(f_{i+1}-f_i)^2=\sum_{i\ge 1} \mu_i a_i(f_{i}-f_{i-1})^2$$
with the minimal domain ${\scr D}^{\min}(D)$ consisting of the functions in the closure of
$\scr K$ with respect to the norm $\|\cdot\|_D$: $\|f\|_D^2=\|f\|^2+ D(f)$.
Next, define
$$\lz_0=\inf\{D(f): \|f\|=1, f\in {\scr K}\}= \inf\{D(f): \|f\|=1, f\in {\scr D}^{\min}(D)\}.$$
From now on, we often write $f_\infty$ or $f(\infty)$ as the limit
of $f$ at infinity provided it exists. In the definition of $\lz_0$,
it is natural to add the boundary condition $f_\infty=0$ but this
can be ignored since on the one hand, for each $f\in \kq$, we have
$f_{\infty}=0$, and on the other hand $\kq$ is a core of the
Dirichlet form $\big(D, {\scr D}^{\min}(D)\big)$ (i.e., the form is regular) by
[10; \prp\;6.59]. For a large part of the paper, we are dealing with this minimal
Dirichlet form or the minimal process.

We now make a connection between $\az^*$ and $\lz_0$.
The proofs of the next three propositions are delayed for a moment.

\proclaim{\prp\;1.1} For a general non-ergodic symmetric semigroup $\{P_t\}_{t\ge 0}$ with
Dirichlet form $(D, {\scr D}(D))$\,(not necessarily regular) on $L^2(\mu)$, the parameter $\lz_0$,
$$\lz_0=\inf\{D(f): \|f\|=1,\; f\in {\scr D}(D)\}, \tag 1.8$$
is the largest $\vz$ such that
$$\|P_t f\|\le \|f\|\,e^{-\vz t},\qqd t\ge 0, f\in L^2(\mu). \tag 1.9$$
\endproclaim

It was proved in [2; \thm\;5.3] that for birth--death
processes, under (1.2), the exponentially ergodic convergence rate $\az^*$
coincides with the $L^2$-exponential one, denoted by $\lz_1$:
$$\|P_t f-\pi(f)\|\le \|f-\pi(f)\|e^{-\lz_1 t}\qqd \text{for all $t\ge 0$ and $f\in L^2(\mu)$}, $$
where $\pi(f)=\int f \d \mu/\mu(E)$. For non-ergodic birth--death processes, we have similarly ${\az^*}=\lz_0$,
as mentioned at
the end of [2]. Here is a generalization.

\proclaim{\prp\;1.2} For a general non-ergodic $\mu$-symmetric Markov chain with
Dirichlet form $(D, {\scr D}(D))$, we have ${\az^*}=\lz_0$ defined by (1.8).
\endproclaim

About (1.3), we have the following result.

\proclaim{\prp\;1.3} Let ${\scr D}^{\max}(D)=\{f\in L^2(\mu): D(f)<\infty\}$. Then the Dirichlet form $(D, {\scr D}^{\max}(D))$ is
regular iff $(1.3)$ holds.  In other words, the Dirichlet form corresponding to the rates
$(a_i)$ and $(b_i)$ is unique iff $(1.3)$ holds.
\endproclaim

\prp\;1.2 reduces the study on $\az^*$ to the
first (or principal) eigenvalue $\lz_0$. This is the starting point of
this paper. In the two cases we have discussed so far, the state $0$ is a
reflecting (Neumann) boundary, denoted by code ``N''. For $\lz_1$, since
the process starting from any point will certainly come back, the
infinity may be regarded as a reflecting (Neumann) boundary.
However, for $\lz_0$, the situation is different. As we will prove
in the next section, the corresponding eigenfunction decreases to
zero at infinity. Hence, the infinity may be regarded as an absorbing (Dirichlet)
boundary, denoted by code ``D''. Thus, for the temporary convenience, we
rewrite $\lz_1=\lz^{\text{\rm NN}}$ and $\lz_0=\lz^{\text{\rm ND}}$. Replacing the Neumann
boundary at $0$ by the Dirichlet one (i.e., $b_0=0$), we obtain two
more cases for which we have the decay rates (eigenvalues)
$\lz^{\text{\rm DN}}$ and $\lz^{\text{\rm DD}}$, respectively. The main body of this paper is devoted to
study these four cases. Now, the rate $\az^*$ coincides with, case by case, one of $\lz^{\text{\rm NN}}$, $\lz^{\text{\rm ND}}$,
$\lz^{\text{\rm DN}}$, and $\lz^{\text{\rm DD}}$.
Here are simple examples to show the difference in the different cases.

\proclaim{\xmp s\;1.4} \roster
\item Let $a_i=\dz i$, $b_i=\bz i +\gz$, $\dz>\bz$. Then
$\lz^{\text{\rm NN}}=\dz-\bz$ if $\gz>0$ and so does $\lz^{\text{\rm DN}}$ if $\gz=0$.
\item Let $a_i=i$, $b_i=2(i+\gz)$. Then $\lz^{\text{\rm ND}}=\gz$ if $\gz>0$ and $\lz^{\text{\rm DD}}=1$ if $\gz=0$.\endroster
\endproclaim

The rate in the first example is the difference of the coefficients of leading terms,
independent of $\gz$. This is somehow natural. Surprisingly, the second one is determined by the constant term only
except $\gz=0$ at which case there is a jump from $\lz^{\text{\rm ND}}$ to $\lz^{\text{\rm DD}}$.
Thus, for the convergence rate, the role played by the parameters $(a_i, b_i)$ is
mazed and then one may wonder how far we can go (see \thm\;1.5 below for a preliminary answer).

The main body of the paper is devoted to the quantitive study of the convergence rate.
For this, our key result (variational formulas) plays a full power.
For those readers who are interested only in the qualitative criteria
and basic estimates, here is a short statement.

\proclaim{\thm\;1.5\,(Criterion and basic estimates)} Let (1.3) hold.
Then in spite of $b_0>0$ or $b_0=0$, the exponential convergence rate $\az^*$ defined in (1.5)
for the unique process is positive
\roster
\item iff $\dz^{(4.4)}<\infty$ in the case of $\sum_i \mu_i<\infty$; and otherwise,
\item iff $\dz^{(3.1)}<\infty$,
\endroster
where
$$\dz^{(4.4)}=\sup_{n\ge 1}\sum_{i=1}^n \frac{1}{\mu_i a_i} \sum_{j=n}^\infty \mu_j,\qqd
\dz^{(3.1)}=\sup_{n\ge 0}\sum_{i=0}^n\mu_i \sum_{j=n}^\infty \frac{1}{\mu_j b_j}.$$
More precisely, we have the basic estimate $\dz^{-1}/4\le \az^*\le \dz^{-1}$, where the constant $\dz$ is
equal to $\kz^{(6.13)}$ or $\kz^{(7.5)}$ according to $b_0>0$ or $b_0=0$, respectively:
$$\big(\kz^{(6.13)}\big)^{-1}=\inf_{m>n\ge 0}\bigg[\bigg(\sum_{i=0}^{n}\mu_i\bigg)^{-1}
+\bigg(\sum_{i=m}^\infty \mu_i\bigg)^{-1}\bigg]\bigg(\sum_{j=n}^{m-1} \frac{1}{\mu_j b_j}\bigg)^{-1},$$
$$\big(\kz^{(7.5)}\big)^{-1}=\inf_{m\ge n\ge 1}\bigg[\bigg(\sum_{i=1}^n\frac{1}{\mu_i a_i}\bigg)^{-1}
+ \bigg(\sum_{i=m}^{\infty}\frac{1}{\mu_i
b_i}\bigg)^{-1}\bigg]\bigg(\sum_{j=n}^m \mu_j\bigg)^{-1}.$$
Here, the superscript of $\kz^{(7.5)}$, for instance, means that it is in the case studied in Section 7 and
the constant is given in (7.5).
\endproclaim

The proof of \thm\;1.5 and its extension are given in Section 7.
The more general qualitative results are presented in Section 8 and
in Summary 9.12 for the killing case.

To have an impression about the progress made in the paper, let us have a look at the new points
made in the well-developed case, Section 6.
\roster
\item The uniqueness condition (1.2) is replaced by using the maximal process.
\item The more complete dual variational formulas are presented in \thm\,6.1.
\item Even though the criterion \thm\;6.2 is known before, the upper bound in its improvement (\crl\;6.4)
is newly added so that the ratio of the bounds is now no more than 2 as shown by a group of examples. Moreover,
a new criterion (\crl\;6.6) which has been expected naturally (in view of \thm\;4.2) for a long time, is now presented.
\item A more effective sequence for the upper estimate given in \thm\;6.3 is introduced to replace the original one.
The monotonicity of the approximating sequences are proved here for the first time.
\endroster
We now return to prove the propositions above.
\demo{\prf\;of \prp\;$1.1$} Replace by $\vz_{\max}$ the largest exponential rate in
(1.9). Then we have $\vz_{\max}\ge 0$ because of the contractivity of the semigroup
in every $L^p$-space $(p\ge 1)$. We need to show that $\lz_0=\vz_{\max}$. The proof of $\lz_0\ge
\vz_{\max}$ is easier since by an elementary property of the
Dirichlet form and (1.9), we have for every $f$ with $\|f\|=1$,
$$D(f)=\lim_{t\downarrow 0}\uparrow \frac{1}{t} (f-P_t f,\, f)
  \ge \lim_{t\downarrow 0}\frac{1}{t} (1-e^{-\vz_{\max} t})=\vz_{\max},\tag 1.10$$
where $\lim\! \uparrow$ means an increasing limit.
Hence, we have $\lz_0\ge \vz_{\max}$.
To prove $\vz_{\max}\ge \lz_0$, assume that $\lz_0>0$. Otherwise, the
assertion is trivial. Noticing that $D(f)=(-\ooz f,\, f)$ for the generator $\ooz$ of $\{P_t\}$
on $L^2(\mu)$ and for every $f\in {\scr D}(\ooz) $, we have
$$\frac{\d}{\d t}\|P_t f\|^2=2(P_t f,\; \ooz P_t f)=-2D(P_t f).\tag 1.11$$
Next, since $P_t f\in {\scr D}(D)$ for each $f\in L^2(\mu)$, by the
definition of $\lz_0$, we have
$$-2D(P_t f)\le -2 \lz_0\|P_t f\|^2.$$
Thus, $\|P_t f\|\le \|f\|e^{-\lz_0 t}$ for all $t\ge 0$ and $f\in
{\scr D}(\ooz)$, and then for all $f\in L^2(\mu)$ since the density of ${\scr D}(\ooz)$ in
$L^2(\mz)$ and the contractivity of the semigroup $\{P_t\}_{t\ge 0}$.
The assertion now follows since $\vz_{\max}$ is the
largest rate. \qed\enddemo

\demo{\prf\;of \prp\;$1.2$} The proof for ${\az^*}\ge \lz_0$ is rather easy. Simply applying \prp\;1.1 to the indicator
function $f=\dbl_{\{k\}}$, we obtain
$$p_{ik}(t)\le \sqrt{\mu_k/\mu_i}\, e^{-\lz_0 t}.$$
Note that this also provides a non-trivial estimate of the constant in (1.5).

To prove that $\lz_0\ge {\az^*}$, we may assume that ${\az^*}>0$.
One may follow the proof of [12; proof of part (4) of Theorem 8.13]. In the last
part of the original proof, we have
$$\|P_t f\|^2=(f, P_{2t} f)\le \|f\|_\infty^2 e^{-2{\az^*} t}\sum_{i,j\in\text{supp}(f)}\mu_i C_{ij}$$
for every bounded $f$ with compact support.
Here, we have used the assumption that
$p_{ij}(t)\le C_{ij} e^{-{\az^*} t}$.
\qed\enddemo

\demo{\prf\;\prp\;$1.3$} Since the $Q$-matrix is conservative,
by [10;  {\lmm\;} 6.52 and \thm\;6.61], $(D, {\scr D}^{\max}(D))$
is a Dirichlet form and is indeed the maximal one. Note that in the
conservative case, every $Q$-process (in particular, the semigroup
generated by a Dirichlet form) satisfies the backward Kolmogorov's
equation by [10; \thm\;1.15\,(1)].

(a) Let (1.3) hold. Then the Dirichlet form should be regular.
Otherwise, we have two different birth--death semigroups generated
by $(D, {\scr D}^{\max}(D))$ and the minimal Dirichlet form $(D,
{\scr D}^{\min}(D))$, respectively. They satisfy first the backward
and then also the forward Kolmogorov's equations by [10;
\thm\;6.16]. This is impossible since condition (1.3) is the
uniqueness criterion for the process satisfying  the Kolmogorov's
equations simultaneously, due to Karlin and McGregor (1957a,
\thm\;15) (cf. Hou et al. (2000, \thm\;6.4.6\,(1); 1994,
\thm\;12.7.1)). Note that criterion (1.3) is equivalent to the
uniqueness for the process satisfying one of the Kolmogorov
equations since every symmetric process as well as the minimal one satisfies both of the equations. This is the
reason why (1.3) is weaker than (1.2).

(b) Next, let (1.3) fail. Then we have $\sum_i\mu_i<\infty$ and
$\sum_i (\mu_i b_i)^{-1}<\infty$. Moreover, (1.2) fails. Note that
the birth--death $Q$-matrix has at most a single exit boundary, and
there is precisely one if (1.2) fails. Besides, the non-trivial
(maximal) exit solution $z_\lambda$ is bounded from above by 1. In
view of [10;  \prp\;6.56], there are infinitely many Dirichlet
forms. The minimal one is regular but not the maximal one $(D, {\scr
D}^{\max}(D))$. \qed\enddemo

Actually, \prp\;1.3 is a particular case of a result we will study
at the end of Section 9 (\thm\;9.22).

The remainder of the paper is organized as follows. In the next two sections, we study
$\lz^{\text{\rm ND}}$. Sections 4, 6 and 7 are devoted to $\lz^{\text{\rm DN}}$,
$\lz^{\text{\rm NN}}$ and $\lz^{\text{\rm DD}}$, respectively. By exchanging N and D, we formally obtain a dual of $\lz^{\text{\rm ND}}$
and $\lz^{\text{\rm DN}}$ \big(resp. $\lz^{\text{\rm NN}}$ and $\lz^{\text{\rm DD}}$\big) which is studied in
Section 5 (resp.\,7). In each case, we present a group of dual variational
formulas for the first (non-trivial) eigenvalue. By using the formulas, we then deduce a criterion
for the positivity of the eigenvalue and an approximating procedure
for estimating the eigenvalue.
The criteria and basic estimates in a quite general setup are given in Section 8.
A closely related topic, having general killings, is studied in Section 9.
In the study of this paper, the author has benefited a great deal from our
previous work and from many authors' contribution. A part of the
contributions is noted in the context. In the ergodic case under (1.2),
a large number of references are given in [10, 12] and the author
apologizes for omitting them here.
At the end of the paper (Section 10), some remarks on the related results,
some open problems or open topics, and so on are discussed. The analog of \thm\;1.5
for one-dimensional diffusions is also included.

\proclaim{Notation 1.6} To be economical, we use the same notation
$\lz_0$, $\dz$, $\kz$, $I$ and $I\!I$ and so on,
from time to time in different sections with similar but
different meaning. To distinguish them if necessary, we write
$\lz_0^{(\#)}$ for instance
to denote the $\lz_0$ defined by formula $(\#)$.
\endproclaim

\head{2. Absorbing (Dirichlet) boundary at infinity: dual variational formulas}\endhead

This section begins with the study on the property of eigenfunction of $\lz_0$. It
is fundamental in our analysis and has been studied several times
before, see, for instance, [3;  Lemma 4.2]; [4; proofs of
\thm s 3.2 and 3.4]; Chen, Zhang and Zhao (2003, Section 2); Shao
and Mao (2007, Proposition 3.1). The main body of this section is devoted to prove a group
of variational formulas (\thm\;2.4 and \prp\;2.5). Their applications are given in the
next section.

Fix an integer $N$: $1\le N\le \infty$, and let $E=\{k\in {\Bbb Z}_+: 0\le k <N+1\}$.
Throughout the paper, the infinite case that $N=\infty$ is more essential but the finite
case that $N<\infty$ is also included which may be meaningful in matrices theory.
To avoid the confusion of these two
cases in reading the paper, one may read the infinite case first and then go back to
check the modification for the finite case. Besides, note that when $N<\infty$, neither $(1.2)$
nor $(1.3)$ is needed.

Let us start at a general situation. Consider the operator
$\ooz$ corresponding to the birth--death $Q$-matrix
with birth rates $b_i$, death rates $a_i$, and killing rates
$c_i\ge 0\,(i\in E)$ as follows.
$$\text{\hskip-2em}\ooz f(i)\!=\!b_i(f_{i+1}\!-\!f_i)\!+\!a_i (f_{i-1}\!-f_i)\!-\!c_i f_i,\;\; i\in E,\; f_{N+1}=0\text{ if }N<\infty.\tag 2.1$$
In other words, when $N<\infty$, the state $N+1$ is an absorbing
(Dirichlet) boundary. When $c_i\not\equiv 0$ for $1\le i<N$, unless otherwise stated,
we assume that $a_0=0$ and $b_N=0$ if $N<\infty$ (the other
$a_i$ and $b_i$ are positive), otherwise, simply replace the original
$c_0$ and $c_N$ by $a_0+c_0$ and $b_N+c_N$, respectively. Now,
since $a_0=0$, $f_{-1}$ is free in the last formula. The first
eigenvalue $\lz_0$ is now defined by
$$\lz_0=\inf\{D(f): \|f\|=1, f\in {\scr K}\},\tag 2.2$$
where
$$D(f)=\sum_{0\le i<N} \mu_i b_i (f_{i+1}-f_i)^2+ \sum_{i\in E}\mu_i c_i f_i^2,\qqd
 f_{N+1}=0\text{ if }N<\infty.\tag 2.3$$
We say that $g$ is an ``eigenfunction'' of $\lz\in {\Bbb R}$, if $g$ satisfies the ``eigenequation'':
$$\ooz g =-\lz g,\qqd  g_{N+1}=0\text{ if }N<\infty.\tag 2.4$$
Note that the ``eigenvalue'' and ``eigenfunction'' used in this paper
are in a generalized sense rather than the standard ones since here we do not require $g\in L^2(\mu)$.

\proclaim{\prp\;2.1}
\roster
\item Every eigenfunction $g$ of $\lz\in {\Bbb R}$ satisfies
$$\mu_k b_k(g_k-g_{k+1})=\sum_{i=0}^k (\lz-c_i)\mu_i  g_i,
\qqd k\in E,\; g_{N+1}=0\text{ if }N<\infty. \tag 2.5$$
\item If $\lz_0>0$, then $c_i\not\equiv 0\,(0\le i\le N)$ whenever $N<\infty$, and the non-zero
eigenfunction $g$ of $\lz_0$ is either positive or negative on
$E$.
\item The non-zero eigenfunction $g$ of $\lz=0$ is either positive
and nondecreasing, or negative and nonincreasing on $E$. Furthermore, let $g>0$ for
instance. Then $g_{k+1}>g_k$ for all $k: i\le k <N$ whenever $c_i>0$.
\endroster
\endproclaim

\demo{\prf}
(a) Recall the eigenequation
$$\ooz g(i)=b_i (g_{i+1}-g_i) +a_i (g_{i-1}-g_i)-c_i g_i=-\lz g_i,\qqd i\in E, \tag 2.6$$
or more generally, the Poisson equation
$$b_i (g_i-g_{i+1}) -a_i (g_{i-1}-g_i)=h_i,\qqd i\in E,\;
\text{$g_{N+1}=0$ if $N<\infty$}, \tag 2.7$$
for a given function $h$.
Multiplying both sides by $\mu_i$, we get
$$\mu_{i} b_{i} (g_i-g_{i+1})-\mu_{i-1} b_{i-1} (g_{i-1}-g_i)= \mu_i h_i, \qqd i\in E. \tag 2.8$$
When $i=0$, the second term on the left-hand side is set to be zero.
Making a summation over $i$, we obtain
$$\mu_k b_k (g_k-g_{k+1})=\sum_{i=0}^k \mu_i h_i, \qqd k\in E,\;
\text{$g_{N+1}=0$ if $N<\infty$}. \tag 2.9$$
With $h_i=(\lz-c_i)g_i$, this gives us the first assertion of the proposition.

(b) To prove the second assertion, note that $\lz_0=0$ if $c_i\equiv 0$
$(0\le i\le N<\infty)$ in which case both $0$ and $N$ are reflecting and
the process is ergodic. Now, since $\lz_0>0$,
one may assume that $g_0\ne 0$, otherwise $g_i\equiv 0$ by induction.
Next, replacing $g$ by $g/g_0$ if necessary, we can assume that $g_0=1$.
If $g$ is not positive, then there would exist a $k_0\in E$, $k_0\ge 1$ such that $g_i>0$ for $i<k_0$
and $g_{k_0}\le 0$. We then modify $g$ from $k_0$:  set $\tilde g_i=g_i$ for $i<k_0$
and $\tilde g_i=0$ for $i\ge k_0+1$. By choosing a suitable value $\vz>0$ at $k_0$,
the new function $\tilde g\in {\scr K}$
gives us $D(\tilde g)/\|\tilde g\|^2<\lz_0$,
which is a contradiction to the definition of $\lz_0$.
Hence, $g$ does not change its sign.

We are now going to specify $\vz$. Note that
$$\align
(-\ooz {\tilde g})(k_0-1)
&=-b_{k_0-1} (\vz - g_{k_0-1})+ a_{k_0-1}(g_{k_0-1}-g_{k_0-2})+c_{k_0-1}g_{k_0-1}\\
&=  (-\ooz g)(k_0-1)+ b_{k_0-1}(g_{k_0}-\vz)\\
&=\lz_0 g_{k_0-1}+ b_{k_0-1}(g_{k_0}-\vz)\\
&<\lz_0 g_{k_0-1}
\endalign$$
since $\vz>0\ge g_{k_0}$. Note also that
$$(-\ooz {\tilde g})(k_0)=-b_{k_0} (0-\vz)+ a_{k_0}(\vz-g_{k_0-1})+ c_{k_0}\vz
=\vz (a_{k_0}+b_{k_0}+ c_{k_0})-a_{k_0}g_{k_0-1}.$$
Next, since $D(f)=(f, -\ooz f)$ for every $f\in {\scr K}$ and for each $i$,
$\ooz f(i)$ depends on three points $i$ and $i\pm 1$ only,
we obtain
$$\align
D\big(\tilde g\big)&=\sum_{0\le i\le k_0-2}\mu_i g_i (-\ooz g)(i)
+\mu_{k_0-1}{g}_{k_0-1}(-\ooz {\tilde g})(k_0-1)
+\mu_{k_0}{\tilde g}_{k_0}(-\ooz {\tilde g})(k_0)\\
&<\lz_0 \sum_{i=0}^{k_0-1}\mu_i g_i^2
+\vz\mu_{k_0} [\vz (a_{k_0}+b_{k_0}+ c_{k_0})-a_{k_0}g_{k_0-1}].\endalign$$
Because
$$\|\tilde g\|^2= \sum_{i=0}^{k_0-1}\mu_i g_i^2+\mu_{k_0} \vz^2,$$
for $D(\tilde g)/\|\tilde g\|^2<\lz_0$, it suffices that
$$\vz [\vz (a_{k_0}+b_{k_0}+ c_{k_0})-a_{k_0}g_{k_0-1}]<\lz_0 \vz^2.$$
Equivalently,
$$\vz(a_{k_0}+b_{k_0}+ c_{k_0}-\lz_0)< a_{k_0}g_{k_0-1}.$$
This clearly holds for sufficiently small $\vz>0$.

(c) If $\lz=0$, then (2.5) becomes
$$\mu_k b_k(g_{k+1}-g_k)=\sum_{i=0}^k c_i \mu_i  g_i,
\qqd k\in E,\; g_{N+1}=0\text{ if }N<\infty. \tag 2.10$$
Clearly, if $g_0=0$, then $g_i\equiv 0$ by induction. Without loss of generality,
assume that $g_0=1$. By (2.10) and induction, it follows that $g_{k+1}-g_k\ge 0$ for all
$i\in E$. Actually, $g_{k+1}>g_k$ for all $k$: $i\le k<N$ provided $c_i>0$.
\qed\enddemo

In view of (2.5), the eigenfunction $g$ may not be monotone if $c_i\not\equiv 0$.

{\it For the remainder of this section, we assume that $c_i=0$ for $i<N$
but $c_N>0$ if $N<\infty$. However, to simplify our notation, set
$c_i\equiv 0$ but let $b_N>0$  if $N<\infty$. In view of the
definition of the state space $E$, the point $N+1$ is regarded as a
Dirichlet boundary. From now on in the paper, when we talk about $\lz_0^{(2.2)}$, it is defined
by $(2.2)$ but in the present setting.}

\proclaim{\prp\;2.2}
\roster
\item Let $g$ be a non-zero eigenfunction of $\lz_0>0$. Then $g$ is either positive or negative.
\item Let $g$ be a positive eigenfunction of $\lz >0$. Then $g$ is strictly decreasing.
Furthermore,
$$\sum_{k=n}^N \frac{1}{\mu_k b_k} \sum_{i=0}^k \mu_i g_i=\frac{g_n-g_{N+1}}{\lz},
\qqd n\in E.\tag 2.11$$
In particular,
$$\sum_{n=0}^N \mu_n g_n \, \nu[n, N]
=\sum_{n=0}^N \nu_n \sum_{k=0}^n\mu_k g_k
=\frac{g_0-g_{N+1}}{\lz}<\infty,\tag 2.12$$
where $\nu [\ell, m]=\sum_{k=\ell}^m \nu_k$, $\nu_k=(\mu_kb_k)^{-1}$.
Moreover, if (1.2) holds, then $g_\infty:=\lim_{N\to\infty}g_N=0$.
\item Let $\lz_0=0$. Then $N=\infty$ and the eigenfunction $g$ must be a
constant function.\endroster
\endproclaim

\demo{\prf}
(a) The first assertion follows from \prp\;2.1\,(2).

(b) Let $\lz >0$. Since $g>0$, by (2.5) with $c_i\equiv 0$, it follows that
$g_i$ is strictly decreasing in $i$.
By (2.5) again, we have
$$g_n-g_{N+1}=\sum_{k=n}^N (g_k-g_{k+1})=
 \lz \sum_{k=n}^N \frac{1}{\mu_k b_k}\sum_{i=0}^k \mu_i g_i
=\lz \sum_{i=0}^N \mu_i g_i\, \nu[i\vee n, N].$$
We obtain formula (2.11) and then (2.12).
If $g_\infty>0$, then by condition (1.2), the left-hand side of (2.11) is bounded below by
$$g_\infty \sum_{k=n}^\infty \frac{1}{\mu_k b_k} \sum_{i=0}^k \mu_i=\infty  \tag 2.13$$
which is a contradiction since the right-hand side of (2.11) is bounded from the above
by $g_0/\lz<\infty$. Therefore, we must have $g_\infty=0$.

With some additional work, condition (1.2) for $g_\infty=0$ will be removed (see
\prp\;2.5 below).

(c) We now prove the last assertion of the proposition. When $N<\infty$, it is
well known that $\lz_0>0$. Now, let $\lz_0=0$ and then $N=\infty$. By (2.6) with $c_i\equiv 0$,
we have
$$g_{i+1}-g_i=\frac{a_i}{b_i}(g_i-g_{i-1}),\qqd i\ge 0.$$
From this and induction, it follows that $g_n=g_0$ for all $n\ge 1$ since $a_0=0$.
\qed\enddemo

We remark that for finite state space with absorbing at
$N+1<\infty$, Proposition 2.2 was actually proved in  [4;
proof d) of Theorem 3.4] with a change of the order of the state
space. Next, when $N=\infty$ and $\lz_0>0$, in contrast with the ergodic case
where $g\in L^1(\mu)$ (cf. [12; \prp\;3.5]), here one may
have $g\notin L^2(\mu)$ and then $g\notin L^1(\mu)$. However, $g\in
L^1(\nu)$ since $g_n$ is strictly decreasing and $\sum_n
\nu_n<\infty$, which is a consequence of \thm\;3.1 below.

\proclaim{\crl\;2.3} Let $\lz_0>0$. Then $\lim_{i\to\infty}P_t f(i)=0$ for all $t\ge 0$ and $f\in\kq$.
\endproclaim

\demo{\prf} It suffices to show that $\lim_{i\to\infty}\sum_{k=1}^n p_{ik}(t)=0$. We now prove a stronger
conclusion: $\lim_{i\to\infty}P_t g(i)=0$ for all $t\ge 0$, where $g>0$ with $g_0=1$ is the eigenfunction
of $\lz_0$. Since $g$ is bounded, by using the well-known fact that
$$e^{-\lz_0 t}g_i= P_t g(i),\qqd t\ge 0,$$
the conclusion now follows from Propositions 2.2 and 2.5\,(2) below.\qed
\enddemo

For a specialist who does not want to know many details, at the first reading,
one may have a glance at the remainder of this section and the next section, especially \prp\;2.7, and then go to
Section 4 directly. From here to the end of the next section, we are dealing
with a case which is a dual of the one studied in Section 4. However, for the reader
who is unfamiliar with this topic, it is better just to follow the context since
we present everything in detail in these two sections.
A large part of the details in
Sections 4 and 6 are omitted since they are supposed to be known.

To state the main results of this section, we need some notation.
First, we define two operators as follows.
$$I_i(f)= \frac{1}{\mu_i b_i (f_i-f_{i+1})} \sum_{j\le i} \mu_j f_j,
\quad\quad I\!I_i(f)= \frac{1}{f_i} \sum_{j= i}^N
 \frac{1}{\mu_j b_j}\sum_{k\le j} \mu_k f_k. \tag 2.14$$
They are called an operator of single sum (integral)
or double sum, respectively. Here for the first operator, we use a convention: $f_{N+1}=0$ if $N<\infty$.
The second operator can be alternatively expressed as
 $$I\!I_i(f)=\frac{1}{f_i} \sum_{k\in E} \mu_k f_k\, \nu[i\vee k, N],
 \qqd \nu[\ell, m]=\sum_{i=\ell}^m \nu_i,\qd \nu_i=\frac{1}{\mu_i b_i}. \tag 2.15$$
 Next, define a difference operator $R$ as follows.
 $$\!R_i(v)\!=\!a_{i}\big(1\!-v_{i-1}^{-1}\big)\!+\! b_i(1\! -\! v_{i}),\;\; i\!\in\! E,\;v_{-1}\!>\!0\text{ is free},\;
   v_N\!:=0\text{ if }N\!<\!\infty.\tag 2.16$$
The domain of the operators $I\!I$, $I$ and $R$ are defined, respectively, as follows.
$$\align
&{\scr F}_{I\!I}=\{f\!: f>0\text{ on }E\},\\
&{\scr F}_{I}=\{f\!: f>0 \text{ on }E \text{ and is strictly decreasing}\},\\
&{\scr V}_1=\!\big\{v\!: \text{for all $i\,(0\le i\!<\!N)$},\, v_i \!\in\! (0, 1) \text{ if }{\tsize\sum_j}\nu_j\!<\!\infty
\text{ and }v_i \!\in\! (0, 1] \text{ if }{\tsize\sum_j}\nu_j\!=\!\infty\big\}\!.
\endalign$$
These sets are used for the lower estimates. For the upper estimates, we need some
modifications of them as follows.
$$\align
{\widetilde{\scr F}}_{I\!I}&=\big\{f: f>0 \text{ up to some }m: 1\le m<N+1\text{ and then vanishes}\big\},\\
{\widetilde{\scr F}}_I&=\big\{f: f \text{ is strictly decreasing on some interval } [n, m]\,(0\le n<m<N+1),\\
&\text{\hskip3.3em} f_i=f_n \text{ for }i\le n,\; f_m>0,\text{ and } f_i=0 \text{ for }i> m\big\},\\
{\widetilde{\scr V}}_1&=\cup_{m=1}^{N-1}\big\{v:
a_{i+1} (a_{i+1}+b_{i+1})^{-1}\!<\!v_i\!<\!
1-a_i\big(v_{i-1}^{-1}-1\big) b_i^{-1}\\
&\text{\hskip5em}\;\text{ for } i=0,1, \ldots, m-1 \text{ and }
v_i=0 \text{ for }i\ge m\big\}.
\endalign
$$
Here and in what follows, to use the above operators on these modified sets, we adopt
the usual convention $1/0=\infty$. Besides, the operator $I\!I$ should be generalized as follows:
$$I\!I_i (f)= \frac{1}{f_i}\, \sum_{i\le j\, \in\supp(f)}
 \frac{1}{\mu_j b_j}\sum_{k\le j} \mu_k f_k, \qqd i\in \supp(f).\tag 2.17$$
From now on, we should remember that $I\!I_{\bullet} (f)$ is defined on $\supp(f)$ only.
Fortunately, we need only to consider the following two cases: either $\supp(f)=\{0,1,\ldots, m\}$
for a finite $m$ or $\supp(f)=E$.

To avoid the heavy notation, we now split our main result of this section into a theorem and a
proposition below.

\proclaim{\thm\;2.4} The following variational formulas hold for $\lz_0$ defined by (2.2).
\roster
\item Difference form:
$$\inf_{v\in {\widetilde{\scr V}}_1}\,\sup_{i\in E}\,R_i(v)
=\lz_0
=\sup_{v\in {\scr V}_1}\,\inf_{i\in E}\, R_i(v).$$
\item Single summation form:
  $$\inf_{f\in {\widetilde{\scr F}}_I}\,\sup_{i\in E}I_i(f)^{-1}
    =\lz_0=\sup_{f\in {\scr F}_I}\,\inf_{i\in E}I_i(f)^{-1}.$$
\item Double summation form:
  $$\inf_{f\in {\widetilde{\scr F}}_{I\!I}}\,\sup_{i\in\supp(f)}I\!I_i (f)^{-1}
    =\lz_0=\sup_{f\in {\scr F}_{I\!I}}\,\inf_{i\in E}I\!I_i(f)^{-1}.$$
 \endroster
Moreover, the supremum on the right-hand side of the above three formulas
can be attained.
\endproclaim

The next result extends the domain of $\lz_0$ or adds some additional sets of test functions
for the operators $I$ and $I\!I$, respectively. Roughly speaking, a larger
set of test functions provides more freedom in practice and a smaller one is helpful for
producing a better estimate.

\proclaim{\prp\;2.5} \roster
\item We have
$$\lz_0=\inf\{D(f): \|f\|=1,\; f_{N+1}=0\},\tag 2.18$$
where $f_\infty:=\lim_{N\to\infty}f_N$ in the case of $N=\infty$.
\item When $\lz_0>0$, the eigenfunction $g$ satisfies $g_{N+1}=0$.
\item Moreover, we have
$$\align
\lz_0&=\inf_{f\in {\widetilde{\scr F}}_I'}\sup_{i\in E} I_i(f)^{-1} \tag2.19\\
&=\inf_{f\in {\widetilde{\scr F}}_{I\!I}\cup {\widetilde{\scr F}}_{I\!I}'}\,\sup_{i\in\supp(f)}I\!I_i (f)^{-1},\tag 2.20\\
\inf_{f\in {\widetilde{\scr F}}_I}\,\sup_{i\in\supp(f)}\,&I\!I_i
(f)^{-1}
    =\lz_0=\sup_{f\in {\scr F}_{I}}\,\inf_{i\in E}I\!I_i(f)^{-1},\tag 2.21
\endalign$$
where
$$\align
{\widetilde{\scr F}}_I'
&=\big\{f: f \text{ is strictly decreasing and positive up to some } m: 1\le m<N+1 \\
 &\qqd\qd\text{ and then vanishes}\big\}\subset {\widetilde{\scr F}}_I,\\
 {\widetilde{\scr F}}_{I\!I}'&=\big\{f: f>0 \text{ on $E$ and }fI\!I(f)\in L^2(\mu)\big\}.
 \endalign$$
Besides, the supremum $\sup_{f\in {\scr F}_{I}}$ in (2.21) can also be attained.
\endroster\endproclaim

The condition ``$f_{N+1}=0$'' in (2.18) explains the meaning of
``absorbing (Dirichlet) boundary at infinity'' used in the title of
this and the next sections.

Among the different groups of variational forms, the difference form is the simplest
one in the practical computations. For instance, when $N=\infty$, by choosing
$v_i\equiv c <1$, we obtain the following simple
 lower estimate:
$$\lz_0\ge \inf_{i\in E}\,[b_i(1- c)-a_i(c^{-1}-1)].$$
This is non-trivial and is indeed sharp for a linear model (\xmp\;3.5,
$c=1/2$). The difference form of the variational formulas will be
used in Section 5 to deduce a dual representation of $\lz_0$. In
general, the estimates produced by the operator $R$ can be improved
by using the operator $I$ and further improved by using $I\!I$. The
price is that more computation is required successively. The single
summation form of the variational formulas enables us to deduce a
criterion for $\lz_0>0$ (\thm\;3.1). Whereas the double summation
form of the variational formulas enables us to deduce an
approximating procedure to improve step by step the lower and upper
estimates of $\lz_0$ (\thm\;3.2).

Next, we mention that when $N=\infty$, for the upper estimates (the
left-hand side of the formulas given in Theorem 2.4 or the formula
given in (2.20)), the truncating procedure or the condition
``$fI\!I(f)\in L^2(\mu)$'' cannot be removed. For instance, the
formally dual formula $\inf_{0<v\le 1}\, \sup_{i\in E}R_i(v)$
of the lower estimate $\sup_{0<v\le 1}\, \inf_{i\in E}R_i(v)\,[=\sup_{v\in {\scr V}_1}\, \inf_{i\in E}R_i(v)]$
is not an upper bound of $\lz_0$, and is indeed trivial. To see
this, simply take $\bar v_i\equiv 1\,(i<\infty)$. Then $R_i(\bar
v)\equiv 0$ and so
$$\inf_{0<v\le 1}\, \sup_{i\in E}R_i(v)\le \sup_{i\in E}R_i(\bar v)=0.$$
More concretely, take $b_i\equiv 2$ and $a_i\equiv 1$. Then for ${\bar v}_i\equiv c<1$, we have
$$\inf_{v\in {\scr V}_1}\, \sup_{i\in E}R_i(v)\le \inf_{c<1}\, \sup_{i\in E} R_i(\bar v)= 2\,\inf_{c<1} (1-c)=0,$$
but $\lz_0=\big(\sqrt{2}-1\big)^2$ as will be seen in the next section (\xmp\;3.4).
Therefore, the quantity $\inf_{0<v\le 1}\, \sup_{i\in E}R_i(v)$,
as well as $\inf_{v>0}\, \sup_{i\in E}R_i(v)$, has no use for an
upper estimate of $\lz_0$.

\demo{\prf s of \thm\;$2.4$ and \prp\;$2.5$}

\flsh{\bf Part I}. Recall that $\lz_0^{(\#)}$ denotes the one
defined by the formula $(\#)$. In particular, {the notation $\lz_0$
used from now on in this section is $\lz_0^{(2.2)}$}.

To prove the lower estimates, we adopt the following circle argument:
$$\align
\lz_0&\ge\lz_0^{(2.18)}\ge \sup_{f\in {\scr F}_{I\!I}}\inf_{i\in E}I\!I_i(f)^{-1}
=\sup_{f\in {\scr F}_{I}}\inf_{i\in E}I\!I_i(f)^{-1}=\sup_{f\in {\scr F}_I}\inf_{i\in E}I_i(f)^{-1}\\
&\ge \sup_{v\in {\scr V}_1}\inf_{i\in E}\, R_i(v)
\ge \lz_0.  \tag 2.22\endalign  $$
Clearly, $\lz_0^{(2.18)}=\lz_0$ if $N<\infty$.
However, the identity is not trivial in the case of
$N=\infty$. Besides, we will show that each supremum in (2.22) can be attained; and furthermore
the eigenfunction $g$ satisfies $g_{N+1}=0$ whenever $\lz_0>0$.
\medskip

(a) Prove that $\lz_0\ge\lz_0^{(2.18)}\ge \sup_{f\in {\scr F}_{I\!I}}\,\inf_{i\in E}I\!I_i(f)^{-1}$.
\medskip

When $N=\infty$, the first inequality is trivial since
$$\{\|f\|=1, f\in {\scr K}\}\subset \{\|f\|=1,\; f_\infty=0\}.$$

The proof of the second inequality is parallel to the first part of
the proof of [4; \thm\;2.1]. Let $g$ satisfy $g_{N+1}=0$ and
$\|g\|=1$, and let $(h_i)$ be a positive sequence. Then by a good use
of the Cauchy-Schwarz inequality, we obtain
$$ \align
1 &= \sum_i \mu_i g_i^2\qd\text{(since $\|g\|=1$)}\\
&= \sum_i \mu_i \bigg( \sum_{j=i}^N(g_j- g_{j+1})\bigg) ^2\qd\text{(since $g_{N+1}=0$)}\\
&\le  \sum_i \mu_i  \sum_{j=i}^N \frac{(g_{j+1}-g_j ) ^2 \mu_j b_j}
{h_j}\sum_{k=i}^N \frac{h_k}{ \mu_k b_k}.\endalign$$
Exchanging the order of the first two sums on the right-hand side, we get
$$\align
1 &\le \sum_j \mu_j b_j (g_{j+1}-g_j ) ^2 \frac{1}{h_j}
\sum_{i \le  j}\mu_i \sum_{k=i}^N \frac{h_k}{ \mu_k b_k}\\
&\le  D(g)\, \sup_{j\in E} \frac{1}{h_j}
\sum_{i \le  j}\mu_i \sum_{k=i}^N \frac{h_k}{ \mu_k b_k}\\
&=: D(g)\, \sup_{j\in E} H_j.
 \endalign $$
We mention that the right-hand side may be infinite but we do not care at the moment.
Now, let $f\in  {\scr F}_{I\!I} $ satisfy $c:=\sup_{j\in E} I\!I_j(f) < \infty
$ and take $h_j = \sum_{i \le  j}\mu_i f_i$. Then $h_j\le c f_j/v_j<\infty$
for all $j$. By the proportional property,
we have
$$ \sup_{j\in E} H_j\le   \sup_{j\in E}
 \frac{1}{f_j}\sum_{k=j}^N \frac{h_k}{\mu_k b_k}
= \sup_{j\in E} \frac{1}{f_{j}}
\sum_{k=j}^N \frac{1}{\mu_k b_k} \sum_{i \le k}\mu_i f_i
=\sup_{j\in E} I\!I_j(f)<\infty. $$
Combining these facts together, we obtain
$\lz_0^{(2.18)}\ge\inf_{j\ge 0} I\!I_j(f)^{-1}$
whenever $\sup_{j\in E} I\!I_j(f) < \infty$. The inequality is
trivial if $\sup_{j\in E} I\!I_j(f) = \infty$ and so it holds for all $f\in  {\scr F}_{I\!I} $.
By making the supremum with respect to $f\in  {\scr F}_{I\!I} $,
we obtain the required assertion.
\medskip

(b) Prove that $\sup\limits_{f\in {\scr F}_{I\!I}}\,\inf\limits_{i\in E}I\!I_i(f)^{-1}
=\sup\limits_{f\in {\scr F}_{I}}\,\inf\limits_{i\in E}I\!I_i(f)^{-1}
= \sup\limits_{f\in {\scr F}_I}\,\inf\limits_{i\in E}I_i(f)^{-1}$.
\medskip

Let $f\in {\scr F}_I\subset {\scr F}_{I\!I}$. Without loss of generality, assume that $\sup_{i\in E}I_i(f)<\infty$.
By using the proportional property, we obtain
$$\align
\sup_{i\in E} I\!I_i (f)&=\sup_{i\in E}\frac{1}{f_i}\sum_{j= i}^N\frac{1}{\mu_j b_j}\sum_{k\le j}\mu_k f_k\\
&\le\sup_{i\in E}\sum_{j= i}^N\frac{1}{\mu_j b_j}\sum_{k\le j}\mu_k f_k
\bigg/\sum_{j= i}^N(f_j-f_{j+1})\qd\text{(since $f_{N+1}\ge 0$)}\\
&\le \sup_{i\in E}\frac{1}{f_i-f_{i+1}}\bigg(\frac{1}{\mu_i b_i}\sum_{k\le i}\mu_k f_k\bigg)\qd\text{(note that $f_i>f_{i+1}$)}\\
&= \sup_{i\in E} I_i (f)<\infty. \tag 2.23
\endalign$$
Making the infimum with respect to $f\in {\scr F}_I$, we get
$$\inf_{f\in {\scr F}_I}\sup_{i\in E} I\!I_i (f)\le \inf_{f\in {\scr F}_I}\sup_{i\in E} I_i (f).$$
Since
${\scr F}_I\subset {\scr F}_{I\!I}$, the left-hand side is bounded below
by $\inf_{f\in {\scr F}_{I\!I}}\sup_{i\in E} I\!I_i (f)$.
We have thus proved that
$$\sup_{f\in {\scr F}_{I\!I}}\,\inf_{i\in E}I\!I_i(f)^{-1}
\ge \sup_{f\in {\scr F}_{I}}\,\inf_{i\in E}I\!I_i(f)^{-1}
\ge \sup_{f\in {\scr F}_I}\,\inf_{i\in E}I_i(f)^{-1}.$$

There are two ways to prove the inverse inequality. The first one is longer but contains
a useful technique.
Let $f\in {\scr F}_{I\!I}$ with $c:=\sup_{i\in E}I\!I_i(f)<\infty$. Set
$$g_i=\sum_{j= i}^N \frac{1}{\mu_j b_j} \sum_{k\le j}\mu_k f_k
=\sum_{j= i}^N \nu_j \sum_{k\le j}\mu_k f_k>0,\qqd i\in E,\; g_{N+1}:=0\text{ if }N<\infty.$$
Then $g_i$ is strictly decreasing in $i$,
$g_i<g_0\le c f_0<\infty$ for all $i$.
Hence, $g\in {\scr F}_I$. Noticing that
$$g_i-g_{i+1}=\sum_{j= i}^N \nu_j \sum_{k\le j}\mu_k f_k-
\sum_{j= i+1}^N \nu_j \sum_{k\le j}\mu_k f_k=\nu_i \sum_{k\le i}\mu_k f_k$$
(here and in what follows, $\sum_{k=i}^j$ means $\sum_{i\le k<j+1}$ and $\sum_{\emptyset}=0$
by the standard convention), we have
$$\align
\ooz g(i)&=b_i(g_{i+1}-g_i)+a_i (g_{i-1}-g_i)\\
&=-b_i \nu_i \sum_{k\le i}\mu_k f_k + a_i \nu_{i-1} \sum_{k\le i-1}\mu_k f_k\\
&=-\frac{1}{\mu_i} \sum_{k\le i}\mu_k f_k + \frac{a_i}{\mu_{i-1} b_{i-1}} \sum_{k\le i-1}\mu_k f_k\\
&=-f_i,\qqd 1\le i<N.
\endalign$$
Actually, this holds also for $i=0$ and $i=N$ if $N<\infty$. Applying (2.7) to $h=f$, by (2.9), it follows that
$$ \mu_k b_k (g_k -g_{k+1})=\sum_{j\le k}\mu_j g_j f_j/g_j
\ge \sum_{j\le k}\mu_j g_j \inf_{i\in E} I\!I_i (f)^{-1},\qqd k\in E.$$
That is,
$$\sup_{i\in E} I\!I_i (f)\ge \frac{1}{\mu_k b_k (g_k -g_{k+1})}\sum_{j\le k}\mu_j g_j
=I_k(g),\qqd k\in E. $$
Making the supremum with respect to $k$, we obtain
$$\inf_{k\in E} I_k (g)^{-1}\ge \inf_{i\in E} I\!I_i (f)^{-1}, $$
and hence,
$$\sup_{g\in {\scr F}_I}\inf_{k\in E} I_k (g)^{-1}\ge \inf_{i\in E} I\!I_i (f)^{-1}.$$
This lower bound becomes trivial if $\sup_{i\in E}I\!I_i(f)=\infty$, and hence, the inequality holds
for all $f\in {\scr F}_{I\!I}$.
Making the supremum with respect to $f\in {\scr F}_{I\!I}$, we obtain
$$\sup_{g\in {\scr F}_I}\inf_{k\in E} I_k (g)^{-1}\ge
\sup_{f\in {\scr F}_{I\!I}}\inf_{i\in E} I\!I_i (f)^{-1}.  $$
We have thus proved the required assertion.

The second proof is to show that
$$\sup_{f\in {\scr F}_I}\inf_{i\in E} I_i(f)^{-1}\ge
\lz_0$$
and thus completes a smaller circle argument. To do so,
without loss of generality, assume that $\lz_0>0$. Let $g>0$ be the
eigenfunction of $\lz_0$. Applying (2.9) to $h=\lz_0 g$, we obtain
$I_i(g)=\lz_0^{-1}$ for all $i\in E$, and hence, $\inf_{i\in E}
I_i(g)^{-1}=\lz_0$. Noticing that $g\in {\scr F}_I$ by \prp\;2.2, the
assertion is now obvious.
\medskip

(c) Prove that
$\sup_{f\in {\scr F}_{I\!I}}\inf_{i\in E} I\!I_i(f)^{-1}\ge \sup_{v\in {\scr V}_1}\inf_{i\in E} R_i (v)$.
\medskip

Note that by a change
of the sequence $\{v_i\}_{i=0}^{N-1}$:
$$u_i=v_0 v_1 \cdots v_{i-1},\qqd i\in E,\; v_{-1}>0\text{ is free},\; v_N:=0\text{ if }N<\infty,$$
the quantity $R_i(v)$ becomes
$$a_{i}\bigg(1-\frac{u_{i-1}}{u_i}\bigg)+b_{i}\bigg(1-\frac{u_{i+1}}{u_i}\bigg), \qqd i\in E,\; u_{-1}>0\text{ is free},\;
u_{N+1}:=0\text{ if }N<\infty.$$
To save our notation, we use $R_i(u)$ to denote this quantity. Clearly, $\{u_i\}$
is positive and $v_i\le 1$ for all $i$ mean that $\{u_i\}$ is non-increasing.

Before moving further, we prove that if $\inf_{i\in E} R_i(u)>0$ for a positive
sequence $u=(u_i)$, then $u_i$ must be strictly decreasing in $i$.
To do so, let
$$f_i=(a_{i}+b_i) u_i-a_i {u_{i-1}} -b_{i} {u_{i+1}}.$$
Then $f_i=u_i R_i(u)>0$ for all
$i\in E$ by assumption, and so $f\in {\scr F}_{I\!I}$. Noticing that
$$\mu_k f_k=\mu_{k+1} a_{k+1}(u_k-u_{k+1})-\mu_{k} a_{k}(u_{k-1}-u_{k}),$$
we obtain
$$0<\sum_{k\le j}\mu_k f_k=\mu_{j+1} a_{j+1}(u_{j}-u_{j+1})=\mu_{j} b_{j}(u_{j}-u_{j+1}).$$
Hence, $u_i$ is strictly decreasing in $i$ (equivalently, $v_i:=u_{i+1}/u_i<1$).
This proves the required assertion.
The reason of using ${\scr V}_1$ rather than $\{v: v_i>0, 0\le i<N\}$ should be clear now.

We now return to our main assertion. For this, without loss of generality, assume that
$\inf_{i\in E} R_i(u)>0$ for a given strictly decreasing $u=(u_i)$. Otherwise, the assertion is trivial.
From the last formula, we obtain
$$0<\sum_{j= i}^N \nu_j \sum_{k\le j}\mu_k f_k
=\sum_{j= i}^N(u_j-u_{j+1})=u_i-u_{N+1}\le u_i.$$
Therefore,
$$0< R_i(u)=\frac{f_i}{u_i}\le {f_i}\bigg({\sum_{j= i}^N \nu_j \sum_{k\le j}\mu_k f_k}\bigg)^{-1}
=I\!I_i (f)^{-1}, \qqd i\in E.$$
It follows that
$$\inf_{i\in E} R_i(u)\le \inf_{i\in E}I\!I_i (f)^{-1}
\le  \sup_{f\in {\scr F}_{I\!I}}\inf_{i\in E} I\!I_i(f)^{-1}.$$
The assertion now follows by making the supremum with respect to $u$.
\medskip

(d) Prove that $\sup_{v\in {\scr V}_1}\,\inf_{i\in E}\, R_i(v)\ge \lz_0$.
\medskip

Assume that $\lz_0>0$ for a moment (in particular, if $N<\infty$). Then by \prp\;2.2, the
corresponding eigenfunction $g$ (with $g_0=1$) of $\lz_0$ is positive and
strictly decreasing. From the eigenequation
$$-\ooz g(i)=\lz_0 g_i,\qqd i\in E,\;g_{N+1}:=0\text{ if }N<\infty,$$
it follows that
$$a_{i}\bigg(1-\frac{g_{i-1}}{g_i}\bigg)+b_{i}\bigg(1-\frac{g_{i+1}}{g_i}\bigg)={\lz_0},\qqd i\in E.$$
Let $v_i=g_{i+1}/g_i$. Then $v_i\in (0, 1)$ for all $i<N$ and so $v=(v_i)\in {\scr V}_1$.
Moreover, $R_i(v)=\lz_0$ for all $i\in E$.
Therefore, we certainly have $\sup_{v\in {\scr V}_1}\,\inf_{i\in E}\, R_i(v)\ge \lz_0$,
as required.

It remains to prove that $\sup_{v\in {\scr V}_1}\,\inf_{i\in E}\, R_i(v)\ge 0$ when $N=\infty$.
First, let $\sum_{i=0}^\infty \nu_i<\infty$. Choose a positive $f$ such that
$$\sum_{k=0}^\infty \mu_k f_k \fz_k<\infty, \qqd \fz_k:=\sum_{j=k}^\infty \nu_j.$$
Define
$$h_i=\sum_{j=i}^\infty \nu_j \sum_{k\le j} \mu_k f_k, \qqd i\ge 0.$$
Then
$$h_i=\sum_{k=0}^\infty \mu_k f_k \fz_{i\vee k}\le \sum_{k=0}^\infty \mu_k f_k \fz_k<\infty.$$
Set $\bar\nu_i=h_{i+1}/h_i\,(i\ge 0)$. Then ${\bar\nu}\in {\scr V}_1$ since $h_i$ is strictly
decreasing.
A simple computation shows that $R_i(\bar v)=I\!I_i(f)^{-1}>0$ for all $i\ge 0$. Hence,
$\sup_{v\in {\scr V}_1}\,\inf_{i\ge 0}\, R_i(v)\ge 0$. Next, let $\sum_i \nu_i=\infty$
and set $\bar v_i\equiv 1$. Then $R_i(\bar v)\equiv 0$ and so the same conclusion holds.

The proof of the last paragraph indicates the reason why in ${\scr V}_1$ we define ``$v_i\in (0, 1)$''
and ``$v_i \in (0, 1]$'' separately according to ``$\sum_i \nu_i<\infty$'' or ``$\sum_i \nu_i=\infty$''.
Although we have known from proof (c) that for $\inf_i R_i (v)> 0$, it is necessary that $v<1$ but
this condition may not be sufficient for $\inf_i R_i (v)\ge 0$. The extremal ${\bar v}_i\equiv 1$ is
used only in the case of $\sum_i \nu_i=\infty$ in which
we indeed have $\lz_0=0$ (cf. \thm\;3.1 below).

We have thus completed the proof of circle (2.22).

(e) We now prove that each supremum in (2.22) can be attained.
The case that $\lz_0=0$ is easier since
$$0=\lz_0\ge \inf_{i\in E}I\!I_i(f)^{-1}\ge 0
\qd\text{and} \qd 0=\lz_0\ge \inf_{i\in E}I_i(f)^{-1}\ge 0$$
for every $f$ in the corresponding domain, as an application of (2.22). Similarly,
the conclusion holds for the operator $R$  as seen from proof (d): noting that in the degenerated case
that $\sum_i \nu_i=\infty$, we have $\lz_0=0$ and then $v_i\equiv 1$ by \prp\;2.2\,(3).

Next, we consider the case that $\lz_0>0$ with eigenfunction $g$: $g_0=1$. Then
for the operator $R$, the supremum is attained at $v_i=g_{i+1}/g_i$ as seen from
the first paragraph of proof (d). For the operator $I$, it is attained at $f=g$ as an application of
\prp\;2.1 with $c_i\equiv 0$: $I_i(g)\equiv \lz_0^{-1}$. At the same time,
in view of part (2) of \prp\;2.2, we have
$I\!I_i(g)\equiv \lz_0^{-1}$ whenever $g_{N+1}=0$.

It remains to rule out the possibility that $g_{N+1}>0$. Otherwise, by part (2)
of \prp\;2.2 again, we have
$N=\infty$ and
$$M_i:= \sum_{j\ge i}\nu_j \sum_{k\le j}\mu_k\in (0, \infty).$$
Let $\tilde g=g-g_\infty$. Then $\tilde g\in {\scr F}_{I\!I}$.
Noting that
$$\align
\sum_{j\ge i}\nu_j \sum_{k\le j}\mu_k {\tilde g}_k
&= \sum_{j\ge i}\nu_j \sum_{k\le j}\mu_k {g}_k-g_\infty M_i\\
&=\frac{g_i-g_\infty}{\lz_0}-g_\infty M_i\qd\text{(by (2.11))},\endalign$$
we obtain
$$\sup_{i\ge 0} I\!I_i(\tilde g)
=\sup_{i\ge 0} \bigg[\frac{1}{\lz_0} - \frac{g_\infty M_i}{g_i-g_\infty}\bigg]
=\frac{1}{\lz_0}-g_\infty \inf_{i\ge 0}\frac{M_i}{g_i-g_\infty}.$$
By using the proportional property and (2.5), it follows that
$$ \inf_{i\ge 0}\frac{M_i}{g_i-g_\infty}\ge
 \inf_{i\ge 0}\frac{\nu_i\sum_{k\le i}\mu_k}{g_i-g_{i+1}}
=\frac{1}{\lz_0}.$$
Thus, we get
$$\sup_{i\ge 0} I\!I_i(\tilde g)\le \frac{1}{\lz_0} (1-g_{\infty})<\frac{1}{\lz_0}.$$
Hence, $\inf\limits_{i\ge 0} I\!I_i(\tilde g)^{-1}> \lz_0$, which is a contradiction to
proof (a):
$\lz_0 \ge \inf\limits_{i\in E} I\!I(\tilde g)^{-1}.$
We have thus proved that $g_\infty=0$ whenever $\lz_0>0$.
Note that this paragraph uses \prp\;2.2 and proof (a) only.
\medskip

\flsh{\bf Part I\!I}. Next, to prove the upper estimates, we adopt the following circle argument:
$$\align
\lz_0&\le \inf_{f\in {\widetilde{\scr F}}_{I\!I}\cup {\widetilde{\scr F}}_{I\!I}'}\,\sup_{i\in \supp(f)}I\!I_i (f)^{-1}\tag 2.24\\
&\le \inf_{f\in {\widetilde{\scr F}}_{I\!I}}\,\sup_{i\in \supp(f)}I\!I_i (f)^{-1}\tag 2.25\\
&=\inf_{f\in {\widetilde{\scr F}}_I}\,\sup_{i\in \supp(f)}I\!I_i (f)^{-1}= \inf_{f\in {\widetilde{\scr F}}_I}\,\sup_{i\in E}I_i(f)^{-1}\tag 2.26\\
&\le \inf_{f\in {\widetilde{\scr F}}_I'}\,\sup_{i\in E}I_i(f)^{-1}\tag 2.27\\
&\le \inf_{v\in {\widetilde{\scr V}}_1}\,\sup_{i\in E}\,R_i(v)\tag 2.28\\
&\le \lz_0.\tag 2.29
\endalign$$
Since inequalities (2.25) and (2.27) are obvious,
we need only to prove (2.24), (2.26), (2.28) and (2.29).
\medskip

(f) Prove that $\lz_0\le \inf_{f\in {\widetilde{\scr F}}_{I\!I}\cup
{\widetilde{\scr F}}_{I\!I}'}\,\sup_{i\in \supp(f)}I\!I_i (f)^{-1}$.
\medskip

We remark that  in the particular case that the eigenfunction $f$ is in $L^2(\mu)$, then the
function $g:=fI\!I(f)$ is nothing but just $f/\lz_0\in L^2(\mu)$.
Hence, the infimum in (2.24) is attained at this $f\in
{\widetilde{\scr F}}_{I\!I}'$ and the equality sign in (2.24) holds.

We now consider the general case. Let $f\in {\widetilde{\scr
F}}_{I\!I}$. Then there exists an $m$ such that $f_i>0$ for $i\le m$
and $f_i=0$ for $i>m$. Set $g=\dbl_{\supp(f)}f I\!I (f)$. That is,
$$g_i=
\cases
\sum_{j= i}^m \nu_j \sum_{k\le j}\mu_k f_k,\quad & i\le m\\
0, \quad & i\ge m+1.
\endcases$$
Clearly, $g\in L^2(\mu)$ and
$$g_i-g_{i+1}=\cases
\nu_i \sum_{k\le i}\mu_k f_k, \qd & i\le m\\
0, \qd &   i\ge m+1.
\endcases$$
We now have
$$D(g)=\sum_{i\le m} \mu_i b_i (g_{i+1}-g_i)^2
=\sum_{i\le m} (g_i-g_{i+1})\sum_{k\le i}\mu_k f_k
=\sum_{k\le m} \mu_k f_k \sum_{k\le i\le m} (g_i-g_{i+1}).$$
Since $g_{m+1}=0$, we get
$$D(g)=\sum_{k\le m} \mu_k f_k g_k
\le \sum_{k\le m} \mu_k g_k^2 \max_{0\le i\le m}(f_i/g_i)
= \|g\|^2 \sup_{i\in\supp(f)} I\!I_i (f)^{-1}.$$
Dividing both sides by $\|g\|^2\in (0, \infty)$, it follows that
$$\lz_0\le D(g)/\|g\|^2\le \sup_{i\in\supp(f)} I\!I_i(f)^{-1},\qqd f\in {\widetilde{\scr F}}_{I\!I}. \tag 2.30$$
For $f\in {\widetilde{\scr F}}_{I\!I}'$, the same conclusion clearly
holds if $N<\infty$. When $N=\infty$, since $g\in L^2(\mu)$ by assumption,
we have $0<g<\infty$. As a tail sequence of a convergent series (which sum equals $g_0$), we certainly have
$g_i\downarrow g_\infty=0$ as $i\uparrow\infty$. Hence, the same
proof replacing $m$ with $\infty$, plus the fact that
$\lz_0=\lz_0^{(2.18)}$ proved in Part I, shows that
$$\lz_0=\lz_0^{(2.18)}\le \sup_{i\in \supp(f)} I\!I_i (f)^{-1},\qqd f\in {\widetilde{\scr F}}_{I\!I}'. $$
Combining this with (2.30), we prove the required assertion.

The proof indicates the reason why the truncating procedure is used for the upper estimates
since in general the eigenfunction $g$ may not belong to $L^2(\mu)$ as shown by \prp\;2.2.
\medskip

(g) Prove that
$$\inf_{f\in {\widetilde{\scr F}}_{I\!I}}\,\sup_{i\in
\supp(f)}I\!I_i(f)^{-1} \!=\inf_{f\in {\widetilde{\scr
F}}_I}\,\sup_{i\in \supp(f)}I\!I_i(f)^{-1} \!=\inf_{f\in
{\widetilde{\scr F}}_I}\,\sup_{i\in E}I_i(f)^{-1}.$$
\medskip

Let $f\in {\widetilde{\scr F}}_I$. Then there exist $n<m$ such that
$f_i=f_{i\vee n}\dbl_{\{i\le m\}}$, $f_m>0$, and $f$ is strictly
decreasing on $[n, m]$. Clearly, we have
$$\min_{i\le m} I\!I_i
(f)=\min_{n\le i\le m} I\!I_i  (f)\qd\text{ and }\qd\inf_{i\in E} I_i (f)=
\min_{n\le i\le m} I_i (f)$$
since, by assumption, $1/0=\infty$. By the
proportional property, first we have
$$\align
\min_{n\le i\le m} I\!I_i  (f)
&=\min_{n\le i\le m} {\sum_{j= i}^m \nu_j\sum_{k\le j}\mu_k f_k}\bigg/{\sum_{j= i}^m (f_j-f_{j+1})}\\
&\ge \min_{n\le i\le m} \frac{1}{\mu_i b_i ( f_i-f_{i+1})}\sum_{k\le i}\mu_k f_k\\
&= \min_{n\le i\le m} I_i (f),\endalign$$
and then
$$\sup_{f\in {\widetilde{\scr F}}_{I\!I}}\,\inf_{i\in \supp(f)} I\!I_i  (f)
\ge\sup_{f\in {\widetilde{\scr F}}_I}\,\inf_{i\in \supp(f)} I\!I_i
(f) \ge \sup_{f\in {\widetilde{\scr F}}_I}\inf_{i\in E} I_i (f)$$
since ${\widetilde{\scr F}}_I\subset {\widetilde{\scr F}}_{I\!I}$.

As in proof (b), there are two ways to prove the inverse
inequality. First, let $f\in {\widetilde{\scr F}}_{I\!I}$. As in
proof (f), set $g=\dbl_{\supp(f)}f I\!I (f)$. Clearly,
$g\in {\widetilde{\scr F}}_I'\subset{\widetilde{\scr F}}_I$ and moreover,
$$\align
b_i(g_{i+1}-g_i)+ a_i(g_{i-1}-g_i)&=-\frac{1}{\mu_i}\sum_{k\le i}\mu_k f_k
 +\frac{a_i}{\mu_{i-1}b_{i-1}}\sum_{k\le i-1}\mu_k f_k\\
 &=-\frac{1}{\mu_i}\sum_{k\le i}\mu_k f_k+\frac{1}{\mu_i}\sum_{k\le i-1}\mu_k f_k\\
 &=-f_i, \qqd i\le m.
\endalign
$$
When $i=0$, the second term on the left-hand side disappears since $a_0=0$. It follows that
$$\mu_i b_i(g_{i+1}-g_i)+ \mu_i a_i(g_{i-1}-g_i)=-\mu_i f_i,\qqd i\le m,$$
and furthermore,
$$\mu_k b_k (g_k -g_{k+1})=\sum_{j\le k}\mu_j g_j f_j/g_j
\le \sum_{j\le k}\mu_j g_j \max_{0\le i\le m} I\!I_i  (f)^{-1},\qqd k\le m.$$
That is,
$$\min_{0\le i\le m} I\!I_i  (f)\le \frac{1}{\mu_k b_k (g_k -g_{k+1})}\sum_{j\le k}\mu_j g_j=I_k(g),\qqd k\le m. $$
Making the infimum with respect to $k$, we obtain
$$\max_{0\le k\le m} I_k (g)^{-1}\le \max_{0\le i\le m} I\!I_i  (f)^{-1}.$$
One may rewrite $\max_{0\le k\le m}$ as $\sup_{k\in E}$ on the
left-hand side since $I_k(g)=\infty$ for all $k\ge m+1$. Since $g\in
{\widetilde{\scr F}}_I'\subset{\widetilde{\scr F}}_I$, we now have
$$\inf_{g\in {\widetilde{\scr F}}_I}\,\sup_{k\in E} I_k (g)^{-1}
\le\inf_{g\in {\widetilde{\scr F}}_I'}\,\sup_{k\in E} I_k
(g)^{-1}\le
 \sup_{i\in \supp(f)} I\!I_i  (f)^{-1}.$$
Next,  making the infimum with respect to $f\in {\widetilde{\scr
F}}_{I\!I}$, we obtain
$$\inf_{g\in {\widetilde{\scr F}}_I}\,\sup_{k\in E} I_k (g)^{-1}
\le\inf_{g\in {\widetilde{\scr F}}_I'}\,\sup_{k\in E} I_k
(g)^{-1}\le
 \inf_{f\in {\widetilde{\scr F}}_{I\!I}}\,\sup_{i\in \supp(f)} I\!I_i  (f)^{-1}.$$

The second proof for the inverse inequality is to show that
$$\inf_{f\in {\widetilde{\scr F}}_I'}\,\sup_{i\in E} I_i(f)^{-1}\le
\lz_0.$$
For this, recall the definition
$$\lz_0=\inf\{D(f):\|f\|=1,\; f_i=0 \text{ for all } i> \text{some }m: 1\le m<N+1\}.$$
Because of
$$\align
&\{\|f\|=1, f_i=0 \text{ for all } i> m: 1\le m<N+1\}\\
&\qd \subset \{\|f\|=1,\; f_i=0 \text{ for all } i> m+1: 1\le m<N+1\},\endalign$$
it is clear that
$$\lz_0^{(m)}:=\inf\{D(f):\|f\|=1,\; f_i=0 \text{ for all } i> m: 1\le m<N+1\}\;\downarrow \lz_0\;
$$
as $ m\uparrow N $. Note that $\lz_0^{(m)}$ is just the first
eigenvalue of the Dirichlet form $(D, {\scr D}(D))$ restricted to
$\{0, 1, \ldots, m\}$ with Dirichlet (absorbing) boundary at $m+1$.
Now, let $g=g^{(m)}$ be the eigenfunction of $\lz_0^{(m)}>0$ with
$g_0=1$. Extend $g$ to the whole space by setting $g_i=0$ for all
$i>m$. By using \prp\;2.2, it follows that $g\in {\widetilde{\scr
F}}_I'$ with $\supp(g)=\{0, 1, \ldots, m\}$. Furthermore, by (2.9)
with $h=\lz_0 g$, we have $I_i(g)^{-1}=\lz_0^{(m)}>0$ for all $i\le
m$, and hence,
$$\sup_{i\in E} I_i(g)^{-1}=\sup_{i\le m} I_i(g)^{-1}=\lz_0^{(m)}.$$
Thus,
$$\lz_0^{(m)} = \sup_{i\in E} I_i(g)^{-1}
\ge \inf_{f\in {\widetilde{\scr F}}_I'\!,\; \supp(f)=\{0, 1, \ldots, m\}}\;\sup_{i\in E} I_i(f)^{-1}
\ge \inf_{f\in {\widetilde{\scr F}}_I'}\;\sup_{i\in E} I_i(f)^{-1}.$$
The assertion now follows by letting $m\to N$.
\medskip

(h) Prove that $\inf_{f\in {\widetilde{\scr F}}_{I\!I}}\,\sup_{i\in
\supp(f)} I\!I_i (f)^{-1}\le \inf_{v\in {\widetilde{\scr
V}}_1}\sup_{i\in E} R_i (v)$.
\medskip

Let $u$ with $\supp(u)=\{0,1,\ldots, m\}$ be given such that
$v_i:=u_{i+1}/u_i\in {\widetilde{\scr V}}_1$. Then, the constraint
$$v_i< 1- a_i\big(v_{i-1}^{-1}-1\big) b_i^{-1}, \qqd 0\le i\le m,\; v_m=0,$$
is equivalent to $\min_{0\le i\le m} R_i (v)>0$, and the constraint
$$v_i>a_{i+1} (a_{i+1}+b_{i+1})^{-1}, \qqd 0\le i\le m-1,$$
comes from the requirement that $v_i>0$ for all $i<m$. Since the case of
$i=m$ in the first constraint is contained in the second one,
we obtain the constraint described in ${\widetilde{\scr V}}_1$.
In particular, we have
$$a_1 (a_1+b_1)^{-1}<v_0<1-a_0\big(v_{-1}^{-1}-1\big)=1$$
and so $v_0\in (0, 1)$. By induction, we have $v_i\in (0, 1)$ for all $i<m$. The
existence of such a $u$ is guaranteed since $m<\infty$, as will be shown in
proof (i) below. Now, let
$$f_i=
\cases (a_i+ b_i) u_i -a_i u_{i-1}-b_i u_{i+1}, &\qd i\le m,\\
0, & \qd i> m.\endcases
$$
Then by assumption, $f_i/u_i=R_i(u)>0$ for $i\le m$. Hence, $f\in
{\widetilde{\scr F}}_{I\!I}$. Next, we have
$$0<\sum_{k\le j}\mu_k f_k=\mu_{j} b_{j}(u_{j}-u_{j+1}), \qqd j\le m.$$
Hence,
$$\sum_{j= i}^m \nu_j \sum_{k\le j}\mu_k f_k=u_i-u_{m+1}=u_i>0,\qqd i\le m. $$
Therefore, we obtain
$$R_i(u)=\frac{f_i}{u_i}= {f_i}\bigg/{\sum_{j= i}^m \nu_j \sum_{k\le j}\mu_k f_k}
=I\!I_i  (f)^{-1}, \qqd i\le m$$
and then
$$\sup_{i\in E} R_i(u) =\max_{i\le m}R_i(u)= \sup_{i\in \supp(f)} I\!I_i (f)^{-1}\ge
\inf_{f\in {\widetilde{\scr F}}_{I\!I}}\,\sup_{i\in \supp(f)} I\!I_i
(f)^{-1}.$$ To be consistent with the convention of $R_i(v)$, here we
adopt the convention: $R_i(u)=-\infty$ for all $i>m$. The assertion
now follows by making the infimum with respect to $u$.
\medskip

(i) Prove that $\inf_{v\in {\widetilde{\scr V}}_1}\,\sup_{i\in
E}\,R_i(v)\le \lz_0$.
\medskip

As in the last part of proof (g), denote by $g$ (with $g_0=1$) the eigenfunction of $\lz_0^{(m)}>0$.
Then $\supp(g)=\{0,1,\ldots, m\}$, and $g$ is strictly decreasing on $\{0, 1, \ldots, m\}$
by part (2) of \prp\;2.2.
The definition of $g$ gives us
$$b_i (g_i-g_{i+1}) - a_i (g_{i-1}-g_i)=\lz_0^{(m)} g_i,\qqd i\le m,\; g_{m+1}=0.$$
That is,
$$a_{i}\bigg(1-\frac{g_{i-1}}{g_i}\bigg)+b_{i}\bigg(1-\frac{g_{i+1}}{g_i}\bigg) =\lz_0^{(m)},\qqd i\le m.$$
Let $v_i=g_{i+1}/g_i$ for $i\le m$ and $v_i=0$ for $i> m$. Then
$v_i\in (0, 1)$ for $i\in \{0,1,\ldots, m-1\}$, and
$R_i(v)=\lz_0^{(m)}$ for all $i\le m$. It is now easy to see that
$v\in {\widetilde{\scr V}}_1$. We have thus constructed a $u\,(=g)$
required in proof (h). Clearly $R_i(v)=-\infty$ for all $i>m$.
Therefore,
$$\align
\lz_0^{(m)}&=\max_{0\le i\le m}\, R_i(v)\\
&\ge \inf_{v\in {\widetilde{\scr V}}_1:\;\supp(v)=\{0,1,\ldots, m-1\}}\;\max_{0\le i\le m}\, R_i(v)\\
&\ge \inf_{v\in {\widetilde{\scr V}}_1:\;\supp(u)=\{0,1,\ldots,n\}\text{ for some }n\ge 0}\;\sup_{i\in E}\, R_i(u)\\
&= \inf_{v\in {\widetilde{\scr V}}_1}\,\sup_{i\in E} R_i(v).
\endalign$$
Letting $m \to N$, we obtain the required assertion.

We have thus completed the circle argument of (2.24)--(2.29) and then the proofs of
\thm\;2.4 and \prp\;2.5 are finished.\qed\enddemo

Before moving further, we mention a technical point in the proof above.
Instead of the approximation with finite state space used in Part II of the above proof,
it seems more natural to use the truncating procedure for the eigenfunction $g$.
However, the next result shows that this procedure is not practical in general.

\proclaim{Remark 2.6} Let $g\ne 0$ be the eigenfunction of $\lz_0>0$ and define
$g^{(m)}=g \dbl_{\le m}$. Then
$$\min_{i\,\in\, \supp(g^{(m)})} I\!I_i\big(g^{(m)}\big)
=\frac{1}{\lz_0 }\bigg[1-\frac{g_{m+1}}{g_{m}}\bigg].$$
In particular, the sequence
$\big\{\min_{i\in \supp(g^{(m)})} I\!I_i\big(g^{(m)}\big)\big\}_{m\ge 1}$ may not converge to $\lz_0^{-1}$
as $m\uparrow \infty$.\endproclaim

\demo{\prf} Note that
$$\align
\min_{i\in \supp(g^{(m)})} I\!I_i\big(g^{(m)}\big)
&=\min_{0\le i\le m}\frac{1}{g_i}\sum_{j= i}^m \nu_j \sum_{k\le j}\mu_k g_k\\
&=\min_{0\le i\le m}\frac{1}{g_i}
\sum_{j= i}^m \nu_j \frac{\mu_j b_j (g_j-g_{j+1})}{\lz_0}\qd\text{(by (2.5))}\\
&=\min_{0\le i\le m}\frac{1}{\lz_0 g_i}(
g_i-g_{m+1})\\
&=\frac{1}{\lz_0 }\bigg[1-\frac{g_{m+1}}{g_{m}}\bigg].
\endalign$$
This proves the main assertion. For Example 3.4 in the next section, we have
$$\lim_{m\to\infty}\bigg(1-\frac{g_{m+1}}{g_m}\bigg)
=1-\sqrt{\frac{a}{b}}\,<1,$$
and so
$$\lim_{m\to\infty}\,\min_{i\,\in\, \supp(g^{(m)})} I\!I_i\big(g^{(m)}\big)< {\lz_0 }^{-1}.\qed$$
\enddemo

To conclude this section and also for later use, we introduce a variational formula of $\lz_0$
in a different difference form.

\proclaim{\prp\;2.7} On the set ${\scr V}:=\{v: v_i>0,\; 0\le i<N\}$, redefine
$$R_i(v)=a_{i+1}+ b_i-a_i/v_{i-1} - b_{i+1} v_{i},\qqd i\in E,\;\text{$v_{-1}>0$ is free},$$
where $a_{N+1}=b_{N+1}=0$ and $v_N$ is free if $N<\infty$. Then
\roster
\item we have
$$\sup_{v\in {\scr V}}\,\inf_{i\in E}\, R_i(v)\ge \lz_0.  \tag 2.31$$
The equality sign holds once $\sum_{i=0}^N\mu_i=\infty$. In this case,
we indeed have
$$\lz_0=\sup_{v\in {\scr V}}\,\inf_{i\in E}\, R_i(v)
=\sup_{v\in {\scr V}_*}\,\inf_{i\in E}\, R_i(v),
$$
where
$${\scr V}_*=\{v: v_{i-1}>a_i/b_i,\; 0\le i<N+1\}.$$
\item In general, we have
$$\lz_0=\sup_{v\in {\scr V}_*}\,\inf_{i\in E}\, R_i(v), \tag 2.32$$
and the supremum in (2.32) can be attained.
\endroster
\endproclaim

\demo{\prf}
(a) First, we prove that $\sup_{v\in {\scr V}_*}\,\inf_{i\in E}\, R_i(v)\ge 0$.
Given a positive, non-increasing $f$, $f_{N+1}=0$ if $N<\infty$,
define
$$u_i=(\mu_i b_i)^{-1}\sum_{j\le i}\mu_j f_j\in (0, \infty),\qqd
i<N+1.$$
Then
$$b_i u_i -a_i u_{i-1}=f_i>0, \qqd i\in E,\; u_{-1}>0\text{ is free}.$$
This implies that $(v_i:=u_{i+1}/u_i: i<N)\in {\scr V}_*$.
As before, we also use
$$R_i(u):=a_{i+1}+b_i-a_i\frac{u_{i-1}}{u_i} -b_{i+1}\frac{u_{i+1}}{u_i}, \qqd i\in E$$
instead of $R_i(v)$. Clearly,
$$R_i(u)=\frac{f_i-f_{i+1}}{u_i}\ge 0, \qqd i\in E.$$
Hence $\inf_{i\in E}\, R_i(u)\ge 0$ and the required assertion is now obvious.

(b) By (a), without loss of generality, assume that $\lz_0>0$. Then by \prp\;2.2,
the corresponding eigenfunction
$g$ of $\lz_0$ is positive and strictly decreasing.
With $u_i:=g_i-g_{i+1}>0\,(i\in E)$, the eigenequation
$$-\ooz g(i)=b_i u_i-a_i u_{i-1}=\lz_0 g_i,\qqd i\in E,\qd g_{N+1}:=0\text{ if }N<\infty,$$
gives us $(v_i:=u_{i+1}/u_i: i<N)\in {\scr V}_*$. Next, by
making a difference of $-\ooz g(i)$ and $-\ooz g(i+1)$ and noting that
$\ooz g(N+1)$ is setting to be zero if $N<\infty$, we obtain
$$(a_{i+1}+b_i) u_i -b_{i+1} u_{i+1}-a_i u_{i-1}=\lz_0 u_i,\qqd i\in E.$$
Thus, we have $R_i(v)=\lz_0$ for all $i\in E$.
Therefore, (2.31) holds.

(c) To prove the equality sign in (2.31) whenever $\sum_{i=0}^N \mu_i=\infty$,
in view of Part I of the proofs of Theorem 2.4 and \prp\;2.5 and (b), it
suffices to show that
$$\sup_{f\in {\scr F}_I}\,\inf_{i\in E}I_i(f)^{-1}
\ge \sup_{v\in {\scr V}}\,\inf_{i\in E}\, R_i(v).$$

In view of (a), without loss of generality, assume that
$\inf_{i\in E}R_i(u)>0$ for a given $u>0$. Define $f_i=b_i u_i-a_i u_{i-1}$ for $i\in E$, $f_{N+1}=0$ of
$N<\infty$. Then it is clear that
$${(f_i-f_{i+1})}/{u_i}=R_i(u)>0,\qqd i\in E.\tag 2.33$$
Hence, $f$ is strictly decreasing.

We now prove that $f\in {\scr F}_I$ whenever $\sum_i \mu_i=\infty$. First, we have
$$\align
\sum_{k\le i}\mu_k f_k=&\sum_{k\le i}\mu_k (b_k u_k-a_k u_{k-1})
=\sum_{k\le i}(\mu_k b_k u_k-\mu_{k-1}b_{k-1} u_{k-1})
=\mu_i b_i u_i>0,\\
&\text{\hskip 16em} i\in E.\tag 2.34\endalign$$
In particular, $f_0>0$. If $f_{k_0}\le 0$ for some $k_0\ge 1$, then $f_{k_0+1}<0$ and
$$\sum_{k_0+1\le i\le n} \mu_i f_i< f_{k_0+1} \sum_{k_0\le i\le n} \mu_i \to -\infty \qqd\text{as } n\to\infty$$
since $\sum_{i=0}^N \mu_i=\infty$. This implies that
$$\sum_{k_0+1\le i\le n} \mu_i f_i \to -\infty \qqd\text{as } n\to\infty.$$
Now, by (2.33), we would get
$$0<\mu_n b_n u_n=\sum_{i\le n} \mu_i f_i=\sum_{i\le k_0} \mu_i f_i+
\sum_{k_0+1\le i\le n} \mu_i f_i \to -\infty \qqd\text{as } n\to\infty,$$
which is impossible. Therefore, $f>0$ and then $f\in {\scr F}_I$.

Combining (2.33) with (2.34), we obtain that
$$R_i(u)= \frac{f_i-f_{i+1}}{u_i}
= {\mu_i b_i(f_i-f_{i+1})}\bigg/{\sum_{k\le i}\mu_k f_k}
=I_i(f)^{-1},\qqd i\in E.$$ Hence, we have first
$$\inf_{i\in E}\,R_i(u)=\inf_{i\in E}\,I_i(f)^{-1}\le \sup_{f\in {\scr F}_I}\inf_{i\in E}\,I_i(f)^{-1},$$
and then
$$\sup_{u>0}\,\inf_{i\in E}\,R_i(u)\le
\sup_{f\in {\scr F}_I}\inf_{i\in E}\,I_i(f)^{-1},$$
as required. We have thus proved the equality in (2.31) under $\sum_i \mu_i=\infty$.

Actually, we have proved in the last paragraph that
$(f_i=)b_i u_i-a_i u_{i-1}>0$ for all $i\in E$ and so $(v_i:=u_{i+1}/u_i)\in {\scr V}_*$ whenever $\inf_{i\in
E}R_i(u)>0$. This means that the set ${\scr V}\setminus {\scr V}_*$ is
useless since for each $v\in {\scr V}\setminus {\scr V}_*$,
we have $\inf_{i\in E}R_i(u)\le 0$.
Now, because of ${\scr V}_*\subset{\scr V}$ and (a), using the equality in (2.31), we obtain the
last assertion of part (1).

(d) To prove part (2) of the proposition,
note that the inequality ``$\le$'' is proved in (b). For
the inverse inequality, recalling that the main body in proof (c) is to show that
the function $f_i\,(i\in E)$ defined there is positive, this is now automatic due to the
definition of ${\scr V}_*$. The equality sign in (2.32) has already checked in proofs
(a) and (b) in the cases $\lz_0=0$ and $\lz_0>0$, respectively.
\qed\enddemo

\proclaim{Remark 2.8}{\rm For the equality in (2.31), the condition
$\sum_i \mu_i=\infty$ cannot be removed. For instance, consider the
ergodic case for which $\sum_i \mu_i<\infty$ but $\lz_0=0$ by
\thm\;3.1 below and so (2.31) is trivial. However, as proved in
[3;  \thm\;1.1] (cf. \thm\;6.1 below), the left-hand side of
(2.31) coincides with another eigenvalue (called $\lz_1$) which can
be positive. In this case, the equality in (2.31) fails. This also
explains the reason for the use of ${\scr V}_*$.}
\endproclaim

\proclaim{Remark 2.9}{\rm The test sequences with the same notation $(v_i)$ used in \thm\;2.4 and \prp\;2.7 are usually
different. Corresponding to the eigenfunction $(g_i)$ of $\lz_0$, the sequence constructed in
proof (d) of \thm\;2.4 and \prp\;2.5 is $v_i=g_{i+1}/g_i$, but the one constructed in proof (b)
of \prp\;2.7 is
$$v_i=\frac{g_{i+1}-g_{i+2}}{g_i-g_{i+1}}=\frac{1-g_{i+2}/g_{i+1}}{g_i/g_{i+1}-1}.$$
Thus, the mapping from the first sequence to the second one is as follows:
$$(v_i)_{0\le i<N}\to \bigg(\frac{1-v_{i+1}}{v_i^{-1}-1}\bigg)_{0\le i<N}, \tag 2.35$$
where on the right-hand side, $v_{N}$ is set to be zero if $N<\infty$.
}\endproclaim

\head{3. Absorbing (Dirichlet) boundary at infinity: criterion, approximating procedure and examples}\endhead

This section is a continuation of the last one. As applications of the variational formulas
given in the last section, a criterion
for the positivity of $\lz_0$ and an approximating procedure for $\lz_0$ are presented. The section is ended
by a class of examples and then the study on the first case of our
classification is completed.

\proclaim{\thm\;3.1 (Criterion and basic estimates)} The decay rate $\lz_0>0$ iff $\dz<\infty$, where
$$\dz=\sup_{n\in E}\mu[0, n]\, \nu[n, N]
=\sup_{n\in E}\sum_{j=0}^n\mu_j\,\sum_{k=n}^N \frac{1}{b_k\mu_k}.\tag 3.1$$
More precisely, we have $(4\dz)^{-1}\le \lz_0\le \dz^{-1}$. In particular, when $N=\infty$,
we have $\lz_0=0$ if
the process is recurrent (i.e., $\nu[1, \infty)=\infty$) and $\lz_0>0$ if the process is explosive
(i.e., condition (1.2) does not hold).
\endproclaim

\demo{\prf} (a) Let $\fz_n= {\sum_{j=n}^N\nz_j}=:\nu[n, N]$,\; $\nz_j=(b_j\mu_j)^{-1}$.
To prove the lower estimate, without loss of generality, assume that $\fz_0<\infty$.
Otherwise, $\dz=\infty$ and so the estimate is trivial. Next,
let $M_n=\mu[0, n]:=\sum_{k=0}^n \mu_k$.
By using the summation by parts formula
$$\sum_{k=0}^n x_k y_k=X_n y_n -\sum_{k=0}^{n-1}X_k(y_{k+1}-y_k),
\qqd X_n:=\sum_{j=0}^n x_j,\tag 3.2$$
in viewing the definition of $\dz$ and using the decreasing property of $\fz$, we get
$$\align
\sum_{j=0}^n \mu_j \sqrt{\fz_j}&= M_n\sqrt{\fz_n}+ \sum_{k=0}^{n-1}M_k \big(\sqrt{\fz_k}-\sqrt{\fz_{k+1}}\,\big)\\
&\le \frac{\dz}{\sqrt{\fz_n}}+{\dz} \sum_{k=0}^{n-1}
\frac{\sqrt{\fz_k}-\sqrt{\fz_{k+1}}}{\fz_k}.
\endalign$$
Noting that
$$ \big(\sqrt{\fz_k}-\sqrt{\fz_{k+1}}\,\big)/\fz_k\le 1/\sqrt{\fz_{k+1}} - 1/\sqrt{\fz_k},$$
we obtain
$$\sum_{j=0}^n \mu_j \sqrt{\fz_j}\le\frac{ 2\dz}{\sqrt{\fz_n}}.$$
Therefore,
$$I_n\big(\sqrt{\fz}\,\big)\le \frac{1}{\mu_n b_n \big(\sqrt{\fz_n}-\sqrt{\fz_{n+1}}\,\big)}\cdot\frac{2\dz}{\sqrt{\fz_n}}
=\frac{2\dz}{\sqrt{\fz_n}}\big(\sqrt{\fz_n}+\sqrt{\fz_{n+1}}\,\big)\le
4\dz.$$ By part (2) of Theorem 2.4, we have $\lz_0\ge
(4\dz)^{-1}$.

(b) Next, fix arbitrarily $n<m$ and let $f_i=\nu [i\vee n, m]\,
\dbl_{\{i\le  m\}}$. Then $f\in {\widetilde{\scr F}}_I$. To compute
$I_i(f)$, note that when $i<n$ or $i>m$, we have $f_i-f_{i+1}=0$ but
$\sum_{j\le i}\mu_j f_j\ge \mu_0 f_0>0$; and when $n\le i\le m$, we
have $f_i-f_{i+1}=\nu_i=(b_i\mu_i)^{-1}$. Hence, we have
$$I_i(f)
=\cases
\mu[0, n]\, \nz[n, m]+ \sum_{n+1\le j\le i}\mu_j\, \nz[j, m], &\quad\text{$n\le i\le m$},\\
\infty\;\text{(by convention, $1/0=\infty$)}, &\quad\text{otherwise.}
\endcases$$
Clearly, $I_i(f)$ achieves its minimum at $i=n$,
$$\inf_{i\in E}I_i(f)=  \mu[0, n]\, \nz[n, m].$$
Since $n, m\,(n<m)$ are arbitrary, by letting $m\to N$ and making the supremum in $n$, it follows that
$$\sup_{f\in {\widetilde{\scr F}}_I}\inf_{i\in E} I_i(f)
\ge  \sup_{n\in E}\mu[0, n]\, \nu[n, N]=\dz.$$
By using part (2) of Theorem 2.4 again, we obtain $\lz_0\le \dz^{-1}$. Note that in this proof,
we do not preassume that $\dz<\infty$.

(c) The particular assertion for the recurrent case is obvious. The explosive case is also easy since
$$\infty> \sum_{i=0}^\infty \mu_i\, \nu[i, \infty)> \sum_{i=0}^n\mu_i\, \nu[i, \infty)>\mu[0, n]\,\nu[n, \infty)$$
for all $n$, and so $\dz<\infty$.\qed\enddemo

The next result is parallel to [7; Theorem 2.2], and is a
typical application of parts (2) and (3) of \thm\;2.4. It provides
us a way to improve step by step the estimates of $\lz_0$. In view of
\thm\;3.1, the result is meaningful only if $\dz<\infty$.

\proclaim{\thm\;3.2 (Approximating procedure)}
Write $\nu_j=(\mu_j b_j)^{-1}$ and
$\fz_i=\nu[i, N]:=\sum_{j= i}^N \nu_j,\; i\in E$.
\roster
\item When $\fz_0<\infty$, define $f_1\!=\sqrt{\fz}$, $f_n\!=f_{n-1} I\!I(f_{n-1})$ and
$\dz_n\!=\sup_{i\in E} I\!I_i(f_n).$ Otherwise, define $\dz_n\equiv \infty$.
Then $\dz_n$ is decreasing in $n$ (denote its limit by $\dz_\infty$) and
$$\lz_0\ge \dz_{\infty}^{-1}\ge \cdots \ge \dz_1^{-1}\ge (4 \dz)^{-1},$$
where $\dz$ is defined in Theorem 3.1.
\item For fixed $\ell, m\in E$, $\ell<m$ and $m\ge 1$, define
$$\align
f_1^{(\ell, m)}&= \nu[\cdot\vee \ell, m]\,\dbl_{\le m},\\
f_n^{(\ell, m)}&= \dbl_{\le m}\,f_{n-1}^{(\ell, m)} I\!I \big(f_{n-1}^{(\ell, m)}\big),\qd n\ge 2,
\endalign
$$
where $\dbl_{\le m}$ is the indicator of the set $\{0, 1, \ldots, m\}$, and then define
$$\dz_n'=\sup_{\ell, m:\, \ell < m}\,\min_{i\le m} I\!I_i \big(f_{n}^{(\ell, m)}\big).$$
Then $\dz_n'$ is increasing in $n$ (denote its limit by $\dz_\infty'$) and
$$\dz^{-1}\ge {\dz_1'}^{-1}\ge \cdots\ge {\dz_\infty'}^{\!\!-1}\ge \lz_0.$$
Next, define
$${\bar\dz}_n
=\sup_{\ell,  m:\,\ell< m}\frac{\big\|f_n^{(\ell, m)}\big\|^2}{D\big(f_n^{(\ell, m)}\big)},\qqd n\ge 1.$$
Then ${\bar\dz}_n^{-1}\ge \lz_0$, ${\bar\dz}_{n+1}\ge \dz_n'$ for all $n\ge 1$ and $\bar\dz_1=\dz_1'$.
\endroster
\endproclaim

As the first step of the above approximation, we obtain the following improvement of
\thm\;3.1.

\proclaim{\crl\;3.3\,(Improved estimates)} We have
$$\dz^{-1}\ge {\dz_1'}^{-1}\ge \lz_0\ge \dz_1^{-1}\ge (4\dz)^{-1}, \tag 3.3$$
where
$$\align
\text{\hskip-6em}\dz_1&=\sup_{i\in E}\frac{1}{\sqrt{\fz_i}}\sum_{k\in E} \mu_k \fz_{i\vee k}\sqrt{\fz_k}\\
\text{\hskip-6em}&=\sup_{i\in E}\bigg[\sqrt{\fz_i}\sum_{k=0}^{i}\mu_k\sqrt{\fz_k}+
\frac{1}{\sqrt{\fz_i}}\sum_{i+1\le k<N+1}\mu_k \fz_k^{3/2}\bigg]. \tag 3.4\\
\text{\hskip-6em}\dz_1'\!&=\sup_{\ell\in E} \frac{1}{\fz_{\ell}}\!\sum_{k\in E} \mu_k \fz_{k\vee\ell}^2
\!\!=\sup_{\ell\in E}\bigg[{\fz_{\ell}}\, \mu [0, \ell]
\!+\!\frac{1}{\fz_{\ell}}\!\sum_{k=\ell +1}^N \mu_k \fz_{k}^2\bigg]\in [\dz,\, 2\dz].\text{\hskip-3em} \tag 3.5
\endalign$$
\endproclaim

\demo{\prf s of \thm\;$3.2$ and \crl\;$3.3$} (a) First, we prove part (1) of Theorem 3.2.
Noting that if $\fz_0=\infty$, then $\dz=\infty$ and $\dz_n=\infty$ for all $n\ge 1$,
the assertion becomes trivial in view of \thm\;3.1. Thus, we can assume that $\fz_0<\infty$.

By (2.23), we have
$$\dz_1=\sup_{i\in E} I\!I_i(f_1)\le \sup_{i\in E} I_i(f_1).$$
Proof (a) of Theorem 3.1 shows that the last one is bounded
from above by $4\dz$. This
gives us the lower bound of $\dz_1^{-1}$ as required.

We now prove the monotonicity of $\{\dz_n\}$. By induction, assume that $f_n<\infty$
and $\dz_n<\infty$. Then $f_{n+1}<\infty$. Note that
$$\align
\sum_{j\le k}\mu_j f_{n+1}(j)&=\sum_{j\le k}\mu_j f_n(j)f_{n+1}(j)/f_{n}(j)\\
&\le \sup_{i\in E}I\!I_i(f_n)\sum_{j\le k}\mu_j f_{n}(j)\\
&=\dz_n \sum_{j\le k}\mu_j f_n(j).\endalign$$
Multiplying both sides by $\nz_k$ and making a summation of $k$ from $i$ to $N$, by (2.14),
it follows that
$$f_{n+2}(i)\le \dz_n f_{n+1}(i).$$
Because $\dz_n<\infty$ and $ f_{n+1}(i)<\infty$, we obtain
$f_{n+2}<\infty$ and $I\!I_i (f_{n+1})\le \dz_n<\infty.$ Now, making
the supremum over $i$, we obtain $\dz_{n+1}\le \dz_n<\infty$.

We have thus proved part (1) of \thm\;3.2.

(b) To prove the monotonicity of $\dz_n'$ given in part (2) of Theorem 3.2,
we use the proportional property twice:
$$\align
\min_{i\le m} \big[ f_{n+1}^{(\ell, m)}\big/ f_{n}^{(\ell, m)}\big](i)
&=  \min_{i\le m} \sum_{j= i}^m\nu_j\sum_{k\le j} \mu_k f_{n}^{(\ell, m)}(k)
\bigg/ \sum_{j= i}^m\nu_j\sum_{k\le j} \mu_k f_{n-1}^{(\ell, m)}(k)\\
&\ge \min_{i\le m} \sum_{k\le i} \mu_k f_{n}^{(\ell, m)}(k)
\bigg/ \sum_{k\le i} \mu_k f_{n-1}^{(\ell, m)}(k)\\
&\ge \min_{i\le m} f_{n}^{(\ell, m)}(i)
\big/ f_{n-1}^{(\ell, m)}(i).
     \endalign $$
This implies that $\delta_{n+1}^{\prime}\ge \delta_n^{\prime}$.

By part (2) of Theorem 2.4, we also have $\dz_n'\le \lz_0^{-1}$ for all $n\ge 1$.
The assertion that ${\bar\dz}_n\le \lz_0^{-1}$ is obvious. Next, let $f=f_n^{(\ell, m)}$.
Then $g:=\dbl_{\supp(f)}f I\!I (f)=f_{n+1}^{(\ell, m)}$. As a consequence of (2.30),
we obtain ${\bar\dz}_{n+1}\ge \dz_n'$.

We have thus proved part (2) of Theorem 3.2 except the last assertion that
${\bar\dz}_1=\dz_1'$.

(c) We now prove (3.4) and $\dz_1'\ge \dz$. By (2.15), we have
$$\align
f_{n+1}(i)&=\sum_{k\in E} \mu_k f_n(k)\,\nu[i\vee k, N]\\
&=\sum_{k\in E} \mu_k f_n(k)\,\fz_{i\vee k}\\
&=\fz_i \sum_{k= 0}^{i} \mu_k f_n(k)+ \sum_{i+1\le k<N+1}\mu_k\,\fz_{k}  f_n(k). \tag 3.6
\endalign$$
In particular, with $f_1=\sqrt{\fz}$, we get
$$\align
f_{2}(i)&=\fz_i \sum_{k= 0}^{i} \mu_k \sqrt{\fz_k}+ \sum_{i+1\le k<N+1} \mu_k \fz_k^{3/2}.\tag 3.7
\endalign$$
From this, we obtain (3.4).

To prove $\dz_1'\ge \dz$, we need some preparation. As an analog of (3.6),  we have
$$\align
f_{n+1}^{(\ell, m)}(i)
&=\dbl_{\{i\le m\}} \sum_{k\le m} \mu_k f_n^{(\ell, m)}(k)\,\nz [i\vee k, m]. \tag 3.8
\endalign$$
In particular,
$$f_{2}^{(\ell, m)}(i)
=\dbl_{\{i\le m\}} \sum_{k\le m} \mu_k \,\nz [k\vee \ell, m]\,\nz [i\vee k, m].\tag 3.9$$
Since the right-hand
side is decreasing in $i$ for $i\le \ell<m$, $f_{1}^{(\ell, m)}(i)=f_{1}^{(\ell, m)}(\ell)$ for all $i\le \ell$,
and $f_{1}^{(\ell, m)}(i)=0$ for $i>m$, it follows that
$$\align
\min_{i\le m} I\!I_i \big(f_{1}^{(\ell, m)}\big)&=\min_{\ell\le i\le m} I\!I_i \big(f_{1}^{(\ell, m)}\big)\\
&=\min_{\ell\le i\le m} \frac{1}{\nz [i, m]} \sum_{j= i}^m \nu_j \sum_{k\le j} \mu_k\,\nz [k\vee \ell, m]\, \dbl_{\{k\le m\}}\\
&=\min_{\ell\le i\le m} \sum_{j= i}^m \nu_j \sum_{k\le j} \mu_k\,\nz [k\vee \ell, m]
  \bigg/\sum_{j= i}^m \nu_j\\
&\ge \min_{\ell\le i\le m} \sum_{k\le i}\mu_k\,\nz [k\vee \ell, m]\qd\Big[=\inf_{i\in E} I_i\big(f_{1}^{(\ell, m)}\big)\Big].
\endalign$$
Here in the last step, we have used the proportional property. Since the sum on the
right-hand side is increasing in $i$,
it is clear that
$$\min_{\ell\le i\le m} \sum_{k\le i}\mu_k\,\nz [k\vee \ell, m]=\nz [\ell, m]\sum_{k\le \ell}\mu_k
=\mu[0, \ell]\,\nu [\ell, m].$$
We have thus proved that
$$\dz_1'=\sup_{\ell< m}\,\min_{i\le m} I\!I_i \big(f_{1}^{(\ell, m)}\big)
\ge \sup_{\ell< m}\,\mu[0, \ell]\,\nu [\ell, m]
\ge \sup_{\ell\in E}\, \mu[0, \ell]\,\fz_{\ell}
=\dz.$$
A different proof of this is given in proof (d) below.

(d) We now compute $\dz_1'$. Note that by (3.9), we have
$$\align
f_2^{(\ell, m)}(i)
&\!=\dbl_{\{i\le m\}}\sum_{k=0}^m \mu_k \,\nz [k\vee \ell, m]\, \nz [i\vee k, m].
\endalign$$
Since $f_2^{(\ell, m)}(i)$ is decreasing in $i$ and
$f_1^{(\ell, m)}(i)$ is a constant on $\{0,1,\ldots, \ell\}$, it is clear that
$\min_{ i\le m}f_2^{(\ell, m)}(i)\big/f_1^{(\ell, m)}(i)=\min_{\ell\le i\le m}f_2^{(\ell, m)}(i)\big/f_1^{(\ell, m)}(i)$.
Besides, when $\ell\le i\le m$, we have
$$f_2^{(\ell, m)}(i)=\nz [\ell, m] \nz [i, m]\sum_{k=0}^{\ell} \mu_k
+ \nz [i, m] \sum_{\ell+1\le k\le i} \mu_k \nz [k, m]
+\sum_{i+1\le k\le m} \mu_k \nz [k, m]^2.$$
It follows that
$$\align
\min_{i\le m}\! \frac{f_2^{(\ell, m)}(i)}{f_1^{(\ell, m)}(i)}
&\!=\!\min_{\ell\le i\le m}\!\bigg[\nz [\ell, m]\!\sum_{k=0}^{\ell} \mu_k
\!+\!\!\! \sum_{k=\ell+1}^{i}\! \mu_k\nz [k, m]
\!+\!\frac{1}{\nz [i, m]}\!\sum_{k=i+1}^{m}\! \mu_k\nz [k, m]^2\bigg].
\endalign$$
We show that the sum on the right-hand side is increasing
in $i$. That is,
$$\align
\sum_{\ell+1\le k\le i}& \mu_k\nz [k, m]
+\frac{1}{\nz [i, m]}\sum_{k=i+1}^{m} \mu_k\nz [k, m]^2\\
&\le \sum_{k=\ell+1}^{i+1} \mu_k\nz [k, m]
+\frac{1}{\nz [i+1, m]}\sum_{i+2\le k\le m} \mu_k\nz [k, m]^2,\qqd \ell\le i\le m-1.\endalign$$
Collecting the terms, this is equivalent to
$$\frac{1}{\nz [i, m]} \mu_{i+1} \nz [i+1, m]^2\le \mu_{i+1} \nz [i+1, m]
+\bigg(\frac{1}{\nz [i+1, m]}-\frac{1}{\nz [i, m]}\bigg)\sum_{k=i+2}^{m} \mu_k\nz [k, m]^2.$$
Now, the conclusion becomes obvious because by the decreasing property of $\nz [i, m]$ in $i$,
the first term is controlled by the second, and the last one is nonnegative. We have thus obtained that
$$\min_{\ell\le i\le m} \frac{f_2^{(\ell, m)}(i)}{f_1^{(\ell, m)}(i)}
=\nz [\ell, m] \sum_{k=0}^{\ell} \mu_k
+\frac{1}{\nz [\ell, m]}\sum_{k=\ell+1}^m \mu_k \nz [k, m]^2.  \tag 3.10$$
As will be seen soon that the right-hand side is increasing in $m\,(>\ell)$, hence, we obtain
$$\dz_1'=\sup_{\ell<m}\min_{\ell\le i\le m} \frac{f_2^{(\ell, m)}(i)}{f_1^{(\ell, m)}(i)}
=\sup_{\ell\in E}\bigg[\fz_{\ell} \sum_{k=0}^{\ell} \mu_k
+\frac{1}{\fz_\ell}\sum_{\ell+1\le k<N+1} \mu_k\fz_k^2\bigg]. \tag 3.11$$
From this, it follows once again that $\dz_1'\ge \dz$. We now turn to prove the monotone property:
$$\align
\mu[0,\ell] &\,\nu [\ell, m+1]+\frac{1}{\nu [\ell, m+1]}\sum_{i=\ell+1}^{m+1} \mu_i\, \nu [i, m+1]^2\\
&\ge \mu[0,\ell]\, \nu [\ell, m]+\frac{1}{\nu [\ell,
m]}\sum_{i=\ell+1}^m \mu_i\, \nu [i, m]^2.
\endalign$$
Equivalently,
$$\mu[0,\ell]\, \nu_{m+1}+\frac{\mu_{m+1}}{\nu [\ell, m+1]}\,\nu_{m+1}^2
+ \sum_{i=\ell+1}^m \mu_i \bigg(\frac{\nu [i, m+1]^2}{\nu [\ell,
m+1]}- \frac{\nu [i, m]^2}{\nu [\ell, m]}\bigg)\ge 0.$$ This becomes
obvious since the term in the last bracket is positive:
$$\frac{\nu [i, m+1]^2}{\nu [i, m]^2}=\bigg(1+\frac{\nu_{m+1}}{\nu [i, m]}\bigg)^2
> 1+\frac{\nu_{m+1}}{\nu [\ell, m]}=\frac{\nu [\ell, m+1]}{\nu [\ell, m]},\qqd \ell\le i\le m.$$

(e) To show that $\dz_1'\le 2 \dz$, assume $\dz<\infty$.
By using the summation by parts formula (3.2) with $x_k=\mu_k$,
$X_k=\sum_{j=0}^k \mu_j$, and $y_k=\fz_{k\vee i}^2$, we get
$$\align
\sum_{k=0}^M \mu_k \fz_{k\vee i}^2
&=\fz_M^2 X_M +\sum_{k=0}^{M-1} X_k \big[\fz_{k\vee i}^2-\fz_{(k+1)\vee i}^2\big]\\
&=\fz_M^2 X_M +\sum_{k=i}^{M-1} X_k \big[\fz_{k}^2-\fz_{k+1}^2\big]\\
&=\fz_M^2 X_M +\sum_{k=i}^{M-1} X_k \nu_k (\fz_{k} +\fz_{k+1})\\
&<\fz_M^2 X_M + 2\sum_{k=i}^{M-1} X_k \nu_k \fz_{k}\\
&\le \dz \fz_M + 2\dz \sum_{k=i}^{M-1} \nu_k \qd(\text{since }X_k\fz_k\le \dz),\qqd i< M<N+1.
\endalign$$
If $N=\infty$, letting $M\to N$, it follows that
$$\sum_{k=0}^N \mu_k \fz_{k\vee i}^2\le 2 \dz \fz_i.$$
The same conclusion holds in the case that $N<\infty$ since
$$\dz \fz_N + 2\dz \sum_{k=i}^{N-1} \nu_k< 2 \dz \sum_{k=i}^{N} \nu_k=2 \dz\fz_i.$$
Hence,
$$\dz_1'=\sup_{i\in E}\frac{1}{\fz_i}\sum_{k=0}^N \mu_k \fz_{k\vee i}^2\le 2\dz.$$

(f) Now, it remains to compute $\bar\dz_1$. Since $f_1^{(\ell, m)}(i)=\nu[i\vee \ell, m] \,\dbl_{\{i\le m\}}$,
we have
$$\big\| f_1^{(\ell, m)}\big\|^2=\sum_i\mu_i\, \nu[i\vee \ell, m]^2 \dbl_{\{i\le m\}}
=\mu[0,\ell]\, \nu [\ell, m]^2 +\sum_{i=\ell+1}^m \mu_i\, \nu [i,
m]^2,$$ and
$$\align
D\big(f_1^{(\ell, m)}\big)
&=\sum_{i} \mu_i b_i \big(f_1^{(\ell, m)}(i+1)-f_1^{(\ell, m)}(i)\big)^2\\
&=\sum_{i=\ell}^{m-1} \mu_i b_i \big(\nu [i+1, m]-\nu [i, m]\big)^2+ \mu_m b_m \nu_m^2\\
&=\sum_{i=\ell}^{m-1} \nu_i+\nu_m\\
&=\nu [\ell, m].
\endalign
$$
Thus,
$$\frac{\big\| f_1^{(\ell, m)}\big\|^2}{D\big(f_1^{(\ell, m)}\big)}=
\mu[0,\ell]\, \nu [\ell, m]+\frac{1}{\nu [\ell, m]}\sum_{i=\ell+1}^m
\mu_i\, \nu [i, m]^2.$$
Hence, we have returned to (3.10). Since the right-hand side is increasing in $m$ as we have
seen in the proof of (3.11), we obtain
$$\bar\dz_1=\sup_{\ell<m}\frac{\big\| f_1^{(\ell, m)}\big\|^2}{D\big(f_1^{(\ell, m)}\big)}=\dz_1'.\qed$$
\enddemo

To conclude this section, we present some examples to illustrate the power of our results.
The first one is standard having constant rates.
\proclaim{\xmp\;3.4} Let $b_i\equiv b>0\,(i\ge 0)$, $a_i\equiv a>0\,(i\ge 1)$, $b>a$. Then
\roster
\item
$\lz_0=\big(\sqrt{a}-\sqrt{b}\,\big)^2$ with eigenfunction $g$:
$$g_n=\bigg(\frac{a}{b}\bigg)^{n/2}\bigg(n+1-n\sqrt{\frac{a}{b}}\,\bigg),\qqd n\ge 0,\qqd
g\notin L^1(\mu)\cup L^2(\mu).$$
\item $\dz=b (b-a)^{-2}$, $\dz_1'=(a+b)(b-a)^{-2}={\bar\dz}_1>\dz$, and
$\dz_1=\lz_0^{-1}$ which is exact. Note that $\dz_1/\dz_1'<2$ whenever $a\ne b$ and $\lim_{b\to a}\dz_1/\dz_1'=2$.
When $a=b$, we have $\lz_0=\dz_1^{-1}={\dz_1'}^{-1}=0$.
\endroster
\endproclaim

\demo{\prf} (a) First, we have
$\mu_n=(b/a)^n$, $n\ge 0$. Hence,
$$\sum_n \mu_n=\infty,\qqd \sum_n \mu_n g_n \ge \sum_n \mu_n g_n^2\ge \sum_n 1=\infty.$$
Next, since
$$\nu_i=\frac{1}{\mu_i b_i}=\frac{1}{b}\bigg(\frac{a}{b}\bigg)^i,$$
we have
$$\fz_{\ell}=\sum_{i\ge \ell} \nu_i= \frac{1}{b-a}\bigg(\frac{a}{b}\bigg)^{\ell}$$
and then
$$\sum_{i=0}^\infty  \mu_i\sum_{k=i}^\infty \frac{1}{b_k \mu_k}=\sum_{i=0}^\infty  \mu_i\fz_i=\infty.$$
Hence, (1.2) holds. It is easy to check that (2.12) holds:
$$\sum_{n=0}^\infty \mu_n g_n \nu[n, \infty)=\frac{1}{b-a}\sum_{n=0}^\infty
\bigg(\frac{a}{b}\bigg)^{n/2}\bigg(n+1-n\sqrt{\frac{a}{b}}\,\bigg)=\frac{1}{\lz_0}.$$

(b) To study $\lz_0$, according to (a), the Dirichlet form is regular and so the condition ``$f\in {\scr K}$''
in the definition of $\lz_0$ can be ignored. Thus,
$$\lz_0=\inf_{\|f\|=1} D(f)=b\inf_{\|f\|=1} \sum_{i\ge 0}\mu_i (f_{i+1}-f_i)^2.$$
It suffices to consider the case that $b=1$. Write $\gz=b/a>1$. Then we have
$$
g_k=\gz^{-k/2}(k+1-k \gz^{-1/2}),\qd \mu_k=\gz^k, \qd \nu_k=\gz^{-k}, \qd \fz_k= \gz^{-k+1}/(\gz-1),$$
and the required quantities are reduced to
$$\lz_0=\frac{(\sqrt{\gz}-1)^2}{\gz},\qd \dz=\frac{\gz^2}{(\gz-1)^2},\qd \dz_1=\frac{\gz}{(\sqrt{\gz}-1)^2},
\qd \dz_1'=\frac{\gz (\gz+1)}{(\gz-1)^2}.$$

Now, to prove part (1) of \xmp\;3.4, write $\xi=(\sqrt{\gz}-1)^2 \gz^{-1}$ for distinguishing with $\lz_0$.
Since $(g, \xi)$ satisfies the eigenequation, applying anyone of the variational formulas for the lower
estimate given in Theorem 2.4 with $f_i=g_i$ or
$$v_i=\frac{g_{i+1}}{g_i}= \sqrt{\frac{a}{b}}\bigg(1+
\frac{1-\sqrt{a/b}}{1+i (1-\sqrt{a/b}\,)}\bigg)
=\gz^{-1/2}\bigg(1+\frac{1-\gz^{-1/2}}{1+i (1-\gz^{-1/2})}\bigg),$$ it follows that
$\lz_0\ge \xi$. We have seen that the equality sign holds once $g\in L^2(\mu)$. Unfortunately, we
are now out of this case. Therefore, we need to show that $\lz_0\le
\xi$. To do so, one may use the truncated function of $g$:
$g_i^{(m)}=g_i \dbl_{\{i\le m\}}$. Then by the Stolz theorem, we
have
$$\lz_0\le \lim_{m\to\infty}\frac{D(g^{(m)})}{\|g^{(m)}\|^2}
 =\lim_{m\to\infty}\bigg[b_m+\frac{\mu_{m-1}b_{m-1}}{\mu_m}\bigg(1-\frac{2 g_{m-1}}{g_m}\bigg)
 \bigg]. \tag 3.12$$
The last limit equals $\xi$. Alternatively, noting that the leading order of $g\notin L^2(\mu)$ is $\gz^{-k/2}$, one
may adopt the test function $f_i=z^{-i/2}$ for $z>\gz$. Then $f\in L^2(\mu)$. The required
assertion follows by computing $D(f)/\|f\|^2$ and then letting $z\downarrow \gz$.
This proof benefits very much from the explicitly known expression of $\lz_0$.

(c) The computation of $\dz$ is easy:
$$\dz=\sup_{n\ge 0}\,\fz_n\! \sum_{j=0}^n \mu_j
=\frac{1}{(\gz-1)^2}\sup_{n\ge 0} \gz^{-n+1} \big(\gz^{n+1}-1\big)
=\frac{\gz^2}{(\gz-1)^2}. $$

(d) To compute $\dz_1$,
by (3.7), we have
$$\align
f_{2}(i)&=\fz_i \sum_{k= 0}^i \mu_k \sqrt{\fz_k}+ \sum_{k=i+1}^\infty \mu_k \fz_k^{3/2}\\
&=\frac{1}{(\gz -1)^{3/2}}\bigg\{\gz^{-i+1}\sum_{k=0}^i \gz^{k/2+1/2}
+\sum_{k\ge i+1} \gz^{-k/2+3/2}\bigg\}\\
&=\frac{\gz^{-i/2+3/2}}{(\gz -1)^{3/2}(\sqrt{\gz}-1)}\big(\sqrt{\gz}-\gz^{-i/2}+1\big).
\endalign$$
Therefore, we obtain
$$\dz_1=\sup_{i\ge 0}\frac{f_2(i)}{f_1(i)}=\frac{\gz}{(\gz-1)(\sqrt{\gz}-1)} \big(\sqrt{\gz}+1\big)
=\frac{\gz}{(\sqrt{\gz}-1)^2}=\frac{1}{\lz_0}.$$
Noting that even if neither $f_1$ nor $f_2$ is the eigenfunction, we still obtain the sharp
estimate.

(e) To compute $\dz_1'$, by (3.5), we have
$$\align
\dz_1'&=\sup_{\ell\in E}\bigg[\fz_{\ell}\sum_{k=0}^{\ell}
\mu_k+\frac{1}{\fz_{\ell}}\sum_{k\ge \ell+1} \mu_k \fz_k^2\bigg]\\
&=\frac{1}{\gz-1}\sup_{\ell\in E}\bigg[\gz^{-\ell+1} \sum_{k=0}^\ell
\gz^k+ \gz^{\ell+1}
\sum_{k\ge \ell+1} \gz^{-k}\bigg]\\
&=\frac{1}{(\gz-1)^2}\sup_{\ell\in E}\big[ \gz^2-\gz^{-\ell +1}+\gz\big]\\
&=\frac{\gz (\gz+1)}{(\gz-1)^2}.\qed
\endalign
$$
\enddemo

The next example is a typical linear model for which, interestingly, we have
a very simple and common eigenfunction. Moreover, the eigenvalue $\lz_0$ is determined
by the constant term $2\gz$ in the rates, but not the difference of the coefficients of
the leading term $i$, as in the ergodic case (cf. Example 6.8 below).

\proclaim{\xmp\;3.5} Let $b_i=2(i+\gz)\,(i\ge 0)$, $\gz>0$, $a_i=i\,(i\ge 1)$. Then
\roster
\item $\lz_0=\gz $,
$g_n=2^{-n}$ for $ n\ge 0$, and $g\in L^2(\mu)$.
\item When $\gz=1$, we have $\dz=\log 2\approx 0.69$, $\dz_1'\approx 0.84$, and $\dz_1\approx 1.09$.
Then $\dz_1/\dz_1'\approx 1.3<2$.
\endroster
\endproclaim

\demo{\prf} The uniqueness condition (1.2) is trivial since the birth rates are linear:
$$\sum_{k=0}^ \infty \frac{1}{b_k \mu_k} \sum_{i=0}^k \mu_i
=\sum_{k=0}^ \infty\bigg[\frac{1}{b_k}+ \frac{1}{b_k \mu_k} \sum_{i=0}^{k-1} \mu_i\bigg]
\ge \sum_{k=0}^ \infty\frac{1}{b_k}=\infty. $$

(a) Because
$$\mu_0=1,\qqd \mu_n=\frac{2^n \gz (1+\gz)\cdots (n-1+\gz)}{n!},\qqd n\ge 1,$$
it follows that $\mu_n> \gz 2^n/n$ and so $\sum_n\mu_n=\infty.$ Next, since
$$\mu_n b_n=\frac{2^{n+1} \gz (1+\gz)\cdots (n+\gz)}{n!}>\gz 2^{n+1},$$
we have
$\sum_n (\mu_n b_n)^{-1}<\infty$. Furthermore, we have
$$\sum_{n=0}^\infty \mu_n g_n^2=\sum_{n=0}^\infty \frac{2^{-n} \gz (1+\gz)\cdots (n-1+\gz)}{n!}.$$
The ratio test tells us $g\in L^2(\mu)$. Since $\lz_0$ is explicit
and $g\in L^2(\mu)$, it is simple to check that $(g_n)$ is the
eigenfunction of $\lz_0$. Hence, the proof of part (1) is done. For
this example, the sequence $(v_i)$ takes a simple form: $v_i\equiv
1/2$.

(b) When $\gz=1$, we have $\lz_0=1$,
$$\mu_i=2^i,\qqd \mu_i b_i= (i+1) 2^{i+1},\qqd \fz_i= \sum_{k\ge i+1} \frac{1}{2^{i} i},\qqd i\ge 0.$$
In particular, $\fz_0=\log 2$, $\fz_1=\log 2-1/2$.
Numerical computations show that the supremum in the definition of $\dz$, $\dz_1'$ and $\dz_1$
are attained at $0$, $0$ and $1$, respectively, and moreover,
$$\align
\dz&=\fz_0 \mu_0=\fz_0=\log 2\approx 0.69,\\
\dz_1'&=\fz_0 \mu_0+\frac{1}{\fz_0}\sum_{k\ge 1}\mu_k\fz_k^2
  =\log 2+\frac{1}{\log 2}\sum_{k\ge 1} 2^k \fz_k^2\approx 0.84,\\
{\dz_1}&=\sqrt{\fz_1}\,( \mu_0\sqrt{\fz_0}+\mu_1\sqrt{\fz_1}\,)+
  \frac{1}{\sqrt{\fz_1}}\sum_{k\ge 2}\mu_k\fz_k^{3/2}\\
  &=2\log 2-1+\frac{1}{2}\sqrt{(2\log 2) (2 \log 2-1)}
  +\sqrt{\frac{2}{2\log 2-1}}\sum_{k\ge 2} 2^k \fz_k^{3/2}\\
  &\approx 1.09.
\endalign
$$
We have thus proved part (2) of the conclusion.\qed\enddemo

The next example is often used in the study of convergence rates. For which, the first
eigenfunction is unknown but $\lz_0$ can still be computed.

\proclaim{\xmp\;3.6} Let $b_i=(i+1)^2$ and $a_i=i^2$. Then
$\dz=\pi^2/6\approx 1.64$, $\dz_1'\approx 2.19$, and $\dz_1=4$
which is sharp ($\lz_0=1/4$). Besides, $\dz_1/\dz_1'\approx 1.83<2$.
\endproclaim

\demo{\prf}
(a) Since $\mu_i\equiv 1$, $\nu_i= (i+1)^{-2}$, we have $\mu[0, i]=i+1$ and
$\fz_i=\sum_{j\ge i+1} j^{-2}$. For $\dz$ and $\dz_1'$, the supremum is attained
at $0$, therefore,
$$\dz= \fz_0=\sum_{k\ge 1}\frac{1}{k^2}=\frac{\pi^2}{6},$$
and
$$\dz_1'=\frac{1}{\fz_0}\sum_{k=0}^\infty {\fz_k^2}\approx 2.19.$$

(b) For $\dz_1$, the supremum is attained at $\infty$ and is equal
to $4$. By \crl\;3.3, this means that $\lz_0\ge 1/4$. This can be
also deduced by part (1) of \thm\;2.4 with $v_i=1-(2 i+4)^{-1}$ for
which the minimum of $R_i(v)$ is attained at $i=0$ and $i=\infty$.
It is even more simpler to use $v_i=1-(2 i +3)^{-1}$. Next, it is
known that $\lz_0\le 1/4$ (cf. \xmp\;5.5 below), hence, the estimate
is sharp. A direct proof for the upper estimate goes as follows.
Since the lower estimate is sharp, it indicates to use the test
function
$$f_i=\bigg(\sum_{j=i}^\infty \frac{1}{(j+1)^2}\bigg)^{1/2}\sim\frac {1}{\sqrt{i+1}}.$$
However, the last function is not in $L^2(\mu)$, and so one needs an approximating
procedure. Now, a carefully designed test function is the following:
$$f^{(\az)}_i=\frac{1}{\sqrt{(i+1)\az^{i+1}}},\qqd \az>1.$$
Then
$$\align
\mu\big(f^{(\az)\,2}\big)&=\sum_{i=0}^\infty  \frac{1}{(i+1)\az^{i+1}}
  =\sum_{i=1}^\infty  \frac{1}{i\az^{i}}=\log[\az (\az -1)^{-1}]<\infty,\\
D\big(f^{(\az)}\big)&=\sum_{i=0}^\infty (i+1)^2 \bigg[\frac{1}{\sqrt{(i+2) \az^{i+2}}}
  - \frac{1}{\sqrt{(i+1) \az^{i+1}}}\bigg]^2\\
&=\sum_{i=1}^\infty \frac{i^2}{\az^i} \bigg[\frac{1}{\sqrt{(i+1) \az}}
  - \frac{1}{\sqrt{i}}\bigg]^2\\
&=\sum_{i=1}^\infty \frac{i}{(i+1)\az^{i+1}} \frac{[(i+1)\az -i]^2}{[\sqrt{(i+1) \az}+\sqrt{i}\,]^2}\\
&\le\frac 1 4 \sum_{i=1}^\infty \frac{1}{(i+1)\az^{i+1}} [(i+1)\az -i]^2\\
&=\frac 1 4  (2+\log[\az (\az -1)^{-1}]).
\endalign $$
The required assertion now follows from
$$\lz_0\le \frac{2+\log[\az (\az -1)^{-1}]}{4 \log[\az (\az -1)^{-1}]}\to \frac 1 4 \qd\text{as}\;\; \az\downarrow 1.
 \qed$$
\enddemo

The last example below does not satisfy the non-explosive condition (1.2).

\proclaim{\xmp\;3.7} Let $b_i=(i+1)^4$ and $a_i=i(i-1/2)(i^2+3 i +3)$. Then
$\sum_i \mu_i<\infty$, $\sum_i \nu_i<\infty$, $\lz_0=1/2$,
$\dz\approx 1.83$, $\dz_1'\approx 1.9$, and $\dz_1\approx 2$.
Moreover, $\dz_1/\dz_1'\approx 1.05<2$.
\endproclaim

\demo{\prf} A simple computation shows that
$$\mu_i=\frac{i!^3}{\prod_{k=1}^i (k-1/2)(k^2+3 k+3)},\qqd
\nu_i=\frac{\prod_{k=1}^i (k-1/2)(k^2+3 k+3)}{(i+1) (i+1)!^3}.$$
From this, it follows that $\sum_i \mu_i<\infty$ and $\sum_i \nu_i<\infty$,
as an application of the typical Kummer's test: for a positive sequence
$\{x_n\}$, $\sum_n x_n$ converges or diverges according to $\kz>1$ or $\kz<1$,
respectively, where
$$\kz=\lim_{n\to\infty} n \bigg(\frac{x_n}{x_{n+1}}-1\bigg).   \tag 3.13$$

For each of $\dz$, $\dz_1'$ and $\dz_1$, the supremum is attained at $0$.

To see that $\lz_0=1/2$, first we check that $R_i(v)\equiv 1/2$ for
$$v_i=1-\frac{1}{2(i+1)}.$$
This gives us $\lz_0\ge 1/2$ by part (1) of \thm\;2.4. Since the corresponding
eigenfunction $g$,
$$g_i= \prod_{k=0}^{i-1} v_k=\frac{(2 i -1)!}{2^{2 i-1}i (i-1)!^2},\qd i\ge 1, \qqd g_0=1,$$
decreases strictly to $0$ and $\sum_i \mu_i<\infty$, we have
$g\in L^2(\mu)$. Now, because $-\ooz g=\lz_0 g$, $g_{\infty}=0$, and $D(f)=-(g, \ooz g)$,
it follows that $\lz_0= 1/2$ by (2.18).
\qed\enddemo

\head{4. Absorbing (Dirichlet) boundary at origin and reflecting (Neumann) boun\-dary at infinity}\endhead

This section deals with the second case of the boundary conditions.
The process has state space $E=\{i: 1\le i< N+1\}\,(N\le \infty)$,
birth rates $b_i>0$ but $b_N=0$ if $N<\infty$, and death rates $a_i>0$.
The rate $a_1>0$ is regarded as a killing from 1. Define
$$\lz_0=\inf\{D(f)/\mu(f^2): f\ne 0,\; D(f)<\infty\}, \tag 4.1$$
where $\mu(f)=\sum_{k\in E} \mu_k f_k$, and
$$\aligned
D(f)=\sum_{k\in E}\mu_k a_k(f_k-f_{k-1})^2,\qqd f_0:=0,\\
\mu_1=1,\qqd \mu_k=\frac{b_1\cdots b_{k-1}}{a_2\cdots a_k},\qd 2\le k<N+1.
\endaligned$$
The constant $\lz_0^{(4.1)}$ describes the optimal constant
$C=\lz_0^{-1}$ in the following {\it weighted Hardy inequality}:
$$\mu\big(f^2\big)\le C D(f),\qqd f_0=0$$
(cf. [9]). In other words, we are studying the
discrete version of the weighted Hardy inequality in this section.
To save the notation, in this and the subsequent
sections, we use the same notation $\lz_0$, $I$, $I\!I$, $R$ and so
on as in Section 2. Each of them plays a similar role but may have different meaning
in different sections.

To study $\lz_0$, as in Section 2, we need some parallel notation
originally introduced in [3, 7]:
$${ I}_i(f)= \frac{1}{\mu_{i}  a_{i} (f_{i}-f_{i-1})} \sum_{j=i}^N \mu_j f_j,
\qquad {I\!I}_i(f)= \frac{1}{f_i} \sum_{j=1}^i
 \frac{1}{\mu_j  a_j}\sum_{k=j}^N \mu_k f_k. $$
Here, for the first operator, we adopt the convention: $f_0=0$. The second one can be re-written as
 $${I\!I}_i(f)=\frac{1}{f_i} \sum_{k=1}^N \mu_k f_k\, \nu[1, i\wedge k],
 \qqd \nu[\ell, m]=\sum_{j=\ell}^m \nu_j,\qqd \nu_j=\frac{1}{\mu_{j} a_{j}}. $$
Next, define
$$R_i(v)=a_{i}\big(1-v_{i-1}^{-1}\big)+ b_i(1 - v_{i}),\qqd i\in E,\; v_0:=\infty
$$
($v_N$ is free if $N<\infty$ since $b_N=0$) and
$$\align
{{\scr F}}_{I\!I}&=\{f\!: f_i>0\; \text{for all } i\in E\},\\
{{\scr F}}_{I}&=\big\{f: f>0 \text{ and is strictly
increasing on }E\big\},\\
{\scr V}_1&=\{v: v_i >1\text{ for all }i\in E \}.
\endalign
$$
The modifications of ${{\scr F}}_{I\!I}$ and ${{\scr F}}_{I}$ are as
follows:
$$\align
{\widetilde{\scr F}}_{I\!I}&=\{f: \text{ there exists $m\in E$}
\text{ such that }\text{$f_i=f_{i\wedge m}>0$ for $i\in E$}\},\\
{\widetilde{\scr F}}_I&=\big\{f:
\text{ there exists  $m\in E$}\text{ such that }f_i=f_{i\wedge m}>0\text{ for $i\in E$\text{ and $f$ is} }\\
&\qquad\quad\;\text{ strictly increasing in $\{1,\ldots,
m\}$}\big\}.
\endalign
$$
Here, we use again the convention: $1/0=\infty$. Note that for the localization, $f$
is stopped at $m$ rather than vanishing after $m$ used in Sections 2 and 3. This is due to the fact that the
Neumann boundary is imposed at $m$ but not the Dirichlet one. Besides, for the
operator $I\!I$ here, the restriction on $\supp(f)$ used in Section
2 is no longer needed. Finally, define a local operator $\widetilde
R$ (depending on $m$) acting on
$$\align
{\widetilde{\scr V}}_1&=\cup_{m\in E}\big\{v: 1<v_i<1+a_i
\big(1-v_{i-1}^{-1}\big) b_i^{-1} \text{ for } i= 1, 2,
\ldots, m-1\\
&\qquad\qquad\qquad\text{ and } v_i\!=\!1 \text{ for }i\ge m\big\}
\endalign
$$
by replacing $a_m$ with ${\tilde a}_m:=\mu_m a_m\big/\sum_{k=m}^N
\mu_k$ in $R_i(v)$ for the same $m$ as in ${\widetilde{\scr V}}_1$.
Again, the change of $a_m$ is due to the Neumann boundary at $m$.
Note that if $v_i=1$ for all $i\ge m$, then ${\widetilde
R}_i(v)=R_i(v)=0$ for all $i> m$.

Before stating our main results in this section, we mention an
exceptional case that $\sum_i \mu_i=\infty$. On the one hand, by
choosing $f_0=0$ and $f_i=1$ for $i\ge 1$, it follows that
$$D(f)=\mu_1 a_1<\infty, \qqd \mu(f^2)=\sum_{i\ge 1}\mu_i=\infty$$
and so $\lz_0=0$. On the other hand, if $\sum_{i=1}^N\mu_i<\infty$, then
for every $f$ with $\mu(f^2)=\infty$, by setting $f^{(m)}=f_{\cdot\wedge m}\in L^2(\mu)$,
we get
$$\align
\infty>D\big(f^{(m)}\big)&=\sum_{i=1}^m \mu_i a_i (f_i-f_{i-1})^2\;\uparrow\; D(f)\qqd
\text{as } m\to\infty,\\
\infty>\mu\big(f^{(m)\,2}\big)&\ge\sum_{i=1}^m \mu_i f_i^2 \to \infty=\mu(f^2)\qqd
\text{as } m\to\infty.
\endalign$$
In words, for each non-square-integrable function $f$,
both $\mu(f^2)$ and $D(f)$ can be approximated by a sequence of square-integrable ones.
Hence, we can rewrite $\lz_0$ as follows:
$$\lz_0=\inf\{D(f): \mu(f^2)=1\}. \tag 4.2$$
In this case, as will be seen soon but not obvious, we also have
$$\lz_0=\inf\big\{{D}(f)\!: \mu(f^2)\!=\!1, f_i\!=\!
f_{i\wedge m}\text{ for some }m\!\in\! E \text{ and all }i\!\in\!
E\big\},\tag 4.3$$ Besides, we mention that the Dirichlet eigenvalue
$\lz_0$ is independent of $ b_0\ge 0$ (cf. [4; \thm\;3.4] or
[12; \thm\;3.7]).

For a large part of the paper, we do not use the uniqueness
condition (1.2) (note that a change of a finite number of the rates
$a_i$ and $b_i$ does not interfere in the uniqueness). Under (1.2), the process
is ergodic iff $\sum_i \mu_i<\infty$ (see [10; \thm\;4.45\,(2)], for instance).
If (1.2) fails but $N=\infty$,
then the decay rate for the minimal process is delayed to Section 7.
In (2.2), the condition ``$f\in {\scr K}$'' means that we deal with the
minimal process. This condition is removed in (4.2). It means that
we are in this section dealing with the maximal process in the sense
that the domain ${\scr D}^{\max}(D)$ of $D$ ignored in (4.2) is
taken to be the largest one: $\{f\in L^2(\mu): D(f)<\infty\}$
(that is the maximal process described at the beginning of Section 6
but killed at $1$). When
$N=\infty$, even though there is now a killing at $1$ (i.e.,
$a_1>0$), the regularity for (or the uniqueness of) the Dirichlet
form is still equivalent to (1.3):
$$\sum_{k=1}^ \infty \bigg(\frac{1}{b_k \mu_k}+ \mu_k\bigg)
= \infty  \tag "(1.3)'"
$$
since a modification of a finite number of rates does not change the
regularity (cf. \thm\;9.22 for further information). In this section
and Section 6, starting from any point in $E$, even though the
process can visit every larger state, it will come back in a
finite time. In this sense, the point infinity is regarded as a
reflecting boundary.

It is the position to finish the comparison of (4.1) and (4.2).
We have seen that $\lz_0^{(4.1)}=\lz_0^{(4.2)}$ once $\sum_i \mu_i<\infty$. We now claim that they
can be different otherwise. To see this, note that on the one hand, $\lz_0^{(4.1)}=0$
if $\sum_i \mu_i=\infty$, as proved above. On the other hand,
once (1.3)' holds \big(in particular, if $\sum_i \mu_i=\infty$, then\big)
by \prp\;1.3, $\lz_0^{(4.2)}$ coincides with
$$\inf\big\{D(f): f\in {\scr K}\!,\; \mu\big(f^2\big)=1\big\},$$
which is the one used in (7.1) below and can often be non-zero.
Thus, in general, $\lz_0^{(4.2)}\ge \lz_0^{(4.1)}$ and they
can be different. As will be seen in \thm\;7.1\,(2), in the special case that both of the series in (1.3)'
are divergent, we have $\lz_0^{(4.2)}= \lz_0^{(7.1)}=0$.

\proclaim{\thm\;4.1} Assume that $\sum_{i=1}^N\mu_i<\infty$. Then
the following variational formulas hold for $\lz_0$ defined by one of (4.1)---(4.3).
\roster
\item Difference form:
$$\align
&\inf_{v\in {\widetilde{\scr V}}_1}\,\sup_{i\in E}\,{\widetilde
R}_i(v) =\lz_0 =\sup_{v\in {\scr V}_1}\,\inf_{i\in
E}\,R_i(v).\endalign$$
\item Single summation form:
  $$\inf_{f\in {\widetilde{\scr F}}_I}\,\sup_{i\in E} I_i(f)^{-1}
    =\lz_0=\sup_{f\in{\scr F}_I}\,\inf_{i\in E} I_i(f)^{-1},$$
\item Double summation form:
  $$  \align
\lz_0&=\sup_{f\in{\scr F}_{I\!I}}\,\inf_{i\in E}{I\!I}_i(f)^{-1}
  =\sup_{f\in{\scr F}_{I}}\,\inf_{i\in E}{I\!I}_i(f)^{-1},\\
  \lz_0&=\inf_{f\in {\widetilde{\scr F}}_I}\,\sup_{i\in E}{I\!I}_i(f)^{-1}
  =\inf_{f\in {\widetilde{\scr F}}_{I\!I}}\,\sup_{i\in E}{I\!I}_i(f)^{-1}
  =\inf_{f\in {\widetilde{\scr F}}_{I\!I}\cup{\widetilde{\scr F}}_{I\!I}'}\,\sup_{i\in E}{I\!I}_i(f)^{-1},
  \endalign$$
  where
  ${\widetilde{\scr F}}_{I\!I}'=\big\{f: f_i>0 \text{ for all $i\in E$ and }f{{I\!I}}(f)\in L^2(\mu)\big\}$.
\endroster
\endproclaim

The next result was proved in [6] except the exceptional case
that $\sum_i \mu_i=\infty$ in which case $\lz_0=0$ (and $\dz=\infty$) and so
the assertion is trivial. See also \crl\;5.2 below. Note that $(\nu_j)$ below
is different from (2.15).

\proclaim{\thm\;4.2 (Criterion and basic estimates)} The rate $\lz_0$ defined by (4.1)
\big(or equivalently by (4.2) provided $\sum_{i\in E} \mu_i<\infty$\big)
is positive iff $\dz<\infty$, where
$$\dz=\sup_{n\in E} \nu [1, n]\,\mu[n, N]
=\sup_{n\in E} \sum_{i=1}^n \frac{1}{{\mu_{i}} { a_{i}}}
\sum_{j=n}^N{\mu_{j}}.
\tag4.4$$
More precisely, we have $(4\dz)^{-1}\le \lz_0\le {\dz}^{-1}$.
In particular, we have $\lz_0=0$ if $\sum_{i\in E} \mu_i=\infty$ and $\lz_0>0$ if
either $N<\infty$ or (1.3)' fails.
\endproclaim

\proclaim{\thm\;4.3 (Approximating procedure)}
Assume that $\sum_{i=1}^N \mu_i\!<\!\infty$ and $\dz<\infty$. Write $\fz_0=0$,
$\fz_i=\nu [1, i]:=\sum_{j=1}^i (\mu_j a_j)^{-1}$, $i\in E$.
\roster
\item Define $f_1= \sqrt{ \varphi }$, $f_n= f_{n-1} {I\!I}(f_{n-1})$ and
$ \delta_n=\sup_{i\in E}{I\!I}_i (f_n)$.
Then $ \delta_n$ is decreasing in $n$ and
$$ \lambda_0 \ge \delta_\infty^{-1}\ge\cdots
\ge {\delta_1}^{-1}  \ge (4 \delta )^{-1}.$$
\item For fixed $m \in E$, define
$$\align
f_1^{(m)}&= \varphi (\cdot \wedge
m),\\
f_n^{(m)} &=\big[f_{n-1}^{(m)} {I\!I}\big(f_{n-1}^{(m)}\big)\big](\cdot \wedge m),\qqd n\ge 2
\endalign$$ and then define
$\delta_n^{\prime} =
\sup_{m \in E} \inf_{i \in E}\,{I\!I}_{i} \big(f_n^{(m)}\big).$
Then $ \delta_n^{\prime}$ is increasing in
$n$ and
$${ \delta}^{-1}\ge \delta_1^{\prime\,-1}\ge\cdots \ge
\delta_\infty^{\prime\,-1} \ge \lambda_0.$$
Next, define
$${\bar{\dz}}_n
=\sup_{m\in E}\frac{\mu\big(f_n^{(m)\,2}\big)}{D\big(f_n^{(m)}\big)},\qqd n\in E.
$$
Then ${\bar{\dz}}_n^{-1}\ge \lz_0$, ${\bar\dz}_{n+1}\ge \dz_n'$ for all $n\ge 1$ and
${\bar{\dz}}_1 =\dz_1'$.
\endroster
\endproclaim

As the first step given in \thm\;4.3, we obtain the following improvement of
\thm\;4.2.

\proclaim{\crl\;4.4\,(Improved estimates)} For the rate $\lz_0$ defined by (4.1)
\big(or equivalently by (4.2) provided $\sum_{i\in E} \mu_i<\infty$\big), we have
$${ \delta}^{-1}\ge {\delta}_1^{\prime\,-1}\ge \lambda_0
\ge  \delta_1^{-1}  \ge (4 \delta )^{-1},$$
where
$$\align
\text{\hskip-8em}\dz_1&=\sup_{i\in E} \frac{1}{\sqrt{\fz_i}}\sum_{k\ge 1} \mu_k \fz_{i\wedge k}\sqrt{\fz_k}\\
\text{\hskip-8em}&=\sup_{i\in E}\bigg[\frac{1}{\sqrt{\fz_i}}\sum_{1\le k<i}
\mu_k \fz_{k}^{3/2}+ \sqrt{\fz_i}\sum_{k=i}^N
\mu_k\sqrt{ \fz_{k}}\, \bigg],\tag 4.5\\
\text{\hskip-8em}\dz_1'&\!=\!\sup_{m\in E}
\frac{1}{\fz_{m}}\!\sum_{k=1}^N \mu_k {\fz}_{k\wedge
m}^2\!\! =\!\sup_{m\in E}\!\bigg[\frac{1}{\fz_{m}}\!\sum_{k=1}^{m-1} \mu_k
{\fz}_{k}^2 \!+\!{\fz_{m}} \mu[m, N] \bigg]\!\in
[\dz,\, 2\dz].\text{\hskip-4em} \tag 4.6
\endalign$$
\endproclaim

\demo{\prf\; of \thm\;$4.1$} Note that
$${ I}_i(f)= \frac{1}{\mu_{i}  a_{i} (f_{i}-f_{i-1})} \sum_{j=i}^N \mu_j f_j
=\frac{1}{\mu_{i-1}  b_{i-1} (f_{i}-f_{i-1})} \sum_{j=i}^N \mu_j
f_j.$$ Hence, ${ I}_i(f)$ coincides with ${I}_{i-1}(f)$ used in
[3, 4, 6, 7, 12], whenever $ b_0>0$. The same change
is made for the operator $I\!I(f)$ in this section.

Throughout this proof, we use $\lz_0=\lz_0^{(4.3)}$ to denote the one given in
(4.3).
Similar to the proofs of \thm\;2.4 and \prp\;2.5, we adopt the following circle arguments:
$$\align
\text{\hskip-4em}\lz_0&\ge \lz_0^{(4.2)}\tag 4.7\\
&\ge  \sup_{f\in{\scr F}_{I\!I}}\,\inf_{i\in E}{I\!I}_i(f)^{-1}
=\sup_{f\in{\scr F}_{I}}\,\inf_{i\in E}{I\!I}_i(f)^{-1}
=\sup_{f\in{\scr F}_I}\,\inf_{i\in E} I_i(f)^{-1}\tag 4.8\\
&\ge \sup_{v\in {\scr V}_1}\,\inf_{i\in E}\,R_i(v)\tag 4.9\\
&\ge \lz_0\tag 4.10\endalign$$ and
$$\align
{\text{\hskip-4em}\lz_0} &\le  \inf_{f\in{\widetilde{\scr
F}}_{I\!I}\cup{\widetilde{\scr F}}_{I\!I}'}\,
   \sup_{i\in E}{I\!I}_i(f)^{-1}\tag 4.11\\
&\le\! \inf_{f\in{\widetilde{\scr F}}_{I\!I}}\,\sup_{i\in
E}{I\!I}_i(f)^{-1}\! =\!\inf_{f\in{\widetilde{\scr
F}}_I}\,\sup_{i\in E}{I\!I}_i(f)^{-1}
\!=\!\inf_{f\in {\widetilde{\scr F}}_I}\,\sup_{i\in E}{I}_i(f)^{-1}\text{\hskip-2em}\tag 4.12\\
&\le \inf_{v\in {\widetilde{\scr V}}_1}\,\sup_{i\in E}\,{\widetilde
R}_i(v)\tag 4.13\\
&\le \lz_0\tag 4.14
\endalign$$
Assertion (4.7) is obvious. The following assertions are proved
in [4; \thm\;3.3], or [12; \S 3.8] and [7; \S 2] (see
also the remark given in the next paragraph):
$$\align
&\text{\hskip-5em} \sup_{f\in{\scr F}_I}\,\inf_{i\in E} I_i(f)^{-1}=
\sup_{f\in{\scr F}_{I}}\,\inf_{i\in E}{I\!I}_i(f)^{-1}=
    \sup_{f\in{\scr F}_{I\!I}}\,\inf_{i\in E}{I\!I}_i(f)^{-1}\le\lz_0^{(4.2)}.\text{\hskip-1em}\tag 4.15\\
&\text{\hskip-5em}\inf_{f\in {\widetilde{\scr
F}}_I}\,\sup_{i\in E} I_i(f)^{-1} =\inf_{f\in {\widetilde{\scr
F}}_I}\,\sup_{i\in E}{I\!I}_i(f)^{-1}
    =\inf_{f\in {\widetilde{\scr F}}_{I\!I}}\,\sup_{i\in E}{I\!I}_i(f)^{-1}.\tag 4.16
\endalign$$
In particular, we have known (4.8) and (4.12) since the inequality in (4.12) is trivial.
It remains to prove (4.9)--(4.11), (4.13) and (4.14).

In [7; \S 2] and [12; \S 3.8], only the ergodic case under condition (1.2) is
considered. But for (4.15) and (4.16), one does not need (1.2). Actually, one can now follow the proofs
of \thm\;2.4 and \prp\;2.5 with a little change. For instance, to
prove the last inequality in (4.15), following proof (a) of
\thm\;2.4 and \prp\;2.5, let $g$ satisfy $\|g\|=1$ and $g_0=0$. Then
$$ \align   1 &= \sum_{i\in E} \mu_i g_i^2\qd\text{(since $\|g\|=1$)}\\
&= \sum_{i\in E} \mu_i \bigg( \sum_{j=1}^i(g_j- g_{j-1})\bigg) ^2\qd\text{(since $g_0=0$)}\\
&\le  \sum_{i\in E} \mu_i  \sum_{j=1}^i \frac{(g_{j}-g_{j-1} )^2 \mu_j a_j}
{h_j}\sum_{k=1}^i \frac{h_k}{ \mu_k a_k}.\endalign$$
Exchanging the order of the first two sums on the right-hand side, we get
$$\align
1 &\le \sum_{j\in E} \mu_j a_j (g_{j}-g_{j-1})^2 \frac{1}{h_j}
\sum_{i=j}^N\mu_i \sum_{k=1}^i \frac{h_k}{ \mu_k a_k}\\
&\le  D(g) \sup_{j\in E} \frac{1}{h_j}
\sum_{i =j}^N\mu_i \sum_{k=1}^i \frac{h_k}{ \mu_k a_k}\\
&=: D(g) \sup_{j\in E} H_j.
 \endalign $$
The next step is to choose $h_j=\sum_{i=j}^N\mu_i f_i$ for a given $f\in {\scr F}_{I\!I}$
with $\sup_{j\in E}I\!I_j(f)\!<\infty$. From these, it should be clear what change is required
in order to prove (4.15) and (4.16).

We now begin to work on the additional part of the proof.
\medskip

(a) Prove that $\sup_{f\in{\scr F}_{I\!I}}\,\inf_{i\in E}{I\!I}_i(f)^{-1}
\ge \sup_{v\in {\scr V}_1}\,\inf_{i\in E}\,R_i(v)$.
\medskip

As in proof (c) of \thm\;2.4 and \prp\;2.5, we use $R_i (u)$,
$$R_i (u)= a_{i}\bigg(1-\frac{u_{i-1}}{u_i}\bigg)+b_{i}\bigg(1-\frac{u_{i+1}}{u_i}\bigg),\qqd i\in E,$$
($u_{N+1}\text{ is free if }N<\infty$ since $b_N=0$),
instead of $R_i (v)$, where $u_i>0$ for $i\in E$ and $u_0=0$. Then
 $v_i>1\,(i\in E)$ means that $u_{i+1}>u_i> 0$, and $v_i=1$
for $i\ge m$ means that $u_i=u_{i\wedge m}>0$.

Without loss of generality, assume that
$\inf_{i\in E}R_i(u)>0$ for a given strictly increasing $u$ with $u_0=0$.
Define $f_i=( a_i+ b_i) u_i- a_i u_{i-1}- b_i u_{i+1} \,\big[=u_i {R}_i(u)\big]$
for $i\in E$ and $f_0=0$. Then by assumption,
$${f_i}/{u_i}=R_i(u)>0,\qqd i\in E.$$
Hence, $f\in {\scr F}_{I\!I}$. Next, since
$$\align
0<\mu_k f_k&=\mu_k  a_k (u_k-u_{k-1})
-\mu_{k+1} a_{k+1} (u_{k+1}-u_k)\endalign$$
and the strictly increasing property of $u_i$ in $i$, it follows that
$$0<\sum_{k=j}^N\mu_k f_k\le \mu_j  a_j (u_{j}-u_{j-1}),$$
and so
$$u_i=\sum_{j=1}^i (u_j-u_{j-1})\ge \sum_{j=1}^i\nu_j \sum_{k=j}^N\mu_k f_k>0.$$
We obtain
$$R_i(u)
= {f_i}/{u_i}\le {I\!I}_i(f)^{-1}, \qqd i\in E.$$
Therefore, we have first
$$\inf_{i\in E}\, R_i(u)
\le\inf_{i\in E}\,{I\!I}_i(f)^{-1}
\le \sup_{f\in {\scr F}_{I\!I}}\inf_{i\in E}\,{I\!I}_i(f)^{-1},$$
and then
$$\sup_{v\in {\scr V}_1}\,\inf_{i\in E}\, R_i(v)\le
\sup_{f\in {\scr F}_{I\!I}}\inf_{i\in E}\,{I\!I}_i(f)^{-1},$$
as required.
\medskip

(b) Prove that $\sup_{v\in {\scr V}_1}\,\inf_{i\in E}\,R_i(v)\ge \lz_0$.
\medskip

First, we show that $\sup_{v\in {\scr V}_1}\,\inf_{i\in E}\,R_i(v)\ge 0$.
For a given positive $f\in L^1(\mu)$, let
$u =f I\!I(f)$. Then $u_{i+1}/u_i>1$ and $R_i(u)=f_i/u_i>0$ for all $i\in E$.
With $(v_i=u_{i+1}/u_i)\in {\scr V}_1$,
this implies $\inf_{i\in E}R_i(v)\ge 0$ and then the required assertion follows.

Alternatively, since $a_1>0$, the eigenfunction is still strictly increasing
when $\lz_0=0$ by part (3) of \prp\;2.1.
Hence the proof in the case of $\lz_0=0$ can be combined
into the next paragraph, and then the last paragraph can be omitted.

By assumption, we have $\sum_{i\in E} \mu_i<\infty$. When $\lz_0>0$,
it was proved in proof (d) of [12; \thm\;3.7] that the
eigenfunction of $\lz_0^{(4.2)}$ is strictly increasing. Even though
${\lz}_0$ could formally be bigger than $\lz_0^{(4.2)}$, the same
proof still works for the eigenfunction $g$ of $\lz_0$ since the
modified function $\bar g$ used there satisfies ${\bar g}_i={\bar
g}_{i\wedge n}$ for some $n$. Having this at hand, the proof is just
a use of the eigenequation:
$$-\ooz g(i):= - b_i(g_{i+1}-g_i)+  a_i (g_i-g_{i-1})
=\lz_0 g_i,\qqd i\in E,\; g_0:=0$$
($g_{N+1}^{}\text{ is free if }N<\infty$ since $b_N=0$).
With $v_i:=g_{i+1}/g_i>1$ for $i<N$,
this gives us $v\in {{\scr V}}_1$ and
$R_i(v)\equiv \lz_0$, and so the assertion follows.

We have thus completed the circle argument of (4.7)---(4.10).
\medskip

(c) Prove that $\lz_0 \le  \inf_{f\in{\widetilde{\scr F}}_{I\!I}
\cup{\widetilde{\scr F}}_{I\!I}'}\,\sup_{i\in
E}{I\!I}_i(f)^{-1}$.
\medskip

In the original proof of [7; Theorem 2.1], when $N=\infty$,
from the estimate
$$D(g)\le \mu(g^2)\, \sup_{i\in E} {{I\!I}}_i (f)^{-1}$$
for $f\in {\widetilde{\scr F}}_{I\!I}$ and $g:=[f{I\!I}(f)](\cdot\wedge m)$ to
conclude that $\lz_0\le D(g)/\mu(g^2)$, one requires an additional
condition $g\in L^2(\mu)$, provided $m=\infty$ is allowed. This is
the reason why the set ${\widetilde{\scr F}}_{I\!I}'$ in part (3) of
Theorem 4.1 is added. Anyhow, with the modified conditions, the same
proof gives us the required assertion (cf. proof (f) of Theorem
2.4 and \prp\;2.5).
\medskip

(d) Prove that $\inf_{f\in{\widetilde{\scr F}}_{I\!I}}\,\sup_{i\in
E}{I\!I}_i(f)^{-1} \le \inf_{v\in {\widetilde{\scr
V}}_1}\,\sup_{i\in E}\,{\widetilde R}_i(v)$.
\medskip

Given $u$ with $u_0=0$ and $u_i=u_{i\wedge m}$ for all $i\in E$ so
that $(v_i:=u_{i+1}/u_i)\in {\widetilde{\scr V}}_1$, let
$$f_i=
\cases
( a_i+ b_i) u_i- a_i u_{i-1}- b_i u_{i+1},&\qd i\le m-1\\
 {\tilde a}_m (u_m- u_{m-1}), &\qd i\ge m.
\endcases
$$
It is simple to check that $f_0=0$,
$$f_i/u_i={\widetilde
R}_i(u)>0\text{ for } i\in \{1,\ldots, m\} \text{ and }f_i=f_m
\text{ for }i>m,$$ and so $f\in {\widetilde{\scr F}}_{I\!I}$.
Moreover, since
$$\align
\sum_{k=j}^{m-1}\mu_k f_k&=\mu_j a_j( u_j-u_{j-1})- \mu_m a_m(
u_m-u_{m-1})\\
&=\mu_j a_j( u_j-u_{j-1})- f_m \sum_{k=m}^N \mu_k,\endalign$$
we get
$$0<\sum_{k=j}^N\mu_k f_k=\mu_j a_j( u_j-u_{j-1}).$$
It follows that
$$0<u_i=\sum_{j=1}^i (u_j-u_{j-1})=
 \sum_{j=1}^i\nu_j \sum_{k=j}^N\mu_k f_k,\qqd i\in \{1,\ldots,
m\},$$
and then ${\widetilde R}_i(u)=f_i/u_i= {I\!I}_i(f)^{-1}$ for
$i\in \{1,2,\ldots, m\}$. Therefore, we have
$$\max_{1\le i\le m}{\widetilde
R}_i(u)=\! \max_{1\le i\le m} {I\!I}_i(f)^{-1}\! \ge\!\! \inf_{f\in
{\widetilde{\scr F}}_{I\!I},\,f_i=f_{i\wedge m}} \max_{1\le i\le m}
{I\!I}_i(f)^{-1}\!\ge\! \inf_{f\in {\widetilde{\scr F}}_{I\!I}}
\sup_{i\in E} {I\!I}_i(f)^{-1}\!,$$
and then
$$\inf_{v\in {\widetilde{\scr V}}_1}\sup_{i\in E}\,{\widetilde
R}_i(v)\ge \inf_{f\in {\widetilde{\scr F}}_{I\!I}} \sup_{i\in E}\,
{I\!I}_i(f)^{-1}.$$
\medskip

(e) Prove that $\inf_{v\in {\widetilde{\scr V}}_1}\,\sup_{i\in
E}\,{\widetilde R}_i(v)\le \lz_0$.
\medskip

Recall the definition of $\lz_0$:
$$\lz_0=\inf\big\{{D}(f): \mu(f^2)=1, f_i=f_{i\wedge m}
\text{ for some }m\in E\text{ and all }i\in E\big\}.$$
Clearly, we have
$${\lz}_0^{(m)}:=\inf\big\{{D}(f): \mu(f^2)=1,\,
f_i=f_{i\wedge m}\text{ for all }i\in E\big\}\;\downarrow \lz_0\; \text{ as } m\uparrow N.$$
We now explain the meaning of $\lz_0^{(m)}$ as follows. Let
$$\align
{\tilde \mu}_i&=\mu_i,\; 1\le i<m,\qd {\tilde \mu}_m=\sum_{i=m}^N \mu_i,\\
{\tilde a}_i&=a_i,\; 1\le i<m,\qd {\tilde a}_m=\mu_m a_m/{\tilde \mu}_m,\\
&{\widetilde D}(f)=\sum_{i=1}^m {\tilde \mu}_i {\tilde a}_i (f_i-f_{i-1})^2.\tag 4.17
\endalign
$$
Then ${\tilde \mu}_i{\tilde a}_i=\mu_i a_i$ for $i=1,\ldots, m$,
${\widetilde D}(f)=D(f)$ and ${\tilde\mu}(f^2)=\mu(f^2)$ for
every $f$ with $f=f_{\cdot\wedge m}$. Thus, $\lz_0^{(m)}$ is just
the first eigenvalue of the local Dirichlet form $\big({\widetilde
D}, {\scr D}\big({\widetilde D}\big)\big)$ having the state space
$\{1, \ldots, m\}$, with Dirichlet (absorbing) boundary at $0$ and
Neumann (reflecting) boundary at $m$. Let $g$ ($g_0=0$) be the
eigenfunction of the local first eigenvalue ${\lz}_0^{(m)}$. Extend
$g$ to the whole space by setting $g_i=g_{i\wedge m}$. Next, set
$u_i=g_i$ for $i<N$. Then
$${\widetilde
R}_i(u)= \cases \lz_0^{(m)}>0,
&\qd i\in\{1,\ldots, m\},\\
0, &\qd i>m.
\endcases\tag 4.18$$
Furthermore, for $v_i\!:=\!u_{i+1}/u_i$, we have $v_0\!=\!\infty$, $v_i\!>\!1$
on $\{1,\ldots, m-1\}$, and $v_i=1$ for $i\ge m$. Thus, by (4.18), it
is easy to check that $v=(v_i)\!\in {\widetilde{\scr V}}_1$. Therefore,
$$\align
\lz_0^{(m)}&=\max_{1\le i\le m}\, {\widetilde
R}_i(v)\\
&\ge \inf_{v\in{{\widetilde{\scr V}}_1}:\;v_i=1 \text{ for }i\ge
m}\,\max_{1\le i\le m}\,{\widetilde
R}_i(v)\\
&\ge \inf_{v\in{{\widetilde{\scr V}}_1}:\;v_i=1 \text{ for }i \,\ge
\text{ some }n> 1}\,\sup_{i\in E}\,{\widetilde
R}_i(v)\\
&= \inf_{v\in {\widetilde{\scr V}}_1}\sup_{i\in E}{\widetilde
R}_i(v).
\endalign$$
The assertion now follows by letting $m\to N$.
\qed\enddemo

\demo{\prf\; of \thm\;$4.3$}

(a) We remark that the sequence $\big\{f_n^{(m)}\big\}_{n\in E}$ is
clearly contained in ${\widetilde{\scr F}}_I$. But the
modified sequence used in [7; \thm\;2.2],
$${\tilde f}_1^{(m)}=\varphi (\cdot \wedge m),\qqd
{\tilde f}_n^{(m)} = {\tilde f}_{n-1}^{(m)}(\cdot \wedge m)\,
{I\!I}\big({\tilde f}_{n-1}^{(m)}(\cdot \wedge m)\big),\quad n\ge 2,$$ is
usually not contained in ${\widetilde{\scr F}}_{I\!I}$. However,
$$\align
\delta_n^{\prime} &=
\sup_{m\in E}
\inf_{i\in E} {I\!I}_i\big( f_n^{(m)} \big)\\
&=\sup_{m\in E}
\min_{1\le i\le m} {I\!I}_i\big(f_n^{(m)} \big)\\
&= \sup_{m\in E}
\min_{1\le i\le m} {I\!I}_i\big({\tilde f}_n^{(m)} (\cdot\wedge m)\big)\\
&=\sup_{m\in E}
\inf_{i\in E} {I\!I}_i\big({\tilde f}_n^{(m)} (\cdot\wedge m)\big).
\endalign$$
Here in the last step, we have used the convention $1/0=\infty$.
Hence, these two sequences produce the same $\{\dz_n'\}$.

(b) The approximating procedure given in \thm\;4.3 is mainly a copy
of [7; Theorem 2.2] (cf. the proof of \thm\;3.2). For
later use, here we review the proof of part (1). From [6;
proof of \thm\;3.5], we have known that
$$I_j (f_1)=\frac{1}{\mu_j a_j (f_1(j)-f_1(j-1))}
\sum_{k\ge j}\mu_k f_1(k)\le 4\dz,\qqd j\ge 1.$$
Hence (Alternatively, by the proportional property),
$$f_2(i)=
\sum_{j=1}^i \frac{1}{\mu_j a_j }\sum_{k\ge j}\mu_k f_1(k)
\le 4\dz\sum_{j=1}^i (f_1(j)-f_1(j-1)) = 4\dz
f_1(i).$$ This gives us the assertion $\dz_1=\sup_{i\ge 1} I\!I_i
(f_1)\le 4\dz$.

To prove the monotonicity of $\{\dz_n\}$ and $\{f_n\}\subset L^1(\mu)$,
we adopt induction. As we have just seen,
$$\dz_1=\sup_{i\ge 1} \frac{f_2(i)}{f_1(i)}\le 4\dz.$$
This means that $f_1\in L^1(\mu)$ (or equivalently, $f_{2}<\infty$)
and $\dz_1<\infty$ since $\dz<\infty$ by assumption. Assume that
$f_n\in L^1(\mu)$ (or equivalently, $f_{n+1}<\infty$) and
$\dz_n<\infty$. Then
$$\sum_{k\ge j}\mu_k f_{n+1}(k)=\sum_{k\ge j}\mu_k f_n(k) [f_{n+1}(k)/f_n(k)]
\le \dz_n\sum_{k\ge j}\mu_k f_{n}(k).$$
Multiplying both sides by $\nu_j$ and making summation from $1$ to $i$, it follows that
$$f_{n+2}(i)\le \dz_n f_{n+1}(i),\qqd i\ge 1.$$
Since $f_{n+1}<\infty$ and $\dz_n<\infty$ by assumption, we have $f_{n+2}<\infty$, and
$$\frac{f_{n+2}(i)}{f_{n+1}(i)}=I\!I_i(f_{n+1})<\infty,\qqd i\ge 1.$$
This proves not only $f_{n+1}\in L^1(\mu)$ but also $\dz_{n+1}\le \dz_n<\infty$.

The assertion that ${\bar\dz}_n^{-1}\ge \lz_0$ is obvious by (4.2).
Similar to proof (b) of \thm\;3.2, the assertion
${\bar\dz}_{n+1}\ge \dz_n'$ is a consequence of the last part of the
proof of [7; \thm\;2.1].\qed\enddemo

\demo{\prf\; of \crl\;$4.4$}

(a) The degenerated case that $\sum_i \mu_i=\infty$ is trivial since
$\lz_0^{(4.1)}=0$ and $\dz=\dz_1=\dz_1'=\infty$. The
main assertion of \crl\;4.4 is a consequence of \thm\;4.3. Here, we
consider (4.6) only since the proof of (4.5) is easier. Note that
$${I\!I}_i\big(f_1^{(m)}\big)=\frac{1}{\varphi_{i\wedge m}}
\sum_{j=1}^i\frac{1}{\mu_j  a_j}\sum_{k=j}^N\mu_k \varphi_{k\wedge m}.$$
The right-hand side is clearly increasing in $i$ for $i\ge m$ and is
decreasing (not hard to check) in $i$ when $i\le m$. Hence, ${I\!I}_i\big(f_1^{(m)}\big)$
achieves its minimum at $i=m$. Then, by exchanging the order of the
sums, it follows that the minimum is equal to
$$\frac{1}{\varphi_m}\sum_{k=1}^N \mu_k \varphi_{k\wedge m}^2.$$
This observation is due to Sirl,  Zhang and Pollett (2007). We have thus proved the first equality
in (4.6).

Next, following the proof of [6; Theorem 3.5], we have
$${D}\big(f_1^{(m)}\big)=\sum_{i=1}^m \mu_i a_i (\fz_i-\fz_{i-1})^2
=\fz_m,$$
and
$$\mu\big(f_1^{(m)\,2}\big)=\sum_{k=1}^N \mu_k \varphi_{k\wedge m}^2.$$
Combining these facts together, it follows
that ${\bar\dz}_1=\dz_1'$.

(b) Finally, we prove the estimates in (4.6).
The lower estimate of $\dz_1'$ is rather easy since
$$\frac{1}{\varphi_m}\sum_{k=1}^N \mu_k \varphi_{k\wedge m}^2
\ge \frac{1}{\varphi_m}\sum_{k=m}^N \mu_k \varphi_{k\wedge m}^2=\fz_m \sum_{k=m}^N \mu_k.$$
For the upper estimate, use the summation by parts formula:
$$\sum_{k=1}^N \mu_k \varphi_{k\wedge m}^2
=\sum_{k=1}^m\big[ \varphi_{k}^2- \varphi_{k-1}^2\big]
\sum_{j=k}^N \mu_j
=\sum_{k=1}^m\frac{\varphi_{k}+ \varphi_{k-1}}{\mu_k a_k}
\sum_{j=k}^N \mu_j.
$$
It follows that
$$\frac{1}{\varphi_m}\sum_{k=1}^N \mu_k \varphi_{k\wedge m}^2
< \frac{2}{\varphi_m}\sum_{k=1}^m\frac{1}{\mu_k a_k}\bigg[\varphi_{k}
\sum_{j=k}^N \mu_j\bigg]
\le \frac{2\dz}{\varphi_m}\sum_{k=1}^m\frac{1}{\mu_k a_k}
=2\dz.$$
The estimate now follows by making the supremum with respect to $m\in E$.\qed\enddemo

\head{5. Dual approach}\endhead

This section is devoted to the duality of the processes studied in the previous sections, as well as
a duality to be used in the next two sections. Again, the section is ended by a class of examples.

Suppose that we are given a birth--death process with state space $E=\{i: 0\le i<N+1\}\,(N\le \infty)$,
birth rates $b_i>0$ ($b_0>0$, especially) but $b_{N}\ge 0$ if $N< \infty$, and death rates $a_i>0$ but $a_0=0$.
The case that $b_N>0$ is used in this section while the case of $b_N=0$ is
for use in Section 7.
Define a dual chain with state space $\widehat E=\{i: 1\le i<N'+1\}$ and with rates as follows:
$$\hat b_0=0, \qd \hat b_i=a_i, \qd \hat a_i=b_{i-1},\qqd i\in {\widehat E},\tag 5.1$$
where $a_{N+1}=b_{N+1}=0$ if $N<\infty$ by convention and
$$ N'=
{\cases
N,&\qd N<\infty\text{ and } b_N=0,\\
N+1,&\qd N<\infty\text{ and }  b_N>0,\\
\infty, &\qd N=\infty.
\endcases}$$
The dual process with rates $\big({{\hat a}_i, {\hat b_i}}\big)$ has an absorbing at $0$. When $N<\infty$,
for the dual process,
the state $N+1$ is absorbing if $b_{N}= 0$ (then ${\hat a}_{N+1}=0$ but ${\hat b}_N>0$);
otherwise, it is a reflecting boundary since ${\hat a}_{N+1}=b_N>0$.
In a word, the absorbing boundary is dual to the reflecting one and vice versa.
This dual technique goes back to Karlin and McGregor (1957b, \S 6). Next, define
$${\hat\mu}_1=1,\qqd {\hat\mu}_n
=\frac{{\hat b}_1\cdots {\hat b}_{n-1}}{{\hat a}_2\cdots {\hat a}_{n}},\qqd 2\le n< N'+1.\tag 5.2$$
When $N<\infty$ and $b_N>0$, then ${\hat a}_{N+1}>0$, and so ${\hat\mu}_n$ can be defined up to
$n=N+1$. Otherwise, it can be defined up to $n=N$ only. It is now easy to check
(noticing the difference of $(\nu_j)$ and $(\hat\nu_j)$) that
$$\hat\mu_n=\frac{b_0}{\mu_n a_n}=b_0 \nu_{n-1},\qd \hat\nu_n:=\frac{1}{\hat\mu_n\hat a_n}
=\frac{1}{b_0}{\mu_{n-1}}, \qqd 1\le n < N'+1. \tag 5.3$$
Actually, the rates $\big({\hat a}_i, {\hat b}_i\big)$ in (5.1) are determined by the transform
given in (5.3): ${\hat \mu}_n=b_0 \nu_{n-1}$
and ${\hat\nu}_n=\mu_{n-1}/b_0$. From this, it follows that
$$\gather \mu_n=b_0 {\hat\nu_{n+1}}={\hat a_1} \hat \nu_{n+1},\qd
\nu_n=\frac{1}{b_0}\hat\mu_{n+1}=\frac{1}{\hat a_1} \hat\mu_{n+1}, \qqd 0\le n< N',\\
\mu_N={\hat a}_1\big({\hat\mu}_N {\hat b}_N\big)^{-1}\qqd \text{if }\; N<\infty \text{ and } b_N=0,
\tag 5.4\endgather$$
and so
$$\aligned
&\sum_{n= 1}^{N'}\frac{1}{{\hat \mu}_n {\hat a}_n}=\sum_{n= 1}^{N'}{\hat \nu}_n
\!=\frac{1}{b_0}\sum_{n= 0}^{N'-1}\mu_{n},\\
&\sum_{n= 1}^{N'} {\hat\mu}_n=b_0\sum_{n= 1}^{N'} {\nu}_{n-1}
=b_0\sum_{n= 0}^{N'-1} \frac{1}{\mu_n b_n}. \endaligned\tag 5.5$$
Note that by (5.1),
$$a_{i+1}+ b_i-\frac{a_i}{v_{i-1}} - b_{i+1} v_{i}
={\hat b_{i+1}}+{\hat a_{i+1}}-\frac{{\hat b_i}}{v_{i-1}}-{\hat a_{i+2}}\, v_i.$$
By a change of the variables $(v_i)\in {\scr V}$:
$$v_i=\frac{{\hat b_{i+1}}}{{\hat a_{i+2}}}\, \hat v_{i+1},
\tag 5.6$$ or
$${\hat v}_i=\frac{{\hat a}_{i+1}}{{\hat b}_i}\, v_{i-1}
=\frac{b_i}{a_i}\, v_{i-1},\tag 5.7
$$ we get
$$a_{i+1}+ b_i-\frac{a_i}{v_{i-1}} - b_{i+1} v_{i}
={\hat a_{i+1}}\bigg(1-\frac{1}{{\hat v_{i}}}\bigg)+
{{\hat b_{i+1}}}(1-{{\hat v}_{i+1}}).$$
Since $b_0>0$, $v_{-1}>0$ but $a_0=0$, from (5.7), it is clear that we should set
$\hat v_0=\infty$. Next, by (5.7) again,
$$v_{i-1}> \frac{a_i}{b_i}\Longleftrightarrow {\hat v}_i>1.$$
It remains to examine the boundary condition on the right-hand side when $N<\infty$.
\roster
\item First, let $b_N=0$. Then $v=(v_i>0: 0\le i<N-1)$, $v_{-1}$ and $v_{N-1}$ are free.
The dual state space is ${\widehat E}=\{1,2, \ldots, N\}$. The dual test sequence is
${\hat v}=\big({\hat v}_i>0: 1\le i< N\big)$, ${\hat v}_N=0$.
\item Next, let $b_N>0$. Then $v=(v_i>0: 0\le i<N)$, $v_{-1}$ and $v_{N}$ are free.
The dual state space is ${\widehat E}=\{1,2, \ldots, N+1\}$ with reflecting at $N+1$.
Hence, ${\hat v}=\big({\hat v}_i>0: 1\le i< N+1\big)$,
${\hat v}_{N+1}=0$.
\endroster

We have thus proved the following result.

\proclaim{\prp\;5.1} For the dual processes defined above, the following identities hold:
$$\align
&\sup_{v}\inf_{0\le i<N'}\bigg[a_{i+1}+ b_i-\frac{a_i}{v_{i-1}} - b_{i+1} v_{i}\bigg]\\
&\qd =\sup_{{\hat v}}\inf_{1\le i<N'+1}\bigg[{\hat a_{i}}\bigg(1 -
\frac{1}{{\hat v}_{i-1}}\bigg)+ {\hat b_{i}}(1-{\hat v}_{i})\bigg],\tag 5.8
\endalign$$
where $v=(v_i>0: 0\le i<N'-1)$ with free $v_{-1}$, and ${\hat v}=\big({\hat v}_i>0: 1\le i< N'\big)$ with
${\hat v}_0=\infty$, $v_{N'-1}$ is free and ${\hat v}_{N'}=0$ if $N<\infty$;
$$\align
&\sup_{v\in {\scr V}_*}\inf_{0\le i<N+1}\bigg[a_{i+1}+ b_i-\frac{a_i}{v_{i-1}} - b_{i+1} v_{i}\bigg]\\
&\qd =\sup_{{\hat v}\in {\scr V}_1}\inf_{1\le i<N+2}\bigg[{\hat a_{i}}\bigg(1 -
\frac{1}{{\hat v}_{i-1}}\bigg)+ {\hat b_{i}}(1-{\hat v}_{i})\bigg]\tag 5.9\endalign$$
in the case that $b_N>0$ if $N<\infty$,
where ${\scr V}_*$ is given in \prp\;2.7, and ${\scr V}_1$ is defined in \thm\;4.1
replacing $N$ by $N+1$ when $N<\infty$.

In these formulas, $a_{N+1}=b_{N+1}=0$ if $N<\infty$ by convention.
\endproclaim

\proclaim{\crl\;5.2} Given rates $(a_i, b_i)$ as in Section 2 (then $b_N>0$ if $N<\infty$),
let $\lz_0=\lz_0^{(2.2)}$ and define $\dz$ by (3.1). Next,
define the dual rates $\big(\hat a_i, \hat b_i\big)$ as above.
Correspondingly, we have $\hat\lz_0$ and $\hat\dz$ defined by (4.1)
and (4.4) replacing $N$ by $N+1$ if $N<\infty$, respectively, in terms of the dual rates. Then we have
$\lz_0=\hat\lz_0$ and $\dz=\hat\dz$.
\endproclaim

\demo{\prf} Having relationship (5.9) at hand, the assertion
that $\lz_0=\hat\lz_0$ follows by a combination part (2) of
\prp\;2.7 and part (1) of Theorem 4.1, provided $\sum_i \hat\mu_i
<\infty$.

Next, by (4.4), (5.3), and (3.1), we have
$$\align
{\hat\dz} &=\sup_{1\le i <N+2} \sum_{j=1}^i \frac{1}{{\hat\mu_{j}} {\hat
a_{j}}}\sum_{j= i}^{N+1}{\hat\mu_{j}}\\
&=\sup_{1\le i <N+2} \sum_{j=1}^i \hat\nu_j\sum_{k= i}^{N+1}{\hat\mu_k}\\
&=\sup_{1\le i <N+2} \sum_{j=1}^i \frac{{\mu_{j-1}}}{b_0}\sum_{k= i}^{N+1}
b_0 \nu_{k-1}\\
&=\sup_{0\le i <N+1} \sum_{j=0}^i \mu_{j} \sum_{k= i}^N\frac{1}{\mu_{k} b_{k}}\\
&=\dz.
\endalign$$
This proves that $\dz=\hat\dz$. In particular, if $\sum_i \hat\mu_i
\big(=\sum_i \nu_i\big) =\infty$, then by \thm\;3.1 and \crl\;4.4, we get
$\lz_0=\hat\lz_0=0$. We have thus completed the proof of $\lz_0=\hat\lz_0$.\qed
\enddemo

As will be seen in \thm\;7.1\,(2), in the degenerated case that $\sum_i \mu_i=\infty$
and $\sum_i (\mu_i b_i)^{-1}=\infty$, the dual of the process studied in Section 2 also goes
to the one studied in Section 7.

Before moving further, let us discuss the duality used here. Very recently, Chi Zhang
provides us a nice explanation which leads to a deeper understanding of the duality (5.1).
Consider a simple example as follows:
$${Q=\pmatrix
 -b_0 & b_0 & 0 & 0 \\
 a_1 & -a_1-b_1 & b_1 & 0 \\
 0 & a_2 & -a_2-b_2 & b_2 \\
 0 & 0 & a_3 & -a_3-b_3
\endpmatrix},\qqd a_i, b_i>0.$$
Introduce an invertible matrix:
$${M\!=\!\pmatrix
\mu_0 b_0 & -\mu_0 b_0 & 0 & 0\\
0 & \mu_1 b_1 & -\mu_1 b_1 & 0\\
0 & 0 & \mu_2 b_2 & -\mu_2 b_2\\
0 & 0 & 0 & \mu_3 b_3
\endpmatrix}\Longrightarrow
{M^{-1}\!=\!\pmatrix
\frac{1}{\mu_0 b_0} & \frac{1}{\mu_1 b_1} & \frac{1}{\mu_2 b_2} & \frac{1}{\mu_3 b_3}\\
0 & \frac{1}{\mu_1 b_1} & \frac{1}{\mu_2 b_2} & \frac{1}{\mu_3 b_3}\\
0 & 0 & \frac{1}{\mu_2 b_2} & \frac{1}{\mu_3 b_3}\\
0 & 0 & 0 & \frac{1}{\mu_3 b_3}
\endpmatrix}\!.$$
Then
$$\align
MQM^{-1}&=
{\pmatrix
 -a_1-b_0 & a_1 & 0 & 0 \\
 b_1 & -a_2-b_1 & a_2 & 0 \\
 0 & b_2 & -a_3-b_2 & a_3 \\
 0 & 0 & b_3 & -b_3
\endpmatrix}\\
&= {\pmatrix
-{\hat a}_1-{\hat b}_1 & {\hat b}_1 & 0 & 0 \\
 {\hat a}_2 & -{\hat a}_2-{\hat b}_2 & {\hat b}_2 & 0 \\
 0 & {\hat a}_3 & -{\hat a}_3-{\hat b}_3 & {\hat b}_3 \\
 0 & 0 & {\hat a}_4 & -{\hat a}_4
\endpmatrix}\\
&={\widehat Q}.\endalign$$
Hence, the dual matrix ${\widehat Q}$ is just the classical similar transformation of $Q$ and so they
have the same spectrum. In particular, the eigenequation $Q g=-\lz_0 g$ ($g\ne 0$)
is transferred into
$${\widehat Q} (M g)=(MQ M^{-1})(Mg)=\lz_0 M g={\hat\lz}_0 \big(M g\big).$$
Hence, the eigenfunction $g$ of $\lz_0$ is transformed to ${\hat g}=M g$ of ${\hat\lz}_0=\lz_0$.
Correspondingly, the test function $f$ is transformed to ${\hat f}=M f$.
From this, it should be clear that all the operators $R$ and $\widehat R$,
$I$ and $\hat I$, $I\!I$ and $\widehat{I\!I}$ are
closely related to each other and then so are the variational formulas.

Having these facts at hand, one can simplify a part of the previous proofs.
However, we prefer to keep all the details here since they are needed when
we go to the more general situation, so called the Poincar\'e-type inequalities (Section 8), or
can be used as a reference for studying the continuous case. For the Poincar\'e-type inequalities,
the current duality seems not available.

By \crl\;5.2, we have two ways to estimate $\lz_0=\hat\lz_0$:
using either the rates $(a_i, b_i)$ or $\big(\hat a_i,\hat b_i\big)$. The corresponding
formulas for $\dz_1'$, $\hat\dz_1'$, $\dz_1$ and $\hat\dz_1$ are collected in Tables 5.1 and
5.2.

\medskip

\centerline{{\bf Table 5.1}: Expressions of $\dz=\hat\dz$, $\dz_1'$, $\hat\dz_1'$, $\dz_1$ and $\hat\dz_1$
in terms of the rates $(b_i, a_i)$:}
\nopagebreak
\vskip-2em
$$\align
&\text{\hskip-5em}\dz=\hat\dz=\sup_{0\le i<N+1}\mu[0, i]\, \nu[i, N]
=\sup_{0\le i<N+1}\sum_{j=0}^i\mu_j \sum_{k=i}^N \nu_k,\tag 5.10\\
&\text{\hskip-5em}\dz_1'=\sup_{0\le i<N\!+\!1}\! \frac{1}{\nu[i,N]}\!\sum_{k=0}^N \mu_k \nu[k\!\vee i,N]^2\\
&\text{\hskip-4em}=\sup_{0\le i<N\!+\!1}\!\!\bigg[\mu[0,i]{\nu[i,N]}
+\frac{1}{\nu[i,N]}\sum_{k=i +1}^N \mu_k \nu[k,N]^2\bigg],\tag 5.11\\
&\text{\hskip-5em}{\hat\dz_1}'= \sup_{0\le i<N+1} \frac{1}{\mu[0, i]}\!\sum_{k=0}^N
\nu_k \mu[0, k\!\wedge i]^2\\
&\text{\hskip-4em}=\sup_{0\le i<N+1}\!\bigg[\mu[0,i]{\nu[i,N]}\!+\!\frac{1}{\mu[0, i]}\!\sum_{k= 0}^{i-1}
\nu_k \mu[0, k]^2\bigg],\tag 5.12\\
&\text{\hskip-5em}\dz_1=\sup_{0\le i<N+1}\frac{1}{\sqrt{\nu[i,N]}}\sum_{k=0}^N \mu_k
\nu[i\vee k,N]\sqrt{\nu[k,N]}\\
&\text{\hskip-4em}=\!\sup_{0\le i<N+1}\!\bigg[\sqrt{\nu[i,N]}\,\sum_{k=0}^{i}\mu_k\sqrt{\nu[k,N]}+\!
\frac{1}{\sqrt{\nu[i,N]}}\sum_{k=i+1}^N\!\!\mu_k \nu[k,N]^{3/2}\bigg]\!,
\text{\hskip-2em} \tag 5.13\\
&\text{\hskip-5em}{\hat\dz_1}= \sup_{0\le i<N+1} \frac{1}{\sqrt{\mu[0, i]}}\sum_{k=0}^N \nu_k
\mu[0, k\wedge i] \sqrt{\mu[0, k]}\\
&\text{\hskip-4em}=\sup_{0\le i<N+1}\bigg[\frac{1}{\sqrt{\mu[0, i]}}\sum_{k= 0}^{i-1}
\nu_k \mu[0, k]^{3/2}+ \sqrt{\mu[0, i]}\, \sum_{k=i}^N \nu_k \sqrt{\mu[0, k]}\bigg].
\text{\hskip-1em}\tag 5.14\\
\endalign
$$
\newpage

\centerline{{\bf Table 5.2}: Expressions of $\dz=\hat\dz$, $\dz_1'$,
$\hat\dz_1'$, $\dz_1$ and $\hat\dz_1$ in terms of the rates $(\hat
b_i,\hat a_i)$:}
\nopagebreak
\vskip-2em
$$\align
&\text{\hskip-5em}\dz=\hat\dz=\sup_{1\le i<N+1} \hat\nu [1, i]\,\hat\mu[i, N]
  =\sup_{1\le i<N+1} \sum_{k=1}^i {\hat\nu_{k}} \sum_{j= i}^N {\hat\mu_{j}},\tag 5.15\\
&\text{\hskip-5em}\dz_1'=\sup_{1\le i<N+1} \frac{1}{\hat\mu[i,N]}\sum_{k= 1}^N\hat \nu_k\hat \mu[k\vee i,N]^2\\
&\text{\hskip-4em}=\sup_{1\le i<N+1}\!\bigg[{\hat\mu[i,N]}\hat\nu[1,i]
+\frac{1}{\hat\mu[i,N]}\!\sum_{k=i +1}^N\hat \nu_k\hat \mu[k,N]^2\!\bigg], \tag 5.16\\
&\text{\hskip-5em}\hat\dz_1'=\sup_{1\le i<N+1}
\frac{1}{\hat\nu[1,i]}\sum_{k=1}^N\hat \mu_k {\hat\nu}[1,k\!\wedge i]^2\\
&\text{\hskip-4em}=\sup_{1\le
i<N+1}\bigg[{\hat\mu[i,N]}\hat\nu[1,i]+\frac{1}{\hat\nu[1,i]}\sum_{k=1}^{i-1}\hat
\mu_k {\hat\nu}[1,k]^2
\bigg]. \tag 5.17\\
&\text{\hskip-5em}\dz_1=\sup_{1\le i<N+1} \frac{1}{\sqrt{\hat\mu[i,N]}}\sum_{k= 1}^N\hat \nu_k\hat \mu[k\vee i,N]\sqrt{\hat\mu[k,N]}\\
&\text{\hskip-4em}=\!\sup_{1\le i<N+1}\!\bigg[\sqrt{\hat\mu[i,N]}\,\sum_{k=1}^{i}\hat \nu_k
\sqrt{\hat\mu[k,N]}
+\!\frac{1}{\sqrt{\hat\mu[i,N]}}\sum_{k=i +1}^N\!\!\!{\hat \nu}_k\,{\hat \mu}[k,\!N]^{3/2}\bigg]\!,
\text{\hskip-2em} \tag 5.18\\
&\text{\hskip-5em}\hat\dz_1=\sup_{1\le i<N+1} \frac{1}{\sqrt{\hat\nu[1,i]}}\sum_{k=
1}^N \hat\mu_k \hat
\nu[1,k\wedge i]\sqrt{\hat\nu[1,k]}\\
&\text{\hskip-4em}=\sup_{1\le
i<N+1}\bigg[\frac{1}{\sqrt{\hat\nu[1,i]}}\sum_{k=1}^{i-1} \hat\mu_k
\hat\nu[1,k]^{3/2}+ \sqrt{\hat\nu[1,i]}\,\sum_{k=i}^N \hat\mu_k\sqrt{
\hat\nu[1,k]}\, \bigg],\text{\hskip-1em}\tag 5.19
\endalign
$$

The next four examples are dual of Examples 3.4--3.7, respectively.

\proclaim{\xmp\;5.3} For Example 3.4, we have ${\hat a}_i\equiv b\,(i\ge 1)$,
${\hat b}_i\equiv a\,(a>0)$, $b\ge a$. Then ${\hat \dz}=\dz=b (a-b)^{-2}$,
$\hat\dz_1'=\dz_1'=(a+b)(a-b)^{-2}$, and
$\hat\dz_1=\dz_1=\lz_0^{-1}=\big(\sqrt{a}-\sqrt{b}\,\big)^{-2}$. In
particular, if we take $\hat a_i=4$ and $\hat b_i=1\,(i\ge 1)$, then
$\hat\lz_0=1$,
$$\align
&\hat\dz_1'=5/9= 0.{\dot 5}, \qd\hat\dz_2'= 0.{6\dot 4}, \qd\hat\dz_3' \approx 0.71,
\qd\hat\dz_4' \approx 0.755,\qd\hat\dz_5' \approx 0.79;\\
&\bar{\hat\dz}_1=0.{\dot 5}, \qd\bar{\hat\dz}_2\approx 0.71,
\qd\bar{\hat\dz}_3\approx 0.79, \qd\bar{\hat\dz}_4\approx 0.835
\qd\bar{\hat\dz}_5\approx 0.8647.\endalign$$ Thus, $\hat\dz_n'$ and
$\bar{\hat\dz}_n$ are increasing and close to $\hat\lz_0^{-1}$ as
$n\uparrow$.
\endproclaim

\demo{\prf} To compute $\hat\dz_1'$ and $\hat\dz_1$, we
use Table 5.1. For simplicity, write $\gz=b/a>1$.
Then
$$\mu_k=\gz^k,\qqd \mu[0, i]=\frac{\gz^{i+1}-1}{\gz-1},\qqd \nu_k=\frac{1}{b} \gz^{-k}.$$

(a) Note that
$$\align
&\frac{1}{\mu[0, i]}\sum_{k= 0}^{i-1}
\nu_k \mu[0, k]^2+ \mu[0, i] \sum_{k=i}^\infty \nu_k\\
&\qd=\frac{1}{b}\bigg[\frac{\gz-1}{\gz^{i+1}-1}
\sum_{k=0}^{i-1}\gz^{-k} \bigg(\frac{\gz^{k+1}-1}{\gz-1}\bigg)^2
+\frac{\gz^{i+1}-1}{\gz-1}\sum_{k\ge i}\gz^{-k}\bigg]\\
&\qd=\frac{1}{b(\gz-1)}\bigg[\frac{1}{\gz^{i+1}-1}
\sum_{k=0}^{i-1}\gz^{-k} \big(\gz^{k+1}-1\big)^2
+\big(\gz^{i+1}-1\big)\sum_{k\ge i}\gz^{-k}\bigg]\\
&\qd=\frac{1}{b(\gz-1)}\bigg[\frac{\gz(1 + \gz)}{\gz-1}-\frac{2 (i+1) \gz}{\gz^{i+1}-1}\bigg].
\endalign$$
Since the second term in the last $[\cdots]$ is negative and $\gz>1$,
the right-hand side attains its supremum at $i=\infty$. By (5.12),
we have thus obtained
$${\hat\dz_1}'=\frac{\gz(1+\gz)}{b(\gz-1)^2}=\frac{a+b}{(a-b)^2}.$$

(b) Next, note that
$$\align
&\frac{1}{\sqrt{\mu[0, i]}}\sum_{k= 0}^{i-1}
\nu_k \mu[0, k]^{3/2}+ \sqrt{\mu[0, i]} \sum_{k=i}^\infty \nu_k \sqrt{\mu[0, k]}\\
&\qd=\frac{1}{b(\gz-1)}\bigg[\frac{1}{\sqrt{\gz^{i+1}\!-1}}
\sum_{k=0}^{i-1}\gz^{-k} \big(\gz^{k+1}\!-1\big)^{3/2}
\!+\!\sqrt{\gz^{i+1}\!-1}\, \sum_{k\ge i}\gz^{-k}\sqrt{r^{k+1}\!-1}\,\bigg]\\
&\qd\le \frac{1}{b(\gz-1)}\bigg[\frac{1}{\sqrt{\gz^{i+1}-1}}
\sum_{k=0}^{i-1}\gz^{(k+3)/2}
 +\sqrt{\gz^{i+1}-1}\, \sum_{k\ge i}\gz^{-k/2+1/2}\,\bigg]\\
&\qd= \frac{1}{b(\gz-1)}\bigg[\frac{1}{\sqrt{\gz^{i+1}-1}}
\frac{\gz^{3/2}(\gz^{i/2}-1)}{\sqrt{\gz}-1}
 +\frac{\gz^{-i/2+1/2}\sqrt{\gz^{i+1}-1}}{1-1/\sqrt{\gz}} \bigg]\\
&\qd\le \frac{1}{b(\gz-1)}\bigg[\frac{\gz}{\sqrt{\gz}-1}
 +\frac{\gz\sqrt{\gz}}{\sqrt{\gz}-1} \bigg]\\
&\qd=\frac{\gz}{b(\sqrt{\gz}-1)^2}\\
&\qd=\frac{1}{(\sqrt{a}-\sqrt{b})^2}.
\endalign$$
By (5.14), this means that $\hat\dz_1\le \hat\lz_0^{-1}$ and so the equality sign
must hold because $\hat\dz_1^{-1}$ is a lower estimate: $\hat\lz_0\ge \hat\dz_1^{-1}$.

(c) We now compute the approximating sequences $\big\{{\hat\dz}_n'\big\}$
and $\big\{\bar{\hat\dz}_n\big\}$ for the upper estimate,
using the dual rate $\big(\hat a_i,\hat b_i\big)$. In the particular case, we have
$$\hat\mu_i=4^{1-i},\qd \hat\nu_i= 4^{i-2},\qd \hat\fz_i=\hat\nu[1,i]=\frac{4^i-1}{12}.$$
The approximating sequences can be computed successively by using the following formulas:
$$\align
f_1^{(m)}(i)&=\frac{4^i-1}{12},\qqd i \in \{1,2,\ldots, m\},\\
f_n^{(m)}(i)&=\frac{1}{3} \bigg\{ \sum_{k=1}^{i-1} (1-4^{-k}) f_{n-1}^{(m)}(k)
+(4^i-1) \sum_{k=i}^{m-1} 4^{-k} f_{n-1}^{(m)}(k)\\
&\qqd\; + \frac 1 3 (4^i-1) 4^{1-m} f_{n-1}^{(m)}(m)\bigg\},\qqd i \in \{1,2,\ldots, m\},\; n\ge 2.
\endalign$$
Then $\hat\dz_n'=\sup_{m\ge 1}\min_{1\le i\le m} f_{n+1}^{(m)}(i)\big/f_{n}^{(m)}(i)$.
For the first five of  $\{\hat\dz_n'\}$, the minimum are all attained at $m$ and so the
computations become easier.

To compute $\bar{\hat\dz}_n$, simply use the formula
$$\bar{\hat\dz}_n=\sup_{m\ge 1}
\frac{\sum_{i=1}^m 4^{1-i} f_n^{(m)}(i)^2 + 3^{-1} 4^{1-m}
f_n^{(m)}(m)^2}{ \sum_{i=1}^m 4^{2-i}
\big(f_n^{(m)}(i)-f_n^{(m)}(i-1)\big)^2},\qqd f_n^{(m)}(0):=0.
\qed$$\enddemo

\proclaim{\xmp\;5.4} For Example 3.5 with $\gz=1$
$\big({\hat b}_i=i,\; {\hat a}_i=2\, i\big)$, we have
${\hat\dz_1}'\approx 0.75< \dz_1'\approx 0.84$ and
${\hat\dz_1}\approx 1.12> {\dz_1}\approx 1.09$. Besides, ${\hat\dz_1}/{\hat\dz_1}'\approx 1.5$.
\endproclaim

\proclaim{\xmp\;5.5} For Example 3.6, we have
${\hat a}_i={\hat b}_i=i^2\,(i\ge 1)$, ${\hat b}_0=0$,
${\hat\dz}_1'=2< \dz_1'\approx 2.19$ and
${\hat\dz}_1= {\dz_1}=4$ which is sharp. Besides, ${\hat\dz}_1/{\hat\dz}_1'=2$.
\endproclaim

\demo{\prf} By \xmp\;3.6 and \crl\;5.2, it follows that
 ${\hat \lz}_0=\lz_0= 1/4$. Here, we present an easier proof for the upper
 estimate.
Note that when ${\hat a}_i={\hat b}_i$ for $i\ge 2$, we have
$${\hat \mu}_1=1,\;\;{\hat \mu}_i=\frac{{\hat b}_1\cdots {\hat b}_{i-1}}{{\hat a}_2\cdots {\hat a}_i}
=\frac{{\hat b}_1}{{\hat a}_i},\;\;i\ge 2;\qqd
{\hat \mu}_i{\hat  b}_i ={\hat b}_1,\;\;i\ge 1. \tag 5.20$$
In the present case, we have
${\hat \mu}_i=i^{-2}\,(i\ge 1)$ and ${\hat \mu}_i {\hat a}_i\equiv 1$.
Let $f_i^{(m)}=\sqrt{i\wedge m}$. Then
$$\align
{\hat \mu}\big(f^{(m)\,2}\big)&= \sum_{i=1}^{m} \frac 1 i + m\sum_{i\ge m+1} \frac{1}{i^2},\\
{\widehat D}\big(f^{(m)}\big)&=\sum_{i=1}^m \big(\sqrt{i}-\sqrt{i-1}\,\big)^2
=\sum_{i=1}^m \frac{1}{\big(\sqrt{i}+\sqrt{i-1}\,\big)^2}
\le 1+\frac 1 4 \sum_{i=1}^{m-1} \frac 1 i.
\endalign
$$
Hence,
$${\hat\lz}_0\le \varliminf_{m\to\infty}\frac{{\widehat D}\big(f^{(m)}\big)}
{{\hat \mu}\big(f^{(m)\,2}\big)}=\frac{1}{4}.
\qed$$\enddemo

\proclaim{\xmp\;5.6} For Example 3.7, we have
${\hat a}_i=i^4\,(i\ge 1)$, ${\hat b}_i=i(i-1/2)(i^2+3i+3)$,
${\hat\lz}_0=\lz_0=1/2$,
${\hat\dz}_1'\approx 1.83 < \dz_1'\approx 1.9$ and
${\hat\dz}_1\approx {\dz_1}\approx 2$. Besides, ${\hat\dz}_1/{\hat\dz}_1'=1.09$.
\endproclaim

\demo{\prf} First, we have
$${\hat\mu}_i=\frac{\prod_{k=1}^{i-1}(k-1/2)(k^2+3 k +3)}{i\, i!^3},\qd
{\hat\nu_i}=\frac{(i-1)!^3}{\prod_{k=1}^{i-1}(k-1/2)(k^2+3 k +3)},\qd i\ge 1.
$$
By (5.5) and Example 3.7, we have $\sum_i {\hat\mu}_i<\infty$ and $\sum_i {\hat\nu}_i<\infty$,
and so the minimal dual process is explosive (but here we are dealing with the maximal one).
The sharp lower bound can be deduced from part (1)
of \thm\;4.1 with the dual test sequence
$${\hat v}_i=1+\frac{1}{i (i^2+3 i +3)},\qqd i\ge 1.$$
From this, it follows that the corresponding eigenfunction
$${\hat g}_i=\prod_{k=1}^{i-1}{\hat v}_k,\qd i\ge 2, \qd  {\hat g}_1=1,$$
increases strictly to a finite limit since $\sum_{i\ge 1} i^{-1} (i^2+3 i +3)^{-1}<\infty$.
The sequence $({\hat v}_i)$ comes from the one computed in Example 3.7 plus a use of (2.35) and (5.7).\qed\enddemo

The precise value of $\lz_0$ for the next example is unknown. Its eigenfunction is
non-polynomial. It is interesting to compare this example with the ergodic one given
in \S 6 for which $\lz_1=2$, as well as the one with rates $a_i=i+1$ and
$b_i=i^2\,(i\ge 1)$ given in \S 7 for which $\lz_0=2$.

\proclaim{\xmp\;5.7} Let $\hat b_0=0$,
$\hat b_i=i+2\,(i\ge 1)$ and $\hat a_i=i^2$. It is the dual of the process
studied in \S 2 with rates $a_i=i+2\,(i\ge 1)$ and $b_i=(i+1)^2\,(i\ge 0)$.
Then ${\hat\lz}_0\in (0.395, 0.399)$,
${\hat\dz_1}'\approx 2.37< \dz_1'\approx 2.48$ and
${\hat\dz_1}\approx 2.63 > {\dz_1} \approx 2.61$.
Besides, ${\hat\dz_1}/{\hat\dz_1}'\approx 1.1$.
\endproclaim

It is interesting that for all of \xmp s 5.3--5.7, we have
${\hat\dz_1}'\le \dz_1'$ and ${\hat\dz_1}\ge {\dz_1} $ which then means that
\crl\;3.3 is more effective than \crl\;4.4.
The effectiveness of the bounds $\dz_1$ and $\dz_1'$ given in
\crl\;4.4 was also checked by Sirl, Zhang and Pollett (2007) for
some models from practice.

\proclaim{\rmk\;5.8}{\rm It is now a suitable position to mention a method for the
numerical computation of $\lz_0$
defined in \S 4. The idea is meaningful in the other cases.
From proof (b) of \thm\;4.1, it follows that there is a sequence $(v_i: v_i>1, 1\le i<N)$ such that
$$R_i(v)=a_i(1-v_{i-1}^{-1})+b_i(1- v_i)=\lz_0,\qqd v_0=\infty,\;v_N=0\text{ if }N<\infty.$$
Hence, we have
$$\cases
v_1-1=(a_1-\lz_0) b_1^{-1},\\
v_i-1=\big[ a_i (1-v_{i-1}^{-1})-\lz_0\big] b_i^{-1},\qd 2\le i <N.
\endcases\tag 5.21$$
In other words, replacing $v_i-1$ by $u_i$, when $z=\lz_0$, the equation
$$\cases
u_1=(a_1-z) b_1^{-1},\\
u_i=\big[{a_i u_{i-1}} {(1+u_{i-1})^{-1}}-z\big] b_i^{-1},\qd 2\le i<N,
\endcases\tag 5.22$$
has a positive solution $(u_i=u_i(z))_{1\le i<N}$. Thus, one may use the maximal $z$
so that (5.22) has a positive solution as an approximation of $\lz_0$ (based on part (1)
of \thm\;4.1). In this
way, we obtain the approximation of ${\hat\lz}_0$ given in \xmp\;5.7.}
\endproclaim

\head{6. Reflecting (Neumann) boundaries at origin and
infinity (ergodic case)}\endhead

We now turn to studying the first non-trivial eigenvalue in the ergodic case. Let $E=\{i:
0\le i <N+1\}\,(N\le \infty)$, $ b_0>0$, $b_N=0$ if $N<\infty$,
$$\lz_1=\inf\big\{D(f): \mu(f)=0,\; \mu(f^2)=1\big\}, \tag 6.1$$
where $\mu (f)=\int f \d \mu$,
$$D(f)=\sum_{0\le i<N} \mu_i b_i (f_{i+1}-f_i)^2=
\sum_{1\le i<N+1} \mu_i a_i (f_{i}-f_{i-1})^2 \tag 6.2$$
with domain ${\scr D}^{\max}(D)=\{f\in L^2(\mu): D(f)<\infty\}$. In (6.1), we presume that
$$\sum_{i=0}^N\mu_i<\infty.  \tag 6.3$$
Then the Dirichlet form $(D, {\scr D}^{\max}(D))$ has a trivial eigenvalue
$\lz_0=0$ with constant eigenfunction $\dbl$, and here we are
working on the next ``eigenvalue'' $\lz_1$ of $(D, {\scr D}^{\max}(D))$. If
(6.3) does not hold, then $\dbl\notin L^2(\mu)$ and so $\lz_1$
is not meaningful. Moreover, by (1.3) and \prp\;1.3, the
Dirichlet form is unique. In this case, the corresponding process is
explosive, or zero-recurrent, or transient. The decay rate is
described by $\lz_0$ which has already been treated in Sections 2
and 3. Hence, throughout this section, we assume (6.3).

Note that condition (6.3) plus (1.2) means that the unique process
is ergodic. When $N=\infty$ and (1.2) fails, the minimal process was treated in Sections 2
and 3, and in this section, we are dealing with the maximal process (cf.
[10;  \prp\;6.56]) as
in Section 4, it is indeed the unique honest reversible process.
Denote by $Q=(q_{ij})$ the birth--death $Q$-matrix. Then under
(6.3), the maximal process $P_{ij}^{\max}(\lz)$ (Laplace transform)
can be expressed as
$$P_{ij}^{\max}(\lz)=P_{ij}^{\min}(\lz)+\frac{{z}_i(\lz)\, \mu_j\, {z}_j (\lz)}{\lz \sum_k \mu_k\, {z}_k(\lz)},
\qqd i, j\in E,\qd \lz>0,$$
where for each fixed $j$, $\{P_{ij}^{\min}(\lz): i\in E\}$
is the minimal solution to the equations
$$x_i=\sum_{k\ne i}\frac{q_{ik}}{\lz+q_i} x_k+\frac{\dz_{ij}}{\lz+q_i},\qqd i\in E,$$
and $(z_i(\lz): i\in E)$ is the maximal solution to the equation
$$\cases
(\lz I -Q)u=0,\\
0\le u\le 1,
\endcases \qqd \lz>0$$
(cf. [10;  \prp\;6.56]). According to a result due to Z.K.
Wang (1964) (cf. Wang and Yang (1992, \S 6.8, \thm\;2)): if $N=\infty$ and (1.2)
fails, then every honest process (may be non-symmetric) is ergodic
and so is the maximal one. Certainly, within the symmetric context,
by using (1.4), it is
easy to check directly the ergodicity of the maximal process.

Here, we mention a technical point. If (6.3) fails, then as mentioned before, by
(1.3), there is precisely one symmetrizable process (Dirichlet form)
which is nothing but the minimal one. Thus, if (1.2) also fails,
then the unique process must be explosive and so there is no honest
symmetrizable process. This is a different point to the reversible
case (i.e., (6.3) holds) for which there exists exactly one honest
reversible process as just mentioned above.

We use the same notation $ I$, ${I\!I}$, ${\scr F}_I$, ${\scr F}_{I\!I}$,
${\widetilde{\scr F}}_I$ and ${\widetilde{\scr F}}_{I\!I}$ defined
in Section 4 with an addition ``$f_0=0$'' in the last four sets,
but redefine $R$ and ${\scr V}$ as follows:
$$\align
&R_i (v)= a_{i+1}+  b_i- a_i/v_{i-1} -  b_{i+1} v_{i},\qd 0\le i<N,\\
&\text{\hskip17em}\; v_{-1}\!>\!0\text{ is free and so is $v_{N-1}$ if }N\!<\!\infty,\\
&{\scr V}=\{v: v_i >0 \text{ for all }i: 0\le i<N-1\}.
\endalign$$
The local operator ${\widetilde R}$ is modified from $R$, replacing
$a_m$ by ${\tilde a}_m:= \mu_m a_m\big/\sum_{k=m}^N \mu_k$ for $v$
with $\supp(v)=\{0, 1, \ldots, m-2\}$ in the set
$$\align
&{\widetilde{\scr V}}=\bigcup_{m=2}^{N}\bigg\{v\!:
\frac{a_{i+1}}{a_{i+2}+b_{i+1}}<v_i<
\frac{a_{i+1}+b_i-a_i/v_{i-1}}{b_{i+1}}\;\;\text{for }i=0, 1, \ldots,  m-2\\
&\qquad\qquad\qquad   \text{ and }v_i=0 \text{ for }i\ge m-1\Big\}.
\endalign$$

\proclaim{\thm\;6.1} Under (6.3), the following variational formulas for $\lz_1$
hold. \roster
\item Difference form:
$$\inf_{v\in {\widetilde{\scr V}}}\,\sup_{0\le i<N}\,{\widetilde R}_i(v)
=\lz_1 =\sup_{v\in {\scr V}}\,\inf_{0\le i<N}\,R_i(v).$$
\item Summation form:
  $$\inf_{f\in {\widetilde{\scr F}}_I \cup {\widetilde{\scr F}}_I'}\,\sup_{1\le i\in E} I_i(\bar f)^{-1}
    =\lz_1=\sup_{f\in{\scr F}_I}\,\inf_{1\le i\in E} I_i(\bar f)^{-1},$$
    where
$$\gather
{\widetilde{\scr F}}_I'=\{f\in L^2(\mu): f_0=0, f \text{ is strictly increasing}\},\\
   \text{$\bar f=f - \pi(f)$,\qd $\pi=\mu/Z$\qd and \qd $Z={\tsize\sum_{i\in E}}\,\mu_i$.}\endgather$$
\endroster
\endproclaim

The use of $\bar f$ in the last line is based on the property
$\bar f = \overline{f+c}$ for every constant $c$ and so we can fix $f_0$ to be $0$.

\demo{\prf\; of \thm\;$6.1$} In the ergodic case under (1.2), the assertion
 $$\inf_{f\in {\widetilde{\scr F}}_I \cup {\widetilde{\scr F}}_I'}\,\sup_{1\le i\in E} I_i(\bar f)^{-1}
    \ge \lz_1$$
was proved in [7; \thm\;2.3] (but in the case that
$k=\infty$ in the original proof, one requires the
$L^2$-integrability condition included in $\widetilde{\scr F}_I'$,
as was pointed out in proof (c) of Theorem 4.1). The proof
remains the same in the present general situation with an obvious
modification when $N<\infty$. Next, in the ergodic case  under (1.2), the
following result
$$\lz_1=\sup_{v\in {\scr V}}\,\inf_{0\le i<N}\,R_i(v)
=\sup_{f\in{\scr F}_I}\,\inf_{1\le i\in E} I_i(\bar f)^{-1}  \tag
6.4$$ is just [3;  \thm\;1.1]. In the present general
situation, the proof for the second equality in (6.4) needs a slight
change only (cf. [3;  \lmm\;2.1]). To prove the first equality in (6.4), we claim that
$$\align
\lz_1&=  \inf\big\{D(f):\mu\big(\big|f-\pi(f)\big|^2 \big)=1,\;
f_i=f_{i\wedge m} \text{ for some }m\in E,\; m\ge 1\big\}\\
&=:{\tilde\lz}_1.\tag 6.5
\endalign$$
To see this, first it is clear that ${\tilde\lz}_1\ge \lz_1$.
Next, the proof of [4; \thm\;3.2] gives us
$${\lz}_1\ge \sup_{f\in{\scr F}_I}\,\inf_{1\le i\in E} I_i(\bar f)^{-1},$$
and furthermore, the equality sign with ${\lz}_1$ replaced by ${\tilde\lz}_1$ holds.
Once again, the
key point for the last statement is to show that the eigenfunction
of ${\tilde\lz}_1$ is strictly
increasing. For this, the original proof needs only a modification
replacing ${\lz}_1$ by ${\tilde\lz}_1$
(as indicated in proof (b) of \thm\;4.1). Therefore, (6.4) holds in the
present general situation.

Now, we need only to show that
\roster
\item"$\bullet$" $\inf_{f\in {\widetilde{\scr F}}_I}\,\sup_{1\le i\in E} I_i(\bar f)^{-1}
    \le\inf_{v\in {\widetilde{\scr V}}}\,\sup_{0\le i<N}\,{\widetilde R}_i(v)$, and
\item"$\bullet$" $\inf_{v\in {\widetilde{\scr V}}}\,\sup_{0\le i<N}\,{\widetilde R}_i(v)
    \le \lz_1$.
\endroster
\medskip

(a) Prove that $\inf_{f\in {\widetilde{\scr F}}_I}\,\sup_{1\le i\in
E} I_i(\bar f)^{-1}
    \le\inf_{v\in {\widetilde{\scr V}}}\,\sup_{0\le i<N}\,{\widetilde R}_i(v)$.
\medskip

As before, write ${\widetilde R}(u)$ instead of ${\widetilde R}(v)$.
Given $u$ with $\supp(u)=\{0,1,\ldots, m-1\}$ so that
$(v_i)\in{\widetilde{\scr V}}$, where $v_i=u_{i+1}/u_i>0$ for $i<
m-1$  and $v_i=0$ for $i\ge m-1$, let
$$f_i=
\cases
 a_i u_{i-1}- b_i u_i, \qd & i< m,\\
{\tilde a}_m u_{m-1}, & i\ge m.
\endcases
$$
Since the constraint in $\widetilde{\scr V}$ is equivalent to $\min_{i\le
m-1}{\widetilde R}_i (v)>0$, it is easy to check that
$$(f_{i+1}-f_i)/u_i={\widetilde R}_i(u)>0
\text{ for } i< m \text{ and } f_i=f_{i\wedge m},$$ and so $f+b_0 u_0\in
{\widetilde{\scr F}}_I$. Moreover, since
$$\align
\sum_{k=i}^N\mu_k f_k&=\sum_{k= i}^{m-1}\mu_k f_k+f_{m}\sum_{j=m}^N\mu_j\\
&=\sum_{k= i}^{m-1}\mu_k( a_k u_{k-1}- b_k u_k)+f_{m}\sum_{j=m}^N\mu_j\\
&=\mu_i  a_i u_{i-1}-\mu_m a_m u_{m-1}+{\tilde
a}_{m}u_{m-1}\sum_{j=m}^N\mu_j\\
&=\mu_i  a_i u_{i-1}, \qqd i\le m-1,
\endalign$$
we get
$$\mu (f)=\sum_{k=0}^N\mu_k f_k=\mu_0  a_0 u_{-1}=0,$$
and so
$$\align
&\sum_{k=i}^N\mu_k\bar f_k= \mu_i  a_i u_{i-1},\qqd i\le m-1,\\
&\sum_{k=m}^N\mu_k\bar f_k= f_{m}\sum_{k=m}^N\mu_k={\tilde a}_m
u_{m-1}\sum_{k=m}^N\mu_k=\mu_m a_m u_{m-1}.
\endalign$$ It follows that
$$u_{i-1}=\frac{1}{\mu_i  a_i} \sum_{k=i}^N\mu_k \bar f_k, \qqd i\le m.$$
Hence,
$${\widetilde R}_{i-1}(u)=\frac{f_i-f_{i-1}}{u_{i-1}}=  I_i(\bar f)^{-1}, \qqd 1\le i\le m.$$
Therefore, we have
$$\max_{0\le i< m}\,{\widetilde R}_i(u)\!=\! \max_{1\le i\le m}\, I_i(\bar f)^{-1}
\!\ge\!\! \inf_{f\in {\widetilde{\scr F}}_I,\;f_i=f_{i\wedge m}}\, \max_{1\le i\le m}\, I_i(\bar f)^{-1}
\!\ge\! \inf_{f\in {\widetilde{\scr F}}_I}\, \sup_{1\le i\in E}\,
I_i(\bar f)^{-1}\!\!,$$ and then
$$\inf_{v\in {\widetilde{\scr V}}}\,\sup_{0\le i<N}\,{\widetilde R}_i(v)\ge
\inf_{f\in {\widetilde{\scr F}}_I}\, \sup_{1\le i\in E}\, I_i(\bar
f)^{-1}. $$ Here, we have used the fact that ${\widetilde R}_i(v)=-\infty$ for
$i\ge m-1$ if $\supp(v)=\{0, 1, \ldots, m-2\}$ and in the last step,
we have returned to the original notation ${\widetilde R}(v)$
instead of ${\widetilde R}(u)$.
\medskip

(b) Prove that $\inf_{v\in {\widetilde{\scr V}}}\,\sup_{0\le
i<N}\,{\widetilde R}_i(v) \le \lz_1$.
\medskip

Because of
$$\big\{ \mu(f)=0,\; \mu\big(f^2\big)=1,\;
f_i=f_{i\wedge m}\big\} \subset
\big\{ \mu(f)=0,\; \mu\big(f^2\big)=1,\; f_i=f_{i\wedge
(m+1)}\big\},$$ by (6.5), it is clear that
$$\lz_1^{(m)}:=\inf\big\{D(f):\mu(f)=0,\; \mu\big(f^2\big)=1,\;
f_i=f_{i\wedge m}\big\} \;\downarrow \lz_1\; \text{ as } m\uparrow
N.$$ Actually, this is a special case of an approximation result
given in [2; \thm\;4.2 and \crl\;4.3] or [10; \thm\;9.20
and \crl\;9.21]. Note that $\lz_1^{(m)}$ is just the first
non-trivial eigenvalue of the local Dirichlet form $\big({\widetilde
D},$ $ {\scr D}\big({\widetilde D}\big)\big)$ defined by (4.17)
replacing the Dirichlet boundary at $0$ by the Neumann one (having
the state space $\{0, 1, \ldots, m\}$), with Neumann (reflecting)
boundary at $m$. Denote by $g$ the first eigenfunction of
$\lz_1^{(m)}$ and extend it to the whole space by setting
$g_i=g_{i\wedge m}$. Now, if we set $u_i=g_{i+1}-g_i$ for $i\in E$,
then $u_i>0$ for $i\le m-1$, $u_i=0$ for $i\ge m$, and furthermore,
$${\widetilde R}_i(u)=\lz_1^{(m)}>0 \qqd\text{ for all } i\le m-1.$$
Moreover, by the definition of $g$, we have
$$\cases
 b_i u_i-  a_i u_{i-1}=-\lz_1^{(m)} g_i,\qqd i\le m-1,\\
 {\tilde a}_m u_{m-1}=\lz_1^{(m)} g_m.
\endcases$$
Making a difference of this with the one
replacing $i$ by $i+1$, we get ${\widetilde R}_i(u)=\lz_1^{(m)}$ for all $i\le
m-1$ (From this, the reason should be clear why in the definition
of ${\widetilde{\scr V}}$, we use ``$v_i=0$ for $i\ge m-1$'' rather than
``$v_i=v_{i\wedge m}$''). Thus,
$$\align
\lz_1^{(m)}&=\max_{0\le i< m}\, {\widetilde R}_i(u)\\
&\ge \inf_{u:\;\supp(u)=\{0,1,\ldots, m-1\};\; (v_i=u_{i+1}/u_i)\in {\widetilde{\scr V}}}\,\max_{0\le i< m}\, {\widetilde R}_i(u)\\
&\ge \inf_{u:\;\supp(u)=\{0,1,\ldots, n\}\text{ for some }n\ge 0,\;n<N;\; (v_i=u_{i+1}/u_i)\in {\widetilde{\scr V}}}\,\sup_{0\le i<N}\, {\widetilde R}_i(u)\\
&= \inf_{v\in {\widetilde{\scr V}}}\sup_{0\le i<N} {\widetilde
R}_i(v).
\endalign$$
Here in the last step, we have returned to the original notation
${\widetilde R}(v)$ instead of ${\widetilde R}(u)$.
Letting $m \to N$, we obtain the required assertion.
\qed\enddemo

With the same rates $(a_i, b_i)$ here but endow with the Dirichlet
boundary at $0$, we return to the situation studied in Section 4.
The next result, taken from [7; \thm\;2.2] and [6;
\thm\;3.5], is a comparison of $\lz_1$ with the quantities $\lz_0$,
$\dz$, $\dz_1$ and $\dz_1'$ given in Section 4. See also \crl\;6.6
below for an improvement.

\proclaim{\thm\;6.2 (Criterion and basic estimates)} Under (6.3), $\lz_1>0$ iff $\dz<\infty$.
More precisely, we have
$$\frac{1}{4\,\dz}\le\frac{1}{\dz_1}\le \lz_0
\le \lz_1\le {\lz_0}{Z}\le \frac{Z}{\dz_1'}
\le \frac{Z}{\dz}. \tag 6.6$$
\endproclaim

The next two results are mainly taken from [7; \thm\;2.4]
with an addition on the monotonicity of $\{\ez_n\}$ and $
\{{\ez}_n'\}$.

\proclaim{\thm\;6.3 (Approximating procedure)} Let $(6.3)$ hold  and
$\dz<\infty$. Write $\fz_0=0$, $\fz_i=\sum_{0\le j \le i-1} (\mu_j
b_j)^{-1}=:\nu[0, i-1]\, (1\le i\in E)$, $\bar f=f - \pi(f)$,
$\pi=\mu/Z$, and $Z=\sum_{k\in E}\mu_k=:\mu[0, N]$. \roster
\item Define $f_1=\sqrt{\fz}$,
$f_n =\bar f_{n-1} {I\!I}(\bar f_{n-1})$ and $ \eta_n=\sup_{1\le
i\in E}{I}_i (\bar f_n ).$ Then $\ez_n$ is decreasing in $n$ and
$\lz_1\ge \eta_n^{-1} \ge (4 \delta )^{-1}$ for all $n\ge 1$.
\item For fixed $m\!\in\! E\!: m \!\ge\! 1$, define
$$f_1^{(m)}=\fz_{\cdot\wedge m},\qqd
f_n^{(m)} =\big[\bar f_{n-1}^{(m)}
{I\!I}\big(\bar f_{n-1}^{(m)}\big)\big](\cdot \wedge m),\qqd n\ge 2,
$$ and then define
$$\ez_n'=\sup_{1\le m\in E}\,\inf_{1\le i\in E}{I}_i (\bar f_n ),\qqd
{\bar{\ez}}_n =\sup_{1\le m\in E} \frac{\mu\big({\bar
f}_n^{(m)\,2}\big)}{D\big(f_n^{(m)}\big)},\qqd n\ge 1. $$
Then $\ez_n'$ is increasing in $n$ and ${\ez_n'}^{\!-1}\ge {\bar{\ez}}_n^{-1}\ge \lz_1$ for all $n\ge 1$.
\endroster
\endproclaim

The notation ``$\bar f_{n-1} {I\!I}(\bar f_{n-1})$'' used in the theorem may
have $0/0$ but it should not cost any confusion. Note that here we use the same
$(\nu_j)$ as in (2.15). In other words, when $b_0>0$, we use (2.15). But for its
dual, it is more convenient to use ${\hat\nu}_j=(\hat\mu_j \hat a_j)^{-1}$ as in
Section 4 since $b_0=0$. This is consistent with the notation used in Section 5.

As a consequence of \thm\;6.3, we have the following improvement of \thm\;6.2.

\proclaim{\crl\;6.4\,(Improved estimates)} Let $(6.3)$ hold. Then we have
$$(4\dz)^{-1}\le\ez_1^{-1}\le\lz_1\le {\bar\ez}_1^{-1}, \tag 6.7$$
where
$$\align
\text{\hskip-6em}\ez_1&=\sup_{1\le i\in E}\big(\sqrt{\fz_i}+\sqrt{\fz_{i-1}}\,\big)
\bigg[\psi_i -\psi_1\frac{\mu[i, N]}{\mu[0, N]}
\bigg],\qd
\psi_i:=\sum_{j=i}^N\mu_j \sqrt{\fz_j},\text{\hskip-3em}\tag 6.8\\
\text{\hskip-6em}{\bar{\ez}}_1 &=\sup_{1\le m\in E}
\frac{1}{\fz_m}\bigg[\sum_{1\le k\in E}\mu_k\fz_{k\wedge m}^2-\frac{1}{Z}
\bigg(\sum_{1\le k\in E}\mu_k\fz_{k\wedge m}\bigg)^2\bigg]\\
\text{\hskip-6em}&=\sup_{1\le m\in E}\bigg\{
 \frac{1}{\fz_m}\bigg[\sum_{1\le k\le m-1}\mu_k\fz_k^2-\frac{1}{\mu [0, N]}
\bigg(\sum_{1\le k\le m-1}\mu_k\fz_k\bigg)^2\bigg]\\
\text{\hskip-6em}&\qqd\qqd\qd +\frac{\mu [m, N]}{\mu [0, N]}\bigg[\fz_m
{\mu [0, m-1]}- 2\sum_{1\le k\le m-1}\mu_k\fz_k\bigg]
\bigg\}.\tag 6.9\endalign
$$
\endproclaim

\demo{\prf\; of \thm\;$6.3$}

{\it Part $1$}. We prove that
$\{f_n\}\subset L^1(\mu)$ in three steps. This was missed in the original paper [7].
Certainly, we need only to consider the case that $N=\infty$.

(a) First, we show that the functions $\{h_n\}$,
$$h_0(i)\equiv 1,\qd i\in E,\qqd
h_n(i)=\sum_{j=1}^i \frac{1}{\mu_j a_j} \sum_{k=j}^\infty \mu_k h_{n-1}(k),\qd i\ge 1,\;n\ge 1,$$
are all in $L^1(\mu)$. Clearly, $h_1$ (and then $h_n$ for $n\ge 2$)
may increase to infinity if the minimal process is recurrent
which is the main problem we need to handle. The required assertion says that
even though $h_n$ can be unbounded but is still in $L^1(\mu)$. For this,
to distinguish with $\{f_n\}$ used in \thm\;6.3, let $\{{\tilde f}_n\}$ be
the sequence defined in part (1) of \thm\;4.3:
$$\align
{\tilde f}_1(i)&=\bigg(\sum_{k=1}^i \nu_{k-1}\bigg)^{1/2}\,=f_1(i),\qqd i\ge 1,\;\; \nu_{j-1}:=\frac{1}{\mu_j a_j},\\
{\tilde f}_n(i)&=\sum_{j=1}^i \nu_{j-1} \sum_{k=j}^\infty \mu_k {\tilde
f}_{n-1}(k),\qqd i\ge 1,\; n\ge 2,\\
{\tilde f}_n(0)&=0,\qqd n\ge 1.
\endalign
$$

From proof (b) of \thm\;4.3, we have seen that
$${\tilde f}_2(i)=
\sum_{j=1}^i \frac{1}{\mu_j a_j }\sum_{k\ge j}\mu_k {\tilde f}_1(k)
\le  4\dz {\tilde f}_1(i).$$
Because ${\tilde f}_1(i)\ge {\tilde f}_1(1)=a_1^{-1/2}$
for $i\ge 1$, this gives us
$$h_1(i)\le 4\dz \sqrt{a_1}\,{\tilde f}_1(i),\qqd i\ge 1.$$
By induction, it follows that
$$h_n\le \sqrt{a_1}\,(4\dz)^n {\tilde f}_1,\qqd n\ge 1.$$
This proves that $h_n\in L^1(\mu)$ for all $n\ge 1$ since ${\tilde f}_1\in L^1(\mu)$
as mentioned in proof (b) of \thm\;4.3, due to
the assumption $\dz<\infty$.

(b) Next, we study the relation between $\{f_n\}$ and $\{{\tilde f}_n\}$.
By definition, we have
$$\align
f_2(i)&=\sum_{j=1}^i \nu_{j-1} \sum_{k=j}^\infty \mu_k {\bar f}_{1}(k)
={\tilde f}_2(i)-h_1(i) \pi (f_{1}),\qqd i\ge 1,\\
f_3(i)&=\sum_{j=1}^i \nu_{j-1} \sum_{k=j}^\infty \mu_k {\bar f}_{2}(k)
={\tilde f}_3(i)-h_2(i)\pi(f_1)-h_1(i) \pi (f_{2}),\qqd i\ge 1,\\
f_4(i)&
={\tilde f}_4(i)-h_3(i)\pi(f_1)-h_2(i)\pi(f_2)-h_1(i) \pi (f_{3}),\qqd i\ge 1.
\endalign$$
Successively, we obtain
$$
f_n={\tilde f}_n- \sum_{k=1}^{n-1}\pi({f}_{k}) h_{n-k}, \qqd n\ge
2.$$

(c) Since $f_1={\tilde f}_1\in L^1(\mu)$ as shown in proof (b) of \thm\;4.3.
Now, to show that $\{f_n\}\subset L^1(\mu)$, by (a) and (b), it suffices to prove
that $\{{\tilde f}_n\}\subset L^1(\mu)$. This is done in proof (b) of \thm\;4.3.

{\it Part $2$}. We now prove the monotonicity of $\{\ez_n\}$ in two steps.
Since ${\bar f}_n$ values both positive and negative or even zero, the proportional property used in
the proof of the monotonicity of $\{\dz_n\}$ is currently not available. To overcome this difficulty, a finer technique
is needed.

(d) Because
$$\mu_i a_i \big[f_n(i)-f_n (i-1)\big]= \sum_{k=i}^N \mu_k {\bar f}_{n-1} (k),\qqd n\ge 2,$$
by the definition of $I(f)$, we obtain
$$\ez_n= \sup_{1\le i\in E} \sum_{j=i}^N \mu_j {\bar f}_{n} (j)\bigg/\sum_{k=i}^N \mu_k {\bar f}_{n-1} (k),
\qqd n\ge 2. \tag 6.10$$
Since the denominator is positive, the assertion that $\ez_n\le \ez_{n-1}$ is equivalent to
$$\sum_{j=i}^N \mu_j \big[{\bar f}_{n} (j)-\ez_{n-1} {\bar f}_{n-1}(j)\big]\le 0,\qqd i\in E. $$
That is,
$$\text{\hskip-2em}\ez_{n-1}\pi(f_{n-1})-\pi(f_n)\le \frac{1}{\mu[i, N]}\sum_{j=i}^N\mu_j \big[\ez_{n-1} f_{n-1}(j)-f_n(j)\big],\qd i\in E.
\tag 6.11$$
Let us observe the meaning of this inequality: the left-hand side is the infimum (attained at $i=0$) of the right-hand side.

The monotonicity of $\{\ez_n\}$ now follows once we show that the right-hand side
of (6.11) is luckily increasing in $i$, or equivalently,
$${\mu[i, N]}\sum_{j=i+1}^N\mu_j \big[\ez_{n-1} f_{n-1}(j)-f_n(j)\big]
\ge {\mu[i+1, N]}\sum_{j=i}^N\mu_j \big[\ez_{n-1} f_{n-1}(j)-f_n(j)\big].$$
By removing the common term
$${\mu[i+1, N]}\sum_{j=i+1}^N\mu_j \big[\ez_{n-1} f_{n-1}(j)-f_n(j)\big]$$
in both sides, it is enough to check that
$$\text{\hskip-1em}\ez_{n-1}\!\!\!\sum_{j= i+1}^N\!\mu_j \big[f_{n-1}(j)\!-\!f_{n-1}(i)\big]
\ge \!\!\!\sum_{j= i+1}^N \!\mu_j\big[f_{n}(j)\!-\!f_{n}(i)\big],\qd
i\in E,\; n\ge 2. \tag 6.12$$ First, let $n\ge 3$. Then by the
definition of $f_n$ and (6.10), we have
$$\align
f_{n}(j)-f_{n}(i)&=\sum_{s=i+1}^j \nu_{s-1} \sum_{k=s}^N \mu_k {\bar f}_{n-1}(k)\\
&\le \ez_{n-1} \sum_{s=i+1}^j \nu_{s-1} \sum_{k=s}^N \mu_k {\bar f}_{n-2}(k)\\
&= \ez_{n-1}\big[ f_{n-1}(j)-f_{n-1}(i)\big].
\endalign$$
This certainly implies (6.12) in the case of $n\ge 3$, regarded as an application of
the proportional property. Next, let $n=2$. Then by the definition of
$f_2$ and $\ez_1$, we have
$$\align
f_{2}(j)-f_{2}(i)&=\sum_{s=i+1}^j \nu_{s-1} \sum_{k=s}^N \mu_k {\bar f}_{1}(k)\\
&\le \ez_{1} \sum_{s=i+1}^j \big[f_1(s)-f_1(s-1)\big]\\
&= \ez_{1}\big[ f_{1}(j)-f_{1}(i)\big].
\endalign$$
This also implies (6.12) in the case of $n=2$. We have thus proved that $\ez_n\le \ez_{n-1}$
for all $n\ge 2$.

{\it Part $3$}. To prove the monotonicity of $\{\ez_n'\}$, for each fixed $m$, as a dual argument
(exchanging ``$\sup$'' and ``$\le$'' with ``$\inf$'' and ``$\ge$'', respectively) of the above proofs
(d) and (e), we have
$$\inf_{1\le i\in E} I_i\big({\bar f}_n^{(m)}\big)\le \inf_{1\le i\in E} I_i\big({\bar f}_{n+1}^{(m)}\big).$$
Then the assertion follows by making supremum with respect to $m$.

{\it Part $4$}.  The proof of ${\bar\ez}_n\ge \ez_n'$ is given in \lmm\;6.5 below.
\qed\enddemo

In practice, using ${\bar{\ez}}_n$ rather than $\ez_n'$
is based on the following result.

\proclaim{\lmm\;6.5} For every non-decreasing, and non-constant function
$f$ satisfying $f\in L^1(\mu)$ and $D(f)<\infty$, we have
$$\frac{\mu \big(\bar f^2\big)}{{D}(f)}\ge \inf_{1\le i\in E} I_i (\bar f).$$
Similarly, for every nonnegative, non-decreasing, and non-zero function
$f$ satisfying $f\in L^1(\mu)$ and $D(f)<\infty$, we have
$$\frac{\mu \big(f^2\big)}{{D}(f)}\ge \inf_{1\le i\in E} I_i (f).$$
\endproclaim

\demo{\prf} (a) Since $f$ is not a constant, we have $\mu(\bar f^2)>0$
and $D(f)>0$.
Moreover, since $f\in L^1(\mu)$ is also non-decreasing, we claim that
$$\infty> \sum_{k=i}^N \mu_k \bar f_k>0\qqd\text{for all }i\in E,\qqd i\ge 1.$$
Actually, the non-decreasing sequence $\{\bar f_k\}$, starting at $\bar f_0<0$
(since $f$ is non-trivial) and
having mean zero, should be positive for all large enough $k$. Thus, if
$\sum_{k=i_0}^N  \mu_k \bar f_k\le 0$ for some $i_0: 1\le i_0\in E$,
then we would have $\bar f_{i_0}<0$
\big(otherwise $\bar f_i\ge 0$ for all $i\ge i_0$ and then
$\sum_{k=i_0}^N  \mu_k \bar f_k>\mu_j \bar f_j>0$ for large enough $j$\big). This implies that
$$\sum_{k\le i_0-1} \mu_k \bar f_k\le \bar f_{i_0}\sum_{k\le i_0-1} \mu_k<0, $$ and
furthermore,
$$0=\mu(\bar f)=\sum_{k=i_0}^N  \mu_k \bar f_k+ \sum_{k\le i_0-1} \mu_k \bar f_k
\le \sum_{k\le i_0-1} \mu_k \bar f_k <0,$$
which is a contradiction. Because of the assertion we have just proved and using the convention that
$1/0=\infty$, it follows that
$\inf_{1\le i\in E} I_i (\bar f) \in [0, \infty)$.

Let $\gz=\inf_{1\le i\in E} I_i (\bar f)$. Then we have
$$-\sum_{k\le i-1}\mu_k \bar f_k
=\sum_{k=i}^N\mu_k \bar f_k\ge \gz\, \mu_i a_i (\bar f_i-\bar f_{i-1})$$
first for those $i$ with $f_i> f_{i-1}$ and then for all $i: 1\le i\in E$.
Multiplying both sides by $\bar f_i-\bar f_{i-1}\ge 0$, we obtain
$$-(\bar f_i-\bar f_{i-1})\sum_{k\le i-1}\mu_k \bar f_k\ge \gz\,
\mu_i  a_i (\bar f_i-\bar f_{i-1})^2,
\qqd i\in E,\; i\ge 1.$$
Making a summation over $i$ from 1 to $m$, it follows that
$$-\sum_{i=1}^m (\bar f_i-\bar f_{i-1})\sum_{k\le i-1}\mu_k \bar f_k
\ge \gz \sum_{i=1}^m\mu_i a_i (\bar f_i-\bar f_{i-1})^2.$$
Noticing that the mean of
$\bar f$ equals zero and exchanging the order of the sums, the left-hand side is equal to
$$\align
-\sum_{k=0}^m\mu_k \bar f_k \sum_{i=k+1}^m (\bar f_i-\bar f_{i-1})
&=-\sum_{k=0}^m\mu_k \bar f_k (\bar f_m-\bar f_k)\\
&=-\bar f_m\sum_{k=0}^m\mu_k \bar f_k+\sum_{k=0}^m\mu_k \bar f_k^2\\
&=\sum_{k=m+1}^N\mu_k \bar f_k\bar f_m+\sum_{k=0}^m\mu_k \bar f_k^2.
\endalign$$
As mentioned in the last paragraph, $\bar f_m>0$ first for some $m$ and then for all large
enough $m$ since $\bar f$ is non-decreasing, the right-hand side is controlled,
for large enough $m$, from above by
$$\sum_{k=m+1}^N\mu_k \bar f_k^2+ \sum_{k=0}^m\mu_k \bar f_k^2=\mu(\bar f^2).$$
With the assumption $D(f)<\infty$ in mind, the required assertion now follows
immediately by passing the limit as $m\to N$.

(b) For the second assertion, since $f\in L^1(\mu)$ is nonnegative and non-zero,
we have
$$\infty>\sum_{j=i}^N\mu_j f_j>0\qqd \text{for all } i\in E.$$
Now, if $f_0=0$, then there is an $i_0$ such that $f_{i_0-1}=0$ but
$f_{i_0}>0$ and so $I_{i_0}(f)<\infty$. If $f_0>0$ and
$\inf_{1\le i\in E} I_i (f)=\infty$, then $f$ should be a positive constant,
and hence, $D (f)=0$. In this case, the assertion is trivial since
$\mu(f^2)>0$.
Therefore, we may assume that $\gz:=\inf_{1\le i\in E} I_i (f)<\infty$.
We now have
$$\sum_{k=i}^N\mu_k  f_k\ge \gz \mu_i a_i (f_i- f_{i-1}),
\qqd i\in E,\; i\ge 1.$$
Hence,
$$\sum_{i=1}^m (f_i-f_{i-1})\sum_{k=i}^N\mu_k  f_k
\ge \gz \sum_{i=1}^m\mu_i a_i (f_i-f_{i-1})^2.$$
Exchanging the order of the sums, the left-hand side is equal to
$$\sum_{k=1}^N\mu_k  f_k\sum_{i=1}^{k\wedge m}(f_i-f_{i-1})
=\sum_{k=1}^N\mu_k  f_k (f_{k\wedge m}-f_{0})
\le \sum_{k=1}^N\mu_k  f_k f_{k\wedge m}\le \mu(f^2).$$
Combining this with the last inequality, we have obtained the required assertion.
\qed\enddemo

Having the comparison of ${\bar\ez}_n\ge \ez_n'$ (Lemma 6.5) in
mind, one may expect a parallel result for $\dz_n'$ and
${\bar\dz}_n$ defined in \thm\;4.3. All the examples we have ever
computed support the conjecture that ${\bar \dz}_n\ge \dz_n'$, however,
there is still no proof. In general, we have ${\bar\dz}_{n+1}\ge
\dz_n'$ only as stated in \thm\;4.3. Note that $\dz_n'$ is defined
by using $I\!I(f_n)$ rather than $I(f_n)$. If we redefine $\dz_n'$
by using $I(f_n)$ as in [7; \thm\;2.2], denoted by
${\tilde\dz}_n'$ for a moment, then by the second assertion of \lmm\;6.5, we do
have ${\bar\dz}_n\ge {\tilde\dz}_n'$. Besides, by the theorem just
quoted, we also have $\dz_n'\ge {\tilde\dz}_n'\ge \dz_{n-1}'$. This
remark is also meaningful for those $\dz_n'$ and ${\bar\dz}_n$
defined in Section 3.

Note that the factor of the upper and lower bounds of $\lz_1$ given
in \thm\;6.2 is $4 Z>4$. The next result has a factor $4$ only. A simple comparison
of $\kz$ below and $\dz^{(4.4)}$ shows that it is not easy to find such a result. Its
proof is delayed to the next section.

\proclaim{\crl\;6.6\,(Criterion and basic estimates)} Let $(6.3)$ hold.
Then we have $\kz^{-1}/4\le \lz_1\le
\kz^{-1}$, where
$$\kz^{-1}=\inf_{0\le n< m< N+1}\bigg[\bigg(\sum_{i=0}^{n}\mu_i\bigg)^{-1}
+\bigg(\sum_{i=m}^N \mu_i\bigg)^{-1}\bigg]\bigg(\sum_{j=n}^{m-1} \frac{1}{\mu_j b_j}\bigg)^{-1}. \tag 6.13$$
Furthermore, we have
$$\dz_L\wedge \dz_R\ge \kz\ge  Z^{-1} \dz_L,$$
where $Z=\sum_{i=0}^{N}\mu_i$,
$$\dz_L\!=\!\sup_{1\le n< N+1} \sum_{i=1}^n\frac{1}{\mu_i a_i}
\sum_{j=n}^N \mu_j=\dz^{(4.4)},\qqd \dz_R\!=\!\sup_{0\le m< N}\sum_{j=0}^m \mu_j
\sum_{k=m}^{N-1}\frac{1}{\mu_k b_k}.$$
\endproclaim

In the case that the minimal process is ergodic, since
$$1<Z<\infty,\qqd \sum_{j} \frac{1}{\mu_j b_j}=\infty,$$
we have $\dz_R=\infty$ and so the second
assertion of \crl\;6.6 goes back to \thm\;6.2. However, the first assertion of \crl\;6.6 is
clearly finer. An extension of \crl\;6.6 to a more general state space is given in \crl\;7.9
below.

\medskip

Most of the examples below are taken from [10;  Examples 9.27].
The computation of ${\bar{\ez}_1}$, $\ez_1/\bar{\ez}_1$, and $\kz$ is newly
added.

\proclaim{\xmp\;6.7} Let $ b_i=b\,(i\ge 0)$, and $ a_i=a\,(i\ge 1)$,
$a>b$. Then
$$\lz_1=\big(\sqrt{a}-\sqrt{b}\,\big)^2,\qd \dz=\kz=a
(a-b)^{-2},\qd \bar{\ez}_1=\dz_1'=(a+b)/{(a-b)^2},$$ and
$\ez_1=\lz_1^{-1}$ which is sharp. Besides, $\ez_1/\bar{\ez}_1\le
2$, the equality sign holds iff $b=a$. Note that $\lz_1$,
$\ez_1^{-1}$, and $\bar{\ez}_1^{-1}$ all tend to zero as $b\to a$.
Furthermore, $({\bar\ez}_1, \ez_1)\subset (\kz, 4\kz)$.
\endproclaim

\proclaim{\xmp\;6.8} The typical linear model: let $ b_i= \bz_1 i
+\bz_0\,(\bz_0>0,\;\bz_1\ge 0)$, and $a_i\!=\gz_1 i\,(\gz_1\!>\!\bz_1)$ for
$i\!\ge\! 0$. Then $\lz_1\!=\gz_1\!-\!\bz_1$. When $\bz_0\!=\!0$, we have $\lz_0^{(4.2)}\!=\gz_1\!-\!\bz_1$.
\endproclaim

\proclaim{\xmp\;6.9} Let $ b_i= b/(i +1)\,(b>0)$ for $i\ge 0$,
$a_i\equiv a>0 $ for $i\ge 1$. Then
$\lz_1=a-\big(\sqrt{b^2+4ab}-b\big)/2$.
\endproclaim

\proclaim{\xmp\;6.10} Let $ b_i\equiv  b\,(b>0)$ for $i\ge 0$,
$a_i= (i\wedge k)a\,(a>0)$ for $i\ge 1$ and some $k\ge 2$ satisfying
$$k^{-1}\le a/b\le k (k-1)^{-2}.$$ Then
$\lz_1=\big(\sqrt{ak}-\sqrt{b}\big)^2$.
\endproclaim

\proclaim{\xmp\;6.11} Let $ b_0=1$, $ b_i=i$, and $ a_i=2\, i$, $i\ge
1$. Then $\lz_1\ge\lz_0= 1$ but the precise value is unknown.
Moreover,
$$\bar{\ez}_1\approx 0.55,\qd \ez_1\approx 0.9986\qd \text{and}\qd
\ez_1/\bar{\ez}_1\approx 1.82<2.$$ Besides, $\kz\approx 0.4856$ and
so $({\bar\ez}_1, \ez_1)\subset (\kz, 4\kz)$.
\endproclaim

The next one is a continuation of [6; \xmp\;3.10].

\proclaim{\xmp\;6.12} Let $E=\{0, 1\}$. Then $\lz_1 =Z \lz_0={\bar\ez}_1^{-1}=\kz^{-1}$ and
$\lz_0=\dz^{-1}$. Hence, the last upper bound in (6.6) and the one in \crl\;6.6 are sharp
but $\dz^{-1}$ is not an upper bound of $\lz_1$.
\endproclaim

The first lower bound in (6.6) and the one in \crl\;6.6 are sharp for the seventh example
in Table 6.1 below.

\proclaim{\xmp s\;6.13} Here are some additional examples, given in Table 6.1, for which
the quantities $\bar{\ez}_1\le \lz_1^{-1}\le \ez_1$ and $\kz\le
\lz_1^{-1}\le 4\kz$ are compared. For all these examples, we have
$({\bar\ez}_1, \ez_1)\subset (\kz, 4\kz)$ and so the estimates given
in \crl\;6.4 are better than the ones in \crl\;6.6.
\endproclaim


\centerline{{\bf Table 6.1} \quad Exact $\lz_1$ and its estimates for eight examples}
\nopagebreak
\vskip-0.5truecm
$$
  \hfil\vbox{\hbox{\vbox{\offinterlineskip
  \halign{&\vrule#&\strut\;\hfil#\hfil\;&\vrule#&
  \;\hfil#\hfil\;&\vrule#&
  \;\hfil#\hfil\;&\vrule#&
  \;\hfil#\hfil\;\cr
  \noalign{\hrule}
    height2pt&\omit&&\omit&&\omit&&\omit&&\omit&&\omit&\cr
    & ${\pmb {b_i}}\,(i\ge 0)$ && ${\pmb {a_i}}\,(i\ge 1)$ && ${\pmb{\lz_1^{-1}}}$ && $\pmb{{\bar\ez_1}}$
    && $\pmb{{\ez_1}}$&& $\pmb{{\ez_1/\bar\ez_1}}$&& $\pmb{\kz}$&\cr
    height2pt&\omit&&\omit&&\omit&&\omit&&\omit&&\omit&\cr
    \noalign{\hrule}
    height2pt&\omit&&\omit&&\omit&&\omit&&\omit&&\omit&\cr
    & {$i+1$} && {$2 i$} && $1$ && $\approx 0.8$ && $\approx 1.48$&& $\approx 1.85$ && $2/3 $&\cr
    height2pt&\omit&&\omit&&\omit&&\omit&&\omit&&\omit&\cr
    \noalign{\hrule}
    height2pt&\omit&&\omit&&\omit&&\omit&&\omit&&\omit&\cr
    & {$i+1$} && {$\matrix 2 i+3\\a_0=0\endmatrix$} && $1/2$ && $\approx 0.346$&& $\approx 0.638$&& $\approx 1.84$ && $\approx 0.28 $ &\cr
    height2pt&\omit&&\omit&&\omit&&\omit&&\omit&&\omit&\cr
    \noalign{\hrule}
    height2pt&\omit&&\omit&&\omit&&\omit&&\omit&&\omit&\cr
    & {$i+1$} && {$\matrix 2 i\!+\!4\!+\!\sqrt{2}\\a_0=0\endmatrix$} && $1/3$
    && $\approx 0.218$&& $\approx 0.398$&& $\approx 1.83$ && $\approx 0.18 $&\cr
    height2pt&\omit&&\omit&&\omit&&\omit&&\omit&&\omit&\cr
    \noalign{\hrule}
    height2pt&\omit&&\omit&&\omit&&\omit&&\omit&&\omit&\cr
    & $(i+1)^{-1}$ && $1$ && $\matrix 2(3-\sqrt{5})^{-1}\\ \approx 2.618\endmatrix$
    && $\approx 1.92$&& $\approx 3.24$&& $\approx 1.69$ && $\approx 1.6 $&\cr
    height2pt&\omit&&\omit&&\omit&&\omit&&\omit&&\omit&\cr
    \noalign{\hrule}
    height2pt&\omit&&\omit&&\omit&&\omit&&\omit&&\omit&\cr
    & $1$ && $i\wedge 2$ && $\matrix (\sqrt{2}-1)^{-2}\\ \approx 5.8284\endmatrix $
    && $\approx 3$&& $\approx 5.8284$&& $\approx 1.9$ && $2 $ &\cr
    height2pt&\omit&&\omit&&\omit&&\omit&&\omit&&\omit&\cr
    \noalign{\hrule}
    height2pt&\omit&&\omit&&\omit&&\omit&&\omit&&\omit&\cr
    & $i+2$ && $i^2$ && $ 1/2$  && $\approx 0.47$&& $\approx 0.85$&& $\approx 1.81$ && $\approx 0.47 $ &\cr
    height2pt&\omit&&\omit&&\omit&&\omit&\cr
    \noalign{\hrule}
    height2pt&\omit&&\omit&&\omit&&\omit&&\omit&&\omit&\cr
    & {$\matrix i^2\\ b_0=1\endmatrix$} && $i^2$ && $4$
    && $2$&& $\lz_1^{-1}$&& $2$ && $ 1 $&\cr
    height2pt&\omit&&\omit&&\omit&&\omit&&\omit&&\omit&\cr
    \noalign{\hrule}
    height2pt&\omit&&\omit&&\omit&&\omit&&\omit&&\omit&\cr
    & {$\matrix 2+(-1)^i\\b_0\!=\!\dfrac{7\!-\!\sqrt{33}}{2}\endmatrix$}
    && {$\matrix 2 [2\!+\!(-1)^i]\\a_0=0\endmatrix$} && $\matrix (6-\sqrt{33})^{-1}\\ \approx 3.9\endmatrix$
    && $\approx 2.11$&& $\approx 4.21$&& $\approx 2$ && $\approx 1.56 $ &\cr
    height2pt&\omit&&\omit&&\omit&&\omit&&\omit&&\omit&\cr
  \noalign{\hrule}}}}}\hfill$$

\head{7. Bilateral absorbing (Dirichlet) boundaries}
\endhead

This section deals with the fourth case of boundary conditions.
It consists of two parts. The first one is for the ordinary
birth--death processes as studied in the previous sections
and the second one deals with the bilateral birth--death processes
with a more general state space.

First, let us consider the processes with state
space $E=\{i: 1\le i< N+1\}\,(N\le \infty)$ with Dirichlet
boundaries at $0$ ($a_1>0$) and $N+1$ if $N<\infty$. Similar to
Section 2, define
$$\lz_0=\inf\big\{{D}(f): \mu(f^2)=1,\; f \in {\scr K}\big\},\tag 7.1$$
where the symmetric measure $(\mu_i)$ is the same as in Section 4,
$\mu(f)=\sum_{k\in E} \mu_k f_k$, and
$$\aligned
D(f)=\sum_{k\in E}\mu_k a_k(f_k-f_{k-1})^2,\qqd f_0:=0,
\endaligned$$
with domain ${\scr D}^{\min}(D)$.
Clearly, if one changes only the boundary condition at $0$, then
the resulting $\lz_0$ is bigger or equal to $\lz_0^{(2.2)}$. Note that
if (1.3) fails, then the eigenvalues $\lz_0^{(4.2)}$ and
$\lz_0^{(7.1)}$ are different which correspond to the maximal and the minimal
Dirichlet forms, respectively. However, as mentioned in Section 4,
once (1.3) holds, $\lz_0^{(4.2)}$ coincides with $\lz_0^{(7.1)}$. Then there are three
cases. The first one is that $\sum_i \mu_i <\infty$ and
$\sum_i(\mu_i a_i)^{-1}=\infty$. This case is treated in Section 4.
In this section, we are mainly studying the second case that $\sum_i
\mu_i =\infty$ but
$$\sum_{k=1}^N \frac{1}{\mu_k a_k}<\infty. \tag 7.2$$
The third case is that $\sum_i \mu_i =\infty$ and $\sum_i(\mu_i
a_i)^{-1}=\infty$ which is treated in the next theorem.
In this degenerated case, since there is a killing at $1$ (i.e.,
$a_1>0$), the process is transient. Without using duality, by \crl\;7.3 below, we also obtain that $\lz_0=0$.
See the comments right after \crl\;7.3.

\proclaim{\thm\;7.1} \roster\item First, let $(7.2)$ hold.
Define the dual rates $\big(\hat a_i, \hat b_i\big)$ by $(5.1)$ in
the inverse way:
$$\hat b_i=a_{i+1},\qqd \hat a_i=b_i,\qqd 0\le i <N+1,\tag 7.3$$
 and denote by $\hat\lz_1$ the eigenvalue defined in
Section 6 for the dual process. Then we have $\lz_0=\hat\lz_1$.
\item
Next, let (7.2) fail and $\sum_{i\in E} \mu_i =\infty$. Then $\lz_0$
\big(as well as $\lz_0^{(4.2)}$\big) is equal to its dual ${\hat\lz}_0=\lz_0^{(2.2)}=0$.
\endroster
\endproclaim

\demo{\prf} (a) By (7.3), (7.2) and (5.5), we have
$$\sum_{k=0}^N \hat\mu_k<\infty.
\tag 7.4$$ Clearly, the dual process with rates $\big(\hat a_i, \hat
b_i\big)$ has the state space ${\widehat E}=\{i: 0\le i <N+1\}$.
By exchanging $(a_i, b_i, v_i)$ and $\big(\hat a_i, \hat b_i, {\hat v}_i\big)$
in part (1) of Theorem 6.1,
$${\hat\lz}_1=\sup_{{\hat v}}\inf_{0\le i<N}\bigg[{\hat a}_{i+1}+{\hat b}_i-\frac{{\hat a}_i}{{\hat v}_{i-1}}
  -{\hat b}_{i+1}{\hat v}_i\bigg],$$
and in (5.8) with $N'=N$,
$$\sup_{{\hat v}}\inf_{0\le i<N}\bigg[{\hat a}_{i+1}+{\hat b}_i-\frac{{\hat a}_i}{{\hat v}_{i-1}}
  -{\hat b}_{i+1}{\hat v}_i\bigg]
  =\sup_{{v}}\inf_{i\in E} \bigg[{a_{i}}\bigg(1 -
\frac{1}{{v}_{i-1}}\bigg)+ {b_{i}}(1-{v}_{i})\bigg],$$
the first assertion of \thm\;7.1 now follows from the variational formula
given on the right-hand side of (9.2) in Section 9.

(b) Similarly, replacing the use of \thm\;6.1 by \prp\;2.7\,(1), we obtain the second assertion.
In this case, as already mentioned at the beginning of Section 4, we have
$\lz_0=\lz_0^{(4.2)}$. The fact that ${\hat\lz}_0=\lz_0^{(2.2)}=0$ comes from (5.5) and \thm\;3.1.\qed\enddemo

By \thm\;7.1\,(1), all the results obtained in Section 6 can be transformed
into the present setup. For instance, by \crl\;6.4, we obtain the following result.

\proclaim{\crl\;7.2} Under (7.2), we have $(4\dz)^{-1}\le \dz_1^{-1}\le \lz_0\le {\bar\dz}_1^{-1}$, where
$$\align
\dz&= \bigg[\sup_{1\le n<N} \,\mu[1, n]\bigg( \nu[n\!+\!1, N]\!+\!\frac{\,\dbl_{\{N<\infty\}}}{\mu_N b_N}\bigg)\bigg]\bigvee\frac{\mu[1, N]\,\dbl_{\{N<\infty\}}}{\mu_N b_N},\qd\; \nu_k=\frac{1}{\mu_k a_k},\\
\dz_1\!&=\sup_{i\in E}\big(\sqrt{\fz_i}+\sqrt{\fz_{i-1}}\,\big)
  \bigg(\psi_i -\psi_1 \frac{\nu[i+1, N]+(\mu_N b_N)^{-1} \dbl_{\{N<\infty\}}}{\nu[1, N]+(\mu_N b_N)^{-1}\dbl_{\{N<\infty\}}}\bigg),\\
{\bar\dz}_1\!&=\sup_{m\in E} \frac{1}{\fz_m}\bigg[\sum_{k=1}^{N-1}\nu_{k+1} \fz_{k\wedge m}^2+\frac{\fz_m^2}{\mu_N b_N}\,\dbl_{\{N<\infty\}} -\\
&\qqd\qqd- \frac{1}{\nu[1, N]\!+\!(\mu_N b_N)^{-1}\dbl_{\{N<\infty\}}}\bigg(\sum_{k=1}^{N-1}\nu_{k+1}\fz_{k\wedge m}
  +\frac{\fz_m}{\mu_N b_N}\,\dbl_{\{N<\infty\}}\bigg)^{\!\!2}\,\bigg],\\
\fz_i\!&=\mu[1, i],\qqd \psi_i=\sum_{j=i}^{N-1} \nu_{j+1} \sqrt{\fz_j} +\frac{1}{\mu_N b_N} \sqrt{\fz_N}
\,\dbl_{\{N<\infty\}},\qqd i\in E.
\endalign$$
\endproclaim

\demo{\prf} Starting from \crl\;6.4 with its notation, write everything we need in its dual.
First by (5.4), we have
$$\mu_n={\hat a_1} \hat \nu_{n+1},\qd
\nu_n=\frac{1}{\hat a_1} \hat\mu_{n+1}, \qd 0\le n< N,\qqd
\mu_N=\frac{{\hat a}_1}{{\hat\mu}_N {\hat b}_N}\qd \text{if }\; N<\infty.$$
Here, recall that $\nu_n=(\mu_n b_n)^{-1}$ but ${\hat\nu}_n=\big({\hat\mu}_n {\hat a}_n\big)^{-1}$. Then
the constant $\dz$ defined in (4.4) becomes $\dz=\sup_{n\in E} \nu[0, n-1]\,\mu[n, N]$. Moreover, we have
$$\align
\fz_i&=\sum_{j=0}^{i-1}\nu_j=\frac{1}{\hat a_1}{\hat\mu} [1, i],\qqd 1\le i<N+1,\\
\mu[m, n]\!&=\sum_{j=m}^n \mu_j={\hat a}_1 {\hat\nu}[m+1, n+1],\qqd 0\le m\le n<N,\\
\mu[m, N]\!&={\hat a}_1 {\hat\nu}[m\!+\!1, N]\!+\! \mu_N\dbl_{\{N<\infty\}}={\hat a}_1\!
   \bigg[{\hat\nu}[m\!+\!1, N]\!+\! \frac{\dbl_{\{N<\infty\}}}{{\hat\mu}_N {\hat b}_N}\! \bigg],\qqd m\!<\!N,\\
\psi_i&= \sum_{j=i}^{N-1}\mu_j \sqrt{\fz_j}+ \mu_N \sqrt{\fz_N}\,\dbl_{\{N<\infty\}}\\
&=\sqrt{{\hat a}_1}\bigg[\sum_{j=i}^{N-1} {\hat\nu}_{j+1}\sqrt{{\hat\mu}[1, j]}
+ \frac{1}{{\hat\mu}_N {\hat b}_N}\sqrt{{\hat\mu}[1, N]}\,\dbl_{\{N<\infty\}}\bigg].
\endalign$$
Inserting these quantities into (6.8) and (6.9), making a little simplification, and then
ignoring the hat everywhere, we obtain \crl\;7.2.\qed\enddemo

The next result is a criterion for the positivity of $\lz_0$, and is
a particular case of \crl\;8.4 with ${\Bbb B}=L^1(\mu)$ in the next
section. It is not deduced from the last section in terms of duality (\thm\;7.1)
but conversely, it provides an improvement of \thm\;6.2 as shown by
the proof of \crl\;6.6 below.

\proclaim{\crl\;7.3\,(Criterion and basic estimates)} Without condition (7.2), we have $\kz^{-1}/4\le
\lz_0\le\kz^{-1}$, where
$$\kz^{-1}=\inf_{1\le n\le m< N+1}\bigg[\bigg(\sum_{i=1}^n\frac{1}{\mu_i a_i}\bigg)^{-1}
+ \bigg(\sum_{i=m}^{N}\frac{1}{\mu_i
b_i}\bigg)^{-1}\bigg]\bigg(\sum_{j=n}^m \mu_j\bigg)^{-1}. \tag 7.5$$
Furthermore, we have
$$\dz_L\wedge \dz_R\ge \kz \ge
\big(\dbl_{\{S=\infty\}}+ (a_1
S)^{-1} \big) \big(\dz_L\wedge \dz_R\big),$$ where
$$\align
S&=\sum_{i=1}^{N}\frac{1}{\mu_i a_i}+\frac{1}{\mu_N b_N}\dbl_{\{N<\infty\}},\\
\dz_L&=\sup_{1\le n< N+1} \sum_{i=1}^n\frac{1}{\mu_i a_i}
\sum_{j=n}^N \mu_j,\qqd \dz_R=\sup_{1\le m< N+1} \sum_{j=1}^m
\mu_j\sum_{k=m}^{N}\frac{1}{\mu_k b_k}.\endalign$$
\endproclaim

Note that $\dz_L=\dz^{(4.4)}$ and $\dz_R$ almost coincides with $\dz^{(3.1)}$, except for $\dz_R$ there is a shift of the
state space. The second assertion of \crl\;7.3 means that $\lz_0>0$
iff the process goes to either $0$ or $N+1$ exponentially fast. This
is intuitively clear by (7.1). Obviously, we have $\lz_0=\kz^{-1}=0$
if $\sum_i \mu_i=\infty$ and $\sum_j (\mu_j a_j)^{-1}=\infty$ since
then $\dz_L=\dz_R=\infty$. See also \crl\;8.6 below.

\demo{\prf\; of \crl\;$6.6$} For given rates $(a_i, b_i)$ in the setup of Section 6, by (5.3), we have
$$\sum_{i=p}^q {\hat\nu}_i=\frac{1}{b_0}\sum_{i=p}^q \mu_{i-1}=\frac{1}{b_0}\sum_{i=p-1}^{q-1} \mu_{i},\qqd
\sum_{j=p}^q {\hat\mu}_j={b_0}\sum_{j=p}^q \nu_{j-1}={b_0}\sum_{j=p-1}^{q-1} \nu_{j}
$$
and
$$\sum_{i=m}^N \frac{1}{{\hat\mu}_i{\hat b}_i}=\sum_{i=m}^{N-1}{\hat\nu}_{i+1}
+ \frac{1}{{\hat\mu}_N{\hat b}_N}\dbl_{\{N<\infty\}}
=\frac{1}{b_0}\sum_{i=m}^{N-1} \mu_i + \frac{1}{b_0}\mu_N \dbl_{\{N<\infty\}}
=\frac{1}{b_0}\sum_{i=m}^{N} \mu_i.$$
Regarding the process studied in \crl\;7.3 as a dual of the one given in the last section
and then add a hat to each quantity of \crl\;7.3. It follows that
$$\align
{\hat\kz}^{-1}&=\inf_{1\le n\le m<N+1}
\bigg[\bigg(\sum_{i=1}^n {\hat\nu}_i\bigg)^{-1}+\bigg(\sum_{i=m}^N \frac{1}{{\hat\mu}_i{\hat b}_i}\bigg)^{-1}\bigg]
\bigg(\sum_{j=n}^m {\hat\mu}_j\bigg)^{-1}\\
&=\inf_{1\le n\le m<N+1}
\bigg[\bigg(\sum_{i=0}^{n-1} \mu_{i}\bigg)^{-1}+\bigg(\sum_{i=m}^N \mu_i\bigg)^{-1}\bigg]
\bigg(\sum_{j=n-1}^{m-1} \nu_j\bigg)^{-1}\\
&=\inf_{0\le n< m<N+1}
\bigg[\bigg(\sum_{i=0}^{n} \mu_{i}\bigg)^{-1}+\bigg(\sum_{i=m}^N \mu_i\bigg)^{-1}\bigg]
\bigg(\sum_{j=n}^{m-1} \nu_j\bigg)^{-1}\\
&=\kz^{-1}.
\endalign
$$
Next, we have
$${\hat a}_1{\widehat S}=b_0\bigg(\sum_{i=1}^N {\hat\nu}_i + \frac{1}{{\hat\mu}_N{\hat b}_N}\dbl_{\{N<\infty\}}\bigg)
=\sum_{i=0}^{N} \mu_i=Z,$$
${\hat\dz}_L=\dz_R$, and ${\hat\dz}_R=\dz_L$. Since $\sum_i \mu_i<\infty$,
by \thm\;7.1\,(1), we have $\lz_1={\hat\lz}_0$. Thus,
\crl\;6.6 now follows from \crl\;7.3 immediately except for a slight change of the lower bound
in the second assertion. For which, since $Z<\infty$, the term ``$\wedge \dz_R$'' is not needed (cf. Proof of \crl\;8.4).\qed\enddemo

\demo{\prf\;of \thm\;$1.5$}
(a) Condition (1.3) implies that $N=\infty$ and furthermore the uniqueness of the symmetric process
on $L^2(\mu)$ by \prp\;1.3. Now, $\az^*=\lz_1$ or $\lz_0$ by [2; \thm\;5.3] or \prp\;1.2,
respectively.

(b) In the case that $\sum_i \mu_i=\infty$ and $\sum_i (\mu_i b_i)^{-1}=\infty$,
the process is zero-recurrent and so we have $\az^*=0$.
Noting that $\dz^{(3.1)}$, $\dz^{(4.4)}$, $\kz^{(6.13)}$, and $\kz^{(7.3)}$ are all equal to infinity,
the conclusions of the theorem become obvious. Hence, in what follows, we may assume that only one of
$\sum_i \mu_i$ and $\sum_i (\mu_i b_i)^{-1}$ is equal to infinity.

(c) Let $b_0=0$. Then the basic estimate follows from \crl\;7.3.

(d) We now prove the first two parts of the theorem under the assumption that $b_0=0$. In the case that
$\sum_i \mu_i<\infty$ but $\sum_i (\mu_i a_i)^{-1}=\infty$, we have $\kz^{(7.5)}=\dz^{(4.4)}$ which
gives us part (1) of the theorem. Next, if $\sum_i \mu_i=\infty$ but $\sum_i (\mu_i a_i)^{-1}<\infty$,
then for $\dz_L$ and $\dz_R$ given in \crl\;7.3, we have $\dz_L=\dz^{(4.4)}=\infty$ and then
$\kz^{(7.5)}<\infty$ iff $\dz_R<\infty$. Clearly, $\dz_R<\infty$ iff
$\dz^{(3.1)}<\infty$ since $\sum_i (\mu_i a_i)^{-1}<\infty$. This gives us part (2) of the theorem.

(e) Finally, let $b_0>0$. This is a dual case of that $b_0=0$ treated in (c) and (d). By exchanging the measures
$\mu$ and $\nu$, we obtain the remaining conclusions of the theorem.

Actually, part (1) of the theorem is a combination of \thm s 4.2 and 6.2, and part (2) is a combination
of \thm s 3.1 and 7.1.
\qed\enddemo

We are now ready to prove an extension of \thm\;1.5.

\proclaim{\thm\;7.4\,(Criterion and basic estimates)} Without condition (1.3), \thm\;1.5 remains true
provided
\roster
\item the process in part (1) is replaced by the maximal one and $\az^*$ is replaced by
$\az^{\max}$: the largest $\vz$ such that
$$\sum_j |p_{ij}(t)-\pi_j|\le C_i\, e^{-\vz t}, \qqd t\ge 0, \tag 7.6$$
for some $L^{1}(\pi)$-locally integrable function $C_i$ depending on $i$ only; and
\item  the process in part (2) needs no change (i.e., the minimal one).
\endroster
\endproclaim

\demo{\prf} Since $\lz_1$ is equivalent to $\lz_0^{(4.2)}$
(\thm\;6.2) and by duality, $\lz_1=\lz_0^{(7.1)}$ and $\lz_0^{(4.2)}=\lz_0^{(2.2)}$,
it is clear that $\lz_0^{(7.1)}$ is equivalent to $\lz_0^{(2.2)}$. Alternatively, one
can use \crl\;7.3 to arrive at the same conclusion.
Now, part (2) of the theorem follows by \prp\;1.2 for which we do not assume (1.3).
As mentioned in the last proof, part (1) with the original $\az^*$ also follows by
[2] provided (1.3) holds.

Even though in the previous study ([12], for instance), we consider only the
ergodic processes under (1.2), but $\lz_1$ can be actually identified with some
exponentially ergodic convergence rate for more
general ergodic processes (reversible Markov chains, in particular). First,
the fact that the $L^2$-exponential convergence rate is described by $\lz_1$ does
not require the regularity of the Dirichlet form (cf. proof of \prp\;1.1, for instance).
Next, for a Markov process with state space $(X, {\scr X}, \pi)$
and transition probabilities $\{P_t(x, \cdot)\}$, let  ${\tilde\vz}_1$ be the largest $\vz$ such that
$$\|P_t(x, \cdot)-\pi\|_{\text{Var}}\le C(x)\, e^{-\vz t}, \qqd x\in X,\; t\ge 0, \tag 7.7$$
for some $L^{1}(X, \pi)$-locally integrable function $C(x)$ depending on $x$ only.
Then for a reversible process having density $p_t(x, y)$,
we have $\lz_1={\tilde\vz}_1$ provided
$$p_s(\cdot, \cdot)\in L_{\text{\rm loc}}^{1/2}(X, \pi)\qqd\text{for some }\; s>0, \tag 7.8$$
and the set of bounded functions with compact support is dense in $L^2(X, \pi)$.
The outline of the proof is as follows.
\roster
\item "(i)"Prove that ${\tilde\vz}_1\ge \lz_1$.
\item "(ii)"Show that $\|(P_t-\pi)f\|^2\le C_f\,e^{- {\tilde\vz}_1 t}$ for bounded $f$ with compact support.
\item "(iii)"Remove the constant $C_f$ in the last line for fixed $f$.
\item "(iv)"Extend $f$ to $L^{2}(X, \pi)$ and then claim that $\lz_1\ge {\tilde\vz}_1$.
\endroster
By assumption, the last step is obvious. The detailed proof for the first three steps
is given, respectively, in [12]: (8.6), the last formula in \S 8.3 replacing ${\vz}_1$ with ${\tilde\vz}_1$,
and the proof of Lemma 8.12. Actually, this is a small correction to
[12; \thm\;8.13\,(4)]
(i.e., replacing ${\vz}_1$ by ${\tilde\vz}_1$) and its proof.
It is known that ${\tilde\vz}_1>0$ iff ${\vz}_1>0$ (as well as $\vz_2>0$ used in the original proof of the cited theorem).
Hence, the exponential ergodicity is kept but the rates ${\vz}_1\ge {\tilde\vz}_1 \ge \vz_2$ may be different.
By the way, we mention that the change of topology is necessary in many cases. For instance,
the pointwise convergence is natural in the discrete case but is not in the continuous case.
In the ergodic situation, the total variation norm is good enough in general but it is meaningless
in the non-ergodic case.

Having this result at hand, part (1) of the theorem follows since we have $\az^{\max}={\tilde\vz}_1$
in the present context.
\qed\enddemo

We now introduce an interpretation, similar to Section 5, of the duality used in \thm\;7.1.
For the ergodic process with $Q$-matrix,
$${Q=\pmatrix
 -b_0 & b_0 & 0 & 0 \\
 a_1 & -a_1-b_1 & b_1 & 0 \\
 0 & a_2 & -a_2-b_2 & b_2 \\
 0 & 0 & a_3 & -a_3
\endpmatrix},\qqd a_i, b_i>0,$$
we have a simpler transformation matrix
$${M=\pmatrix
\mu_0  & \mu_1  & \mu_2  & \mu_3\\
0 & \mu_1  & \mu_2  & \mu_3 \\
0 & 0 & \mu_2  & \mu_3 \\
0 & 0 & 0 & \mu_3
\endpmatrix}\Longrightarrow
{M^{-1}=\pmatrix
\mu_0^{-1} & -\mu_0^{-1} & 0 & 0\\
0 & \mu_1^{-1} & -\mu_1^{-1} & 0\\
0 & 0 & \mu_2^{-1} & - \mu_2^{-1}\\
0 & 0 & 0 & \mu_3^{-1}
\endpmatrix}.$$
Then
$$\align
MQM^{-1}&=
{\pmatrix
0 & 0 & 0 & 0 \\
 b_0 & -a_1-b_0 & a_1 & 0 \\
 0 & b_1 & -a_2-b_1 & a_2 \\
 0 & 0 & b_2 & -a_3-b_2
\endpmatrix}\\
&= {\pmatrix
0 & 0 & 0 & 0 \\
 {\hat a}_1 & -{\hat a}_1-{\hat b}_1 & {\hat b}_1 & 0 \\
 0 & {\hat a}_2 & -{\hat a}_2-{\hat b}_2 & {\hat b}_2 \\
 0 & 0 & {\hat a}_3 & -{\hat a}_3-{\hat b}_3
\endpmatrix}.\endalign$$
We obtain a process having an absorbing state at $0$ and being killed at
the state $3$. The original trivial eigenvalue with non-zero constant eigenfunction
is transferred into the trivial one with eigenfunction $\dbl_{\{0\}}$.
Our dual matrix ${\widehat Q}$ is now obtained by eliminating the first row and
the first column from the matrix on the right-hand side. The elimination is to
make the symmetrizabi\-lity of ${\widehat Q}$ and at the same time removes the trivial
eigenvalue of the last matrix. Unlike the example given in Section 5 where the
size of the state space stays the same: $\{0, 1, 2, 3\}\to \{1, 2, 3, 4\}$
with a shift for the dual one, here the size of the state space is reduced by one:
$\{0, 1, 2, 3\}\to \{1, 2, 3\}$.

We are now ready to examine some examples.

\proclaim{\xmp s\;7.5}{\rm (1) Let $N=1$. Then the $Q$-matrix is degenerated to be a single
killing $-c$ and so $\lz_0=c=\kz^{-1}$.

(2) Let $N\!=\!2$. Then
$$\lz_0\!=\frac{1}{2}\Big(a_1+a_2+b_1+b_2-\sqrt{(a_1-a_2+b_1-b_2)^2+4 a_2 b_1}\,\Big)\!.$$}
\endproclaim
The next two examples are taken from Chen,  Zhang and Zhao (2003,
\xmp s\;2.2 and 2.3)

\proclaim{\xmp s\;7.6} {\rm
(1) Let $N=2$, $a_1=a_2=1$, $b_1=2$, and $b_2=3$. Then
$\lz_0=2$, and by \crl\;7.2, we have
$${\bar\dz}_1\le
\lz_0^{-1}=0.5\le \dz_1,$$ where
$$\dz_1=\frac{4 + \sqrt{3}}{10}\approx 0.573,\qqd \bar{\dz}_1=\frac{7}{15}=0.4{\dot 6},
\qqd \frac{\dz_1}{\bar{\dz}_1}\approx 1.23.$$
Next, $\kz\le \lz_0^{-1}\le 4 \kz$ with $\kz=3/7$. Obviously,
$\big(\bar{\dz}_1, \dz_1\big)\subset (\kz, 4 \kz)$.

(2) Let $N=2$, $b_1=1$, $b_2=2$,
$$a_1=\frac{2-\vz^2}{1+\vz}, \qqd \vz\in \big[0, \sqrt{2}\,\big),$$
and $a_2=1$. Then
$\lz_0=2-\vz$, and we have
$$\bar{\dz}_1\le \lz_0^{-1}=(2-\vz)^{-1}\le \dz_1,$$
where
$$\align
\dz_1&=\frac{4 + \sqrt{2} + (2 + \sqrt{2}\,) \vz - \vz^2}{8 + 2\, \vz - 3 \vz^2}=
\frac{1}{\lz_0}+ \frac{(1 + \vz) (\sqrt{2} - \vz)}{8 + 2\, \vz - 3 \vz^2},\\
\bar{\dz}_1&=\frac{8 + 6 \,\vz - \vz^2}{16 + 4\, \vz - 6 \vz^2}
=\frac{1}{\lz_0}-\frac{\vz^2}{2(8 + 2\, \vz - 3 \vz^2)}.\endalign$$
Hence,
$$\frac{\dz_1}{\bar{\dz}_1}=2 - \frac{2 (4 - \sqrt{2}\,) (1 + \vz)}{8 + 6\, \vz - \vz^2}
<1.354.$$
Next, $\kz\le \lz_0^{-1}\le 4 \kz$ with
$$\kz=\frac{1}{\lz_0}-\min\bigg\{\frac{1}{8 + 2\, \vz - 3\, \vz^2},\; \frac{\vz^2}{8 - 4\, \vz^2 + \vz^3}\bigg\}.$$
Even though it is not so obvious now but we do
have $\big(\bar{\dz}_1, \dz_1\big)\subset (\kz, 4 \kz)$.
}\endproclaim

\proclaim{\xmp s 7.7} {\rm Because of Theorem 7.1, we can now
transfer [10; Examples 9.27] into the present context, see Table 7.1, by
using $(5.1)$ and $(5.7)$. Here, for the sixth example, we need a
restriction: $1/{k}<{b/a}\le {k}/(k-1)^2\, (k\ge 2)$.}
\endproclaim

\centerline{{\bf Table 7.1} \quad Exact $\lz_0$ for nine examples}

\nopagebreak
\vskip -0.5truecm
$$
  \hfil\vbox{\hbox{\vbox{\offinterlineskip
  \halign{&\vrule#&\strut\;\hfil#\hfil\;&\vrule#&
  \;\hfil#\hfil\;&\vrule#&
  \;\hfil#\hfil\;\cr
    \noalign{\hrule}
    height2pt&\omit&&\omit&&\omit&&\omit&\cr
    & ${\pmb{a_i}}\,(i\ge 1)$ && ${\pmb{b_i}}\,(i\ge 1)$ && ${\pmb{\lz_0}}$ && ${\pmb{v_i}}\,(i\ge 1)$ &\cr
    height2pt&\omit&&\omit&&\omit&&\omit&\cr
    \noalign{\hrule}
    height2pt&\omit&&\omit&&\omit&&\omit&\cr
    & $a$ && $b\,(a<b)$ && $ {\bigl(\sqrt a -\sqrt b\,\bigr)^2}_{}^{}$ && $\sqrt{a/b}$ &\cr
    height2pt&\omit&&\omit&&\omit&&\omit&\cr
    \noalign{\hrule}
    height2pt&\omit&&\omit&&\omit&&\omit&\cr
    & ${\matrix\gz_1 (i\!-\!1)\!+\!\gz_0\!\!\!\\ \gz_0\!>\!0,\gz_1\!\ge\! 0\!\!\endmatrix}$ && $\bz_1 i\,(\bz_1\!>\! \gz_1)\!\!\!$ && ${\bz_1-\gz_1}$ && $\dfrac{\gz_1 i+\gz_0}{\bz_1 i}$ &\cr
    height2pt&\omit&&\omit&&\omit&&\omit&\cr
    \noalign{\hrule}
    height2pt&\omit&&\omit&&\omit&&\omit&\cr
    & {${\matrix i-1+\bz_0\\ \bz_0>0\endmatrix}$} && {$2 (i\!+\!1)\!+\!\bz_0$} && $2$ &&
    $\dfrac{(i+1)(i+\bz_0)}{i[2 (i\!+\!1)\!+\!\bz_0]}$ &\cr
    height2pt&\omit&&\omit&&\omit&&\omit&\cr
    \noalign{\hrule}
    height2pt&\omit&&\omit&&\omit&&\omit&\cr
    & {$i$} && {$ 2 i\!+\!4\!+\!\sqrt{2}$} && $3$
    && $\dfrac{i+1}{2 i\!+\!4\!+\!\sqrt{2}}\bigg[1\!+\!\dfrac{2(i+\sqrt{2}\,)}{i(i\!+\!2\sqrt{2}\!-\!1)}\bigg]$\! &\cr
    height2pt&\omit&&\omit&&\omit&&\omit&\cr
    \noalign{\hrule}
    height2pt&\omit&&\omit&&\omit&&\omit&\cr
    & $\dfrac{a}{i}$ && $b$ && $b\!-\!\dfrac{\sqrt{a^2\!+\!4ab}\!-\!a}{2}\!\!\!$ && $\dfrac{\sqrt{a^2\!+\!4 ab}+a}{2b\, i}$ &\cr
    height2pt&\omit&&\omit&&\omit&&\omit&\cr
    \noalign{\hrule}
    height2pt&\omit&&\omit&&\omit&&\omit&\cr
    & $a$ && $(i\wedge k) b$ && ${\big (\sqrt{bk}-\sqrt{a}\, \big)^2}_{}^{} $ && $\dfrac{1}{i\wedge k}\sqrt{ak/b}$ &\cr
    height2pt&\omit&&\omit&&\omit&&\omit&\cr
    \noalign{\hrule}
    height2pt&\omit&&\omit&&\omit&&\omit&\cr
    & $i+1$ && $i^2$ && $ 2$  && ${i}^{-1}$ &\cr
    height2pt&\omit&&\omit&&\omit&&\omit&\cr
    \noalign{\hrule}
    height2pt&\omit&&\omit&&\omit&&\omit&\cr
    & {$\matrix (i\!-\!1)^2\,(i\!\ge\! 2)\\ a_1>0\endmatrix$} && $i^2$ && $\dfrac 1 4$  && $\dfrac{2 i+1}{2(i+1)}$ &\cr
    height2pt&\omit&&\omit&&\omit&&\omit&\cr
    \noalign{\hrule}
    height2pt&\omit&&\omit&&\omit&&\omit&\cr
    & {$\matrix 2+(-1)^{i-1}\\(i\ge 2)\\a_1\!=\!\dfrac{7\!-\!\sqrt{33}}{2}\endmatrix$}
    && {$2 [2\!+\!(-1)^i]\!\!\!$} && $6-\sqrt{33}$ && $\dfrac{\sqrt{33}+(-1)^i}{8}$ &\cr
    height2pt&\omit&&\omit&&\omit&&\omit&\cr
  \noalign{\hrule}}}}}\hfill$$


We now go to the second part of this section. Consider the birth--death processes with a more general state space
$E=\{i: -M-1<i< N+1\}$, $M$, $N\le \infty$ and with Dirichlet boundaries
at $-M-1$ if $M>-\infty$ and at $N+1$ if $N<\infty$. Its $Q$-matrix
now is $q_{i, i+1}=b_i>0$, $q_{i, i-1}=a_i>0$, and $q_{ij}=0$ if
$|i-j|>1$ for $i, j\in E$. Fix a reference point $\theta\in E$.
Define
$$\gather
\mu_{\theta+n}=\frac{a_{\theta -1} a_{\theta -2}\cdots a_{\theta +n+1}}
 {b_{\theta } b_{\theta-1 }\cdots b_{\theta +n}}, \qqd -M-1-\theta< n\le -2,\\
\mu_{\theta-1 }=\frac{1}{b_{\theta } b_{\theta -1}},\qqd
\mu_{\theta }=\frac{1}{a_{\theta } b_{\theta }},\qqd
 \mu_{\theta+1 }=\frac{1}{a_{\theta } a_{\theta +1}},\tag 7.9\\
 \mu_{\theta+n }=\frac{b_{\theta +1} b_{\theta +2}\cdots b_{\theta +n-1}}
 {a_{\theta } a_{\theta+1 }\cdots a_{\theta +n}}, \qquad  2\le n<N+1-\theta.
\endgather$$
Correspondingly,
$$\gather
\text{\hskip-4em}D(f)=\sum_{-M-1<i\le \uz} \mu_i a_i (f_i-f_{i-1})^2
        +\sum_{\uz\le i <N+1} \mu_i b_i (f_{i+1}-f_{i})^2,\tag 7.10\\
\text{\hskip-4em}  f\in {\scr K},\; f_{-M-1}=0\;\text{if }M<\infty \text{ and }f_{N+1}=0\;\text{if }N<\infty.
\endgather$$

Let us begin with a particular application of \crl\;8.4 to ${\Bbb B}=L^1(\mu)$.

\proclaim{\crl\;7.8\,(Criterion and basic estimates)} Let $\lz_0$ be defined by $(7.1)$ with
the present state space $E$. Then we have $\kz^{-1}/4\le \lz_0\le\kz^{-1}$, where
$$\kz^{-1}=\inf_{m, n\in E:\; m \le n}\bigg[\bigg(\sum_{i=-M}^m\frac{1}{\mu_i a_i}\bigg)^{-1}
\!\!+ \bigg(\sum_{i=n}^{N}\frac{1}{\mu_i
b_i}\bigg)^{-1}\bigg]\bigg(\sum_{j=m}^n \mu_j\bigg)^{-1}. \tag 7.11$$
\endproclaim

By the way, we extend \crl\;6.6 to the present general state space.

\proclaim{\crl\;7.9\,(Criterion and basic estimates)} Let $\sum_{i\in E}\mu_i<\infty$ and define
$\lz_1$ as in (6.1). Then we have $\kz^{-1}/4\le \lz_1\le\kz^{-1}$, where
$$\kz^{-1}=\inf_{m, n\in E:\; m < n}\bigg[\bigg(\sum_{i=-M}^m \mu_i \bigg)^{-1}
\!\!+ \bigg(\sum_{i=n}^{N}\mu_i \bigg)^{-1}\bigg]
\bigg(\sum_{j=m}^{n-1} \frac{1}{\mu_j b_j}\bigg)^{-1}.\tag 7.12 $$
\endproclaim

\demo{\prf} When $M<\infty$, the corollary is simply a modification
of \crl\;6.6 by shifting the left end-point of the state space from
$0$ to $-M$. Thus, when $M=\infty$,
we can choose a sequence $\{M_p\}_{p=1}^\infty$ such that
$M_p\uparrow \infty$ as $p\uparrow \infty$ and then the assertion
holds if $M$ is replaced by $M_p$ for each $p$. In which case, the
corresponding $\lz_1$ is denoted by $\lz_1^{(M_p)}$ for a moment.
Because $\sum_{i\in E} \mu_i<\infty$, following the
proof above (4.2), it follows that
$$\lz_1\!=\inf\big\{D(f): \mu(f^2)=1,\; \mu(f)=0,\; f_i\!=f_{(i\vee m) \wedge n}
\text{ for some } m,n \in E,\; m<n\big\}.$$
Hence, we have $\lz_1^{(M_p)}\downarrow \lz_1 $ as $p\uparrow\infty$.
Similarly, replacing $M$ by $M_p$, we have the notation $\kz^{(M_p)}$. The proof
will be done once we show that
$${\big(\kz^{(M_p)}\big)}^{-1}\downarrow \kz^{-1}
\qqd\text{as } p\uparrow\infty. $$ Obviously, we have
$${\big(\kz^{(M_p)}\big)}^{-1}\downarrow
\qd\text{as } p\uparrow \qd\text{and}\qd
{\big(\kz^{(M_p)}\big)}^{-1}\ge \kz^{-1}. $$ To prove the required
assertion, let $\vz>0$. Then by definition of $\kz$ there exist
$m_0,\;n_0\in E$, $m_0<n_0$ such that
$$\bigg[\bigg(\sum_{i=-M}^{m_0} \mu_i \bigg)^{-1}
+ \bigg(\sum_{i=n_0}^{N} \mu_i \bigg)^{-1}\bigg]\bigg(\sum_{j=m_0}^{n_0-1} \frac{1}{\mu_j b_j} \bigg)^{-1}
\le \kz^{-1}+\vz.$$
Next, since $\sum_i \mu_i<\infty$, for fixed $m_0$, $n_0$ and large enough $M_p\,(-M_p< m_0)$, we have
$$\bigg[\bigg(\sum_{i=-M_p}^{m_0}  \mu_i\bigg)^{-1}
+ \bigg(\sum_{i=n_0}^{N} \mu_i \bigg)^{-1}\bigg]\bigg(\sum_{j=m_0}^{n_0-1} \frac{1}{\mu_j b_j} \bigg)^{-1}
\le \kz^{-1}+2\,\vz.$$
Combining these facts with the definition of $\kz^{(M_p)}$, we obtain
$$\align
\kz^{-1}&\le {\big(\kz^{(M_p)}\big)}^{-1}\\
&= \inf_{-M_p\le m<n<N+1}\bigg[\bigg(\sum_{i=-M_p}^{m}  \mu_i \bigg)^{-1}
+ \bigg(\sum_{i=n}^{N} \mu_i \bigg)^{-1}\bigg]\bigg(\sum_{j=m}^{n-1} \frac{1}{\mu_j b_j} \bigg)^{-1}\\
&\le \kz^{-1}+2\,\vz.\endalign$$
Since $\vz$ is arbitrary, we have proved that
${\big(\kz^{(M_p)}\big)}^{-1}\to \kz^{-1}$
as $p\to\infty. $
\qed\enddemo

For the remainder of this section, we study a splitting technique.
It provides a different tool to study the problem having bilateral Dirichlet boundaries.
This approach is especially meaningful if the duality discussed in Section 5 does not work,
such as in studying the processes on the whole ${\Bbb Z}$ or the Poincar\'e-type inequalities
given in the next section. We remark that Corollaries 7.8 and 7.9 use slightly
the splitting idea only (cf. Proof (b) of \thm\;8.2 below). The idea is splitting the
state space into two parts and then estimating the first (non-trivial) eigenvalue
in terms of the local ones. We have used this
technique several times before: Chen and Wang (1998) with Dirichlet boundary
for the unbounded region, Chen, Zhang and Zhao (2003), as well as Mao and Xia (2009),
with Neumann boundary. The first and the third papers work on a very general setup. Here,
we follow the second one with some addition.

To state our result, we need to construct two birth--death processes on the left- and
the right-hand sides, respectively, for a given birth--death process with rates $(a_i, b_i)$
and state space $E$. Fix a constant $\gz>1$.
\roster
\item "(L)" The process on the left-hand side has state space $E^{\uz -}=\{i: -M-1<i\le \uz\}$, reflects at $\uz$.
Its transition structure is the same as the original one except $a_{\uz}$ is replaced by $\gz a_{\uz}$.
\item "(R)" The process on the right-hand side has state space $E^{\uz +}=\{i: \uz\le i <N+1\}$,
reflects at $\uz$. Its transition structure is again
the same as the original one except $b_{\uz}$ is replaced by $\gz (\gz-1)^{-1} b_{\uz}$.
\endroster
For the process on the right-hand side, the state $\uz$ is a Neumann boundary.
At $N+1$, it is a Dirichlet boundary if $N<\infty$. For this process, the first eigenvalue,
denoted by $\lz_0^{\uz +, \gz}$, has already been
studied in Sections 2 and 3. With a change of the order of the
state space, it follows that the process on the left-hand side has the same boundary
condition, denote by $\lz_0^{\uz -, \gz}$ its first eigenvalue. Note that ignoring a finite number of the
states does not change the positivity of $\lz_0^{\uz \pm, \gz}$, in the qualitative case, we simply denote
them by $\lz_0^{(\pm)}$, respectively.
In general, according to $\sum_{i=\uz}^N (\mu_i b_i)^{-1}$ and/or $\sum_{i=-M}^\uz (\mu_i a_i)^{-1}$,
$\sum_{i=\uz}^N \mu_i $ and/or $\sum_{i=-M}^\uz \mu_i$ being finite or not, there
are eight cases for the processes on ${\Bbb Z}$. For instance,
if $\sum_{i=\uz}^N (\mu_i b_i)^{-1}=\infty$, then $\lz_0^{(+)}=0$ by \thm\;3.1.
Since in this section, we are working on bilateral  Dirichlet boundaries, it is natural
to assume that $\lz_0^{(\pm)}>0$. The other cases may be treated in a parallel way.
For instance, when $\lz_0^{(-)}=0$, it is more natural to consider the process on $[-M, N+1)$ with reflecting
at some finite $-M$ and then pass to the limit as $M\to\infty$ (cf. the proof of \crl\;7.8).
In this case, the eigenfunction should be strictly decreasing once $\lz_0>0$. Hence, there is no
reason to use the splitting technique. Note that the explicit criterion for
$\lz_0^{(\pm)}>0$ is given by \thm\;3.1. We can now state the main result of the second part
of this section as follows.

\proclaim{\thm\;7.10} \roster
\item
In general, the Dirichlet eigenvalue $\lz_0$ of the birth--death
process on \newline $E=\{i: -M-1<i<N+1\}$ satisfies
$$\inf_{\uz\in E}\, \inf_{\gz>1} \big(\lz_0^{\uz -, \gz} \vee \lz_0^{\uz +, \gz}\big)=
\lz_0\ge\!\! \sup_{-M-1\le \uz \le N+1}\, \sup_{\gz>1} \big(\lz_0^{\uz -, \gz} \wedge \lz_0^{\uz +, \gz}\big),\tag 7.13$$
where on the right-hand, when $\uz=-M-1$, define $\lz_0^{\uz -, \gz}=\infty$, and $\lz_0^{\uz +, \gz}$ to
be the first eigenvalue of the original process (independent of $\gz$) reflected at $-M$ if $M<\infty$; when $\uz=N+1$, define $\lz_0^{\uz +, \gz}=\infty$,
and $\lz_0^{\uz -, \gz}$ to be the one reflected at $N$ if $N<\infty$.
\item The second equality in (7.13) replacing $\sup_{-M-1\le \uz \le N+1}$
by $\sup_{\uz \in E }$ also holds provided $\lz_0^{(\pm)}>0$, and moreover,
$$\sum_{i=-M}^\uz \mu_i=\infty\;\text{ if }M=\infty\qd \text{and}\qd
  \sum_{i=\uz}^N \mu_i=\infty\;\text{ if }N=\infty. \tag 7.14$$
\endroster
\endproclaim

\thm\;7.10 was proved in Chen, Zhang and Zhao (2003) for the half-space (i.e., one of $M$ and $N$
is finite), under the hypotheses that $\sum_{i} (\mu_i b_i)^{-1}<\infty$ and $\sum_{i} \mu_i<\infty$
which is essentially the case of having a finite state space.

To prove \thm\;7.10, we need some preparation. First, we couple these two processes on a common state space
${\overline E}=\{i: -M-1<i<N+2\}$. Next, separate the two processes by shifting
the state space $E^{\uz+}$ by one to the right: $1+ E^{\uz+}$. Denote by $({\bar a}_i, {\bar b}_i)$
the rates of the connected process. For this, we need to build a bridge for the processes on the two sides by adding two more rates
${\bar b}_{\uz}=\gz-1$ and ${\bar a}_{\uz+1}=1$. The construction here will become clear once we have a
deeper understanding about the eigenfunction and it will be explained in Part I\!I of the proof
of the theorem. Roughly speaking, there are two possible shapes of the eigenfunction, the
construction enables us to transform one of them to the other so that the splitting with Neumann
boundaries becomes practical.
For which, one needs the parameter $\gz$ as shown in \lmm\;7.12 below.
In detail, we now have
$$ {{\bar a}_i=\cases
a_i,  &\!\!\!\!\!-M-1<i\le \uz-1,\\
\gz a_{\uz}, & i=\uz,\\
1, &i=\uz+1,\\
a_{i-1}, & \uz+2\le i<N+2,
\endcases}\qqd
{{\bar b}_i=\cases
b_i,  &\!\!\!\!\!-M-1<i\le \uz-1,\\
\gz-1 , & i=\uz,\\
\dfrac{\gz b_{\uz}}{\gz -1}, &i=\uz+1,\\
b_{i-1}, & \uz+2\le i<N+2.
\endcases}
$$
Applying (7.9) to the present setup and removing the factor $b_{\uz}
(1-\gz)^{-1}$ (which simplifies the notation but does not change the
ratio ${\overline D}(f)/{\bar\mu}(f^2)$), we obtain
$$\aligned
{\bar\mu}_{i}&= \mu_{i},\qqd\qd -M-1<i\le \uz-1,\\
{\bar\mu}_{\uz}&=\frac{1}{\gz}\, \mu_{\uz},\qqd
{\bar\mu}_{\uz+1}=\frac{\gz-1}{\gz}\, \mu_{\uz},\qqd\qqd\qd\\
{\bar\mu}_{i}&= \mu_{i-1},\qqd \uz+2\le i<N+2.
\endaligned\tag 7.15$$
Then
$${\bar \mu}_i {\bar a}_i=\mu_i a_i,\;\; i\le \uz,\qd {\bar \mu}_{\uz}{\bar b}_{\uz}=\frac{\gz-1}{\gz}\mu_{\uz},
\qd {\bar \mu}_i{\bar b}_i=\mu_{i-1} b_{i-1},\;\; i\ge \uz+1. \tag 7.16$$

The next two results are basic in using the splitting technique.

\proclaim{\lmm\;7.11} Given $f$ on $E$, define
${\bar f}$ on ${\overline E}$ as follows: ${\bar f}_i=f_i $ for
$i\le \uz$ and ${\bar f}_i=f_{i-1} $ for $i\ge \uz+1$. Then we have
${\bar \mu} \big({\bar f}^2\big)=\mu\big(f^2\big)$ and ${\overline D}({\bar f})=D(f)$.
\endproclaim

\demo{\prf} Clearly, we have ${\bar f}_{\uz}={\bar f}_{\uz+1}$. Then
$$\align
{\bar \mu} \big({\bar f}^2\big)&=\sum_{-M-1<i\le \uz-1}\mu_i f_i^2
+\big({\bar \mu}_{\uz}+{\bar \mu}_{\uz+1}\big)f_{\uz}^2 +
\sum_{\uz+2\le i<N+2}\mu_{i-1} f_{i-1}^2 =\mu\big(f^2\big),\\
{\overline D}({\bar f})&=\sum_{-M-1< i \le \uz} {\bar\mu}_i {\bar
a}_i \big({\bar f}_i-{\bar f}_{i-1}\big)^2 +\sum_{\uz+1\le i <N+2}
{\bar\mu}_i {\bar b}_i
\big({\bar f}_{i+1}-{\bar f}_{i}\big)^2\tag 7.17\\
&=\sum_{-M-1< i \le \uz} \mu_i a_i (f_i-f_{i-1})^2+ \sum_{\uz+1\le i
<N+2} \mu_{i-1}b_{i-1}(f_i-f_{i-1})^2\\
&=D(f)\qqd\text{(by (7.10))}.\qed
\endalign$$
\enddemo

\proclaim{\lmm\;7.12} For a given birth--death process with state
space $E$ and rates $(a_i, b_i)$, if its eigenfunction $g$ of $\lz$
satisfies $g_{\uz-1}<g_{\uz}>g_{\uz+1}$ (resp.
$g_{\uz-1}>g_{\uz}<g_{\uz+1}$) for some $\uz\in E$ (of course,
$g_{-M-1}^{}=0$ if $M<\infty$, and $g_{N+1}^{}=0$ if $N=\infty$), let
$$\gz=1+ \frac{b_{\uz}(g_{\uz}-g_{\uz+1})}{a_{\uz}(g_{\uz}-g_{\uz-1})}>1,\tag 7.18$$
and let $\bar g_i=g_i$ for $i\le \uz$, $\bar g_i=g_{i-1}$ for
$i\ge \uz+1$. Then for the $(\bar a_i, \bar b_i)$-process, $\bar
g$ is the eigenfunction of $\bar\lz=\lz$ having the property $\bar
g_{\uz+1}=\bar g_{\uz}$. Furthermore, $\bar g|_{(-M-1, \uz]}$ is the
eigenfunction of ${\bar\lz}$ of the process on the left-hand side reflecting at $\uz$,
and  similarly $\bar g|_{[\uz+1, N+1)}$ is the eigenfunction of the process on the right-hand side reflecting at $\uz+1$.\endproclaim

\demo{\prf} By the construction of $(\bar a_i, \bar b_i)$ and $\bar
g$, we have
$${\overline\ooz}\,{\bar g} (i)=
\cases
\ooz g(i)=-\lz g_i=-{\bar\lz} {\bar
g}_i,\qqd\qqd\;\;\, & i\le \uz-1,\\
\ooz g(i-1)=-\lz g_{i-1}=-{\bar\lz}
{\bar g}_i,\qqd & i\ge \uz+2. \endcases$$
Next, by (7.15), we have
$$\align
{\overline\ooz}\,{\bar g}(\uz)={\bar b}_{\uz} ({\bar g}_{\uz+1}&-
{\bar g}_{\uz}) + {\bar a}_{\uz} ({\bar g}_{\uz-1}- {\bar g}_{\uz})
={\bar a}_{\uz} (g_{\uz-1}-g_{\uz})
=\gz a_{\uz} (g_{\uz-1}-g_{\uz}),\\
{\overline\ooz}\,{\bar g} (\uz+1)&={\bar b}_{\uz+1} ({\bar
g}_{\uz+2}- {\bar g}_{\uz+1}) + {\bar a}_{\uz+1} ({\bar g}_{\uz}-
{\bar g}_{\uz+1})\\
&={\bar b}_{\uz+1} (g_{\uz+1}-g_{\uz})\\
 &=\frac{\gz}{\gz-1} b_{\uz} (g_{\uz+1}-g_{\uz}).
\endalign$$
In the first formula, the term containing ${\bar b}_{\uz}$ vanishes. This is the reason why we can
regard $\uz$ as a reflecting boundary for the process on the left-hand side. Similarly, one can
regard $\uz+1$ as the one for the process on the right-hand side in view of the second formula.
By (7.18), the right-hand sides are the same which is equal to
$$\bigg[1+ \frac{b_{\uz}(g_{\uz}-g_{\uz+1})}{a_{\uz}(g_{\uz}-g_{\uz-1})}\bigg] a_{\uz} (g_{\uz-1}-g_{\uz})
=a_{\uz} (g_{\uz-1}-g_{\uz})+b_{\uz} (g_{\uz}-g_{\uz+1})=-\lz g_{\uz}=-{\bar\lz} {\bar g}_{\uz}.
$$
We have thus proved the lemma.\qed\enddemo

\demo{\prf\; of \thm\;$7.10$. Part I} In this part, we prove \thm\;7.10\,(1)
with the first ``$=$'' replaced by ``$\ge $''.
The proof of this part is relatively easier.
Let $f\in {\scr K}$, $f\ne 0$. Define $\bar f$ as in \lmm\;7.11. For fixed $\uz\in {E}$ and $\gz>1$,
noting that re-labeling the state space does not change $\lz_0^{\uz +, \gz}$,
by (7.17), we have
$${\overline D}({\bar f}) \ge \lz_0^{\uz -, \gz} \sum_{i \le \uz}
{\bar\mu}_i {\bar f}_i^2
+\lz_0^{\uz +, \gz} \sum_{i \ge \uz+1} {\bar\mu}_i {\bar f}_i^2
\ge \big(\lz_0^{\uz -, \gz} \wedge \lz_0^{\uz +,
\gz}\big)\,{\bar\mu} \big({\bar f}^2\big).$$
Hence by \lmm\;7.11,
$$\frac{D(f)}{\mu(f^2)}=\frac{{\overline D}({\bar f})}{{\bar\mu}\big({\bar
f}^2\big)} \ge \lz_0^{\uz -, \gz} \wedge \lz_0^{\uz +, \gz}.$$
Making the supremum with respect to $\gz$ and $\uz$,
it follows that
$$\frac{D(f)}{\mu(f^2)}\ge \sup_{-M-1< \uz < N+1}\, \sup_{\gz>1}\big(\lz_0^{\uz -, \gz} \wedge \lz_0^{\uz +, \gz}\big). $$
At the boundaries, say $\uz=-M-1$ for instance, by (7.10) and the convention, we have
$$D(f)\ge \sum_{-M-1< i< N+1} \mu_i b_i (f_{i+1}-f_i)^2\ge \lz_0^{\uz+, \gz}\, \mu\big(f^2\big)
=\big[\lz_0^{\uz-, \gz}\wedge \lz_0^{\uz+, \gz}\big] \mu\big(f^2\big).$$
Therefore, we indeed have
$$\frac{D(f)}{\mu(f^2)}\ge \sup_{-M-1\le \uz \le N+1}\, \sup_{\gz>1}\big(\lz_0^{\uz -, \gz} \wedge \lz_0^{\uz +, \gz}\big).$$
Making infimum with respect to $f$, we obtain
$$\lz_0=\inf_{f\in {\scr K},\; f\ne 0}\frac{D(f)}{\mu(f^2)}
\ge \sup_{-M-1\le \uz \le N+1}\, \sup_{\gz>1}\big(\lz_0^{\uz -, \gz} \wedge \lz_0^{\uz +, \gz}\big). $$
This proves the (second) inequality in (7.13).

To prove the upper estimate, fix $\uz\in {E}$ and $\gz>1$ again.
As we have seen from the last part of proof (g) of \thm\;2.4 and
\prp\;2.5, if we let $\lz_0^{\uz +, \gz,\, n}$ denote the local eigenvalue with
Neumann boundary at $\uz$ and Dirichlet boundary at $n+1$, then $\lz_0^{\uz +, \gz,\, n}\downarrow
\lz_0^{\uz +, \gz}$ as $n\uparrow \infty$. Thus, for each
$\vz>0$, we have $ \lz_0^{\uz +,
\gz,\, n}< \lz_0^{\uz +, \gz} +\vz$ for large enough $n$. By
\prp\;2.2, we can assume that the corresponding eigenfunction
$g^{(+,\, n)}$ of $\lz_0^{\uz +, \gz,\, n}$ satisfies $g^{(+,\, n)}_{\uz}=1$ and
$g^{(+,\, n)}_{i}=0$ for all $i>n\,(>\uz)$. Similarly, we have $ \lz_0^{\uz
-, \gz,\, m}< \lz_0^{\uz -, \gz} +\vz$ for small enough $-m$, and
moreover, the eigenfunction $g^{(-,\, m)}$ of $\lz_0^{\uz -, \gz,\, m}$
satisfies $g^{(-,\, m)}_{\uz}=1$ and $g^{(-,\, m)}_{i}=0$ for all $i<-m\,(< \uz)$.
Let ${\bar f}$ be defined as above, connecting $g^{(-,\, m)}$ and $g^{(+,\, n)}$. Then
${\bar f}$ has a finite support, ${\bar f}_{\uz}={\bar
f}_{\uz+1}=1$, and moreover by (7.10),
$$\align
{\overline D}({\bar f})&=\sum_{{\overline E}\,\ni i \le \uz}
{\bar\mu}_i {\bar a}_i \big({\bar f}_i-{\bar f}_{i-1}\big)^2
+\sum_{{\overline E}\,\ni i \ge \uz+1} {\bar\mu}_i {\bar b}_i
\big({\bar f}_{i+1}-{\bar f}_{i}\big)^2\\
&= \lz_0^{\uz -, \gz,\, m} \sum_{i \le \uz} {\bar\mu}_i {\bar f}_i^2
+\lz_0^{\uz +, \gz,\, n} \sum_{i \ge \uz+1} {\bar\mu}_i {\bar f}_i^2\\
&\le \big(\lz_0^{\uz -, \gz} \vee \lz_0^{\uz +,
\gz}+\vz\big)\,{\bar\mu} \big({\bar f}^2\big).
\endalign
$$
By \lmm s 7.11 and 7.12, this gives us
$$\lz_0={\bar\lz}_0\le \lz_0^{\uz -, \gz} \vee
\lz_0^{\uz +, \gz}$$
since $\vz$ is arbitrary. Furthermore, we have
$$\lz_0\le \inf_{\uz\in E}\,\inf_{\gz>1}\big(\lz_0^{\uz -, \gz} \vee
\lz_0^{\uz +, \gz}\big)$$
as required.\qed\enddemo

The proof of the equalities in \thm\;7.10 is much harder. For which, we need once again
a deeper understanding of the eigenfunction of $\lz_0$. To have a concrete impression,
we mention that the eigenfunction in \xmp s 7.6\,(2) is $g_0=g_3=0$, $g_1=(1+\vz)g_2$.
Thus, when $\vz=0$, we have $g_1=g_2$. Besides, it is rather easy to see the shape of
eigenfunction $g$ of the examples given in Table 7.1 since $v_i<1$ iff $g_{i+1}<g_i$
for all $i$.

\proclaim{\defn\;7.13}\roster
\item A function $f$ is said to be {\it unimodal} if there exists a finite $k$ such that
$f_i$ is strictly increasing for $i\le k$ and strictly decreasing for $i\ge k$.
\item A function $f$ is said to be a {\it simple echelon} if there exists a $k$ such that
$f_k=f_{k+1}$, $f_i$ is strictly increasing for $i\le k$ and strictly decreasing for $i\ge k+1$.
\endroster\endproclaim

\proclaim{\prp\;7.14} Let $g$ be a positive eigenfunction of $\lz>0$
for a birth--death process. Then $g$ is
strictly monotone, or unimodal, or a simple echelon.
\endproclaim

\demo{\prf} (a)
Let $g_k\ge g_{k+1}$ for some $k$. We prove that $g$ is strictly decreasing for $i\ge k+1$.
To do so, note that
$$b_{k+1}(g_{k+2}-g_{k+1})=-\lz g_{k+1} -a_{k+1} (g_{k}-g_{k+1})\le -\lz g_{k+1}<0.$$
Thus, we have $g_{k+2}<g_{k+1}$. Assume that $g_n< g_{n-1}$ for some $n\ge k+2$. Then the eigenequation
shows that
$$b_{n}(g_{n+1}-g_n)=-\lz g_n -a_n (g_{n-1}-g_n)<-\lz g_n<0.$$
By induction,
this gives us $g_{n+1}<g_n$ for all $n\ge k+1$.

(b) By symmetry, we can handle with the case that $g_k\le g_{k+1}$ for some $k$.
One starts at
$$a_{k}(g_{k-1}-g_{k})=-\lz g_{k} -b_{k} (g_{k+1}-g_{k})\le -\lz g_{k}<0.$$
We obtain $g_{k-1}<g_{k}$ and then $g_{n-1}<g_{n}$ for all $n\le k$ by induction.

(c) By (a) and (b), it follows that there is no local convex part of $g$. Otherwise, there is a
$k$ such that either $g_{k-1}> g_k < g_{k+1}$ or $g_{k-1}> g_k = g_{k+1}<g_{k+2}$ which
contradict what we proved in (a) and (b).

(d) We claim that for every $k$, say $k=0$ for simplicity, the two cases
``$g_{-1}\ge g_0$'' and ``$g_0\le g_1$'' cannot happen at the same time. Otherwise, there are
four situations:
$$g_{-1}=g_0=g_1,\qd g_{-1}>g_0< g_1,\qd
g_{-1}>g_0= g_1,\qd \text{and}\qd g_{-1}=g_0< g_1.$$
The first one cannot happen,
otherwise we have $g_i\equiv 0$. By (c), the second case is impossible.
The last two cases are also impossible by (b) and (a), respectively.

(e) Having these preparations at hand, we are ready to prove the main assertion of the proposition.
Clearly, we need only to consider the case that $g$ is not
strictly monotone.
Choose a starting point, say 0 for instance.
By (d), we have only one possibility: either $g_{-1}\ge g_0$ or $g_0\le g_1$.
Without loss of generality, assume that $g_0\le g_1$.
If $g_0= g_1$, then by (a) and (b), $g$ is a simple echelon.
If $g_0< g_1$, then on the one hand, by (b), $g_i$ is strictly increasing for all $i\le 1$, and on
the other hand, we can find a $k\ge 1$ such that $g_1< g_2<\ldots < g_k \ge g_{k+1}$ since $g$ is not
strictly monotone by assumption. Applying (a) again, it follows that $g$ is strictly
decreasing for all $i\ge k+1$. Hence, $g$ is either unimodal or a simple echelon.
\qed\enddemo

\proclaim{\prp\;7.15} For the birth--death process on ${\Bbb Z}$, the
following assertions hold. \roster
\item The eigenfunction $g$ of $\lz$ satisfies the following
successive formulas:
$$\aligned
\text{\hskip-2em}g_{k+1}&=g_k +
\frac{1}{\mu_k b_k}\bigg[\mu_{\uz}a_{\uz}(g_{\uz}-g_{\uz-1})-\lz\sum_{i=\uz}^k \mu_i g_i\bigg],\qqd k\ge \uz,\\
\text{\hskip-2em}g_{k-1}&=g_k +
\frac{1}{\mu_k a_k}\bigg[\mu_{\uz}a_{\uz}(g_{\uz-1}-g_{\uz})-\lz\sum_{i=k}^{\uz-1} \mu_i g_i\bigg],\qqd k< \uz.
\endaligned\tag 7.19$$
\item If $\lz=0$, then the non-trivial eigenfunction $g$
with $g_{\uz}=1$ for some $\uz\in {\Bbb Z}$ is given by
$$\align
g_n &=\!1+(1-g_{\uz-1})\sum_{j=\uz}^{n-1} \prod_{k=\uz}^j
\frac{a_k}{b_k},\qqd  n\ge \uz,\\
g_n &=\!1-(1-g_{\uz-1})\sum_{j=n}^{\uz-1} \prod_{k=j+1}^{\uz-1}
\frac{b_k}{a_k},\qqd  n< \uz. \endalign$$ In this case, the
function $g$ is either the constant function $\dbl$ or strictly
monotone on ${\Bbb Z}$.
\item If $\lz>0$ and
$$\sum_{i=-\infty}^\uz \mu_i=\sum_{i=\uz}^\infty \mu_i=\infty, \tag 7.20$$
then the non-trivial eigenfunction $g$ of $\lz$ cannot be monotone.
\endroster
\endproclaim

\demo{\prf} (a) Part (1) of the proposition follows from the
eigenequation.

(b) When $\lz=0$, with $u_i:=g_{i+1}-g_i\,(i\in {\Bbb Z})$, the eigenequation
$b_i u_i=a_i u_{i-1}$
gives us
$$u_j=(1-g_{\uz-1})\prod_{k=\uz}^j \frac{a_k}{b_k},\qd j\ge \uz,\qqd
u_j=(1-g_{\uz-1})\prod_{k=j+1}^{\uz-1} \frac{b_k}{a_k},\qd j< \uz.$$
It follows that either $g_i\equiv 1$ or $g$ is strictly monotone on ${\Bbb Z}$.
Now, part (2) of the proposition follows by making a summation of $j$ from $\uz$ to $n-1$
and from $n$ to $\uz-1$, respectively.

(c) Without loss of generality, assume that $g_{\uz}=1$ for some $\uz\in {\Bbb Z}$.
Suppose that $g$ is non-decreasing, then by the first equation in part (1),
we would have
$$\infty>\frac{\mu_{\uz}a_{\uz}(g_{\uz}-g_{\uz-1})}{\lz}\ge
\sum_{k=\uz}^n \mu_k g_k\ge \sum_{k=\uz}^n \mu_k\to \infty\qqd\text{as }n\to\infty.$$
Otherwise, if $g$ is non-increasing, then by the second equation in part (1),
we would have
$$\infty>\frac{\mu_{\uz}a_{\uz}(g_{\uz-1}-g_{\uz})}{\lz}\ge
\sum_{k=n}^{\uz-1} \mu_k g_k\ge \sum_{k=n}^{\uz-1} \mu_k\to \infty\qqd\text{as }n\to -\infty.
$$
We have thus proved part (3) of the proposition.\qed
\enddemo

We remark that \prp\;7.15\;(2) is different from \prp\;2.2
where the eigenfunction of $\lz=0$ must be a constant. Here is a simple example with $\uz=0$:
$a_i=b_i=|i|$ if $i\ne 0$ and $a_0=b_0=1$, then corresponding to $\lz=0$, we have a family
of linear eigenfunctions $\{g_i^{(\gz)}=1+(1-\gz)\, i: i\in {\Bbb Z}\}_{\gz \in{\Bbb R}}$
(normalized at $0$) with one-parameter $\gz$.

\proclaim{\prp\;7.16} Let $(7.14)$ hold and $g$ be a non-zero eigenfunction of $\lz_0>0$. Then
$g$ is either positive or negative on $E$.
\endproclaim

 \demo{\prf} If one of $M$ or $N$ is finite, then the conclusion
 follows from \prp\;2.2\,(1). From now on in the proof, assume that
 $M=N=\infty$.

(a) If the
conclusion of the proposition does not hold, then there is a $k$
(say) such that $g_{k}\le 0$ and either $g_{k-1}>0$ or $g_{k+1}>0$.
By symmetry, assume that $g_{k+1}>0$.

(b) We now prove that $g_i>0$ for all $i\ge k+1$. Given $m, n\in
{\Bbb Z}$ with $m\le n$, denote by $\lz_0^{[m, n]}$ the first
eigenvalue of the process restricted on the state space $\{i: m\le i
\le n\}$ with Dirichlet boundaries at $m-1$ and $n+1$ in the sense
similar to (7.1). If the assertion does not hold, then there is a
$k_0: k_0 > k +1$ such that $g_{k_0}\le 0$. Now, let $\tilde g$
satisfy ${\tilde g}_{k}=0$, ${\tilde
g}_i=g_i$ for $i=k+1, \ldots, k_0-1$, ${\tilde g}_{k_0}=\vz$ for some
$\vz>0$, ${\tilde g}_i=0$ for $i\ge k_0+1$. Note that
$$\align
\big(-\ooz\, {\tilde g}\big)(k+1)&=b_{k+1}\big({\tilde
g}_{k+1}-{\tilde g}_{k+2}\big)
+a_{k+1}\big({\tilde g}_{k+1}-{\tilde g}_{k}\big)\\
&=b_{k+1}(g_{k+1}-g_{k+2})+a_{k+1}(g_{k+1}-g_{k})+a_{k+1}g_{k}\\
&=\lz_0 g_{k+1}+a_{k+1}g_{k}\\
&\le \lz_0 \,{\tilde g}_{k+1},\endalign$$
Because of $\lz_0>0$ and following proof (b) of \prp\;2.1,
we can choose a suitable $\vz>0$ such that
$$\sum_{i=k+1}^{k_0}\mu_i {\tilde g}_i\big(-\ooz {\tilde g}\big)(i)
< \lz_0 \sum_{i=k+1}^{k_0} \mu_i {\tilde g}_i^2.$$
It follows that
$\lz_0^{[k+1, k_0]}<\lz_0$. However, it is obvious that $\lz_0\le
\lz_0^{[k+1, k_0]}$ and so we get a contradiction. We have thus
proved that $g_i>0$ for all $i\ge k+1$.

(c) By (7.14) and proof (c) of \prp\;7.15, $g$ cannot be
non-decreasing since $\lz_0>0$. Hence, there is a $\uz\ge k+2$
such that $g_{k+1}< g_{k+2}<\ldots < g_{\uz}\ge g_{\uz +1}$. In the
case that $g_{\uz}> g_{\uz +1}$, by introducing an additional
point but keeping the same $\lz_0$ as shown in \lmm\;7.12, one can reduce to the case that
$g_{\uz}= g_{\uz +1}$. Hence, one can split the original process
into two as in (L) and (R). Now, starting from $\uz$ at which $g_{\uz}>0$, look at the
process on the left-hand side in the inverse way, one finds the point $k<\uz$ at which
$g_k\le 0$. Applying proof (b) above
to this process, one may get a contradiction. It follows that $g>0$ on
$(-\infty, \uz]\supset (-\infty, k]$.

Therefore, we should have $g>0$ on ${\Bbb Z}$. \qed\enddemo

\demo{\prf\; of \thm\;$7.10$. Part I\!I} We now prove the equalities in (7.13).
By assumption $\lz_0^{(\pm)}>0$ and the second inequality in (7.13),
it follows that $\lz_0>0$. If one of $M$ and $N$ is finite, then the
non-trivial eigenfunction $g$ must be positive by \prp\;2.2\,(1). In this case, it is helpful to include
the boundary into the domain of $g$ for understanding its shape. Then by \prp\;7.14, there are only two possibilities:
\roster
\item "(i)" $g$ is unimodal;
\item "(ii)" $g$ is a simple echelon.
\endroster
Next, if $M=N=\infty$, then by \prp\;7.16, we have $g>0$.
Moreover, by \prp\;7.15, $g$ cannot be monotone. Hence, by \prp\;7.14, $g$ has
again one of shapes (i) and (ii) as above.

We now prove the equalities in (7.13) only in the case that $M=N=\infty$. The proof for the other case
is simpler.

(a) Case (ii). We use the operator $I\!I$ defined in Section 2:
$${\overline {I\!I}}_i^{\uz+, \gz}\big({\bar f}\big)=\frac{1}{\bar f_i} \sum_{j=i}^{N+1} \frac{1}{{\bar \mu}_j {\bar b}_j}
\sum_{k=\uz+1}^j {\bar \mu}_k {\bar f}_k,\qqd \uz+1\le i <N+2.$$
For each ${\bar f}$ satisfying: ${\bar f}_i=f_i$ for $i\le \uz$ and ${\bar f}_i=f_{i-1}$ for $i\ge \uz+1$ for some
$f$ on $E$, by (7.15) and (7.16), we have
$$\align
{\overline {I\!I}}_i^{\uz+, \gz}\big({\bar f}\big)&=\frac{1}{f_{i-1}}\sum_{j=i}^{N+1} \frac{1}{\mu_{j-1} b_{j-1}}
\bigg[{\bar\mu}_{\uz+1} f_{\uz}+\sum_{k=\uz+2}^j \mu_{k-1}f_{k-1} \bigg]\\
&=\frac{1}{f_{i-1}}\sum_{j=i}^{N+1} \frac{1}{\mu_{j-1} b_{j-1}}
\bigg[\bigg(1-\frac{1}{\gz}\bigg){\mu}_{\uz} f_{\uz}+\sum_{k=\uz+1}^{j-1} \mu_{k}f_{k} \bigg]\\
&=\frac{1}{f_{i-1}}\sum_{j=i-1}^{N} \frac{1}{\mu_{j} b_{j}} \sum_{k=\uz}^j \mu_{k}f_{k}
-\frac{\mu_{\uz} f_{\uz}}{\gz f_{i-1}}\sum_{j=i-1}^{N} \frac{1}{\mu_{j} b_{j}}\\
&=\frac{1}{f_{i-1}}\sum_{j=i-1}^{N} \frac{1}{\mu_{j} b_{j}} \sum_{k=\uz+1}^j \mu_{k}f_{k}
+\bigg(1-\frac{1}{\gz}\bigg)\frac{\mu_{\uz} f_{\uz}}{f_{i-1}}\sum_{j=i-1}^{N} \frac{1}{\mu_{j} b_{j}}\\
&\text{\hskip 8em} \uz+1\le i <N+2.  \tag 7.21
\endalign$$
Similarly, we have
$$\align
{\overline {I\!I}}_i^{\uz-, \gz}\big({\bar f}\big)&=
\frac{1}{{\bar f}_i}\sum_{j=-M}^i \frac{1}{{\bar \mu}_j {\bar a}_j}\sum_{k=j}^{\uz} {\bar \mu}_k {\bar f}_k\\
&=\frac{1}{f_i}\sum_{j=-M}^i \frac{1}{\mu_j a_j}\sum_{k=j}^{\uz} \mu_k f_k
-\bigg(1-\frac{1}{\gz}\bigg)\frac{\mu_{\uz} f_{\uz}}{f_i}\sum_{j=-M}^i \frac{1}{\mu_j a_j},\\
&\text{\hskip 8em} -M-1<i\le \uz. \tag 7.22
\endalign$$
Because $g_{\uz}=g_{\uz+1}$, we can regard $\uz$ as a Neumann boundary of the original process restricted on the left--hand side
and at the same time, regard $\uz$ as a Neumann boundary of the original process restricted on the right--hand side.
Because $\lz_0^{(\pm)}>0$, by \prp\;2.5\,(2), we have $g_{\pm\infty}=0$. Hence, by (2.11),
(7.21), and (7.22), we obtain
$$\align
{\overline {I\!I}}_i^{\uz+, \gz}\big({\bar g}\big)&=
\frac{1}{\lz_0}+\bigg(1-\frac{1}{\gz}\bigg)\frac{\mu_{\uz} g_{\uz}}{g_{i-1}}\sum_{j=i-1}^{N} \frac{1}{\mu_{j} b_{j}},\qqd
 \uz+1\le i <N+2,\\
{\overline {I\!I}}_i^{\uz-, \gz}\big({\bar g}\big)&=
\frac{1}{\lz_0}-\bigg(1-\frac{1}{\gz}\bigg)\frac{\mu_{\uz} g_{\uz}}{g_i}\sum_{j=-M}^i \frac{1}{\mu_j a_j},
\qqd -M-1<i\le \uz.
\endalign$$
By \prp\;2.2\,(2), we have
$$\sup_{\uz\le i <N+1}\frac{\mu_{\uz} g_{\uz}}{g_{i}}\sum_{j=i}^{N} \frac{1}{\mu_{j} b_{j}}\le \frac{1}{\lz_0},\qqd
\sup_{-M-1<i\le \uz}\frac{\mu_{\uz} g_{\uz}}{g_{i}}\sum_{j=-M}^i \frac{1}{\mu_j a_j}\le \frac{1}{\lz_0}.$$
Therefore, by the second inequality in (7.13) and \thm\;2.4\,(3), it follows that
$$\align
\lz_0& \ge \sup_{\uz^{\prime}\in E} \,\sup_{\gz>1}\Big[\lz_0^{\uz^{\prime}-, \gz}\wedge \lz_0^{\uz^{\prime}+, \gz}\Big]\\
&\ge \sup_{\gz>1}\Big[\lz_0^{\uz-, \gz}\wedge \lz_0^{\uz+, \gz}\Big]\\
&\ge \sup_{\gz>1}\Big[\Big(\inf_{-M-1<i\le \uz} {\overline {I\!I}}_i^{\uz-, \gz}\big({\bar g}\big)^{-1}\Big)
\wedge \Big(\inf_{\uz+1\le i <N+2}{\overline {I\!I}}_i^{\uz+, \gz}\big({\bar g}\big)^{-1}\Big)\Big]\\
&= \sup_{\gz>1}\, \inf_{\uz+1\le i <N+2}{\overline {I\!I}}_i^{\uz+, \gz}\big({\bar g}\big)^{-1}\\
&=\sup_{\gz>1}\bigg\{\frac{1}{\lz_0}+\bigg(1-\frac{1}{\gz}\bigg)\sup_{\uz\le i <N+1}
\frac{\mu_{\uz} g_{\uz}}{g_{i}}\sum_{j=i}^{N} \frac{1}{\mu_{j} b_{j}}\bigg\}^{-1}\\
&=\lz_0.
\endalign
$$
We have thus proved  in Case (ii) the second equality in (7.13).

To prove the first equality in (7.13), noting the inequality was proved in Part I, we have dually
$$\align
\lz_0& \le \inf_{\uz^{\prime}\in E} \,\inf_{\gz>1}\Big[\lz_0^{\uz^{\prime}-, \gz}\vee \lz_0^{\uz^{\prime}+, \gz}\Big]\\
&\le \inf_{\gz>1}\Big[\lz_0^{\uz-, \gz}\vee \lz_0^{\uz+, \gz}\Big]\\
&\le \inf_{\gz>1}\Big[\Big(\sup_{-M-1<i\le \uz} {\overline {I\!I}}_i^{\uz-, \gz}\big({\bar g}\big)^{-1}\Big)
\vee \Big(\sup_{\uz+1\le i <N+2}{\overline {I\!I}}_i^{\uz+, \gz}\big({\bar g}\big)^{-1}\Big)\Big]\\
&= \inf_{\gz>1}\, \sup_{-M-1<i\le \uz}{\overline {I\!I}}_i^{\uz-, \gz}\big({\bar g}\big)^{-1}\\
&=\sup_{\gz>1}\bigg\{\frac{1}{\lz_0}-\bigg(1-\frac{1}{\gz}\bigg)\sup_{-M-1<i\le \uz}
\frac{\mu_{\uz} g_{\uz}}{g_i}\sum_{j=-M}^i \frac{1}{\mu_j a_j}\bigg\}^{-1}\\
&=\lz_0.
\endalign
$$
However, there is a problem in the second line of the proof. To apply \thm\;2.4\,(3), one
requires that either $g\in L^2(\mu)$ or $g$ is local. Hence, an additional work is required. Anyhow, the
conclusion holds whenever both $M$ and $N$ are finite. We will come back to the proof in proof (c) below.

(b) Case (i). By \lmm\;7.12, this case can be reduced to Case (ii). Actually, the proof becomes easier now.
With $\gz$ given by (7.18), we have
$$\align
{\overline {I\!I}}_i^{\uz+, \gz}\big({\bar g}\big)&\equiv \frac{1}{\lz_0},\qqd \uz+1\le i <N+2,\\
{\overline {I\!I}}_i^{\uz-, \gz}\big({\bar g}\big)&\equiv \frac{1}{\lz_0},\qqd -M-1< i \le \uz.
\endalign
$$
Hence the second equality in (7.13) holds. Moreover, the first equality in (7.13) also holds whenever
both $M$ and $N$ are finite.

(c) To complete the proof for the first equality in (7.13), we need to overcome the unbounded problem.
For this, choose $M_p, N_p\uparrow \infty$ as $p\to \infty$. Denote by $\lz_0^{\uz-, \gz,\, p}$, $\lz_0^{\uz+, \gz,\, p}$
and $\lz_0^{(p)}$, respectively, the quantities $\lz_0^{\uz-, \gz}$, $\lz_0^{\uz+, \gz}$, and $\lz_0$  when
$M$ and $N$ are replaced by $M_p$ and $N_p$.
Note that for a finite state space, we certainly have $\lz_0^{(p)}>0$, its eigenfunction is positive
(by \prp\;2.2\,(1)) and has
properties (i) and (ii) mentioned in the above proof (by \prp\;7.14).
Clearly, for each fixed $\uz$ and $\gz$, we have
$$\lz_0^{\uz\pm, \gz,\, p}\downarrow \lz_0^{\uz\pm, \gz},\qd
\lz_0^{(p)}\downarrow \lz_0 \qd \text{as }p\to\infty.$$
Thus, as proved in (a) and (b), whether we are in Case (i) or (ii), we have for each $p$,
$$ \align
\lz_0^{(p)}&=\inf_{\uz\in [-M_p, N_p]}\, \inf_{\gz>1}\big[\lz_0^{\uz-, \gz,\, p}\vee \lz_0^{\uz+, \gz,\, p}\big]\\
&\ge \inf_{\uz\in [-M_p, N_p]}\,\inf_{\gz>1}\big[\lz_0^{\uz-, \gz}\vee \lz_0^{\uz+, \gz}\big]\\
&\ge \inf_{\uz\in E}\,\inf_{\gz>1}\big[\lz_0^{\uz-, \gz}\vee \lz_0^{\uz+, \gz}\big].\endalign$$
Therefore, by the first inequality in (7.13)  proved in Part I, it follows that
$$\lz_0 \le \inf_{\uz\in E}\,\inf_{\gz>1}\big[\lz_0^{\uz-, \gz}\vee \lz_0^{\uz+, \gz}\big]
\le \lz_0^{(p)}\downarrow \lz_0 \qqd\text{as $p\to\infty$.}$$
We have thus completed the proof of the theorem.\qed\enddemo

Here are remarks about the assumption made in part (2) of \thm\;7.10.
Similar to the upper estimate,
we do have
$$\lz_0^{(p)}= \sup_{\uz\in [-M_p, N_p]}\, \sup_{\gz>1}\Big[\lz_0^{\uz-, \gz,\, p}\wedge \lz_0^{\uz+, \gz,\, p}\Big].$$
The problem is that
$\lz_0^{\uz\pm, \gz,\, p}\downarrow \lz_0^{\uz\pm, \gz}$ as $p\to\infty$ goes to the opposite
direction and the approximating sequences $\{M_p\}$ and $\{N_p\}$ depend on $\uz$ and $\gz$.
Hence, the proof for the upper estimate does not work for the lower one.
Next, to prove the second equality in (7.13), it seems more natural to assume that
$\lz_0^{\uz-, \gz}\wedge \lz_0^{\uz+, \gz}>0$ for some $\uz$ and $\gz$, that is,
$\lz_0^{(-)}\wedge \lz_0^{(+)}>0$, rather than $\lz_0^{(-)}\vee \lz_0^{(+)}>0$ as we made.
However, if one of them is zero, say $\lz_0^{(-)}=0$, then as mentioned
before \thm\;7.10,
we have a single Dirichlet boundary but not the bilateral Dirichlet ones, and the variational
formula takes a different form (i.e., the second inequality in (7.13) at the boundaries).
Condition (7.14) is due to the same reason. In particular, when $M=-1$, for instance, if
$\sum_i \mu_i<\infty$ and $\sum_{i}(\mu_i a_i)^{-1}=\infty$, then $\lz_0^{(+)}=0$ by \thm\;3.1, and we go
back to the case studied in Section 4. In which case, the eigenfunction of $\lz_0$ is strictly increasing.

To conclude this section, we introduce a complement result to [12;
\prp\;5.13] about the principal eigenvalue for general Markov chains.

\proclaim{\prp\;7.17} Let $(q_{ij}: i, j\in E)$ be symmetric with
respect to $(\mu_i)$ on a countable set $E$, not necessarily
conservative (or having killings):
$$d_i:=q_i-\sum_{j\ne i}q_{ij}\ge 0,\qqd q_i:=-q_{ii}\in [0, \infty).$$
Define
$$D(f)=\frac 1 2 \sum_{i,j\in E}\mu_i q_{ij}(f_j-f_i)^2+\sum_{i\in E} \mu_i d_i f_i^2$$
and
$$\lz_0=\inf\big\{D(f): f\text{ has a finite support and }\mu\big(f^2\big)=1\big\}.$$
Then we have $\inf_{i\in E}\, q_{i}\ge \lz_0.$
\endproclaim

\demo{\prf} Without loss of generality, assume that $E={\Bbb Z}_+=\{0, 1,
\ldots\}$. Fix $k\in E$ and take $f={\dbl}_{\{k\}}$. Then
$\mu\big(f^2\big)=\mu_k$ and
$$\align
D(f)&=\sum_{i,j:\; i< j}\mu_i q_{ij}(f_j-f_i)^2+\sum_{i\in E} \mu_i d_i f_i^2\\
&=\sum_{j>k}\mu_k q_{kj}(f_k-f_j)^2+\sum_{i< k}\mu_i q_{ik}(f_i-f_k)^2+\sum_{i\in E} \mu_i d_i f_i^2\\
&=\sum_{j>k}\mu_k q_{kj}+\sum_{i< k}\mu_i q_{ik}+ \mu_k
d_k.\endalign$$ By the symmetry of $\mu_i q_{ij}$, we get
$$D(f)=\sum_{j>k}\mu_k q_{kj}+\sum_{i< k}\mu_k q_{ki}+ \mu_k d_k
=\mu_k \bigg(\sum_{j\ne k}q_{kj} + d_k\bigg) =\mu_k q_k.$$ It
follows that
$$\lz_0\le D(f)\big/\mu\big(f^2\big)=q_k.$$
The assertion now follows since $k\in E$ is arbitrary. \qed\enddemo

\head 8. Criteria for Poincar\'e-type inequalities \endhead

As in [9] for the ergodic case having $N<\infty$ or (1.2),
the results studied in Sections 2, 3 and 7
can be extended to a more general setup, so called Poincar\'e-type
inequalities. In this way, one obtains
various types of stability, not only the $L^2$-exponential one studied
in the other sections of the paper. Here we consider only the criteria
and the basic estimates for the
inequa\-lities. In other words, we extend \thm s 3.1, 4.2, and
6.2 to the general setup with some improvement. At the same
time, we introduce a criterion for the processes studied in Section
7 in this setup. To do so, we need a class of normed linear spaces
$(\Bbb B, \|\cdot\|_{\Bbb B},\, \mu)$ consisting of real Borel
measurable functions on a measurable space $(X, {\scr X},\, \mu)$.
We now modify the hypotheses on the normed linear spaces given in
[12; Chapter 7] as follows. \roster
\item "{(H1)}" In the case that $\mu(X)=\infty$, $\dbl_K\in {\Bbb B}$ for all
   compact $K$. Otherwise,  $1\in {\Bbb B}$.
\item "{(H2)}" If $h\in \Bbb B$ and  $|f|\le h$, then  $f\in \Bbb B$.
\item "{(H3)}" $\|f\|_{\Bbb B}=\sup_{g\in {\scr G}} \int_X |f| g
\d\mu$,
\endroster
where ${\scr G}$, to be specified case by case, is a class of
nonnegative $\scr X $-measurable functions. A typical example
is ${\scr G}=\{\dbl\}$ and then ${\Bbb
B}=L^1(\mu)$. Throughout this section, we assume (H1)--(H3) for
$(\Bbb B, \|\cdot\|_{\Bbb B},\, \mu)$ without mentioning again.

Before moving further, let us mention the following result.

\proclaim{\rmk\;8.1} Without using (1.2), the results in [9] remain true under
the condition $\sum_i \mu_i<\infty$ replacing the original process
with the maximal one if necessary.
\endproclaim

The key reason is that without condition (1.2), the same conclusion holds in Section 4 on
which the cited paper is based on.

In this section, our
state space is $E=\{i: -M-1<i<N+1\}$ ($M, N\le\infty$) as in the
second part of Section 7. The next result is the main one in this
section; it has several corollaries as we have seen in the
last section. Note that the factor $4$ in (8.2) below
is universal, independent of ${\Bbb B}$.

\proclaim{\thm\;8.2} Consider the minimal birth--death
process with Dirichlet boun\-daries at $-M-1$ if $M<\infty$ and at
$N+1$ if $N<\infty$. Assume that
${\scr G}$ contains a locally positive element.
Then the optimal constant ${A_{\Bbb B}}$ in
the {\it Poincar\'e-type inequality}
$$\big\|f^2\big\|_{\Bbb B} \le A_{\Bbb B} D(f),\qqd f\in {\scr D}^{\min}(D), \tag 8.1$$
satisfies
$$B_{\Bbb B}\le A_{\Bbb B}\le 4 B_{\Bbb B}, \qqd \tag 8.2$$
where the isoperimetric constant $B_{{\Bbb B}}$ can be expressed as
follows:
$$\text{\hskip-2em}B_{{\Bbb B}}^{-1}=\inf_{m, n\in E:\; m \le n}\bigg[\bigg(\sum_{i=-M}^m\frac{1}{\mu_i a_i}\bigg)^{-1}
\!\!+ \bigg(\sum_{i=n}^{N}\frac{1}{\mu_i
b_i}\bigg)^{-1}\bigg]\,\|\dbl_{[m, n]}\|_{{\Bbb B}}^{-1}.\tag 8.3$$
In particular, when ${\Bbb B}=L^{p/2}(\mu)\,(p\ge 2)$ (then (8.1) is called the {\it Sobolev-type inequality}), we have
$$B_{p}^{-1}=\inf_{m, n\in E:\; m \le n}\bigg[\bigg(\sum_{i=-M}^m\frac{1}{\mu_i a_i}\bigg)^{-1}
+ \bigg(\sum_{i=n}^{N}\frac{1}{\mu_i
b_i}\bigg)^{-1}\bigg]\bigg(\sum_{j=m}^n \mu_j\bigg)^{-2/p}.\tag 8.4$$
\endproclaim

\demo{\prf} (a) First consider the transient case, in particular
when one of $M$ or $N$ is finite.
We use the proof of [11; \crl\;4.1] or [12;
\crl\;7.5] with a slight modification. In proof (a) there, it was
shown that one can replace ``$f|_K\ge 1$'' by ``$f|_K= 1$'' in
computing the capacity Cap$(K)$ for compact $K$. Without loss of
generality, assume that $f\ge 0$. Otherwise, replace $f$ with $|f|$.
In the proof just mentioned, the condition ``$\sum_i \mu_i<\infty$''
was used so that $\dbl\in {\scr D}(D)$. We cannot use
this assumption now, but for a given nonnegative $f\in {\scr
D}^{\min}(D)\cap {\scr C}_c(E)$, where ${\scr C}_c(E)$ is the set of
continuous functions with compact support, we can simply choose a
nonnegative smooth $h\in {\scr C}_c(E)$ such that $h|_{\supp
(f)}=1$. Then $h\in {\scr D}^{\min}(D)$, $f\wedge h\in {\scr
D}^{\min}(D)$, and so one can use $f\wedge h\in {\scr D}^{\min}(D)$
instead of $f\wedge \dbl$ to arrive at the same conclusion $D(f)\ge
D(f\wedge h)$ as in the original proof (a).

The first step in the original proof (b) shows that one can replace
a finite number of disjointed finite intervals $\{K_i\}$ by the
connected one $[\min \cup_i K_i,\, \max \cup_i K_i]$. This part of
the proof needs no change.

Note that in the original proof, the state space is $\{1, 2,
\ldots\}$ with Dirichlet boundary at $0$. The main body of the
original proof (b) is to find a minimizer (actually unique) $f\in
{\scr C}_c(E)$ for $D(f)$ having the properties $f_0=0$ and
$f|_K=1$. Replacing $N$ with $q$ for the consistence with the
notation used here and let $K=\{m, m+1, \ldots, n\}$ ($1\le m\le n$,
here $m$ and $n$ are exchanged from the original proof). Now, within
the class of $f$: $f_0=0$, $f|_K=1$ and supp$(f)=\{1, \ldots, q\}$
($n\le q<N+1$), the minimal solution is
$$D(f)=\bigg(\sum_{i=1}^m\frac{1}{\mu_i a_i}\bigg)^{-1}
+ \bigg(\sum_{i=n}^{q}\frac{1}{\mu_i b_i}\bigg)^{-1}.  \tag 8.5$$ To
handle with the general state space, one needs to move the original
left-end point $1$ of the state space to somewhere, say $p>-M-1$. In
detail, replace the condition $m\ge 1$ used in defining the compact
set $K$ by $m> -M-1$. At the same time, replace $\{1, \ldots, q\}$
by $\{p, p+1, \ldots, q\}$ with $ -M-1<p\le m$ for the supp$(f)$.
Then the last formula reads as follows:
$$D(f)=\bigg(\sum_{i=p}^m\frac{1}{\mu_i a_i}\bigg)^{-1}
+ \bigg(\sum_{i=n}^{q}\frac{1}{\mu_i b_i}\bigg)^{-1}.$$ In the
original proof, the ergodic condition and (1.2) are mainly used here to remove
the second term on the right-hand side. We now keep it. Since the
right-hand side is increasing in $p$ and decreasing in $q$, by
making the infimum with respect to $f$, it follows that
$$\align
\text{Cap}(K)
&:=\inf\big\{D(f): f\in {\scr D}^{\min}\cap {\scr C}_c(E)\;\text{and}\;f|_K\ge 1\big\}\\
&\,=\bigg(\sum_{i=-M}^m\frac{1}{\mu_i a_i}\bigg)^{-1}\!\! +
\bigg(\sum_{i=n}^{N}\frac{1}{\mu_i b_i}\bigg)^{-1},\qd K\!=\!\{m, m\!+\!1, \ldots, n\}=:[m, n].
\endalign$$
The assertion of the theorem now follows by using
$$B_{\Bbb B}:=\sup_K \frac{\|\dbl_{K}\|_{\Bbb B}}{\text{Cap}(K)}
=\sup_{-M-1< m\le n< N+1}\frac{\|\dbl_{[m, n]}\|_{\Bbb
B}}{\text{Cap}([m, n])}$$ and applying [11; \thm\;1.1] or
[12; \thm\;7.2]. The last result is an extension of Fukushima and
Uemura (2003, \thm\;3.1).

(b) Next, consider the recurrent case:
both $\sum_{i<\uz} (\mu_i a_i)^{-1}$ and $\sum_{i>\uz} (\mu_i b_i)^{-1}$
are diverged. Here is actually a direct proof of the lower estimate in (8.2).
Without loss of generality, assume that the reference point $\uz=0$.
Fix $m'\ge m\ge 0$ and $n'\ge n\ge 0$. Based on the knowledge about the eigenfunction
given in the last section, and similar to proof (b) of \thm\;3.1, define
$$f_i=
\cases\displaystyle
\sum_{k=i \vee n}^{n'} \dfrac{1}{\mu_k b_k} \dbl_{\{i\le n'\}},\qqd i\ge 0,\\
\displaystyle\gz \sum_{k=-m'}^{i\wedge (-m)} \dfrac{1}{\mu_k a_k} \dbl_{\{i\ge -m'\}},\qqd i\le 0,
\endcases$$
where
$$\gz:=\gz (m', m, n, n')= \sum_{k=n}^{n'} \frac{1}{\mu_k b_k} \bigg/
\sum_{k=-m'}^{-m} \frac{1}{\mu_k a_k}.$$
Here, $\gz$ is chosen to make $f$ be a constant on $[-m, n]$.
By (7.10), we have
$$\align
D(f)&=\sum_{i=-m'}^{-m} \mu_i a_i (f_i-f_{i-1})^2+ \sum_{i=n}^{n'} \mu_i b_i (f_{i+1}-f_{i})^2\\
&= \gz^2 \sum_{i=-m'}^{-m} \frac{1}{\mu_i a_i}+ \sum_{i=n}^{n'} \frac{1}{\mu_i b_i}\\
&=\bigg(\sum_{i=n}^{n'} \frac{1}{\mu_i b_i}\bigg)\bigg[1+ \bigg(\sum_{k=n}^{n'} \frac{1}{\mu_k b_k}\bigg)
\bigg(\sum_{k=-m'}^{-m} \frac{1}{\mu_k a_k}\bigg)^{-1}\bigg].
\endalign$$
Moreover,
$$\big\|f^2\big\|_{\Bbb B}\ge \big\|f|_{[-m, n]}^2\big\|_{\Bbb B}=\bigg(\sum_{i=n}^{n'} \frac{1}{\mu_i b_i}\bigg)^2
\|\dbl_{[-m, n]}\|_{\Bbb B}.$$
Hence,
$$A_{\Bbb B}\ge \frac{\|f^2\|_{\Bbb B}}{D(f)}\ge \|\dbl_{[-m, n]}\|_{\Bbb B}\,
\bigg[\bigg(\sum_{k=-m'}^{-m} \frac{1}{\mu_k a_k}\bigg)^{-1}+\bigg(\sum_{j=n}^{n'} \frac{1}{\mu_i b_i}\bigg)^{-1}\bigg]^{-1}.$$
From this, we obtain the lower estimate in (8.2).
Since ${\scr G}$ contains a locally positive element, we have $\|\dbl_{[-m, n]}\|_{\Bbb B}>0$
for large enough $m$ and $n$. Letting $m', n'\to\infty$, by the recurrent assumption,
it follows that $A_{\Bbb B}=\infty$. Besides,
it is obvious that $B_{\Bbb B}=\infty$ in this case and so the first and then the second assertion of the theorem
becomes trivial in the recurrent case. \qed\enddemo

Proof (b) above indicates an easy improvement of the lower bound
of $A_{\Bbb B}$. Use the same $f$ as above, and define
$$\align
h_i^{(m,m',n,n')}&=\bigg[1-\sum_{k=i+1}^{-m}\frac{1}{\mu_k a_k}
\bigg/\sum_{k=-m'}^{-m}\frac{1}{\mu_k a_k}\bigg]^2\dbl_{[-m', -m-1]}(i)\\
&\qd + \bigg[1-\sum_{k=n}^{i-1}\frac{1}{\mu_k b_k}
\bigg/\sum_{k=n}^{n'}\frac{1}{\mu_k b_k}\bigg]^2\dbl_{[n+1,
n']}(i).\endalign$$ Then a simple computation shows that
$$f^2=\bigg(\sum_{i=n}^{n'}\frac{1}{\mu_i b_i}\bigg)^2 \big(\dbl_{[-m, n]}+ h^{(m,m',n,n')}\big).$$
Hence
$$\frac{\|f^2\|_{\Bbb B}}{D(f)}\ge \big\|\dbl_{[-m, n]}+ h^{(m,m',n,n')}\big\|_{\Bbb B}\,
\bigg[\bigg(\sum_{k=-m'}^{-m} \frac{1}{\mu_k
a_k}\bigg)^{-1}+\bigg(\sum_{j=n}^{n'} \frac{1}{\mu_i
b_i}\bigg)^{-1}\bigg]^{-1}.$$
Noting that the right-hand side is
increasing in $m'$ and $n'$, and making a change of the variable
$-m\to m$, we obtain
$$A_{\Bbb B}\ge \sup_{m, n\in E:\, m\le n}\big\|\dbl_{[m, n]}+ h^{(-m,M,n,N)}\big\|_{\Bbb B}\,
\bigg[\bigg(\sum_{k=-M}^{m} \frac{1}{\mu_k
a_k}\bigg)^{-1}+\bigg(\sum_{j=n}^{N} \frac{1}{\mu_i
b_i}\bigg)^{-1}\bigg]^{-1}.$$
Denote by $C_{\Bbb B}$ the right-hand side. Then the conclusion of \thm\;8.2 can be restated as
$B_{\Bbb B}\le C_{\Bbb B}\le A_{\Bbb B}\le 4 B_{\Bbb B}$.
Certainly, this remark is meaningful
in other cases but we will not mention again.

The next result is an easier consequence of \thm\;8.2.

\proclaim{\crl\;8.3} Everything in the premise is the same as in
\thm\;8.2. Then
\roster
\item we have $B_{\Bbb B}\le B_L\wedge B_R$, where
$$ B_{L}=\sup_{n\in {E}}\, \sum_{i=-M}^n \frac{1}{\mu_i a_i} \|\dbl_{[n, N+1)}\|_{\Bbb B},\qqd
B_{R}=\sup_{n\in {E}}\, \sum_{i=n}^{N} \frac{1}{\mu_i b_i} \|\dbl_{(-M-1, n]}\|_{\Bbb B}.$$
The equality sign holds once
$$ S:=\sum_{i=-M}^{N} \frac{1}{\mu_i a_i}+\frac{1}{\mu_N b_N} \dbl_{\{N<\infty\}}=\infty.$$
\item Next, we have $B_{\Bbb B}\ge (B_L\wedge B_R)\dbl_{\{S=\infty\}}+ S^{-1} B$, where
$$B=\sup_{m,n\in E:\, m\le n}\bigg[\bigg(\sum_{i=-M}^m\frac{1}{\mu_i a_i}\bigg)
\bigg(\sum_{k=n}^{N}\frac{1}{\mu_k b_k}\bigg)\|\dbl_{[m, n]}\|_{{\Bbb B}}\bigg].$$
\endroster
\endproclaim

\demo{\prf}  Clearly, by (8.3), we have
$$B_{\Bbb B}^{-1}\ge\! \inf_{ m\le n} \bigg(\sum_{i=-M}^m\frac{1}{\mu_i a_i}
\|\dbl_{[m, n]}\|_{{\Bbb B}}\bigg)^{\!-1}\!\! =\!\inf_{m\in E}
\bigg(\sum_{i=-M}^m\frac{1}{\mu_i a_i} \|\dbl_{[m, N+1)}\|_{{\Bbb
B}}\bigg)^{\!-1}\!\!,$$
and so $B_{\Bbb B}\le B_L.$ The equality sign holds once
$\sum_{i=\uz}^N (\mu_i b_i)^{-1}=\infty$.
Similarly, we have $B_{\Bbb B}\le B_R$.
The equality sign holds once
$\sum_{i=-M}^{\uz} (\mu_i a_i)^{-1}=\infty$.
Hence, $B_{\Bbb B}\le B_L\wedge B_R$ and the equality sign holds once $S=\infty$.

Next, when $S<\infty$, we have
$$B_{\Bbb B}^{-1} \le S\inf_{m\le n}\bigg[\bigg(\sum_{i=-M}^m\frac{1}{\mu_i a_i}\bigg)
\bigg(\sum_{k=n}^{N}\frac{1}{\mu_k b_k}\bigg)\|\dbl_{[m, n]}\|_{{\Bbb B}}\bigg]^{-1}
=S B^{-1}.$$
We have thus proved the corollary.\qed\enddemo

Of course, one can decompose the constant $B$ in \crl\;8.3\,(2). For instance, for fixed
$m_0$, we have
$$B\ge \bigg(\sum_{i=-M}^{m_0}\frac{1}{\mu_i a_i}\bigg)\sup_{m_0\le n<N+1}
\bigg[\bigg(\sum_{k=n}^{N}\frac{1}{\mu_k b_k}\bigg)\|\dbl_{[m_0, n]}\|_{{\Bbb B}}\bigg].$$
The last factor is close to $B_R$ when $m_0$ is negative enough. However, when $m_0\to -M$,
the first term tends to zero since $S<\infty$, unless $M<\infty$. This indicates that bounding $B_{\Bbb B}$ by
$B_L$ and $B_R$ is rather rough, especially in the case that $E={\Bbb Z}$ (cf. \xmp\;8.9 below).
This is a particularly
different point of the processes on the whole ${\Bbb Z}$ or on the half space ${\Bbb Z}_+$
as shown by \crl\;8.4 below.

We now specify \thm\;8.2 and \crl\;8.3 to the half space: either $M$ or $N$ is
finite. This corresponds to the processes studied in the first part
of Section 7 (see \crl\;7.3).

\proclaim{\crl\;8.4} In \thm\;8.2, let $M=-1$. Then we have
$$B_{\Bbb B}^{-1}=\inf_{1\le n\le m< N+1}\bigg[\bigg(\sum_{i=1}^n\frac{1}{\mu_i a_i}\bigg)^{-1}
+ \bigg(\sum_{i=m}^{N}\frac{1}{\mu_i
b_i}\bigg)^{-1}\bigg]\,\|\dbl_{[n, m]}\|_{{\Bbb B}}^{-1}. \tag 8.6$$
Furthermore, we have
$$B_L\wedge B_R\ge B_{\Bbb B}\ge
\big(\dbl_{\{S=\infty\}}+(a_1 S)^{-1}\big)\,(B_L\wedge B_R),$$
where
$$\align
B_L&=\sup_{1\le n< N+1} \sum_{i=1}^n\frac{1}{\mu_i a_i} \|\dbl_{[n,
N+1)}\|_{{\Bbb B}},\qqd B_R=\sup_{1\le m< N+1}
\sum_{k=m}^{N}\frac{1}{\mu_k b_k} \|\dbl_{[1, m]}\|_{{\Bbb B}},\\
S&=\sum_{i=1}^{N}\frac{1}{\mu_i a_i}+\frac{1}{\mu_N b_N}\dbl_{\{N<\infty\}}.\tag 8.7\endalign$$
\endproclaim

\demo{\prf} The first assertion follows from
\thm\;8.2 with $M=-1$ and an exchange of $m$ and $n$ again.
The second one follows from \crl\;8.3 except the last estimate. When $S=\infty$,
we have $B_{\Bbb B}=B_L$. While when $S<\infty$, we have
$$\align
B_{\Bbb B}^{-1} &\le S\inf_{1\le n\le m< N+1}\bigg[\bigg(\sum_{i=1}^n\frac{1}{\mu_i a_i}\bigg)
\bigg(\sum_{k=m}^{N}\frac{1}{\mu_k b_k}\bigg)\|\dbl_{[n, m]}\|_{{\Bbb B}}\bigg]^{-1}\\
&\le a_1 S\inf_{1\le m< N+1} \bigg(\sum_{k=m}^{N}\frac{1}{\mu_k b_k}\|\dbl_{[1, m]}\|_{{\Bbb B}}\bigg)^{-1}\\
&= a_1 S\, B_R^{-1}.
\endalign$$
Therefore,
$$B_{\Bbb B}\ge B_L \dbl_{\{S=\infty\}}+(a_1 S)^{-1} B_R
\ge \big(\dbl_{\{S=\infty\}}+(a_1 S)^{-1}\big)\,(B_L\wedge B_R)$$
as required.\qed\enddemo

When one of $M$ or $N$ is finite and its Dirichlet boundary is
replaced by the Neumann one, the solution becomes simpler. The
next result corresponds to the processes studied in Sections 2 and
3.

\proclaim{\thm\;8.5} Let $M=0$ be the Neumann boundary and assume that
${\scr G}$ contains a locally positive element. Then the
isoperimetric constant $B_{\Bbb B}\!:\!=\!\sup_K {\|\dbl_{K}\|_{\Bbb
B}}/{\text{\rm Cap}(K)}$ can be expressed as
$$B_{{\Bbb B}}=\sup_{0 \le n< N+1} \sum_{i=n}^{N}\frac{1}{\mu_i
b_i}\,\|\dbl_{[0, n]}\|_{{\Bbb B}}.\tag 8.8$$
In particular, for the Sobolev-type inequality, we have
$$B_p=\sup_{0 \le n< N+1} \sum_{i=n}^{N}\frac{1}{\mu_i
b_i} \bigg(\sum_{j=0}^n \mu_j\bigg)^{2/p},\qqd p\ge 2.  \tag 8.9$$
\endproclaim

\demo{\prf} The proof is nearly the same as that of \thm\;8.2 except
one point. In proof (b) of [11; \crl\;4.1] or [12;
\crl\;7.5], to find a minimizer $f$ for $D(f)$, since the constraint
$f_0=0$ and $f_n=1$, $f$ cannot be a constant on $\{0, 1, \ldots,
n\}$. Now, without the constraint $f_0=0$, the minimizer should satisfy
$f_j=1$ for all $j: 0\le j\le n$. Thus, instead of (8.5), the minimal
solution becomes
$$D(f)=\bigg(\sum_{i=n}^{q}\frac{1}{\mu_i b_i}\bigg)^{-1}.  $$
Then the necessary change of the proof of \thm\;8.2 after (8.5)
should be clear.
\qed\enddemo

Applying \thm\;8.5 to ${\Bbb B}=L^1(\mu)$, we return to \thm\;3.1.
Actually, in parallel to [9], one may extend the results in Sections 2 and 3,
\thm\;3.1 in particular, to the present setup of normed linear spaces and then deduce \thm\;8.5.
The next result is obvious, it says that for a null-recurrent process,
the $L^p\,(p\ge 1)$-Sobolev inequality is still not weak enough.

\proclaim{\crl\;8.6} Consider a birth--death process on ${\Bbb Z}_+$.
If $\sum_{i\ge 1} \mu_i=\infty$ and \newline $\sum_{i\ge 1} (\mu_i b_i)^{-1}=\infty,$
then $B_p^{(8.4)} = B_p^{(8.9)}=\infty$ for all $p\ge 2$.
\endproclaim

\proclaim{\rmk\;8.7} {\rm We now compare (8.6) and (8.8) in the particular case that
$\sum_i \mu_i$ $=\infty$. Then the constant $B_{\Bbb B}$ given in
(8.6) becomes
$$B_{\Bbb B}^{(8.6)}=\sup_{1\le m<N+1} \sum_{i=m}^{N}\frac{1}{\mu_i
b_i}\,\|\dbl_{[1, m]}\|_{{\Bbb B}}.$$ Rewrite the constant $B_{\Bbb
B}$ given in (8.8) as
$$B_{{\Bbb B}}^{(8.8)}=\bigg(\sum_{i=0}^{N}\frac{1}{\mu_i
b_i}\,\|\dbl_{\{0\}]}\|_{{\Bbb B}}\bigg)\bigvee\bigg(\sup_{1 \le n<
N+1} \sum_{i=n}^{N}\frac{1}{\mu_i b_i}\,\|\dbl_{[0, n]}\|_{{\Bbb
B}}\bigg).$$
By (H3), we have
$$\|\dbl_{[1, n]}\|_{{\Bbb B}}\le \|\dbl_{[0, n]}\|_{{\Bbb B}}
\le \|\dbl_{\{0\}]}\|_{{\Bbb B}} +\|\dbl_{[1, n]}\|_{{\Bbb B}}.$$
Next, by (H1), we have $\|\dbl_{\{0\}}\|_{{\Bbb B}}<\infty$. It follows that
$B_{\Bbb B}^{(8.6)}<\infty$ iff $B_{\Bbb B}^{(8.8)}<\infty$.
}
\endproclaim

We conclude this section by a simple example to show the role of the Poincar\'e-type inequalities.

\proclaim{\xmp\;8.8} Consider a birth--death process on ${\Bbb Z}_+$ with $\mu_i=(i+1)^{\gz}\,(\gz>1)$
and $b_i\equiv 1$. Then $a_i=i^{\gz}(i+1)^{-\gz}$ and
$$B_p^{(8.9)}=\sup_{n\ge 0} \bigg[\sum_{i=0}^n (i+1)^{\gz}\bigg]^{2/p}\sum_{j=n}^\infty \frac{1}{(j+1)^{\gz}},\qqd p\ge 2.$$
Hence, $B_p^{(8.9)}<\infty$ iff
$$p\ge 2\bigg(1+\frac{2}{\gz -1}\bigg).$$
However, $\dz^{(3.1)}=B_2^{(8.9)}=\infty$ for all $\gz>1$.
\endproclaim

\proclaim{\xmp\;8.9} Let $E={\Bbb Z}$, $b_i\equiv 1$, $\mu_i=e^{i^2}$, and ${\Bbb B}=L^1(\mu)$.
Then for the quantities given in \crl\;8.3, we have $B_L=B_R=\infty$ but $B_{\Bbb B}<\infty$.
\endproclaim

\demo{\prf} Obviously, $B_L=B_R=\infty$. To show that $B_{\Bbb B}<\infty$, since
$$x\vee y \le x+y \le 2 (x\vee y),$$
it suffices to prove that
$$\sup_{m\le n}\bigg[\bigg(\sum_{i=-\infty}^m\frac{1}{\mu_i a_i}\bigg)\bigwedge
\bigg(\sum_{k=n}^{\infty}\frac{1}{\mu_k b_k}\bigg)\bigg]\sum_{j=m}^n \mu_j<\infty.$$
By symmetry, without loss of generality, it is enough to show that
$$\sup_{m\ge n\ge 0}\bigg(\sum_{i=-\infty}^{-m} \frac{1}{\mu_i a_i}\bigg)\sum_{j=-m}^n \mu_j<\infty,$$
or
$$\sup_{m\ge 0}\bigg(\sum_{i=-\infty}^{-m} \frac{1}{\mu_i a_i}\bigg)\sum_{j=-m}^m \mu_j
<\infty.$$
Equivalently,
$$\sup_{m\ge 0}\bigg(\sum_{i=m}^{\infty} \frac{1}{\mu_i b_i}\bigg)\sum_{j=0}^m \mu_j
<\infty. \tag 8.10$$
The assertion now follows by using Conte's inequality:
$$x\bigg(1+\frac{x}{24}+\frac{x^2}{12}\bigg) e^{- 3 x^2/4}< e^{-
x^2}\int_0^x e^{y^2} \le \frac{\pi^2}{8x}\big(1-e^{-
x^2}\big), \qqd x\ge 0$$
and Gautschi's estimate:
$$\align
\frac 1 2 \Big[(x^p+2)^{1/p}&-x\Big] <  e^{x^p}\int_x^\infty e^{-y^p}dy\le C_p
\bigg[\bigg(x^p+\frac{1}{C_p}\bigg)^{1/p} - x\bigg],\qqd x\ge 0,\\
&C_p=\ggz\big(1+{1}/{p}\big)^{p/(p-1)},\qd p>1;\qqd C_2=\pi/4.
\endalign$$
Alternatively, one may check directly that the function under supremum
on the left-hand side of (8.10) is decreasing in $m$\,$(\ge 1)$
and then (8.10) follows easily.\qed
\enddemo

\head{9. General killing} \endhead

In Sections 4 and 7, we have studied the special case having a
killing at $1$ only. We now study the process with general killing,
as described by $(2.1)$ with state space shifted by $1$: $E=\{i:
1\le i <N+1\}$. We use the same symmetric measure $(\mu_i)$ as in
Section 4.

The next preliminary result is quite useful. To which it is more
convenient to use $a_1+c_1$ and $b_N+c_N$ for the killing rates at
boundaries $1$ and $N$ (if $N<\infty$), respectively, rather than
$c_1$ and $c_N$ used in \prp\;2.1. Note that the killing rates in
the next proposition are allowed to be zero identically.

\proclaim{\prp\;9.1} Let $(a_i)$ and $(b_i)$ be positive but $a_1\ge
0$, $b_N\ge 0$ if $N<\infty$, and let $(c_i)$ be nonnegative on $E$.
Define $\lz_0=\lz_0 (a_i, b_i, c_i)$ as follows:
$$\lz_0=\inf\big\{D(f): \mu\big(f^2\big)=1, f\in {\scr K}\big\},$$
where
$$D(f)=\sum_{i\in E}\mu_i b_i (f_{i+1}-f_i)^2+\mu_1 a_1 f_1^2
+\sum_{i\in E}\mu_i c_i f_i^2,\qqd f_{N+1}=0\text{ if }N<\infty.$$
Write $\tilde\lz_0=\lz_0 (a_i, b_i, 0)$ for simplicity. Then we have
\roster
\item $\lz_0 (a_i, b_i, c_i')\ge \lz_0 (a_i, b_i, c_i)$ if $c_i'\ge c_i$ for all $i\in E$.
\item $\lz_0 (a_i, b_i, c_i+c)= \lz_0 (a_i, b_i, c_i)+c$ for constant $c\ge 0$.
\item
$\tilde\lz_0+\sup_{i\in E}c_i\ge \lz_0\ge \tilde\lz_0+\inf_{i\in
E}c_i$ and the equalities hold if $c_i$ is a constant on $E$.
\endroster
\endproclaim

\demo{\prf} Since a change of $\{c_i\}_{i=1}^N$ makes no influence
to $\{\mu_i\}_{i=1}^N$, part (1) is simply a comparison of the
Dirichlet forms on the same space $L^2(\mu)$ with common core ${\scr
K}$. Similarly, one can prove the other assertions. \qed\enddemo

Note that \prp\;9.1 makes a comparison for the killing rates only.
Actually, a more general comparison is available in view of
[3;  \thm\;3.1]. Next, if (1.3) holds, then by \prp\;1.3 and the
remark below (4.3), the Dirichlet is unique, and so the condition
$f\in {\scr K}$ can be ignored in defining $\tilde\lz_0$.

It is worthy to mention that the principal eigenvalue $\lz_0$
studied here can be extended to a more general class of
Schr\"odinger operators. That is, we may replace the nonnegative potential
$(c_i)$ with the one bounded below by a constant: $\inf_i c_i\ge
-M>-\infty$. Then we have $c_i+M\ge 0$ for all $i$ and
$$\lz_0 (a_i, b_i, c_i)=\lz_0 (a_i, b_i, c_i+M)-M\ge -M.$$

Having \prp\;9.1 at hand, all the examples for $\tilde\lz_0$ given
in Sections 3, 5 and 7, can be translated into the case of $\lz_0$
with constant killing rate. For instance, we have the following
example which already shows the complexity of the problem studied in
this section.

\proclaim{\xmp\;9.2} Let $a_i\equiv a>0$ for $i\ge 2$, $b_i \equiv b>0$
and $c_i \equiv c\ge 0$  for $i\ge 1$. \roster
\item If $a_1=a$, or $a_1=0$ but still $a\le b$, then $\lz_0=\big(\sqrt{a}-\sqrt{b}\,\big)^2+c$.
\item If $a_1=0$ and $a>b$, then $\lz_0=c$.
\endroster
\endproclaim

\demo{\prf} By \prp\;9.1, we need only to study ${\tilde\lz}_0. $ In
the last case, since the process is ergodic, we have
${\tilde\lz}_0=0$. Next, we have
${\tilde\lz}_0=\big(\sqrt{a}-\sqrt{b}\,\big)^2$ according different
cases by \roster
\item"(i)" Example 5.3 if $a_1=a$ and $a\ge b$,
\item"(ii)" Example 7.7 if $a_1=a$ and $a\le b$,
\item"(iii)" Example 3.4 if $a_1=0$ and $a\le b$. \qed\endroster\enddemo

From now on, we return to a convention made in Section 2, the rates
$a_1$ and $b_N$ are combined into $c_1$ and $c_N$ if $N<\infty$.
Thus, in \thm\;7.1 for instance, we have $a_1=0$ and $c_1> 0$,
and moreover, $b_N=0$ and $c_N>0$ if $N<\infty$. In
general, we assume that $c_i\not\equiv 0$. Otherwise, we will return
to what we treated in Sections 2 and 3. Define the operator $R$:
$$R_i(v)={a_{i}}\big(1-v_{i-1}^{-1}\big)+{b_{i}} (1 -
 {v}_{i})+c_i,\qqd i\in E,\; v_0=\infty,\;
v_N=0\text{ if }N<\infty,$$ for $v$ in the set ${\scr V}=\{v_i>0:
1\le i<N\}$. Next, define $\widetilde{\scr V}={\scr V}$ if
$N<\infty$. When $N=\infty$, define
$$\align
{\widetilde{\scr V}}=&\bigcup_{m=1}^\infty
\bigg\{v_i: v_i>0\text{ for }i<m,\; v_i=0\text{ for }i\ge m\bigg\}\\
& \bigcup\bigg\{v: v_i>0 \text{ on }E,\; \text{the function }f\!: f_1=1, f_i=\prod_{k=1}^{i-1}v_k\,(i\ge 2)\text{ is in } L^2(\mu)\\
&\qqd\qqd\text{ and satisfies }\ooz f/f \le \ez\text{ on }E\text{
for some constant }\ez\bigg\}. \tag 9.1
\endalign$$
For $v\in {\widetilde{\scr V}}$ with finite support,
$R_{\bullet}(v)$ is also well defined by setting $1/0=\infty$.

\proclaim{\thm\;9.3} Assume that $c_i\not\equiv 0$. For $\lz_0$
defined by $(2.2)$ with state space $E=\{i: 1\le i<N+1\}$, the
following variational formulas hold:
$$\inf_{v\in {\widetilde{\scr V}}}\,
\sup_{i\in E}R_i(v)= \lz_0=\sup_{v\in {\scr V}}\,\inf_{i\in
E}R_i(v).\tag 9.2$$
\endproclaim

\demo{\prf} (a) First, we study the lower estimate. In the case that
$\sum_{k\in E}\mu_k<\infty$, as a particular consequence of
[5; Theorem 1.1], we have
$$\lz_0\ge \sup_{g>0}\,\inf_{i\in E}\frac{-\ooz g}{g}(i)
=\sup_{g>0}\,\inf_{i\in
E}\bigg[{a_{i}}\bigg(1-\frac{{g}_{i-1}}{g_i}\bigg)+{b_{i}}\bigg(1 -
\frac{g_{i+1}}{{g}_{i}}\bigg)+c_i\bigg],\tag 9.3$$ where $g_0:=0$
and $g_{N+1}=0$ if $N<\infty$. The proof remains true when
$\sum_{k\in E}\mu_k=\infty$, simply using $E_m=\{1, 2, \ldots,
m\}$\, ($m<N+1$) instead of the original one. Actually, the
conclusion holds in a very general setup (cf. Shiozawa and Takeda
(2005) and its extension to the unbounded test functions by Zhang
(2007)).

Suppose that $\lz_0>0$ for a moment. Then by \prp\;2.1 (with a shift
by 1 of the state space), the eigenfunction $g$ of $\lz_0$ is
positive. It follows that the first equality sign in (9.3) can be
attained and so does the last equality in (9.2) with
$v_i=g_{i+1}/g_i>0\,(1\le i<N$). Next, if $\lz_0=0$, then $N=\infty$
since $a_i$ and $b_i$ are positive for $i: 2\le i < N$, and
$c_i\not\equiv 0$ (in the case of \thm\;7.1, we have $c_1>0$ and
also $c_N>0$ if $N<\infty$). By setting $v_i\equiv 1$ for $i\in E$,
we get
$$\inf_{i\in E}\bigg[{a_{i}}\bigg(1- \frac{1}{{v}_{i-1}}\bigg)+{b_{i}}(1 -
{v}_{i})+c_i\bigg]\ge \inf_{i\in E} c_i\ge 0. \tag 9.4$$ Hence, the
last term of (9.2) is nonnegative. Therefore, the last equality in
(9.2) is trivial if $\lz_0=0$, in view of (9.3).

(b) Next, we study the upper estimate. We consider only the case
that $N=\infty$. Otherwise, the proof is easier. Given $v\in
{\widetilde{\scr V}}$, let $\gz=\gz(v)=\sup_{1\le i<\infty} R_i(v)$
and as in the definition of ${\widetilde{\scr V}}$, set
$$f_0=0,\; f_1=1,\;
f_i=\prod_{k=1}^{i-1} v_k, \qd i\ge 2. \tag 9.5$$

First, suppose that $\supp(v)=\{1, 2, \ldots, m-1\}$ for a finite
$m$. Then $\supp(f)=\{1, 2, \ldots, m\}$ and
$$\frac{-\ooz f}{f}(i)=R_i(v)\le \gz, \qqd i=1, 2, \ldots, m.$$
Hence,
$$\align
\gz\sum_{k=1}^m \mu_k f_k^2&\ge \sum_{k=1}^m \mu_k f_k (-\ooz f)(k)\\
&=\sum_{k=2}^{m+1} \mu_k a_k f_{k-1}(f_{k-1}- f_k)-\sum_{k=1}^m
\mu_k a_k f_k (f_{k-1}-f_k)
 +\sum_{k=1}^m \mu_k c_k f_k^2\\
&=\sum_{k=1}^m \mu_k \big[a_k (f_{k-1}-f_k)^2+c_k f_k^2\big] + \mu_{m+1}a_{m+1}f_m(f_m-f_{m+1})\\
&= \sum_{k=1}^{m+1} \mu_k \big[a_k (f_{k-1}-f_k)^2+ c_k f_k^2\big].
\endalign$$
We have not only $\gz\ge 0$ (actually $\gz>0$ when $m$ is large
enough since $c_i\not\equiv 0$) but also
$$\lz_0\le \frac{D(f)}{\|f\|^2}\le \gz(v)\tag 9.6$$
for all $v\in {\widetilde{\scr V}}$ with finite support.

(c) Next, we are going to prove (9.6) in the case that $v\in
{\widetilde{\scr V}}$ with $v_i>0$ for all $i\ge 1$. In this case,
the positivity condition of $v$ is not enough for the first equality
in (9.2), as mentioned in Section 2 (above the proofs of \thm\;2.4
and \prp\;2.5). See also the specific situation given in the proof
of Example 9.17 below. This explains why two additional conditions
are included in the second union of the definition of
${\widetilde{\scr V}}$. The condition ``$f\in L^2(\mu)$'' is
essential but not the one ``$\ooz f/f\le \ez$'' since the
eigenfunction $g$ of $\lz_0$ satisfies ``$\ooz g/g=- \lz_0$''. To
prove (9.6), without loss of generality, assume that $\gz=\gz
(v)<\infty$. Otherwise, (9.6) is trivial. Clearly, $\gz\ge
R_1(v)=b_1(1-v_1)+c_1>-\infty$. Note that by assumptions, the
function $f$ possesses  the following properties: \roster
\item"(i)" $f>0$ on $E$.
\item"(ii)" $f\in L^2(\mu)$ and then $P_t f\in L^2(\mu)$, where $P_t=(p_{ij}(t))$ is the minimal semigroup
determined by the Dirichlet form.
\item"(iii)"
$|\ooz f(i)|=\big|\sum_{j}q_{ij}f_j\big|\le \max\{|\ez|, |\gz|\}
f_i$ for all $i\in E. $
\endroster
Here, property (iii) comes from
$$-\ez f\le -\ooz f\le \gz f.$$
Since $(p_{ij}(t))$ satisfies the forward Kolmogorov equation:
$$p_{ij}(t)=\dz_{ij}+\int_0^t \sum_{k} p_{ik}(s)q_{kj}\d s$$
and (i), it follows that
$$P_t f(i)=f_i+ \sum_j \int_0^t\sum_{k} p_{ik}(s)q_{kj} f_j\d s.$$
By (ii), $P_t f(i)<\infty$ and is continuous in $t$. Because of this
and (iii), the order of the last two sums and also the integration
are exchangeable. This leads to
$$P_t f(i)\ge f_i-\gz \int_0^t\sum_{k} p_{ik}(s) f_k\d s=f_i -\gz \int_0^t P_s f(i)\d s,\qqd i\in E, \tag 9.7$$
since by assumption $\ooz f \ge -\gz f$. Therefore, we obtain
$$0< D(f) =\lim_{t\downarrow 0}\frac{1}{t}(f, f-P_t f)
\le \lim_{t\downarrow 0}\frac{\gz}{t}\int_0^t (f, P_s f)\d s=\gz
(v)\|f\|^2<\infty.$$ Here, the first limit is due to (ii) and the
first equality in (1.10), the last inequality comes from (i) and
(9.7). We have thus proved that not only $\gz> 0$ but also $f\in
{\scr D}(D)$ and so we have returned to (9.6). In other words, (9.6)
holds for all $v\in{\widetilde{\scr V}}$. By making infimum with
respect to $v\in{\widetilde{\scr V}}$, we obtain
$$\lz_0\le \inf_{v\in {\widetilde{\scr V}}}\,
\sup_{i\in E}R_i(v).$$

(d) To prove the equality sign in the last formula holds, in view of
proof (b), we have actually proved that for every finite $m$,
$$\align
\lz_0^{(m)}&:=\inf\{D(f): f_0=0,\; f_i=0\text{ for all }i\ge m+1, \; \|f\|=1\}\\
&\le\inf_{v\in {\widetilde{\scr V}}_m} \sup_{1\le i\le m}R_i(v),
\tag 9.8
\endalign$$
where
$${\widetilde{\scr V}}_m=\{v_i: v_i>0\text{ for }i<m,\; v_m=0 \}.$$
Actually, there is a $\bar v\in {\widetilde{\scr V}}_m$ such that
$R_i(\bar v)=\lz_0^{(m)}>0$ for all $i\, (1\le i\le m)$ since
$m<\infty$ and then the equality sign in (9.8) holds. Therefore, the
first equality in (9.2) holds since $\lz_0^{(m)}\downarrow \lz_0$ as
$m\uparrow \infty$. \qed\enddemo

We now begin to study the estimate of $\lz_0$. First, by \prp\;7.17,
we have a simple upper bound:
$$\lz_0\le \inf_{i\in E} (a_i+b_i+c_i).$$
Hence, $\lz_0=0$ whenever $\lim_{n\to\infty}(a_n+b_n+c_n)=0$. The
next result provides us a finer upper bound. It is motivated from
\thm\;3.1.

\proclaim{\prp\;9.4} Let ${\tilde c}_i=c_i-\inf_{i} c_i$. Then
$$\align
\text{\hskip-6em}\lz_0&\le \inf_{i\in E} c_i+\inf_{\ell\in E}
\bigg(\sum_{i=1}^{\ell}\mu_i\bigg)^{-1}\inf_{E\ni m\ge\ell}
\bigg[\bigg(\sum_{k=\ell}^m \frac{1}{\mu_k b_k}\bigg)^{-1}
    +\sum_{i=1}^m\mu_i {\tilde c}_i\bigg]\text{\hskip-3em}\tag 9.9\\
\text{\hskip-6em}&\le  \inf_{i\in E} c_i+ \inf_{\ell\in E}
\bigg(\sum_{i=1}^{\ell}\mu_i\bigg)^{-1}\bigg[\mu_{\ell} b_{\ell}
+\sum_{i=1}^{\ell}\mu_i {\tilde c}_i\bigg]. \tag 9.10\endalign$$
\endproclaim

\demo{\prf} By \prp\;9.1, it is enough to consider the case that
${\tilde c}_i\equiv c_i$, i.e., $\inf_i c_i=0$. Fix $\ell\le m$ and define
$$\fz_i= \fz_i^{(\ell, m)}={\dbl}_{\{i\le m\}}\sum_{k=i\vee \ell}^m \frac{1}{\mu_k b_k},
\qqd i\in E.$$ Then
$$\gather
\mu\big(\fz^2\big)=\sum_{i=1}^\ell \mu_i \fz_{\ell}^2+
\sum_{i=\ell+1}^m \mu_i \fz_i^2\ge \fz_{\ell}^2\sum_{i=1}^\ell \mu_i, \\
D(\fz)=\sum_{k=\ell}^m \frac{1}{\mu_k
b_k}+\fz_{\ell}^2\sum_{i=1}^{\ell}\mu_i c_i +\sum_{i=\ell+1}^m \mu_i
c_i\fz_i^2 \le \fz_{\ell}+\fz_{\ell}^2 \sum_{i=1}^m\mu_i
c_i.\endgather$$ Hence,
$$\frac{D(\fz)}{\mu\big(\fz^2\big)}
\le
\bigg(\sum_{i=1}^{\ell}\mu_i\bigg)^{-1}\bigg[\fz_{\ell}^{-1}+\sum_{i=1}^m\mu_i
c_i\bigg].$$ Because $\fz^{(\ell, m)}\in {\scr K}$, it follows that
$$\align
\lz_0&\le \inf_{\ell\in E}\inf_{E\ni m\ge\ell}\frac{D(\fz)}{\mu\big(\fz^2\big)}\\
&\le \inf_{\ell\in E}
\bigg(\sum_{i=1}^{\ell}\mu_i\bigg)^{-1}\inf_{E\ni m\ge\ell}
\bigg[\bigg(\sum_{k=\ell}^m \frac{1}{\mu_k b_k}\bigg)^{-1}+\sum_{i=1}^m\mu_i c_i\bigg]\\
&\le \inf_{\ell\in E}
\bigg(\sum_{i=1}^{\ell}\mu_i\bigg)^{-1}\inf_{E\ni m\ge\ell}
\bigg[\mu_{\ell} b_{\ell}+\sum_{i=1}^m\mu_i c_i\bigg]\\
&= \inf_{\ell\in E}
\bigg(\sum_{i=1}^{\ell}\mu_i\bigg)^{-1}\bigg[\mu_{\ell}
b_{\ell}+\sum_{i=1}^{\ell}\mu_i c_i\bigg]. \qed\endalign$$
\enddemo

As an immediate consequence of (9.10), we obtain the following
result.

\proclaim{\crl\;9.5} If $\sum_{i=1}^\infty\mu_i=\infty$ and
$$\varlimsup_{m\to\infty} \mz_m b_m\bigg(\sum_{i=1}^m\mu_i\bigg)^{-1}= 0,$$ then
$$\lz_0\le \inf_{i\in E}c_i+\varliminf_{m\to\infty} \sum_{i=1}^{m}\mu_i {\tilde c}_i\bigg/\sum_{i=1}^{m}\mu_i
\le \inf_{i\in E}c_i+\varlimsup_{n\to\infty}{\tilde c}_n.$$
\endproclaim

\demo{\prf} Without loss of generality, assume that ${\tilde
c}_i\equiv c_i$.

By assumptions, it follows that
$$\align
\varliminf_{\ell\to\infty}&\bigg(\sum_{i=1}^{\ell}\mu_i\bigg)^{-1}\bigg[\mu_{\ell} b_{\ell}+\sum_{i=1}^{\ell}\mu_i c_i\bigg]\\
&\le \varlimsup_{\ell\to\infty} \mu_{\ell}
b_{\ell}\bigg(\sum_{i=1}^{\ell}\mu_i\bigg)^{-1} +
\varliminf_{\ell\to\infty}\bigg(\sum_{i=1}^{\ell}\mu_i\bigg)^{-1}\sum_{i=1}^{\ell}\mu_i c_i\\
&=\varliminf_{\ell\to\infty}\bigg(\sum_{i=1}^{\ell}\mu_i\bigg)^{-1}\sum_{i=1}^{\ell}\mu_i
c_i.
\endalign$$
The first inequality now follows from (9.10).

To prove the second inequality, let
$\gz=\varlimsup_{n\to\infty}{c}_n\in [0, \infty]$. Then for every
$\vz>0$, we have $\sup_{k\ge n} c_k\le \gz+\vz$ for large enough
$n$. Hence,
$$\sum_{i=1}^{\ell}\mu_i c_i=\sum_{i=1}^{n}\mu_i c_i+\sum_{i=n+1}^{\ell}\mu_i c_i
\le \sum_{i=1}^{n}\mu_i c_i+(\gz+\vz) \sum_{i=n+1}^{\ell}\mu_i,\qqd
\ell>n.
$$
We have thus obtained
$$\align
\varliminf_{\ell\to\infty}\bigg(\sum_{i=1}^{\ell}\mu_i\bigg)^{-1}\sum_{i=1}^{\ell}\mu_i
c_i &\le
\varlimsup_{\ell\to\infty}\bigg(\sum_{i=1}^{\ell}\mu_i\bigg)^{-1}
\bigg[\sum_{i=1}^{n}\mu_i c_i+(\gz+\vz) \sum_{i=n+1}^{\ell}\mu_i\bigg]\\
&= \gz+\vz\endalign$$ as required.\qed\enddemo

To study the lower estimate of $\lz_0$, we observe that not every
positive sequence $(v_i)$ is useful for the lower estimate given in
(9.2) since one may have $\inf_i R_i (v)<0$. In order for $\inf_i
R_i (v)\ge 0$, it is necessary that
$$0<v_i\le \frac{1}{b_i}\bigg(c_i+a_i+b_i-\frac{a_i}{v_{i-1}}\bigg).$$
From this, we obtain the following necessary condition:
$$\frac{a_{i+1}}{c_{i+1}+a_{i+1}+b_{i+1}}< v_i
\le x_i- \dfrac{y_i}{x_{i-1}-
  \dfrac{y_{i-1}}{x_{i-2}-
  \dfrac{y_{i-2}}{\ddots x_3-
    \dfrac{y_3}{x_2-
     \dfrac{y_2}{x_1
}}}}},$$ where
$$x_i=\frac{c_i+a_i+b_i}{b_i},\qqd y_i=\frac{a_i}{b_i}.$$
However, the condition is clearly not practical. Because of this reason, we are
now going to introduce an alternative variational formula for the
lower estimates.

For a given sequence $(r_i)$, define an operator
$I\!I^r=I\!I^{(r_i)}$ of ``double sum'' on the set of positive
functions $(f_i)$ as follows:
$$I\!I_1^r(f)=0,\qd I\!I_i^r(f)
=\sum_{k=1}^{i-1}\frac{1}{\mu_k b_k}\sum_{j=1}^k r_j\mu_j f_j
=\sum_{j=1}^{i-1} r_j f_j\mu_j \sum_{k=j}^{i-1}\frac{1}{\mu_k b_k},
\qd E\ni i\ge 2. $$ Write $I\!I(f)=I\!I^{\text{\bbt{1}}}(f)$. For a fixed
sequence $(c_i)$, let ${\tilde c}_i=c_i-\inf_{i} c_i$ and define
$${\scr F}=\big\{f>0: f_i<f_1+ I\!I_i^{\tilde c}(f) \text{ for all }E\ni i\ge 2\big\}.\tag 9.11$$
Clearly, if ${\tilde c}_1>0$, then every positive constant function
belongs to ${\scr F}$. Otherwise, every $f>0$ with $f_i<f_1$ for all
$E\ni i\ge 2$ belongs to ${\scr F}$.

\proclaim{\thm\;9.6} Let $I\!I^r$, $(\tilde c_i)$ and $\scr F$ be defined as above. Next, for each fixed $f\in {\scr F}$, define
$$\gather
{{\hskip-4em}\xi=\xi_f=\! \cases\displaystyle\!\!
\inf_{E\ni i\ge 2}\frac{f_1-f_i + I\!I_i^{\tilde c}(f)}{I\!I_i(f)},&N=\infty,\\
\displaystyle\!\! \inf_{E\ni i\ge 2}\frac{f_1-f_i + I\!I_i^{\tilde
c}(f)}{I\!I_i(f)} \bigwedge \frac{\sum_{j=1}^N {\tilde c}_j \mu_j
f_j}{\sum_{j=1}^N \mu_j f_j},\;&N<\infty,
\endcases} \tag 9.12\\
{{\hskip-6em}\zz (\ez, f)=\! \cases\displaystyle\!\! \inf_{E\ni i\ge
2,\, {\tilde c}_i<\ez}
\bigg[{\tilde c}_i+ \frac{(\ez -{\tilde c}_i)f_i}{f_1+I\!I_i^{\tilde c -\ez}(f)}\bigg],&\{E\ni i\ge 2:\, {\tilde c}_i<\ez\}\ne\emptyset,\\
\displaystyle\!\! \ez,&\{E\ni i\ge 2:\, {\tilde
c}_i<\ez\}=\emptyset,
\endcases} {\hskip-3em}\tag 9.13\\
\qqd \ez\in [0, \xi].
\endgather
$$
Then we have
$$\lz_0\ge
\inf_{i\in E} c_i+\zz (\ez, f)\qd\text{and}\qd \ez\ge \zz(\ez,
f),\qqd f\in {\scr F},\; \ez\in [0, \xz].\tag 9.14$$ Moreover, for
fixed $f$, $\zz(\ez, f)$ is increasing in $\ez$ and furthermore,
$$\lz_0=\inf_{i\in E} c_i+\sup_{f\in {\scr F}} \zz (\xi, f).
\tag 9.15$$
\endproclaim

\proclaim{\rmk\;9.7} To indicate the dependence on $(\tilde c_i)$,
rewrite $\zz(\ez, f)$ as $\zz({\tilde c}_i, \ez, f)$. Similarly, we
have $\xi({\tilde c}_i, f)$. Then for each $f\in {\scr F}$ and
constant $\gz\ge 0$, we have a shift property as follows:
$$\xi({\tilde c}_i+\gz, f)= \gz+\xi({\tilde c}_i, f),\qqd
\zz({\tilde c}_i+\gz, \ez+\gz, f)=\gz+\zz({\tilde c}_i, \ez, f).
\tag 9.16$$ Hence, the use of $\inf_{i\in E} c_i$ in \thm\;9.6 is not
essential but only for simplifying the computations. The same
property holds for $(9.10)$ but not for $(9.9)$.
\endproclaim

As will be illustrated later by Examples 9.17 and 9.19, it is not
unusual that $\xi_f>\lz_0$ for some $f\in {\scr F}$. In that case,
we certainly have $\xi_f>\zz (\xi_f, f)$. This means that $\xi_f$
may not be a lower bound of $\lz_0$ and so the use of $\zz(\ez, f)$
in \thm\;9.6 is necessary.

\demo{\prf\; of \thm\;$9.6$} By \prp\;9.1, for simplicity, we
assume that ${\tilde c}_i\equiv c_i$.

(a) First, we prove ``$\lz_0\ge$'' in (9.14). Fix $f\in {\scr F}$.
Then $\xi=\xi_f\ge 0$. Without loss of generality, assume that
$(\xi\ge )\,\ez>0$. Otherwise, the assertion is trivial. Let
$$h_i=f_1+I\!I_i^{c -\ez}(f),\qqd i\in E,\; \ez\in (0, \xz].$$
Since by (9.12),
$$f_1-f_i +I\!I_i^c (f)\ge \ez I\!I_i(f)>0$$
for $E\ni i\ge 2$ and $h_1=f_1>0$, we have $h>0$. Next, define
$v_i=h_{i+1}/h_i$ ($v_0:=\infty$ and $v_N=0$ if $N<\infty$). Then
for $i: 2\le i <N$, since
$$h_i-h_{i+1}=I\!I_i^{c -\ez}(f)-I\!I_{i+1}^{c -\ez}(f)
=\frac{1}{\mu_{i}b_{i}}\sum_{j=1}^{i}(\ez -c_j)\mu_j f_j,$$ we have
$$\align
&a_i\big(1-v_{i-1}^{-1}\big)+b_i(1-v_i)\\
&\qd=\frac{1}{h_i}\big[a_i(h_i-h_{i-1})+ b_i (h_i-h_{i+1})\big]\\
&\qd=\frac{1}{h_i}\bigg[-\frac{a_i}{\mu_{i-1}b_{i-1}}\sum_{j=1}^{i-1}(\ez
-c_j)\mu_j f_j
+\frac{b_i}{\mu_{i}b_{i}}\sum_{j=1}^{i}(\ez -c_j)\mu_j f_j\bigg]\\
&\qd=\frac{(\ez -c_i)f_i}{h_i}.
\endalign$$
This also holds when $i=1$ (noting that $a_1=0$):
$$b_1 (1-v_i)=\frac{b_1}{h_1}(h_1-h_2)=\frac{(\ez -c_1)f_1}{h_1}=\ez -c_1.$$
If $N<\infty$, then at $i=N$, by assumption
$$\ez\le \xi\le \sum_{j=1}^N c_j \mu_j f_j\bigg/\sum_{j=1}^N \mu_j f_j,$$
we get
$$a_N\big(1-v_{N-1}^{-1}\big)+b_N(1-v_N)
= -\frac{a_N}{h_N\mu_{N-1}b_{N-1}}\sum_{j=1}^{N-1}(\ez -c_j)\mu_j
f_j\ge \frac{(\ez -c_N)f_N}{h_N}.
$$
Combining these facts together, we arrive at
$$R_i(v)=c_i +a_i\big(1-v_{i-1}^{-1}\big)+b_i(1-v_i)\ge c_i+ \frac{(\ez
-c_i)f_i}{h_i}, \qqd i\in E. \tag 9.17
$$

We now show that the right-hand side of (9.17) is nonnegative for all $i$ and
so we have ruled out the useless case that $\inf_i R_i(v)< 0$. Since
$h>0$, the assertion is equivalent to
$$c_i h_i \ge (c_i-\ez ) f_i,\qqd i\in E,$$
or
$$c_i\big[f_1-f_i + I\!I_i^c(f)\big]\ge \ez \big[c_i I\!I_i(f)-f_i\big].$$
This is trivial if $c_i I\!I_i(f)\le f_i$ (in particular if $i=1$)
since $f_1-f_i + I\!I_i^c(f)\ge 0$ for all $E\ni i\ge 2$ and $f\in
{\scr F}$. Otherwise, by the definition of $\xz $ and $\ez$, we have
$$f_1-f_i + I\!I_i^c(f)\ge \xz I\!I_i(f)\ge \ez I\!I_i(f)> \ez \big[I\!I_i(f)-f_i/c_i\big],
\qqd E\ni i\ge 2. \tag 9.18$$ We have thus proved the required
assertion.

By \thm\;9.3 and (9.17), we obtain
$$\align
\lz_0&\ge \sup_{f\in {\scr F}}\inf_{i\in E}
\bigg[c_i+ \frac{(\ez -c_i)f_i}{f_1+I\!I_i^{c -\ez}(f)}\bigg]\tag 9.19\\
&=\sup_{f\in {\scr F}}\bigg\{\ez \wedge\inf_{E\ni i\ge 2} \bigg[c_i+
\frac{(\ez -c_i)f_i}{f_1+I\!I_i^{c -\ez}(f)}\bigg]\bigg\}.
\endalign$$
Here, the last line is due to the fact that $I\!I_1^r(f)=0$.

(b) To prove the first assertion of the theorem, we show that for each $i$: $2\le i\in E$,
$$c_i+ \frac{(\ez -c_i)f_i}{f_1+I\!I_i^{c -\ez}(f)}\ge \ez \qd\text{iff}\qd c_i\ge \ez.$$
Clearly, the inequality is equivalent to
$$(\ez -c_i)f_i\ge (\ez -c_i)[f_1+I\!I_i^{c -\ez}(f)].$$
The required assertion then follows since by (9.18), we already have
$$f_i\le f_1+I\!I_i^{c -\ez}(f).$$
As a consequence of the assertion, we have $\ez\ge \zz(\ez, f)$. Now, from (9.19), it
follows that
$$\lz_0\ge \sup_{f\in {\scr F}} \ez\wedge \zz(\ez, f)= \sup_{f\in {\scr F}} \zz(\ez, f).$$
This gives us the first assertion of the theorem.

(c) To prove the monotonicity of $\zz(\ez, f)$ in $\ez$, let
$\ez_1<\ez_2\le \xz$. If $\{E\ni i\ge 2: c_i<\ez_2\}=\emptyset$, then
$\{E\ni i\ge 2: c_i<\ez_1\}=\emptyset$ and so
$$\zz (\ez_2, f)=\ez_2>\ez_1=\zz(\ez_1, f).$$
If $\{E\ni i\ge 2: c_i<\ez_1\}\ne\emptyset$, since
$\{E\ni i\ge 2: c_i<\ez_1\}\subset \{E\ni i\ge 2: c_i<\ez_2\}$,
we need only to show that
$$\frac{\big(\ez_2-c_i\big) f_i}{f_1+I\!I_i^{c-\ez_2}(f)}\ge
\frac{\big(\ez_1-c_i\big) f_i}{f_1+I\!I_i^{c-\ez_1}(f)}\qqd \text{on
}\; \{E\ni i\ge 2: c_i<\ez_2\}\ne \emptyset.$$
Actually, this is enough even if $\{E\ni i\ge 2: c_i<\ez_1\}=\emptyset$ in view of (b).
Now, the required conclusion is trivial on the set $
\{E\ni i\ge 2: \ez_1\le c_i<\ez_2\}$. Hence, it suffices to show that
$$\frac{\ez_2-c_i}{f_1+I\!I_i^{c-\ez_2}(f)}\ge
\frac{\ez_1-c_i}{f_1+I\!I_i^{c-\ez_1}(f)}\qqd \text{on }\; \{E\ni
i\ge 2: c_i<\ez_1\}.$$ A simple computation shows that this is
equivalent to
$$f_1+I\!I_i^c(f)\ge c_i I\!I_i(f),$$
which holds on $\{E\ni i\ge 2: c_i<\ez_1\}$ in view of (9.12) and $\xi>\ez_1$.

(d) To prove (9.15), it suffices to show that the equality in (9.19)
holds for $\ez=\xz$. Noting that the right-hand side of (9.19) is
nonnegative, without loss of generality, we may assume that
$\lz_0>0$. Then, by \prp\;2.1, the eigenfunction $g>0$ of $\lz_0$
satisfies
$$\mu_k b_k(g_k-g_{k+1})=\sum_{j=1}^k (\lz_0-c_j)\mu_j g_j,\qqd k\in E, \; g_{N+1}=0\text{ if }N<\infty.$$
Hence,
$$g_1-g_i=I\!I_i^{\lz_0-c}(g),\qd i\in E,
\qqd \sum_{j=1}^N (\lz_0-c_j)\mu_j g_j=0\qd\text{if }N<\infty$$ and
furthermore, $g\in {\scr F}$. It follows that
$$\align
&\frac{g_1-g_i+I\!I_i^{c}(g)}{I\!I_i(g)}\equiv \lz_0,\qqd E\ni i\ge 2,\\
&\frac{\sum_{j=1}^N c_j \mu_j g_j}{\sum_{j=1}^N \mu_j g_j}=\lz_0\qqd\text{if }N<\infty,\\
&{c}_i+ \frac{(\lz_0-{c}_i)g_i}{g_1+I\!I_i^{{c}-\lz_0}(g)}\equiv
\lz_0,\qqd i\in E.
\endalign$$
Therefore, $\xi_g=\lz_0$, and furthermore, the equality sign in (9.19)
is attained at $(f, \ez )=(g, \lz_0)$. \qed\enddemo

We now make a rough comparison of \thm s 9.6 and 9.3 for the lower
estimate. See also the comment below the proof of \crl\;9.9.

\proclaim{\rmk\;9.8} For a given positive sequence $(v_i)$ such that
$\inf_{i\in E} R_i(v):=\gz_v\ge 0$, corresponding to the sequence
$(f_i)$ and $\xi_f$ defined by (9.5) and (9.12), respectively, we
have $\xi_f\ge \gz_v$.
\endproclaim

\demo{\prf} From the assumption
$$R_i(v)=c_i+a_i\big(1-v_{i-1}^{-1}\big)+b_i (1-v_i)\ge\gz_v=:\gz,\qqd i\in E,$$
it follows that
$$f_k-f_{k+1}\ge \frac{1}{\mu_k b_k}\sum_{j=1}^k (\gz -c_j)\mu_j f_j,
$$
and then
$$f_1-f_i\ge \sum_{k=1}^{i-1}\frac{1}{\mu_k b_k}\sum_{j=1}^k (\gz -c_j)\mu_j f_j
= I\!I_i^{\gz -c} (f),\qqd i\in E.$$ To prove our assertion, without
loss of generality, assume that $\gz>0$. Then it is clear not only
that $f\in {\scr F}$ but also $\xi_f\ge \gz$. \qed\enddemo

As a complement to \rmk\;9.8, it would be nice if we could show that
$$c_i+\frac{(\xi_f-c_i)f_i}{f_1+I\!I_i^{c-\xi_f}(f)}\ge \gz_v\qd\text{on the set $\{E\ni i\ge 2: c_i<\xi_f\}$}.$$
This holds obviously on the subset $\{\gz_v\le c_i < \xi_f\}$, but is
not clear on the subset $\{E\ni i\ge 2: c_i<\gz_v\}$.

The next result is a particular application of \thm\;9.6. It is a
complement of \crl\;9.5. The combination of \prp\;9.4 and \crl\;9.5
with \crl\;9.9 below indicates that when ${\lz}_0(a_i, b_i, 0)$
$=0$, the condition $\lim_{n\to\infty}c_n>0$ is crucial for
$\lz_0(a_i, b_i, c_i)>0$. This is more or less clear in terms of the
Feynman-Kac formula:
$$P_t^c f(x)={\Bbb E}^x\Big[f(X_t)e^{-\int_0^{t\wedge \tz} c_{X_s}^{}\d s}\Big],$$
where $\{P_t^c\}_{t\ge 0}$ is the minimal semigroup generalized by
the operator with rates $(a_i, b_i, c_i)$, $\{X_t\}_{0\le t< \tau}$
is the minimal process with rates $(a_i, b_i)$, and $\tau$ is the
life time of $\{X_t\}$. Note that ${\lz}_0(a_i, b_i, 0)>0$, and hence,
${\lz}_0(a_i, b_i, c_i)>0$ if the uniqueness condition (1.2) fails.
Otherwise, $\tau=\infty$.

\proclaim{\crl\;9.9} Let $\vz\in (0, 1)$. Define
$$\gather
{\text{\hskip-4em}\xi_\vz=\! \cases\displaystyle\!\! \inf_{E\ni i\ge
2}\frac{1-\vz +{\tilde c}_1 z_i+\vz x_i}
{z_i+\vz y_i},&\qd N=\infty,\\
\displaystyle\!\! \inf_{E\ni i\ge 2}\frac{1-\vz +{\tilde c}_1
z_i+\vz x_i} {z_i+\vz y_i}\bigwedge \frac{{\tilde c}_1+\vz\sum_{j=2}^N {\tilde c}_j
\mu_j}{1+\vz\sum_{j=2}^N  \mu_j}, &\qd N<\infty,
\endcases}\tag 9.20\\
{\text{\hskip-4em}\zz_\vz= \cases\displaystyle \inf_{E\ni i\ge 2:\,
{\tilde c}_i<\xi_\vz}&\!\!\!\!\bigg[{\tilde c}_i+
\dfrac{\vz(\xi_\vz-{\tilde c}_i)}{1+ {\tilde c}_1 z_i
+\vz x_i - \xi_\vz (z_i +\vz y_i)}\bigg],\\
&\qqd\qd\qqd\{E\ni i\ge 2:\, {\tilde c}_i<\xi_\vz\}\ne\emptyset,\\
\xi_\vz, &\qqd\qd\qqd\{E\ni i\ge 2:\, {\tilde c}_i<\xi_\vz\}=\emptyset,\endcases}\tag 9.21\\
\endgather$$
where
$$x_i=\sum_{2\le j\le i-1} {\tilde c}_j\,\mu_j\,\nu[j, i-1],\qqd
y_i=\sum_{2\le j\le i-1} \mu_j\,\nu[j, i-1],\qqd z_i=\nu[1, i-1],$$
and $\nu[i, j]=\sum_{i\le k \le j}(\mu_k b_k)^{-1}.$
Then we have $\lz_0\ge \inf_{i\in E} c_i+
\sup_{\vz\in (0, 1)} \zz_\vz.$ The same conclusion holds if
$\xi_\vz$ in (9.21) is replaced by $\ez\in [0, \xi_\vz]$. In
particular, if $\varliminf_{n\to\infty}c_n>0$, then $\lz_0>0$.
\endproclaim

\demo{\prf} (a) The main assertion of the corollary is an
application of \thm\;9.6 to the specific $f\in {\scr F}$: $f_1=1$,
$f_i=\vz\in (0, 1)\,(E\ni i\ge 2)$, for which we have
$$
I\!I_1^r(f)=0,\qd I\!I_i^r(f)={r_1}\!\!\!\sum_{1\le k\le
i-1}\frac{1}{\mu_k b_k} +\vz \!\!\!\sum_{2\le j\le i-1}\!\! r_j
\mu_j\!\!\! \sum_{j\le k \le i-1}\frac{1}{\mu_k b_k},\qd E\ni  i\ge
2.$$ Then (9.20) and (9.21) follows from (9.12) and (9.13),
respectively.

We now prove the particular assertion for which $N=\infty$.

(b) If (1.2) does not hold, then ${\lz}_0(a_i, b_i, 0)>0$ by
\thm\;3.1, and so $\lz_0>0$ by part (3) of \prp\;9.1. Similarly, if
$\inf_i c_i>0$, then we have again $\lz_0>0$. Thus, without loss of
generality, assume that
$$\inf_i c_i=0\text{ and } (1.2)\text{ holds}. $$

(c) With the test function $f$ given in (a), by (9.12), we have
$$\xi_\vz
= \inf_{i\ge 2}\frac{1-\vz+I\!I_i^c(f)}{I\!I_i(f)}.
$$
By assumption, there exist $\gz>0$ and $m\ge 2$ such that
 $c_i> \gz$ for all $i\ge m$. Certainly, we have
$$\xi_\vz\ge \inf_{2\le i\le m}\frac{1-\vz+I\!I_i^c(f)}{I\!I_i(f)}
\bigwedge \inf_{i>m}\frac{I\!I_i^c(f)}{I\!I_i(f)}.$$ For $i>m$, we
have
$$\align
\frac{I\!I_i^c(f)}{I\!I_i(f)} &\ge \sum_{j=m}^{i-1} c_j f_j\mu_j
\sum_{k=j}^{i-1}\frac{1}{\mu_k b_k}
\bigg/\sum_{j=1}^{i-1} f_j\mu_j \sum_{k=j}^{i-1}\frac{1}{\mu_k b_k}\\
&> \vz {\gz} \sum_{j=m}^{i-1}\mu_j \sum_{k=j}^{i-1}\frac{1}{\mu_k
b_k} \bigg/\sum_{j=1}^{i-1}\mu_j \sum_{k=j}^{i-1}\frac{1}{\mu_k
b_k}.\endalign$$ By assumption (1.2), the right-hand side goes to
$\vz\gz>0$ as $i\to \infty$. It follows that
$\inf_{i>m}{I\!I_i^c(f)}/{I\!I_i(f)}>0$, and furthermore, there exists
$\ez\in (0, \gz)$ such that $\xi_\vz> \ez$.

(d) Noting that the set $\{i\ge 2: c_i<\ez\}\subset \{i: 2\le i<m\}$
is finite, by (9.13), we have
$$\zz (\ez, f)=\inf_{i\ge 2,\, c_i<\ez}
\bigg[c_i+ \frac{(\ez -c_i)f_i}{f_1+I\!I_i^{c -\ez}(f)}\bigg] \ge
\min_{i\ge 2,\, c_i<\ez} \frac{(\ez -c_i)f_i}{f_1+I\!I_i^{c
-\ez}(f)}>0.$$ Now, by \thm\;9.6 or proof (a) above, we conclude
that $\lz_0>0$. \qed\enddemo

From proof (c) above, we have seen that when $N=\infty$,
$$\xi_\vz>0\;\;\text{ iff }\;\;
\inf_{i>m}I\!I_i^{\tilde c}(\dbl)/I\!I_i(\dbl)>0\;\;\text{ for all
$m\ge 2$}. \tag 9.22$$ Note that
$$\inf_{i>m}\frac{I\!I_i^{\tilde c}(\dbl)}{I\!I_i(\dbl)}\ge
\inf_{i\ge 1} \sum_{j=1}^i \mu_j {\tilde c}_j\bigg/\sum_{j=1}^i \mu_j$$
and the right-hand side is positive iff
$$\varliminf_{m\to\infty} \sum_{j=1}^m \mu_j {\tilde c}_j\bigg/\sum_{j=1}^m \mu_j>0.
\tag 9.23$$ Thus, \crl\;9.9 is qualitatively consistent with
\crl\;9.5.

In view of \rmk\;9.8, it is not obvious that \thm\;9.6 improves
\thm\;9.3. An easier way to see the improvement is as follows.
Recall that the last assertion of \crl\;9.9 is deduced in terms of
the test function $f$ used in its proof (a). For which, the
corresponding sequence $(v_i)$ is $v_1=1/2$ and $v_i=1$ for all
$i\ge 2$. Inserting this into $R(v)$, we get
$$\inf_{i\ge 1}R_i(v)=(c_1+b_1/2)\wedge (c_2-a_2)\wedge \inf_{i\ge 3} c_i.$$
Thus, for $\inf_{i\ge 1}R_i(v)>0$, it is necessary that $\inf_{i\ge
3} c_i>0$, which is clearly much stronger than the last condition
$\varliminf_{n\to\infty}c_n>0$ used in \crl\;9.9.

As \prp\;9.4, the next result is also motivated from \thm\;3.1.

\proclaim{\crl\;9.10} An explicit lower estimate can be obtained by
\thm\;9.6 using the specific test function $f^{(m)}$:
$$f^{(m)}_i=\bigg(\sum_{j=i\wedge m}^m\frac{1}{\mu_j b_j}\bigg)^{1/2}, \qqd i\in E,$$
where $m$ may be optimized over $\{m\in E: m\ge 2\}$ (or over $E$ if
${\tilde c}_1>0$).\endproclaim

We now show that some special killing (or Schr\"odinger) case can be regarded as a perturbation of the one without killing.
To do so, fix constants $\bz, \gz>0$, and define
$$\gather
{\hat a}_i=b_{i-1},\qd 2\le i<N+1,\qqd {\hat b}_i=a_{i+1},\qd 1\le i<N,\\
{\hat a}_1=\bz,\qd {\hat b}_N=\gz \qd\text{if }N<\infty;\\
{\check a}_i=a_{i+1},\qd 0\le i<N,\qqd {\check b}_i=b_i,\qd 1\le i <N+1. \\
{\check b}_0=\bz,\qd {\check a}_N=\gz\qd\text{if }N<\infty.\tag 9.24\endgather
$$
Note that $\big({\hat a}_i, {\hat b}_i\big)$ and $\big({\check a}_i, {\check b}_i\big)$
are dual each other in the sense of Section 5 but they are clearly different
from $({a}_i, {b}_i )$. Recall that $a_1=0$ and $b_N=0$ by convention. Next, let
$(c_i)$ satisfy
$$c_i\ge
\cases
a_{i+1}-a_i -b_i+  b_{i-1},\qd & 2\le i<N,\\
a_2-b_1+\bz,  &i=1,\\
\gz-a_{N}+b_{N-1}, &i=N<\infty.
\endcases\tag 9.25$$
Note that the right-hand side of (9.25) can be negative. Conversely,
for given rates $\big({\hat a}_i, {\hat b}_i\big)$, the inverse
transform is as follows:
$$\gather
a_i={\hat b}_{i-1},\qd 2\le i<N+1,\qqd b_i={\hat a}_{i+1}, \qd 1\le i <N,\\
 c_i\ge {\hat b}_{i}
-{\hat b}_{i-1}-{\hat a}_{i+1}+{\hat a}_{i},\qqd 1\le i <N+1\;\text{ (or $i\in E$)}. \tag 9.26
\endgather$$

\proclaim{\prp\;9.11} Suppose that the given rates $(a_i, b_i, c_i: i\in
E)$ satisfy (9.25). Define $\lz_0(a_i, b_i, c_i)$ as in \prp\;9.1 without preassuming that
$c_i\ge 0$ for all $i\in E$. Next, define
$\big({\hat a}_i, {\hat b}_i\big)$ and $\big({\check a}_i, {\check b}_i\big)$ by (9.24).
\roster
\item If $\sum_{i=2}^N \mu_i b_{i-1}^{-1}=\infty$, then $\lz_0(a_i, b_i, c_i)\ge {\hat\lz}_0$,  where ${\hat\lz}_0$ is defined by (4.1) with rates $\big({\hat a}_i, {\hat b}_i\big)$.
\item Otherwise, $\lz_0(a_i, b_i, c_i)\ge {\check\lz}_1$,
  where ${\check\lz}_1$ is defined by (6.1) with rates $\big({\check a}_i, {\check b}_i\big)$.
\item The equality sign of the conclusions in parts $(1)$ and $(2)$ holds
provided it does in $(9.25)$.
\endroster
\endproclaim

\demo{\prf} (a) As an application of \prp\;9.1, without loss of generality, we may and will assume that the equality sign for $c_i$ in (9.26) holds. Then, we prove that the equality sign of the conclusions in parts $(1)$ and $(2)$ holds.

Clearly, we have
$${\hat\mu}_1=1,\;\; {\hat\mu}_1 {\hat a}_1=\bz,\;\;
{\hat\mu}_i=\mu_i^{-1},\;\;
{\hat\mu}_i {\hat a}_i=\frac{b_{i-1}}{\mu_i},\qqd
2\le i <N+1.\tag 9.27$$

(b) Recall the operators:
$$\align
\ooz f(i)&=b_i (f_{i+1}-f_i)+a_i(f_{i-1}-f_i)-c_i f_i,\\
{\widehat\ooz} f(i)&={\hat b}_i (f_{i+1}-f_i)+ {\hat a}_i(f_{i-1}-f_i),
\qd f\in {\scr K}, f_0=0, f_{N+1}=0\text{ if }N<\infty.
\endalign$$
Clearly, $\lz_0(a_i, b_i, c_i)$ is the principal eigenvalue of $\ooz$ and the idea is
describing it in terms of the first eigenvalue ${\hat\lz}_{\min}$ of ${\widehat\ooz}$.
Let $U$ be the diagonal matrix with diagonal elements $(\mu_i: i\in E)$. Then
$U^{-1}$ is simply the diagonal matrix with diagonal elements $({\hat\mu}_i: i\in E)$.
For each function $h$ with $h_0=0$ and $h_{N+1}=0$ if $N<\infty$, by (9.27), (9.24) and (9.26), we have
$$\align
\big(\ooz  U^{-1}h\big)(i)&=b_i\big( {\hat \mu}_{i+1}h_{i+1}- {\hat \mu}_{i}h_{i}\big)
+a_i \big( {\hat \mu}_{i-1}h_{i-1}- {\hat \mu}_{i}h_{i}\big)-c_i {\hat \mu}_{i} h_i\\
&={\hat a}_{i+1}\big( {\hat \mu}_{i+1}h_{i+1}- {\hat \mu}_{i}h_{i}\big)
+{\hat b}_{i-1} \big( {\hat \mu}_{i-1}h_{i-1}- {\hat \mu}_{i}h_{i}\big)\\
&\qd -\big({\hat b}_{i}-{\hat b}_{i-1}-{\hat a}_{i+1}+{\hat a}_{i}\big) {\hat \mu}_{i} h_i\\
&=\big({\hat a}_{i+1} {\hat \mu}_{i+1}h_{i+1}
-{\hat b}_{i} {\hat \mu}_{i}h_{i}\big)+\big({\hat b}_{i-1} {\hat \mu}_{i-1}h_{i-1}- {\hat a}_i {\hat \mu}_{i} h_{i}\big)\\
&={\hat \mu}_{i}{\hat b}_{i} ( h_{i+1}-h_{i})+ {\hat \mu}_{i}{\hat a}_{i} ( h_{i-1}-  h_{i} )\\
&={\hat \mu}_{i} {\widehat\ooz}\, h(i)\\
&=\big(U^{-1} {\widehat\ooz}\, h\big)(i),\qqd 2\le i <N.
\endalign
$$
It is easy to check that the identity holds also for $i=1$ and $i=N$, and then for all
$i\in E$. Multiplying $U$ from the left on the both sides,  we obtain
$$U \ooz  U^{-1}={\widehat\ooz}.\tag 9.28$$
Furthermore, we get
$$\langle f, \ooz  g\rangle_{\mu}
= \langle U^{-1}f, (U \ooz  U^{-1})U g\rangle_{\mu}
= \langle U f, (U \ooz  U^{-1})U g\rangle_{\hat\mu}
=\langle {\hat f}, {\widehat\ooz} {\hat g}\rangle_{\hat\mu}$$
for all $f, g\in {\scr K},$
where the mapping $f\to {\hat f}:=Uf$ is an isometry from $L^2(\mu)$ to  $L^2(\hat\mu)$.
Since $f\in {\scr K}$ iff ${\hat f}\in {\scr K}$, it follows that the operators $\ooz $
and ${\widehat\ooz}$ with the same core ${\scr K}$ are isospectral. In particular, $\lz_0(a_i, b_i, c_i)={\hat\lz}_{\min}$.

(c) For assertion (1), since $\sum_{i}\big({\hat \mu}_i {\hat b}_{i}\big)^{-1}=\infty$ by assumption,
it follows that $N=\infty$ and the Dirichlet form corresponding to ${\widehat\ooz}$ is regular by \prp\;1.3.
Hence, the
minimal and the maximal domains of the Dirichlet form are coincided. Therefore,
${\hat\lz}_{\min}$ is equal to $\lz_0^{(4.1)}$ replacing the original rates $(a_i, b_i)$
by $\big({\hat a}_i, {\hat b}_i\big)$.

For assertion (2), since ${\hat a}_1>0$ and $\hat b_N>0$, we come to the setup of Section 7:
${\hat\lz}_{\min}=\lz_0^{(7.1)}$ with $(a_i, b_i)$ replaced
by $\big({\hat a}_i, {\hat b}_i\big)$. Next, because of $\sum_{i}\big({\hat \mu}_i {\hat b}_{i}\big)^{-1}<\infty$, by \thm\;7.1, it turns out ${\hat\lz}_{\min}={\check \lz}_1$
in terms of the dual rates  $\big({\check a}_i, {\check b}_i\big)$ of $\big({\hat a}_i, {\hat b}_i\big)$.
\qed\enddemo

We now summarize our main qualitative result about $\lz_0$. The three
parts given below are obtained by \crl\;9.9, \prp\;9.1 plus
\prp\;1.3, and \crl\;9.5, respectively.

\proclaim{Summary 9.12} We have $\lz_0>0$ whenever $N<\infty$. Next,
let $N=\infty$. Then \roster
\item $\lz_0>0$ if $\varliminf_{n\to\infty} c_n>0$ (in particular, if $\inf_i c_i>0$).
\item $\lz_0=\lz_0(a_i, b_i, c_i)>0$ if $\lz_0\big(a_i, b_i, c_i{\dbl}_{\{1\}}\big)>0$ which
can be checked case by case by

(i) \thm\;3.1 when $c_1=0$ and
$$\sum_{i\ge 1}\frac{1}{\mu_i a_i}<\infty; \tag 9.29$$

(ii) \thm s 7.1 and 6.2 when $c_1>0$ and (9.29) holds;

(iii) \thm\;4.2 when $c_1>0$ but (9.29) fails.
\item $\lz_0=0$ if
$$\varliminf_{m\to\infty}\sum_{i=1}^m \mu_i c_i\bigg/ \sum_{k=1}^m \mu_k=0,
\qd \varlimsup_{m\to\infty}\mu_m b_m
\bigg(\sum_{i=1}^m\mu_i\bigg)^{-1}=0\qd \text{and}\qd\sum_{i}
\mu_i=\infty.$$
\endroster
\endproclaim

\proclaim{Open problem\;9.13 (Explicit criterion for $\lz_0>0$)}{\rm
As will be seen soon in Example 9.18 below, for $\lz_0>0$, it can happen that
$\varliminf_{n\to\infty} c_n=0$. Hence, the simple
condition ``$\varliminf_{n\to\infty} c_n>0$'' in part (1) is
sufficient only but not necessary. Naturally, this condition becomes
necessary for the first one in part (3) for which a sufficient
condition is $\lim_{n\to\infty} c_n=0$. Thus, it is more or less
satisfactory whenever $(c_n)$ has a limit. Otherwise, there is a
gap. In contrast with the first condition in part (3), condition
(9.23) is sufficient for $\xi_\vz>0$ but there is still a distance
to deduce the positivity of $\lz_0$ in view of \crl\;9.9.

Next, since we are dealing with the minimal Dirichlet form, a
general criterion for Hardy-type inequalities (cf. [12;
\thm s\;7.1 and 7.2]) which was successfully used in Section 8,
is also available in the present situation, hence, there is
already a criterion for $\lz_0>0$ in terms of capacity which is
unfortunately not explicit. More seriously, the technique to produce an
explicit result used in [12; pages 134--136] does not work
at the beginning (replacing a finite number of disjointed finite
intervals $\{K_i\}$ by the connected one $[\min \cup_i K_i,\, \max
\cup_i K_i]$) in the present setup. Thus, it is still an unsolved
problem to figure out an explicit criterion for $\lz_0>0$ in the present setup.}
\endproclaim

It is our position to illustrate by examples the application of the
results obtained in this section. First, by using \prp\;9.11 and
(9.26), it is easy to transfer the examples given in Sections 3 and
6 to the present context. However, most of the resulting killing
rates are rather simple. We are now going to construct some new
examples, all of them are out of the scope of \prp\;9.11. In the
most cases, we use simple $(a_i, b_i)$ and pay more attention on
$(c_i)$. Let us begin with the following simplest case.
\proclaim{\xmp\;9.14} Let
$$Q=
\pmatrix
-b_1 - c_1 & b_1\\
a_2 & -a_2 - c_2
\endpmatrix.$$
Then as in Examples 7.5\,(2), we have
$$\lz_0=
\frac 1 2 \Big(a_2 + b_1 + c_1 + c_2 - \sqrt{
   (a_2  + c_2 - b_1 - c_1)^2 + 4 a_2 b_1}\,\Big),$$
with eigenvector
$$g=\bigg(\frac{1}{2 a_2}\Big[a_2  + c_2 - b_1 - c_1 + \sqrt{
 (a_2  + c_2 - b_1 - c_1)^2 + 4 a_2 b_1}\,\Big],\, 1\bigg).
$$
\endproclaim

Even in such a simple case, the role for $\lz_0$ played by the
parameters $a_i$, $b_i$, and $c_i$ is ambiguous. For instance, since
${\tilde c}_1-{\tilde c}_2=c_1-c_2$ $({\tilde c}_k:=c_k-c_1\wedge
c_2)$, one can separate out the constant $c_1\wedge c_2$ from the
above expression of $\lz_0$. However, this obvious separation
property becomes completely mazed for the next example having three
states only.

\proclaim{\xmp\;9.15} Let
$$Q =
\pmatrix
-b_1 - c_1 & b_1 & 0\\
a_2& -a_2 - b_2 - c_2 & b_2\\
0 & a_3 & -a_3 - c_3\endpmatrix.$$ Then
$$\lz_0=-\frac 1 3 \gz_1+2 \sqrt{\frac{-U}{3}}\, \cos \bigg[\frac{1}{3}\,
\text{arc\,cos}\bigg(-\frac{V}{2}\bigg(\frac{-U}{3}\bigg)^{-3/2}\bigg)+\frac{2
\pi}{3}\bigg], \tag 9.30$$ with eigenvector
$$g=\bigg\{\frac{b_1 (a_3+c_3-\lz_0)}{a_3 ({b_1+c_1}-\lz_0)},\,
   1+\frac{{c_3}-\lz_0}{a_3},\, 1\bigg\},$$
where
$$U=\gz_2-\gz_1^2/3,\qqd
V=\gz_3-\gz_1 \gz_2/3+ 2(\gz_1/3)^3,$$ and $\lambda^3 + \gz_1
\lambda^2 +\gz_2 \lambda+\gz_3$ the eigenpolynomial  of $-Q$ with
coefficients:
$$\align
\gz_1&=-{a_2}-{a_3}-{b_1}-{b_2}-{c_1}-{c_2}-{c_3}, \\
\gz_2&= {a_2} {a_3}+{b_1}
   {a_3}+{c_1} {a_3}+{c_2} {a_3}+{b_1} {b_2}+{a_2} {c_1}+{b_2} {c_1}+{b_1}
   {c_2}+{c_1} {c_2}+{a_2} {c_3}\\
   &\qd +{b_1} {c_3}+{b_2} {c_3}+{c_1} {c_3}+{c_2}
   {c_3},\\
\gz_3&=-{a_2} {a_3} {c_1}-{a_3} {b_1} {c_2}-{a_3} {c_1} {c_2}-{b_1}
{b_2}
   {c_3}-{a_2} {c_1} {c_3}-{b_2} {c_1} {c_3}-{b_1} {c_2} {c_3}-{c_1} {c_2}
   {c_3}.
\endalign$$
\endproclaim

\demo{\prf} Since the eigenvalues of $-Q$ are all real, it is easier
to write them down. By using the notation given above, the
eigenvalues of $-Q$ can be expressed as
$$-\frac 1 3 \gz_1+2 \sqrt{\frac{-U}{3}}\, \cos \bigg[\frac{1}{3}\,
\text{arc\,cos}\bigg(-\frac{V}{2}\bigg(\frac{-U}{3}\bigg)^{-3/2}\bigg)+\frac{2
k\pi}{3}\bigg],\qquad k=0, 1, 2. $$ Among them, the minimal one is
$\lz_0$ given in (9.30). Clearly, the solution is indeed rather
complicated in view of the coefficients of the
eigenpolynomial.\qed\enddemo

To see the role played by the killing rate $(c_i)$, in the following
examples, we restrict ourselves to the case that $\lz_0(a_i, b_i,
0)=0$ (and then $N=\infty$). The examples are arranged according the
increasing order of the polynomial rates $(a_i)$ and $(b_i)$.
Actually, all the examples in the paper are either standard or constructed
by using simple rates and simple eigenfunctions. They are used first as a
guidance of the study and then to justify the power of the theoretic results.

In contract to the explosive case (cf. \thm\;3.1), $\lz_0$ can still
be zero for the process having positive killing rate, as shown by
the following example.

\proclaim{\xmp\;9.16} Let $b_1=1$, $a_i=b_i=1$ for $i\ge 2$, and let
$(c_i)$ satisfy $\lim_{n\to\infty} c_n=0.$ Then we have $\lz_0=0$,
even though $c_i$ can be very large locally.
\endproclaim

\demo{\prf} Apply \crl\;9.5.\qed
\enddemo

\proclaim{\xmp\;9.17} Let $a_i=b_i=1$ for $i\ge 2$ and
$c_i=\bz^{-1}(\bz -1)^2\,(\bz>0)$ for $i\ge 1$. Then for every
$a_1\ge 0$ and $b_1>0$, we have $\lz_0=\bz^{-1}(\bz -1)^2$.
\endproclaim

\demo{\prf} Since $\lz_0(a_i, b_i, 0)=0$ and $(c_i)$ is a constant,
this is a consequence of part (3) of \prp\;9.1.

Note that the lower estimates given by \prp\;9.1, \thm\;9.3, and
(9.4) are all sharp for this example. To see this, simply choose
$v_i\equiv 1$ in (9.2) and (9.4). We now consider a more specific
situation: $a_1=0$, $b_1=1-\bz$ and $\bz\in (0, 1)$. If we set
${\bar v}_i\equiv\bz^{-1}$, then it is easy to check that $R_i(\bar
v)\equiv 0$ and so $\inf_{v>0}\sup_{i\ge 1} R_i(v)=0$. This shows
that the truncating procedure used in \thm\;9.3 for the upper
estimate is necessary in the case that the function $f$ defined by
(9.5) does not belong to $L^2(\mu)$, even though $R_i (v)$ is a
constant. In the present case, ${\bar v}_i>1$ for all $i$ and so the
corresponding function $f$ is strictly increasing. Since $\mu_i $ is
a constant for $i\ge 2$, it is clear that $\sum_i\mu_i=\infty$ and
then $f\notin L^2(\mu)$. \qed\enddemo

\proclaim{\xmp\;9.18} Let $a_1=0$, $b_1=5/2$, $a_i=2$ and $b_i=1$
for $i\ge 2$, $c_i=0$ for odd $i$ and $c_i=13/6$ for even $i$. Then
$\lz_0=5/6$. The upper bound provided by (9.9) is approximately $1.03$.
For the lower estimate, \prp\;9.11 is available but not \crl\;9.9.
\endproclaim

\demo{\prf} (a) Let $v_i\equiv 1+(-1)^i/3$. Then it is easy to check
that $R_i(v)\equiv 5/6$. Next, define
$$g_1=1,\qd g_n=\prod_{k=1}^{n-1} v_k,\qqd n\ge 2.\tag 9.31$$
We claim that $g\in L^2(\mu)$ by using Kummer's test. To do so,
note that to study the convergence/divergence of the series $\sum_n
\mu_n g_n^2$, the constant $\kz$ defined by (3.13) takes a simpler
form as follows:
$$\kz=\lim_{n\to \infty} n\bigg(\frac{a_{n+1}}{b_n v_n^2}-1\bigg).  \tag 9.32$$
Now, because $g\in L^2(\mu)$ and
$$- \ooz g /g= R(v)=5/6,$$
we have $\lz_0= 5/6$ by \thm\;9.3. Clearly, this eigenfunction $g$
of $\lz_0$ is not monotone since $g_{i+1}/g_i=v_i= 2/3$ for odd $i$
and $=4/3$ for even $i$.

(b) Next, we study the upper estimates of $\lz_0$. First, we have
$$\mu_1=1,\;\; \mu_i=\frac{5}{2^i},\qd i\ge 2;\qqd
 \mu_i b_i=\frac{5}{2^i},\qd i\ge 1.$$
The upper bound provided by (9.9) is approximately $1.03$.

(c) For a lower estimate, we apply \prp\;9.11\,(2). The modified rates are as follows: ${\check a}_i\equiv 2\,(i\ge 1)$,
${\check b}_1=5/2$, and ${\check b}_i\equiv 1\,(i\ge 2)$. However, $(c_i)$ does not satisfy (9.25)
at $i=2$. We now replace $(c_i)$ by $({\tilde c}_i:= c_i+1/6)$ and choose
${\check b}_0=\bz=2/3$. Then $({\tilde c}_i)$ satisfy (9.25).
With the modified $({\tilde c}_i)$, we are in the ergodic case, and moreover,
${\check \lz}_1=\big(\sqrt{2}-1\big)^2$ with eigenfunction $\check g$: $\check g_0=-1$ and
$$\check g_i= \frac{1}{20}  2^{i/2} \Big[-101+60 \sqrt{2}+  \left(41-25 \sqrt{2}\right)  i\Big],\qqd i\ge 1.$$
Therefore, by \prp\;9.11\,(2), we obtain
$\lz_0(a_i, b_i, {\tilde c}_i)\ge \big(\sqrt{2}-1\big)^2$. Returning to the original  $(c_i)$ by \prp\;9.1\,(2),
we get a rough lower bound as follows:
$$\lz_0=\lz_0(a_i, b_i, {\hat c}_i)-\frac 1 6 \ge \frac{17}{6}-2\sqrt{2}\approx 0.005.$$

Before moving further, let us remark that if only ${\check b}_0$ is changed from $2/3$ to $1/2$, then for the
$\big({\check a}_i, {\check b}_i\big)$-process, we still have
${\check \lz}_1=\big(\sqrt{2}-1\big)^2$ with a similar eigenfunction $\check g$: $\check g_0=-1$ and
$$\check g_i= \frac{1}{10}  2^{i/2} \Big[-67+42 \sqrt{2}+  \left(27-17 \sqrt{2}\right)  i\Big],\qqd i\ge 1.$$
Now, as an application of \prp\;9.11\,(2) with the original $(c_i)$ replacing $c_2=4/3$ by $c_2=3/2$ only,
the resulting $\lz_0$ has a lower bound $\big(\sqrt{2}-1\big)^2\approx 0.17$.

(d) To apply \crl\,9.9, we write $c_i=13 (1+(-1)^i)/12$ and use (9.20) and (9.21):
$$\gather
\text{\hskip-4em}\xi_\vz=\inf_{i\ge 2}\frac{1-\vz + \vz x_i} {z_i+\vz y_i},\;\;
\zz_\vz=\inf_{\text{odd } i\ge 3}
\frac{\vz \,\ez}{1 + \vz x_i  -\ez (z_i+\vz y_i)},\;\; \ez\in [0, \xi_\vz],  \tag 9.33\\
x_i=\frac{13}{72}\big( 2^{2+i}-6\, i -3 -(-1)^i \big), \qd
y_i=2^{i-1}-i,
 \qd z_i=\frac{2^i-2}{5}.\endgather$$
Note that the numerator of $\zz_{\vz}$ given in (9.33) is independent of
$i$ but in the denominator, $x_i$, $y_i$ and $z_i$ all tend to
infinity as $i\to\infty$. To avoid the trivial estimate, one needs to cancel the leading term in $i$ of $\vz x_i
-\ez (z_i+\vz y_i)$ in the denominator. This leads to the following
solution:
$$\ez=\frac{65\, \vz }{9 (2 + 5\,\vz)}.$$
Inserting this into $\vz x_i -\ez (z_i+\vz y_i)$, it follows that the denominator
of $\zz_{\vz}$ in (9.33) becomes
$$1-\frac{65 \left(-2\, i+(-1)^i+3\right) \vz^2+2 \left(78\, i+13 (-1)^i-245\right) \vz-144}{72 (5 \,\vz+2)}.$$
Now, in order to remove the leading term in $i$, the only solution is
$$\vz=78/65>1,$$
which does not belong to the domain of $\vz\in (0, 1)$.
Therefore, the test function used in \crl\;9.9 does not provide enough freedom to cover this example.

Note that without the killing rate, the process with rates $(a_i)$ and $(b_i)$
is exponentially ergodic and so $\lz_0(a_i, b_i, 0)=0$.\qed\enddemo

For the following examples, we assume that $a_i=b_i$ for $i\ge 2$.
Then
$$\mu_1=1,\qd \mu_i={b_1} a_i^{-1},\qd i\ge 2;\qd \mu_i b_i=b_1,\qd i\ge 1.$$
The quantities $\xi_\vz$ and $\zz_\vz$ defined in (9.20) and (9.21),
respectively,
 are now determined by
$$x_i=\sum_{2\le j\le i-1} \frac{i-j}{a_j}{\tilde c}_j,\qqd
y_i=\sum_{2\le j\le i-1} \frac{i-j}{a_j},\qqd z_i=\frac{i-1}{b_1}.$$

\proclaim{\xmp\;9.19} Let $a_1=0$, $b_1=2\bz
(1-\bz)(1-2\bz)^{-1}\,(\bz\in (0, 1/2))$, $a_i=b_i=\bz\, i$ for
$i\ge 2$, $c_i=(1-\bz)^2 (i-1)$ for $i\ge 1$. Then $\lz_0=2\bz
(1-\bz)$. In the special case that $\bz=1/4$, we have $\lz_0=3/8$.
The upper and lower bounds provided by (9.9) and \crl\;9.9 are $3/4$
and approximately $0.274$, respectively.
\endproclaim

\demo{\prf} Let $v_i=\bz (1+i^{-1})$ for $i\ge 1$. Then
$R_i(v)\equiv 2\bz (1-\bz)$. By Kummer's test (cf. (9.32)), the
corresponding function $g$ defined by (9.31) belongs to $L^2(\mu)$.
Hence, the assertion follows from \thm\;9.3. Note that $v_i<1$ for
all $i$, and $g$ is strictly decreasing even though $c_1=0< \lz_0$ and
$c_i>\lz_0$ for all $i>(1+\bz)(1-\bz)^{-1}$ (compare with (2.5)).

As an application of (9.9) with $(\ell, m)=(1, 1)$ or $(3, 4)$, we
obtain
$$\lz_0\le (1- \bz)\bigg\{\frac{2\bz }{1 - 2 \bz} \bigwedge
\frac{23 - 40 \bz +  23 \bz^2}{2(8 -  11\bz)}\bigg\}.$$ To study the
lower bound, for simplicity, we let $\bz=1/4$. Then $\lz_0=3/8$ and
the upper bound in the last formula is $3/4$. Choose
$\vz=\big(\sqrt{409}-5\big)/24$ so that the infimum $\xi_\vz=\big(29
- \sqrt{409}\,\big)/32\approx 0.274$ is attained simultaneously at
$i=2$ and $i=3$. Since $\xi_\vz< c_2$, the set $\{i\ge 2: c_i<
\xi_\vz\}$ is empty. Therefore, the lower bound provided by
\crl\;9.9 is approximately $0.274$.\qed\enddemo

For the following two examples, without the killing rate, the
process is exponentially ergodic and so $\lz_0(a_i, b_i, 0)=0$.

\proclaim{\xmp\;9.20} Let $a_1=0$, $b_1=4/5$, $a_i=b_i=i^2$ for
$i\ge 2$, and
$$c_i=\frac{8}{9}\bigg[\frac{8}{3\, i-8} - \frac{2}{3\, i-4} + 5\bigg],\qqd i\ge 1.$$ Then $\lz_0=4$.
The upper and lower bounds provided by (9.9) and \crl\;9.9 are
$14/3$ and approximately $2.82$, respectively.
\endproclaim

\demo{\prf} The proof is similar as before using
$$v_i=1-\frac{1}{3\,i-4},\qqd i\ge 1.$$
Note that $c_i$ has minimum $0$ at $i=2$. The upper bound provided
by (9.9) with $(\ell, m)=(2, 2)$ is $14/3$. The lower bound produced
by  \crl\;9.9 with $\vz=1$ is $48/17\approx 2.82$. Since $c_1>0$,
the parameter $\vz=1$ is allowed. Then $\xi_\vz=48/11$ is attained
at $i=3$, and $\zz_\vz=48/17$ is attained at $i=2$ (noting that the
set $\{i\ge 2: c_i<\xi_\vz\}$ is a singleton $\{2\}$). \qed\enddemo

\proclaim{\xmp\;9.21} Let $a_1=0$, $b_1=3/2$, $c_1=15$,
$$\align
a_i&=b_i= i\, (i-4^{-1}) (12\, i^2    - 31 i + 27), \qqd i\ge 2,\\
c_i&= i^4 -\frac 1 2\, i^3 -\frac{301}{16}\, i + \frac{227}{8}, \qqd
i\ge 2.
\endalign$$
Then $\lz_0=119/8=14.875$. The upper and lower bounds provided by
(9.9) and \crl\;9.9 are approximately $15.42$ and $13.18$,
respectively.
\endproclaim

\demo{\prf} Note that $c_i$ is convex and has its minimum $0$ at
$i=2$. For
$$v_i=\frac 3 4 - \frac 2 i + \frac{7}{4\, i-1},\qqd i\ge 1,$$
we have $R_i(v)\equiv 119/8$. Note that $v_1>1$ and $v_i<1$ for all
$i\ge 2$. The function $g$ defined by (9.31) is not monotone but is
bounded. Next, since $\mu_i\sim i^{-4}$, we have $g\in L^2(\mu)$.
The assertion now follows from \thm\;9.3.

Clearly, $\inf_{i\ge 1} c_i=11/4$. The upper bound provided by (9.9)
with $(\ell, m)=(2, 4)$ is approximately $15.42$. To get a lower
estimate, we apply \crl\;9.9. Because ${\tilde c}_1>0$, we can
choose $\vz=1$. Then $\xi_\vz=354679/29504$ is attained at $i=4$.
Next, since the set $\{i\ge 2: c_i<\xi_\vz\}$ is a singleton
$\{2\}$, we need only to compute $\zz_\vz$ at $i=2$: $\zz_\vz\approx
10.43$. Thus, the lower bound produced by \crl\;9.9 is approximately
$13.18$. \qed\enddemo

To conclude this section, we return to the uniqueness problem for
birth--death processes with killing of the Dirichlet form as
discussed at the end of Section 1. Certainly, the problem is
meaningful only if $N=\infty$. Recall that for a given $Q$-matrix,
not necessarily conservative (i.e., may have killing), the exit
space ${\scr U}_\lz$ is the set of the solutions $(u_i)$ to the
following equation:
$$\cases
(\lz I -Q)u=0,\\
0\le u\le 1,
\endcases \qqd \lz>0.$$
Note that the dimension of ${\scr U}_\lz$ is independent of $\lz>0$.
By (2.5) replacing $\lz$ with $-\lz$, it follows that the
non-trivial exit solution, if it exists, is unique and is strictly
increasing.

\proclaim{\thm\;9.22 (Uniqueness of the Dirichlet form)} Let
$N=\infty$. \roster
\item The Dirichlet form satisfying the Kolmogorov's equations is unique
if ${\scr U}_\lz=\{0\}$. Equivalently,
$$\sum_{n=1}^\infty \frac{1}{\mu_n b_n}\sum_{k=1}^n \mu_k (1+c_k)=\infty. \tag 9.34$$
\item Let $\sum_{i\in E} \mu_i<\infty$. Then the Dirichlet form is unique iff
$$\sum_{i\in E}\mu_i c_i<\infty\qd\text{and}\qd \sum_{i\in E}\frac{1}{\mu_i b_i}=\infty.$$
\item Let $\sum_{i\in E} \mu_i\!=\!\infty$. Then the Dirichlet form is unique if
$$\text{either }\inf\limits_{i\in {\overset{\bullet}\to{E}}} \sum_{j\in
E}P_{ij}^{\min}(\lz)\!>\!0\qd \text{or}\qd \sum_{i\in E}\mu_i
c_i<\infty$$ holds, where ${\overset{\bullet}\to{E}}=\{i\in E:
c_i>0\}$.
\endroster
\endproclaim

Here are some comments about the theorem. \roster
\item"(i)" Suppose that only a finite number of $c_i$ are non-zero.

Then condition (9.34) is equivalent to (1.2) [Certainly in this
item, we are using the modified (1.2) and (1.3) by removing $0$ from
the state space]:
$$\sum_{n=1}^\infty \frac{1}{\mu_n b_n}\sum_{k=1}^n \mu_k \le
\sum_{n=1}^\infty \frac{1}{\mu_n b_n}\sum_{k=1}^n \mu_k (1+c_k)\le C
\sum_{n=1}^\infty \frac{1}{\mu_n b_n}\sum_{k=1}^n \mu_k ,$$ where
$C=\max_{i: c_i>0}(1+c_i)<\infty$. Hence, condition (9.34) is
stronger than (1.3). In this case, condition (9.34) is even not
needed in part (3) where the last two conditions are automatic.

When $\sum_i\mu_i<\infty$, (1.3) is equivalent to (1.2) which
coincides with (9.34). When $\sum_i\mu_i=\infty$, both (1.3) and
part (3) hold. Therefore, if only a finite number of $c_i$ are
non-zero, then we have
$$\text{criterion } (1.3) \Longleftrightarrow \text{one of parts (2) and (3) holds}.$$
\item"(ii)" When $c_i\ne 0$ for infinite number of $i$, except condition (9.34), an additional condition on
the killing rates $(c_i)$ is required. The condition means that if
$c_i$ increases very fast, then there exist some Dirichlet forms
that do not satisfy the Kolmogorov equations.
\item"(iii)" The second condition in part (3) is the same as the one in part (2). For the first condition in part (3),
it is easy to write down some more explicit sufficient conditions.
This is due to the following fact. Since for each fixed $j$,
$\{P_{ij}^{\min}(\lz): i\in E\}$ is the minimal solution to the
equations
$$x_i=\sum_{k\ne i}\frac{q_{ik}}{\lz+q_i} x_k+\frac{\dz_{ij}}{\lz+q_i},\qqd i\in E,$$
by the linear combination theorem, $\big\{\sum_{j\in
E}P_{ij}^{\min}(\lz): i\in E\big\}$ is the minimal solution to the
equations
$$x_i=\sum_{k\ne i}\frac{q_{ik}}{\lz+q_i} x_k+\frac{1}{\lz+q_i},\qqd i\in E.$$
This minimal solution $(x_i^*)$ can be obtained in the following
way. Let
$$\align
x_i^{(1)}&=\frac{1}{\lz+q_i},\qqd i\in E,\\
x_i^{(n+1)}&=\sum_{k\ne i}\frac{q_{ik}}{\lz+q_i}
x_k^{(n)}+\frac{1}{\lz+q_i},\qqd i\in E, \; n\ge 1.
\endalign
$$
Then $x_i^{(n)}\uparrow x_i^*$ as $n\to\infty$ for every $i\in E$
(cf. [10;  \S 2.1 and \S 2.2]). Hence, for each $n\ge 1$,
$x_i^{(n)}$ is a lower bound of $\sum_{j\in E}P_{ij}^{\min}(\lz)$.
From this discussion, it is clear that the first condition in part
(3) is also a restriction on the growing of the killing rates
$(c_i)$. This is consistent with the second condition there. It is
regretted that we do not know at the moment whether the conditions
in part (3) are necessary or not.
\endroster

\demo{\prf\; of Theorem $9.22$} Part (1) follows from [10; \thm\;3.2] and Chen et.
al. (2005)[1] with a fictitious state $0$. The last cited result is an application of the single birth processes.
Noting that if
$\sum_{i\in E}\mu_i<\infty$ and $\sum_{i\in E}\mu_i c_i<\infty$,
then (9.34) holds iff $ \sum_{i\in E} (\mu_i b_i)^{-1}=\infty$,
hence, part (2) is a special case of [10;  \thm\;6.42]. Next,
noting that the unique exit solution is strictly increasing, when
$\sum_i\mu_i=\infty$, we have ${\scr U}_\lz\cap L^1(\mu)=\{0\}$.
Hence, part (3) is a particular application of [10;
\thm\;6.41]. \qed\enddemo

\head{10. Notes}\endhead

\subhead 10.1\qd Open problems and basic estimates for diffusions\endsubhead

Having seen such a long paper, the reader may feel strange if we claim that the
story is still incomplete even in the context of birth--death processes.
Unfortunately, it is the case.

All of the examples we have done so far show that the following facts hold.
\roster
\item The ratio of the improved upper and lower bounds belongs to $[1, 2]$.
\item The sequence $\{{\bar\dz}_n\}$ is increasing in $n$ and ${\bar\dz}_n\ge \dz_n'$ for all $n$.
\item The sequences $\{{\bar\dz}_n\}$, $\{\dz_n'\}$, and $\{{\dz}_n\}$ all converge to $\lz_0^{-1}$
  as $n\to\infty$.
\item The relation $({\bar\ez}_1, \ez_1)\subset (\kz, 4\kz)$ discussed in Section 6 holds.
\endroster
However, there is still no analytic proof for them. The difficulty for the first
question is that the maximum/minimum of ${\bar\dz_1}$ and $\dz_1$
may locate in different places. In the case that (2) would be true, then the story could be
simplified since we need the first sequence only.
For Questions (2) and (3), the assertions are numerically justified for almost all of the examples in the paper
but the results are not included. We have not worked on Question (3)
hardly enough since one can go ahead only in a finite number of steps in the symbol computation but
the question is certainly meaningful and in the numerical computation,
only in a few steps
one achieves the eigenvalue. For the sequences $\{{\bar\ez}_n\}$ and $\{\ez_n\}$,
we have similar questions as (1) and (3) about, but the corresponding question (2) is answered by \lmm\;6.5.

There is a parallel story for the one-dimensional diffusions. In many cases, one can
easily guess what the result should be, even though there may exist a new difficulty in its proof.
For instance, as a combination of the proofs of \thm\;8.2
and [12; \crl\;7.6], one may prove the following result.

\proclaim{\thm\;10.1} Consider the minimal diffusion on $(-M , N)\,(M, N\le\infty)$ with operator
$$L=a(x) \frac{\d^2}{\d x^2}+ b(x)\frac{\d}{\d x}\qqd \bigg(a(x)>0,\; \frac{b(x)}{a(x)} \text{ is locally integrable\bigg)},$$
and Dirichlet boundaries at $-M$ if $M<\infty$, and at $N$ if $N<\infty$.
Let $C(x)=\int_{\uz}^x b/a$ for some fixed reference point $\uz\in (-M, N)$ and assume additionally
that $e^C/a$ is also locally integrable. Denote by $A_{\Bbb B}$ the optimal constant
in Poincar\'e-type inequality (8.1) with Dirichlet form
$$D(f)=\int_{-M}^N {f'}^2 e^C, \qqd f\in {\scr C}_0^{\infty}(-M, N),$$
Then $A_{\Bbb B}$
satisfies $B_{\Bbb B}\le A_{\Bbb B}\le 4 B_{\Bbb B}$, where
$$B_{\Bbb B}^{-1}=\inf_{ -M<x< y<N}\bigg[\bigg(\int_{-M}^x e^{-C}\bigg)^{-1}+\bigg(\int_{y}^N e^{-C}\bigg)^{-1}\bigg]
\|\dbl_{(x, y)}\|_{\Bbb B}^{-1}.\tag 10.1$$
\endproclaim

By the way, we prove a dual result of \thm\;10.1 for ergodic diffusions. As discussed
in the proof of \thm\;7.5, the exponentially ergodic rate often
coincides with the first non-trivial eigenvalue $\lz_1$ defined
below. Consider a diffusion process with operator $L$ as in
\thm\;10.1, with state space $(-M, N)\,(M, N\le\infty)$ and
reflecting boundaries at $-M$ if $M<\infty$, and at $N$ if
$N<\infty$. For convenience, we define two measures as follows:
$$\align
&\text{Scale measure: } \nu(\d x)=e^{-C(x)}\d x,\qqd
C(x):=\int_{\uz}^x \frac{b}{a},\\
&\qqd \text{where } \uz\in (-M, N) \text{ is a fixed reference
point.}\\
&\text{Speed measure: } \mu(\d x)=\frac{e^{C(x)}}{a(x)}\,\d x.
\endalign
$$
With these measures, the operator $L$ takes a compact form:
$$L= \frac{\d}{\d\mu} \frac{\d}{\d\nu}.  \tag 10.2$$
Next, suppose that $\mu(-M, N)<\infty$, and denote by $\pi$ the
normalized probability measure of $\mu$. Set
$${\scr A}=\{f: f \text{ is absolutely continuous in } (-M, N)\},$$
and define
$$\lz_1=\inf\{D(f): f\in L^2(\mu)\cap {\scr A},\;\pi(f)=0,\; \|f\|=1\},$$
where
$$D(f)=\int_{-M}^{N} a {f'}^2 \d\mu, \qqd f\in {\scr A}.$$
Clearly, in the definition of $\lz_1$, only those $f$ in the set
$\{f\in L^2(\mu)\cap {\scr A}:\; D(f)<\infty\}$ are useful. In other words,
we are here using the maximal Dirichlet form, as in Section 6.

\proclaim{\thm\;10.2} Let $a>0$, $a$ and $b$ be continuous on $[-M,
N]$\,(or $(-M, N]$ if $M=\infty$, for instance). Assume that $\mu(-M,
N)<\infty$. Then for $\lz_1$, we have the basic estimate:
$\kz^{-1}/4\le \lz_1\le \kz^{-1}$, where
$$\kz^{-1}=\inf_{-M<x<y<N}\bigg[
\bigg(\int_{-M}^x \d\mu\bigg)^{-1}+\bigg(\int_y^N
\d\mu\bigg)^{-1}\bigg]\bigg(\int_x^y \d\nu\bigg)^{-1}. \tag 10.3$$
\endproclaim

\demo{\prf} (a) First we show that for the basic estimate,
it suffices to consider the finite $M$ and $N$ with smooth
$a$ and $b$. Since $a$ and $b$ are continuous,
if $M=N=\infty$ for instance, we may choose $M_p,\; N_p\uparrow \infty$
as $p\to\infty$ such that $\uz\in (-M_p, N_p)$ for all $p$. Then, by Chen and Wang (1997, Lemma 5.1),
we have $\lz_1^{(M_p, N_p)}\downarrow \lz_1$ as $p\to\infty$ (This is parallel to the localizing procedure
used in Section 6). At the same time, the isoperimetric
constants ${\kz^{(M_p)}}^{-1}\downarrow \kz^{-1}$ as $p\to\infty$ (cf. proof of \crl\;7.9). Hence, in what follows, we may assume that
$M, N<\infty$. Next, by using the continuity of $a$ and $b$ again, and using a standard smoothing
procedure, we can choose
smooth $a_p$ and $b_p$ such that $a_p\to a$ and $b_p\to b$ (as $p\to\infty$) uniformly on finite intervals, and
furthermore, we can assume that $a_p>0$ on each fixed closed finite interval. Clearly, the corresponding $\kz_p$ converges to
$\kz$ as $p\to\infty$. Therefore, without loss of generality, we assume, unless otherwise stated, that not only
$M, N<\infty$ but also $a$ and $b$ are smooth with $a>0$ on $[-M, N]$.

(b) Recall the following differential form of variational formula for $\lz_1$:
$$\lz_1=\sup_{f\in {\scr F}}\,\inf_{x\in (-M, N)}\bigg[-b' - \frac{a f'' + (a'+b)f'}{f}\bigg](x), \tag 10.4
$$
where
$${\scr F}=\big\{f\in {\scr C}^1[-M, N]\cap {\scr C}^2(-M, N): f(-M)=f(N)=0,\; f|_{(-M, N)}>0\big\}.$$
This is an analog of the variational formula for the lower estimate in \thm\;6.1\,(1).
In the original study by Chen and Wang (1997, (2.3)), the state space is the half-line, not finite, but this is not
essential. It works also for finite state spaces. Besides, it was stated as ``$\ge$'' in (10.4) only.
For ``$=$'', one simply chooses $f=g'$, where $g$ is the eigenfunction of $\lz_1$.
This gives us the boundary condition: $f(-M)=f(N)=0$ since $g'(-M)=g'(N)=0$ by assumption.
Here, one requires that $g\in {\scr C}^3 (-M, N)$ which is satisfied
since we are now in a finite interval having smooth $a$ and $b$. Alternatively, instead of the
original coupling proof, one may use the analytic one which leads to (6.4) for birth--death processes.

(c) We are now going to handle with a more general situation:
$M, N\le \infty$, $a, b\in {\scr C}^1(-M, N)$ and $a>0$ on $(-M, N)$.
Let us define a dual operator $\widehat L$ of $L$. In view of the Karlin and McGregor's construction,
the dual of a birth--death process is simply
an exchange of the scale and speed measures $\hat\mu=\nu$ and $\hat\nu=\mu$ up to a constant (cf. (5.3)).
Thus, in view of (10.2),
the dual operator $\widehat L$, as was introduced by Cox and R\"osler (1983), should be given by
$${\widehat L}= \frac{\d}{\d\hat\mu} \frac{\d}{\d\hat\nu}. \tag 10.5$$
Again, the speed and scale measures $\hat\mu$ and $\hat\nu$ of $\widehat L$ should be expressed as
$$\d {\hat\mu}=\frac{e^{\widehat C}}{\hat a} \d x,\qqd \d {\hat\nu}=e^{-\widehat C} \d x$$
in terms of the coefficients $\hat a$ and $\hat b$ of $\widehat L$ to be determined now.
Because $\hat\mu=\nu$ and $\hat\nu=\mu$, we have
$$\frac{\d\hat\mu}{\d x} \frac{\d\hat\nu}{\d x}=\frac{\d\nu}{\d x} \frac{\d\mu}{\d x}.$$
It follows that $\hat a=a$. Then using the equation $\hat\mu=\nu$, we get
$${\widehat C}= -C + \log \hat a= -C + \log a.$$
Thus, from
$$\frac{\hat b}{\hat a}={\widehat C}'=-\frac{b}{a}+\frac{a'}{a},$$
we get $\hat b=a'-b.$
Therefore, the dual operator $\widehat L$ has the following expression:
$$\widehat L=a(x)\frac{\d^2}{\d x^2}+ \bigg(\frac{\d}{\d x}a(x)-b(x)\bigg)\frac{\d}{\d x}. \tag 10.6$$
For the dual process, the Dirichlet boundary is endowed at $-M$ and $N$ (cf. proof (e) below).
Clearly, the dual operator $\widehat L$ is symmetric on $L^2\big(e^{-C} \d x\big)$.
We remark that the assumption on $a$ and $b$ can be weakened in this paragraph.

(d) Define a Schr\"odinger operator as follows:
$$\align
L_S&= a(x)\frac{\d^2}{\d x^2} + \big(a'(x)+b(x)\big) \frac{\d}{\d x}+ b'(x)\\
&=\frac{\d}{\d x}\bigg(a(x) \frac{\d}{\d x}\bigg)+b(x)\frac{\d}{\d x}+b'(x), \tag 10.7
\endalign$$
with Dirichlet boundaries at $-M$ and $N$ provided they are finite.
Clearly, $L_S$ is symmetric on $L^2\big(e^C \d x\big)$. Denote by $\lz_S$ the
principal eigenvalue of $L_S$:
$$\lz_S=\bigg\{-(f, L_S f)_{L^2(e^C \d x)}: f\in {\scr C}_0^{\infty}(-M, N),\; \int_{-M}^N f^2 e^C=1\bigg\}.$$
In the setup of (b), formula (10.4) becomes
$$\lz_1=\sup_{f\in {\scr F}}\,\inf_{x\in (-M, N)}\frac{- L_S f}{f}(x).$$
This leads to the study on $\lz_S$.

(e) An elementary computation shows that
$$e^C L_S e^{-C}= {\widehat L}. \tag 10.8$$
Note that
$$\int_{-M}^N e^C f L_S g
=\int_{-M}^N e^{-C}(e^C f) \big(e^C L_S e^{-C}\big)(e^C g)
=\int_{-M}^N e^{-C}{\hat f}{\widehat  L} {\hat g},$$
where the mapping $f\to {\hat f}:= e^C f$ is an isometry from $L^2\big(e^C \d x\big)$ to $L^2\big(e^{-C} \d x\big)$,
and that ${\hat f}\in {\scr C}_0^{2}(-M, N)$ iff ${f}\in {\scr C}_0^{2}(-M, N)$.
Since ${\scr C}_0^{2}(-M, N)$ is also a common core of $L_S$ and $\widehat L$ by the assumption on
the coefficients $a$ and $b$,
it follows that the operators $L_S$ and $\widehat L$ with the same core
${\scr C}_0^{\infty}(-M, N)$ are isospectral.
In particular, we have $\lz_S={\hat\lz}_0$.
When $M, N<\infty$, this means that $\widehat L$ has Dirichlet boundaries at $-M$ and $N$ since so does $L_S$.
Now, the basic estimates for $\lz_S$ can be obtained in terms of the ones for ${\hat\lz}_0$, as will be shown in part (f)
below.

To go back to $\lz_1$, noting that by (10.8) again, we also have
$$\frac{-L_S f}{f}=\frac{-(e^C L_S e^{-C})(e^C f)}{e^C f}=\frac{-{\widehat L}{\hat f}}{\hat f}.$$
By (a), we can assume that $M, N<\infty$ and $a>0$ on $[-M, N]$.
From Shiozawa and Takeda (2005) and X. Zhang (2007), it is known that
$${\hat\lz}_0 \ge \sup_{f\in {\scr F}}\;\inf_{x\in (-M, N)} \frac{-{\widehat L}{\hat f}}{\hat f}(x)$$
(i.e., Barta's inequality).
To see that the equality sign holds, simply choose $f$ to be the eigenfunction ${\hat g}$ of ${\hat\lz}_0$.
The fact that ${\hat g}\in {\scr F}$ is guaranteed by the assumptions that $M, N<\infty$, ${\hat a}$ and ${\hat b}$ are continuous,
and ${\hat a}>0$ on $[-M, N]$. This is a standard (regular) Sturm--Liouville eigenvalue problem.
The property ${\hat g}|_{(-M, N)}>0$ is due to the fact that ${\hat \lz}_0$ is the minimal eigenvalue.
We have thus returned to $\lz_1$ from ${\hat\lz}_0$ through $\lz_S$.

(f) For the dual operator $\widehat L$ defined in part (c), applying \thm\;10.1 to ${\Bbb B}=L^1(\hat\mu)$, we obtain
$\hat\kz^{-1}/4\le \hat\lz_0\le \hat\kz^{-1}$, where
$${\hat\kz}^{-1}=\inf_{-M<x<y<N}\bigg[
\bigg(\int_{-M}^x \d\hat\nu\bigg)^{-1}+\bigg(\int_y^N
\d\hat\nu\bigg)^{-1}\bigg]\bigg(\int_x^y \d\hat\mu\bigg)^{-1}.$$
Now, the theorem follows by the dual transform $\hat\mu=\nu$ and $\hat\nu=\mu$.

Finally, the proof of \thm\;10.2 can be summarized as follows:
$$\align
\lz_1& \text{ for general $M, N$ and continuous $a, b$}\\
&\to \lz_1 \text{ for finite $M, N$ and smooth $a, b$}\\
&\qqd \text{ (by approximating and smoothing procedure)}\\
&\to \lz_S \text{ (by coupling method leading to the Schr\"odinger operator)}\\
&\to {\hat\lz}_0 \text{ (by isometry in terms of the dual operator)}\\
&\to \text{ basic estimate of ${\hat\lz}_0$ (by capacitary method: \thm\;10.1)}\\
&\to \text{ basic estimate of }\lz_1 \text{ (by duality)}.\qed
\endalign
$$
\enddemo

Actually, we have also proved the following result (cf. parts (c)--(f) in the last proof)
which is parallel to \prp\;9.11.

\proclaim{\prp\;10.3} Let $M, N\le \infty$, $a, b\in {\scr C}^1(-M,
N)$ and $a>0$ on $(-M, N)$. Then for the Schr\"odinger operator
$L_S$ on $L^2\big(e^C \d x\big)$ having the form (10.7) with
Dirichlet boundaries at $-M$ if $M<\infty$, and at $N$ if
$N<\infty$, we have $\lz_S={\hat\lz}_0$, and furthermore,
${\kz}^{-1}/4\le \lz_S\le{\kz}^{-1},$ where $\kz$ is defined by
(10.3).
\endproclaim

The following simplified estimate of $\kz^{(10.3)}$ is helpful in
practice. Recall that by assumption, $\mu(-M, N)<\infty$. Let
$m(\mu)$ be the median of $\mu$ (i.e., $\mu(-M, m(\mu))=\mu(m(\mu),
N)$). Given $x\in (-M, m(\mu))$, let $y=y(x)$ be the unique solution
to the equation: $\mu(y, N)=\mu(-M, x)$. The A-G inequality
$\az+\bz\ge 2 \sqrt{\az\bz}$ suggests the use of $y(x)$, which then
leads to a simpler bound:
$$\kz^{(10.3)}\ge 2^{-1}\sup_{x\in (-M,\, m(\mu))}\,\mu(-M, x)\, \nu(x, y(x)).$$
We remark that the equality sign here holds in some cases, but the
inequality sign can happen in general. Anyhow, this provides us a
guidance in seeking for the infimum in (10.3). Certainly, the
similar discussion is meaningful for $\kz^{(10.1)}$.

Having \thm s 10.1 and 10.2 at hand, the basic estimates in the other
cases ($\lz^{\text{\rm ND}}$ and $\lz^{\text{\rm DN}}$) mentioned in Section 1 should be clear.

The study on the one-dimensional case provides a comparison tool for the study on
the higher dimensional situation, as we did a lot before. Hence, there is no doubt
for the development in the higher dimensional context.

\subhead 10.2\qd ${\pmb h}$-Transform\endsubhead

In an earlier draft of this paper (roughly speaking, up to \thm\;7.1 plus a part of \thm\;9.3),
the author mentioned an open question: how to handle the case that
(1.3) fails? Then two answers have appeared. The first one is the use of so-called $h$-transform by
Wang (2008a) where the transient case studied in Section 7 is transferred
into the one studied in Section 4. Next, with the help of the duality given in \thm\;7.1,
the ergodic case studied in Section 6 can be also transferred into the one studied in Section 4.
In this way, with a use of \thm\;4.2, Wang obtains a criterion for $\lz_1$ (Section 6) with a
factor 4. To have a taste of this technique, let us quote a particular result here.

\proclaim{\thm\;10.4\,{\rm(Wang (2008a, \thm\;1.2))}} Set $h_i=\sum_{j=i}^N \mu_j$. Then we have
$\dz^{-1}/4\le \lz_1\le \dz^{-1}$, where
$$\dz=\sup_{1\le i< N+1} \bigg(\frac{1}{h_i}-\frac{1}{h_0}\bigg)\sum_{j=i}^N \frac{1}{\mu_i a_i}h_i^2. \tag 10.9$$
\endproclaim
 Comparing this result with \crl\;6.6, the factor 4 is in common but the
isoperimetric constants are quite different.
The advantage here is that only one variable is required in the supremum, but in \crl\;6.6
two variables are needed. The price one has to pay to (10.9) is involving a new quantity $h$.
The natural extension of \crl\;6.6 to the whole line (\crl\;7.9) exhibits an interesting
symmetry of the left and the right half-lines.
Such an extension of \thm\;10.4 with the same factor 4 is unclear to the author.
Along the same line and using [9], Wang then extends the results to
Poincar\'e-type inequalities as well as functional inequalities, see Wang (2008b, c).
Clearly, Wang's papers show that the $h$-transform is a powerful tool and may be
useful in other cases.

While the author's solution to the above open question is the use of the maximal
process as included into this version of the paper. As shown in the paper,
\crl\;6.6 comes from the author's previous general result without using the $h$-transform.
An interesting question in mind is to use the variational formulas in Section 6 to derive
\crl\;6.6 directly. Besides, a direct generalization of Sections 2, 3, and 7 to the
Poincar\'e-type inequalities is still meaningful in practice since the formulas
are quite different (in view of \thm\;10.4) and some of them may be more practicable.

\subhead 10.3\qd Remark on some known results\endsubhead

As mentioned in Section 5, duality (5.1) goes back to
Karlin and McGregor (1957b). The author learned this technique
mainly from van Doorn (1981; 1985) based on which the proof of the
basic result $\lz_1=\az^*$ was done, cf. [2]. It is now known that such
a result holds in a very general setup as indicated in the proof of \thm\;7.4.

We now discuss the situation that (1.3) holds. Then there are three cases.

\roster
\item $\sum_i \mu_i=\infty$ and $\sum_i (\mu_i b_i)^{-1}<\infty$.
\item $\sum_i \mu_i<\infty$ and $\sum_i (\mu_i b_i)^{-1}=\infty$.
\item $\sum_i \mu_i=\sum_i (\mu_i b_i)^{-1}=\infty$.
\endroster

First, let $b_0>0$. In cases (1) or (3), by \thm\;2.4\,(1) and \prp\;2.7\,(1), $\lz_0^{(2.2)}$ is
equal to
$$\gather
\sup_{v\in {\scr V}}\inf_{i\ge 0}\,[a_{i+1}+b_i-a_i/v_{i-1}-b_{i+1}v_i].\tag 10.10\\
{\scr V}:=\{v: v_{-1} \text{ is free},\; v_i>0\text{ for all }i\ge 0\}.
\endgather$$
In case (2), by \thm\;6.1\,(1), $\lz_1$ can be expressed by (10.10).
Thus, in view of \prp\;1.2 and [2; \thm\;5.3], the convergence
rate $\az^*$ can be also expressed by (10.10).

Next, let $b_0=0$. Then in case (2), by \crl\;5.2, \prp\;2.7\,(1), and using (5.8) in an inverse way, it follows
that $\lz_0^{(4.2)}$ is equal to
$$\gather
\sup_{v\in {\scr V}}\inf_{i\ge 1}\,\bigg[a_i\bigg(1-\frac{1}{v_{i-1}}\bigg)+b_{i}(1-v_i)\bigg],\tag 10.11\\
{\scr V}:=\{v: v_{0}=\infty,\; v_i>0\text{ for all }i\ge 1\}.
\endgather$$
In case (1), $\lz_0^{(4.2)}$ is equal to $\lz_0^{(7.1)}$.
By \thm\;7.1\,(1), in terms of \thm\;6.1\,(1) and using (5.8) in an inverse way,
we obtain the same expression (10.11) for $\lz_0^{(7.1)}$.
Finally, in the degenerated case (3), we indeed have $\lz_0^{(4.2)}=\lz_0^{(7.1)}=0$ which can be expressed as (10.11)
by \thm\;7.1\,(2). Hence,
by \prp\;1.2, the convergence rate $\az^*$ can also be expressed by (10.11).
We have thus obtained the following result.

\proclaim{\thm\;10.5\,{\rm (van Doorn (2002))}} Let $(1.3)$ hold. Then the exponential convergence rate
$\az^*$ is given by $(10.10)$ or $(10.11)$, respectively, according to $b_0>0$ or $b_0=0$.
\endproclaim

With a slightly different expression, this result was given in van Doorn (2002)
by the analysis on the extreme zeros of orthogonal polynomials in Karlin and McGregor's representation,
and was actually implied in van Doorn's earlier papers (1985; 1987)
as mentioned in the paper just cited or in [3]. In the last paper, this
result was rediscovered in the study on $\lz_1$, using the coupling methods.
The lower estimate was also obtained by
Zeifman (1991) using a different method in the case that the rates of
the processes are bounded, with a missing
of the equality.

A progress made in the paper is removing Condition (1.3) and even (1.2). In particular,
the situation having finite state spaces is included. This is meaningful
not only theoretically but also in practice since
the infinite situation can be approximated by the finite ones. Besides, when
$b_0>0$ and $\sum_i \mu_i<\infty$,
the duality given by (5.9) is essentially different from (5.8) (cf. Remark 2.8). From the
other point of view, the dual of this case goes to $\lz_0^{(7.1)}$ rather than $\lz_0^{(4.2)}$.
However, we then have to use the maximal process in Section 6, as we did in \thm\;7.1, rather than
the minimal one used in Sections 2 and 3, except using (1.3)
(which is equivalent to (1.2) if $\sum_i \mu_i<\infty$). From analytical point of view,
the use of the maximal process is natural since one looks for the inequality to be
held for the largest class of functions, as illustrated by the weighted Hardy inequality
in Section 4.1.

In van Doorn (2002), some variational formulas of difference form for the upper bound
of $\az^*$ are also presented but we do not use them here.
As far as we know, the criterion
for $\az^*>0$ (\thm\;1.5) has been open for quite a long time; it was answered in the ergodic case
only till [6] in terms of the study on the first non-trivial eigenvalue $\lz_1$. For which, the
criterion was obtained independently by Miclo (1999) based on the weighted Hardy's inequality.
Criterion 3.1 follows from the variational formulas of single summation form (part (2) of \thm\;2.4),
but it is not obvious at all to deduce the criterion from (10.10) (or dually from (10.11)) directly.
More clearly, the variational formula of the difference form
for the lower bound given in (9.3) which is closely related to (10.11) was known for some years
and works in a more general setup, but an explicit criterion for the
killing case is still open (Open Problem 9.13).
Anyhow, having the duality (\crl\;5.2 and \thm\;7.1) at hand, \thm\;1.3 is essentially known from [6],
except the
basic estimates in the ergodic case as well as in the setting of Section 7 is presented here for the first time.
The technique adopted in this paper
depends heavily on the spectral theory, potential theory, and harmonic analysis.
In the transient continuous context, Criterion 3.1 was obtained by
Maz'ja (1985, \S 1.3), as a straightforward consequence of Muckenhoupt (1972).
The discrete version was proved by Mao (2002, \prp\;A.2). In these quoted papers,
the problem in a more general ($L^q$, $L^p$)-setup was done.

In the continuous context, the Hardy-type or Sobolev-type
inequalities (cf. \thm\;10.1) were studied first by P. Gurka and
then by Opic and Kufner (1990, \thm\;8.3). Instead of ${\scr
D}^{\min}(D)$, they considered the following class of functions: the
absolutely continuous functions vanishing at $-M$ and $N$.
This seems not essential in view of $\lz_0^{(2.2)}=\lz_0^{(2.18)}$. With a
different but equivalent isoperimetric constant (i.e., replacing the
sum in (10.1) by maximum ``$\vee$''), they obtained upper and lower
bounds with ratio $2\,\oz^5\approx 22$, where
$\oz=\big(\sqrt{5}+1\big)/2$ is the gold section number. By the way,
we mention that the use of weight functions $w$ and $v$ in (8.6) in
the cited book is formally more general than our setup. One can
first assume that $w$ and $v$ are positive, otherwise replace them
by $w+1/n$ and $v+1/n$, respectively, and then pass to the limit as
$n\to\infty$. Next, it is easy to rewrite $w$ and $v$ as $e^{C}/a$
and $e^{C}$ for some functions $C$ and $a>0$. Note that only $C$ and
$a$ (without using $b$) are needed to deduce the basic estimates in
our proof. Again, in the continuous context, the splitting technique
was also used in \thm\;8.8 of the book just quoted where some basic
estimates were derived in terms of an isoperimetric constant, up to
a factor 8. Their isoperimetric constant is parallel to the
right-hand side of (7.13) replacing $\lz_0^{\uz\pm}$ by the
corresponding $\dz^{(3.1)\pm}$ depending on $\uz$ (certainly,
without using the parameter $\gz$). Our Example 8.9 is an analog of
Examples 6.13 and 8.16 in the quoted book. In contrast with our
probabilistic--analytic proof here, their proof is direct, analytic,
and works in a more general ($L^q$, $L^p$)-setup. We have not seen
the discrete analog of their results in the literature. In the
($L^p$, $L^p$)-sense $(p\ge 1)$, the variational formulas in the
continuous context were obtained in Jin (2006) but it remains open
for the more general ($L^q$, $L^p$)-setup. Even though it is a
typical Sturm-Liouville eigenvalue problem having richer literature,
we are unable to find an analog of \thm\;10.2.

Finally, in computing the examples in the paper, the author uses the software
Mathematica. All the examples were checked by Ling-Di Wang and Chi Zhang using
MatLab. Most of the author's papers cited here can be found in
[8].

\bigskip

\flushpar {{\bf Acknowledgements}} \quad The first version of
the paper dated May 25, 2007.
In June, 2008, the author was invited by Professors Tiee-Jian Wu and Yuan-Shih Chow
to visit the Center for Theoretical Sciences at  Cheng Kung University
and  Taiwan University. The author
was also invited by Professors Mong-Na Lo Haung and Mei-Hui Guo (Sun Yat-sen University),
by Professor Wen-Jang Huang (University of Kaohsiung), by Professor Chen-Hai Tsao
(Dong Hwa University
for the 17$^{\text{th}}$ Conference on Statistics in the Southern area and the Annual Meeting of Chinese
Institute of Probability and Statistics), by Professors Chii-Ruey Hwang, Tzuu-Shuh Chiang,
Yunshyong Chow and Shuenn-Jyi Sheu (Academia Sinica) for visiting
their institutes. The scientific communications, the very warm hospitality of the inviters and the
financial support from the
Center for Theoretical Sciences, and the inviters' universities are deeply acknowledged. The
author has also benefited from several discussions with Professor Shuenn-Jyi Sheu in Beijing.
A part of the results in the paper were reported during the visit.

Partial results of the paper have also been reported in the following conferences:
International Conference on Stochastic Analysis and Related Fields (April 2008, Wuhan),
Workshop on Probability (November 2008, Guangzhou),
Chinese-French Meeting in Probability and Analysis (September 2008, Marseille),
Workshop on IPS in honor of Professor Tom Liggett's 65th Birthday (June 2009, Peking University),
IMS-China International Conference on Statistics and Probability (July 2009, Weihai).
The author acknowledges the organizers Zhi-Ming Ma and Michael R\"ockner, Jia-Gang Ren,
Dominique Bakry and Yves LeJan, Da-Yue Chen, Jia-An Yan, respectively, for their invitation and financial
support.

Because of an honorary professorship from Swansea University, the author
was invited to a workshop and visited University of Strathclyde in UK
in October, 2009. The recent progress was reported at these two universities.
The author appreciates the invitation, the very warm hospitality, and the financial
support by Professors N. Jacob, X.R. Mao, F.Y. Wang, J.L. Wu, C.G. Yuan, and their
institutes. Special thanks are given to the speakers of the celebration workshop for their kindness: T.J. Lyons, M. Davis,
F.Y. Wang, T.S. Zhang, X.R. Mao, and D. Elworthy.

Two drafts of the paper have been reported in a series of seminars
organized by Professors Yong-Hua Mao and Yu-Hui Zhang. Especially,
the drafts are checked  line by line by Zhang's graduate students
Ling-Di Wang and Chi Zhang. The corrections and improvement obtained
from them improved a lot the quality of the paper. Their serious
work is acknowledged. Since May 2007, the author has reported the
results from time to time in his seminars and thanks are given to the
participants for their discussions and comments.


\enddocument